# Towards A Better Understanding of the Leading Digits Phenomena

-------------------------------------------------------------------------------------

## Abstract

That the logarithmic distribution manifests itself in the random as well as in the deterministic (multiplication processes) has long intrigued researchers in Benford's Law. How or why does it appear in two such distinct processes?! In this article it is argued that it springs from one common intrinsic feature of their density curves. On the other hand, the profound dichotomy between the random and the deterministic in the context of Benford's Law is noted here, acknowledging the need to distinguish between them. From its very inception, the field has been suffering from a profound confusion and mixing of these two very different logarithmic flavors, causing mistaken conclusions. One example is Allaart's proof of equality of sums along digital lines, which can only be applied to deterministic processes. Random data lack this equality and consistently show significantly larger sums for lower digits, thus rendering any attempt at test of summation equality irrelevant and futile in the context of forensic analysis regarding accounting and financial fraud detection.

Another digital regularity is suggested here, one that is found in logarithmic as well as non-logarithmic random data sets. In addition, chains of distributions that are linked via parameter selection are found to be logarithmic, either in the limit where the number of the sequences in the chain approaches infinity, or where the distributions generating the parameters are themselves logarithmic.

A new forensic data analysis method in the context of fraud detection is suggested here even for data types that do not obey Benford's Law, and in particularly regarding tax evasion applications. This can also serve as a robust forensic tool to investigate fraudulent fake data provided by the sophisticated and well-educated cheater already aware of Benford's Law, a challenge that would become increasing problematic to tax authorities in the future as Benford's Law becomes almost common knowledge.

## [1] The Leading Digits Phenomena

Leading digits (LD) or first significant digits are the first digits of numbers appearing on the left. Simon Newcomb [N] in 1881 and then Frank Benford [B] in 1938 discovered that low digits lead the first order much more often than high digits in everyday and scientific data and arrived at the exact expression Prob[1st digit is d] = $\text{LOG}_{10}(1 + 1/d)$ known as the logarithmic distribution as the probability that digit d is leading. This is known as Benford's Law. Consideration is also given to second leading digits and third significant digits and so forth for higher order of appearances of digits. The Generalized Benford's Law takes into account all higher orders leading digits, and it will be stated and detailed in later sections. Common observations in the old days before the advent of the hand calculators and computers was that books of tables of logarithms and roots and such were almost always much more worn out at the beginning pertaining to numbers that start with digit 1 or 2 than at the end pertaining to numbers that start with digit 8 or 9, and that wear and tear was progressively less and less severe throughout as value of digit increases. This consistent variation in the use of the pages in almost all such books constituted probably the first hint of the phenomena and may have led Newcomb and Benford to discover and then to investigate the pattern. Here are the exact proportions of Benford's Law for the first digits:
{30.1%, 17.6%, 12.5%, 9.7%, 7.9%, 6.7%, 5.8%, 5.1%, 4.6%}.

Benford's Law is solely a mathematical and statistical fact about how numbers are USED in everyday typical situations expressing physical quantities as well as abstract enteties we fancy contemplating and recording. What was lacking though in the formal statement of Benford' Law concerns a well-defined mathematical and statistical process, an exact sample space for which the law should be applied and examined. Both Benford and Newcomb overlooked this crucial issue. Throughout this article the follwing interpretation is given, assumed and used in order to fascilitate discussion: Benford's Law could be explicitly stated as applying to the totality of all the numbers gathered from everyday and scientific data (in our current modern era), a law about how numbers are being used, occur, and recorded in the physical world or in any abstract manner - when considered in its aggregate. Hence it's also a law about any small subset obtained by sampling from that vast aggreragate – assuming this has been done in a truly random mannar. This definition of the law is not an exact mathematical formulation owing to the fact that this totality is too fluid and constantly changing, and that it can never be measured directy, yet it does enable us to meaningfully assert and state numerous ideas. Only recently this abstract idea has been given a rigrous mathematical foundation by Theodore Hill, creating a specific mathematical model corresponding to real life data and showing that it is logarithmic.

As it happened, the totality of our everyday real data considered as one vast collection of numbers closely mimic a certain random process that is itself in conformity with the logarithmic law (as shown by Ted Hill). It is this abstract entity, the totality of everyday data, a colossal collection of numbers, that Benford's Law in its purest strict form could refer to.

## [2] A Leading Digits Parable

Imagine a simple society in ancient Greece, having very few possessions, with severe limits on the quantities possessed. Believing in 9 Gods, possessing up to 8 olives, 7 dogs, 6 chickens, 5 sheep, 4 oranges, 3 slaves and 2 houses, and having only 1 spouse.

Digital leadership for them depends on the specific topic of their simple conversations. For example, for their number of spouses, they use only the number 1, and for that usage digit 1 leads 100% of the times, and 0% for all other digits. For sheep, being that no one possesses more than five sheep, digits 1 to 5 lead 1/5 each, or equally 20% of the times. A typical conversation may include the statement for example "Yesterday I saw 7 dogs chasing 3 chickens". We wish to aggregate all types of conversations and examine overall resultant LD under the assumption of topic equality. Overall LD distribution by simply calculating the average digital leadership of all the topics digit by digit, as in the following table:

| Digit | Spouse | houses | slaves | oranges | sheep | chicken | dogs | olives | Gods | | Averages |
|---|---|---|---|---|---|---|---|---|---|---|---|
| 1 | 1/1 | 1/2 | 1/3 | 1/4 | 1/5 | 1/6 | 1/7 | 1/8 | 1/9 | -----> | 31.4% |
| 2 | 0 | 1/2 | 1/3 | 1/4 | 1/5 | 1/6 | 1/7 | 1/8 | 1/9 | -----> | 20.3% |
| 3 | 0 | 0 | 1/3 | 1/4 | 1/5 | 1/6 | 1/7 | 1/8 | 1/9 | -----> | 14.8% |
| 4 | 0 | 0 | 0 | 1/4 | 1/5 | 1/6 | 1/7 | 1/8 | 1/9 | -----> | 11.1% |
| 5 | 0 | 0 | 0 | 0 | 1/5 | 1/6 | 1/7 | 1/8 | 1/9 | -----> | 8.3% |
| 6 | 0 | 0 | 0 | 0 | 0 | 1/6 | 1/7 | 1/8 | 1/9 | -----> | 6.1% |
| 7 | 0 | 0 | 0 | 0 | 0 | 0 | 1/7 | 1/8 | 1/9 | -----> | 4.2% |
| 8 | 0 | 0 | 0 | 0 | 0 | 0 | 0 | 1/8 | 1/9 | -----> | 2.6% |
| 9 | 0 | 0 | 0 | 0 | 0 | 0 | 0 | 0 | 1/9 | -----> | 1.2% |

FIG 1: LD Distribution in Ancient Greece

We have tacitly assumed that talking about (or having) 3 olives is just as likely as talking about (or having) 8 olives, so that both come with equal (uniform) probabilities of occurrence. In other words, 2 uniform discrete distributions of occurrences of numbers within each topic itself. For example sheep has uniform distribution for {1, 2, 3, 4, 5} with 1/5 probability, and another one for {6, 7, 8, 9} having 0 probability. The latter default distribution of zero is to be ignored now, and the view is that {6, 7, 8, 9} simply doesn't exist for sheep as part of the probability space.

Taking stock: what we have here is a collection of 9 different uniform discrete distributions where each (per topic) is accorded equal importance. All distributions start at 1 which is the lower bound (LD) for all of them, hence 1 is the anchor orienting, synchronizing and uniting them all in certain positions – so to speak. Length of distributions, that is, their upper bound (UB) is gradually increased from 1 (in the spouse case where LB=UB=1) all the way to 9.

Could this model be relevant at all to Benford's Law and our modern life? The answer is a resounding yes, albeit with some reservation and considerable modification! in spite of the fact that it's quite a restricted model, with 9 as the UB for all distributions as oppose to much higher ranges of values in our modern everyday data, and with no decimals or ratios allowed at all, with negative numbers totally omitted, square roots greatly feared and strictly avoided, and even with the absence of that ubiquitous Normal distribution, so often observed in everyday data!   The remarkable thing is that with a few modifications and improvements, this Greek Parable leads directly (arithmetically that is) to the logarithmic distribution (as will be demonstrated in the next two sections), and with the rationale behind it as an appropriate model representing real life data in some sense or partially still in tact or valid, namely, that that vast aggregate entity we call 'the totality of our everyday data' could – in a way and perhaps only partially - be represented by a huge collection of uniform distributions all sharing that common anchor 1, serving to orient numbers in much the same way as for those ancient Greeks in the parable above!

An important observation here indicating the more universal nature of applicability of the Greek Parable is its total independence on any scale or units of measurements as well as on any number base system such as 10, 16, 2, or e. Since 9 is the upper limit here, there is no need to express say 247 as 2*1+4*10+7*100 for example, and there is no length or weight to measure anything in feet, meters, kilos, or pounds. In other words, we haven't insisted on using the arbitrary base 10 (our fingers) to obtain our result, nor on the French Kilo exclusively to express weights, hence giving us hope that the result thus obtained is quite general! Another section on scale-invariance and base-invariance is devoted to this important aspect of LD.

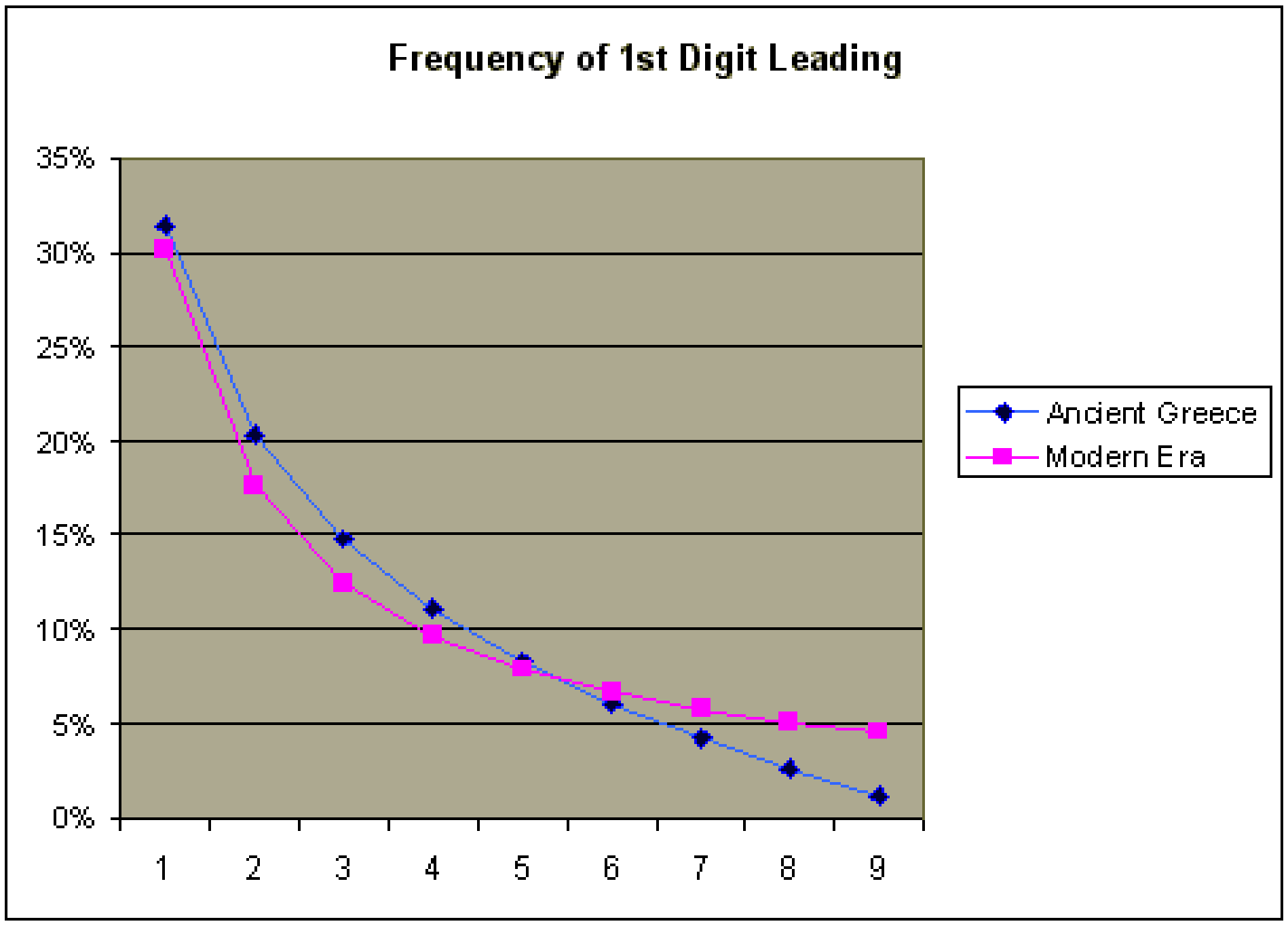

FIG 2: LD Distribution in Ancient Greece vs. Modern Era (the logarithmic)

## [3] Integral Powers of Ten

The abbreviation **IPOT** would stand for an **'integral power of ten'** (singular), namely $\mathbf{10}^{\textbf{integer}}$ (where integer can be negative or zero as well) such as 0.01, 0.1, 0, 1, 10, 100, and so on. The abbreviation **TAIPOT** would stand for **'two adjacent** *(consecutive)* **integral powers of ten'** (plural) namely for $\mathbf{10}^{\textbf{integer}}$ & $\mathbf{10}^{\textbf{integer + 1}}$ such as 10 & 100 , 1 & 10, and so forth, where exponents (log values) differ exactly by unity 1. Integral powers of ten play a crucial role in the understanding of Leading Digits. The value of 9 in the Greek parable above was purposely chosen as that single-digit number closest to 10, as 10 itself is an integral power of ten. Often results and definitions are stated in terms of intervals between TAIPOT numbers, but at times the only thing that is required actually is an integral difference in exponent that is $10^{\textbf{Non-integer}}$ & $10^{\textbf{Non-integer + 1}}$ such as the pair $10^{1.7}$ & $10^{1.7+1}$ namely 50.119 and 501.19.

Beware of calling 0 & 1 TAIPOT numbers, they are not! Calling attention to these two very important numbers demonstrates the unique nature of the interval (0, 1) in LD. An <u>infinite</u> number of IPOT numbers exist in that 'small' interval ! Common (decimal) LOG there varies as in (-∞, 0)! In contrast, on the interval [1, 1000] there exist only four IPOT numbers and decimal LOG varies only from 0 to 3.

Surely 10 was chosen for being the **base** in our number system, and all this generalizes to any other base.

## [4] Simple Averaging Scheme as a Model for Typical Data

Our civilization at the current epoch is certainly much more sophisticated than the society described in the Greek parable. We have use of vastly more topics of conversations than that meager 9-item set, and we have greatly expanded our number system well beyond their limit of 9 whole numbers. We have also advanced beyond their simple count data and use units and scales to record weights, lengths, time etc. Yet, if we were to attempt then to expand on the Greek parable to include the vastly more relevant examples and huge variety of typical speeches and recorded data, we would obviously face a daunting or rather impossible task. How does one go about accounting for, or aggregating LD of say stock prices, expenses at IBM, weights of people, age at death, accidents per year, and so fort, to mention only a tiny portion of everyday data?! What we need here is a limited and well-structured mathematical model that is capable somehow of capturing or representing our typical usage of numbers, to some extent at least. Interestingly, as it turned out, imitating the analysis of the Greek parable – albeit with a much longer range of integers – comes a bit close to the logarithmic! Surely we can't rest there, the resultant discrepancy will have to be ironed out, modifications applied and justified, yes this approach leads to a great deal of deeper understanding here and as such this path will be taken for now.

And so we consider what would happen (collectively) to leading digits distributions of numerous intervals, all starting at 1, made only of the integers, while differing only in their lengths. Those intervals should be made progressively longer by systematically increasing their length, one integer at a time. The plan is then to obtain an aggregate LD distribution representing ALL intervals simply by taking the average of the LD distributions of all intervals, hoping to encounter at the end of the calculating process something quite close to the desired logarithmic distribution, and with another hope in mind, namely that of being able to argue that this scheme is somehow a good representative of real life typical data, (it is not!)

The individual elements of this model are the sets of consecutive integers (intervals) on [1, $N_i$] where each integer within its interval posses equal probability of being picked as in the discrete uniform distribution. Also, each element (an interval) is accorded equal probability of occurring within the entire model itself. This last requirement should be noted carefully, since it implies that LD of the shortest interval say [1, 9] gives as much weight to the overall LD as LD of the longest one say [1, 99999], in spite of the fact that it contains fewer integers by far!!! Note that all intervals start at 1, and we only need then to decide on the starting and ending values of $N_i$ (the upper bounds), in other words, decide on the size of the smallest interval and on the size of the largest interval. In order to decide on this issue, let us first explore the ramification of having a particular N as the highest integer for a SINGLE individual interval, and in particular let's compare the two intervals [1, 3000] and [1, 9999]. If N happens to be a number composed entirely of nines, such as 99, 999, 9999 and so forth, then clearly all nine digits would have equal 1/9 probability of leading, a situation that is quite preferable since this does not bias the model towards any particular digit in the first place. But on the other hand, if N is 3000 say, then digits 1 and digit 2 would have strong advantages and each would lead 1112 numbers out of 3000 total, or 1112/3000, that is 37% each, and 74% for 1 and 2 combined, while all other digits would lead with probability of 111/3000, or a mere 3.7%, a biased situation we may want to avoid. A moment's thought will then convince us that letting the largest interval in the entire scheme fall far from some number made of nines would bias the model towards low digits unfairly. Equally unfairly is to let the smallest interval fall far away from some number made of nines. Hence, it is clear that $N_i$ should vary from a number of nines to another (higher) such number, or else from 1 itself (the LB!) to some large enough number of nines, or, as should be quite obvious now, to an IPOT number such as 10, 100, 1000, and so on. We are not arguing that real life data itself tend to abruptly terminate exactly at 1000, or 100,000 and such, this is obviously not the case, instead we are merely trying to calibrate our mathematical model to avoid any bias or arbitrariness as much as possible.

**How good a representative of real life typical data would such a scheme be?**
Well, numbers in the real world come with incredibly varying upper bounds, besides, they do not increase nicely one integer at a time, and not all of them really start with 1 or even 0. For example, heights of people start perhaps around 1.3 meter with UB around 2.3 meter. Yet, what the scheme imitates quite well is this huge variability in UB as compared with the relative stability of LB which is so typically stuck at 0 or 1!

Is the exclusive usage of integers severely limiting? No! Because fractions can be easily incorporated and attached to each integer in the model without altering the arithmetical result at all. What about the unrealistic assumptions that numbers are uniformly distributed within intervals and that intervals come with equal importance throughout the model? Certainly everyday data frequently come with non-uniform distributions, the Normal distribution is much more common, and of course very few pieces of data (if any) should be accorded equal importance within the totality of everyday data, but it could be argued that this shouldn't favor high digits over low ones or visa versa, as it can work both ways.

The model as described above shall be called "the simple averaging scheme". Below are the computerized results of a few such schemes:

| | Ranges of Upper Bound: | | | | | |
|---|---|---|---|---|---|---|
| Digit | 1 to 9 | 50 to 700 | 10 to 100 | 1 to 100 | 1 to 10,000 | 1 to 100,100 |
| 1 | 31.4% | 28.6% | 24.5% | 25.1% | 24.2% | 24.1% |
| 2 | 20.3% | 20.5% | 18.4% | 18.6% | 18.3% | 18.3% |
| 3 | 14.8% | 15.2% | 14.5% | 14.6% | 14.5% | 14.5% |
| 4 | 11.1% | 11.3% | 11.7% | 11.6% | 11.7% | 11.7% |
| 5 | 8.3% | 8.1% | 9.4% | 9.3% | 9.5% | 9.5% |
| 6 | 6.1% | 5.3% | 7.6% | 7.4% | 7.6% | 7.6% |
| 7 | 4.2% | 3.9% | 6.0% | 5.8% | 6.0% | 6.0% |
| 8 | 2.6% | 3.7% | 4.6% | 4.4% | 4.6% | 4.7% |
| 9 | 1.2% | 3.5% | 3.3% | 3.1% | 3.4% | 3.4% |
| | | Note: All schemes are with lower bound at 1. | | | | |

FIG 3: Simple Averaging Schemes for Various Intervals **(city)**

All this only resembles the logarithmic, but it is not quite close enough. Notice that the Greek Parable is in the table on the left as '1 to 9', and obviously having too limited a range as compared with the rest, which is the cause of its somewhat poorer resemblance to the logarithmic. In the limit where the upper edge of UB is quite large this distribution is known as **Stigler's Law**, which is the last column on the right for one decimal precision. The biased scheme with UB range of 50 to 700 was included here as an illustrative odd case, but even this result is not too different.

What went wrong? Why the discrapancy between this and the logarithmic?
An argument can be raised that as complex as the scheme seems to be, it is still not comlex enough as compared with the totality of everyday data, in other words, it is still too narrow and specific corresponding actually to a particular set of numbers, in a nutshell: **it doesn't average enough things**, hence its resultant LD distribution lacks universality. For a better understanding here let's look into one of the most typical of count data, namely address data, itself an example of real life data.
Street address data, refers specifically to the house number. All addresses start at 1 but street lengths vary, and therefore upper bound – the house with the highest numerical address - is different perhaps for each street, hence only UB is to vary, not LB which has to be fixed at 1. Clearly, our simple averaging scheme would be a perfect representation of the address data if we knew the length of the shortest and the longest streets, and also if we were assured that all street lengths in between are distributed uniformly, namely a fixed and steady increase in length (totally unrealistic). Here, each physical street in the street address data corresponds to an interval [1, Ni] in our abstract model. Another important consideration is that typically upper bounds for street addresses do not vary between IPOT numbers of course. Say we analyse street data for the **city** of Chicago and UB there vary (smoothly and gradually) from 7 to 34098, hence the simple averaging shceme is performed there, aggragate and LD recorded. Same is done for Miami, New York and some other big cities with slighly different UB variability. But the statistician wants to include small towns as well, hence the town of Suffern in NY is included, and there UB varies from 8 to 285. Of all the places in the USA only in one rare and remarkable case, a small town in Pensilvania, that we found UB varying there exactly between 10 and 100. In any case, **Clearly what is needed here is to average results from multiple simple averaging schemes all with different UB layouts!** Hence one can adapt the following point of view: The simple averaging scheme was good enough for just one city in the USA, and the need for a more complex averaging scheme is essentially the need to get the aggragate LD distribution of the whole **country** for all the cities and all the towns. Satisfying this need would be the topic of our next section.

We end this section by noting the remarkable robustness of this scheme albeit in some weak sense. Had we averaged with UB varying between 10 and 90 (instead of the more proper way between 10 and 100) low digits would win even more as digit 9 suffers a great deal here. Had we been overshooting by setting UB to vary between 10 and 110 low digits still win and even more so, as digit 1 would obtain an extra advantage on the interval (100, 110)! So the general trend holds favoring low digits over high ones! This feature gives some kind of stability and consistency to our simple averaging scheme.

## [5] More Complex Averaging Schemes

The next level of complexity is averaging out multiple simple averaging schemes themselves, while gradually varying their focus. This is done by letting LB be fixed (typically at 1), letting UB start also at a fixed point for all schemes (min UB), but allowing it to terminate on different locations. An implicit (and unrealistic) assumption is made if applied to street address, namely that the shortest street (min UB) is the same for ALL cities, although arithmetically this doesn't effect results much. Also, if applied to street address, then min-UB naturally should be fixed at 1 or just slightly above 1 (being the shortest street), while max-UB and all higher parameters should be made much larger than 1 as they stand for the limits of the longest streets, (the table below violates this requierment as it pertains to a generic scheme). For example, consider multiple simple averaging schemes where lower bound is always 1, while upper bound varies from **1 to 100**, **1 to 101**, **1 to 102**, … , **1 to 9,999**, and **1 to 10,000**, and which is then followed by the averaging of all the 9,901 different averages! Few such schemes are performed here, results of which are given below:

| | | | | | | |
|---|---|---|---|---|---|---|
| max UB Highest | 100 | 1,000 | 1,000 | 1,000 | 100 | 4,500 |
| max UB Lowest | 20 | 1 | 100 | 500 | 1 | 800 |
| min UB | 10 | 1 | 10 | 15 | 1 | 30 |
| LB | 1 | 1 | 1 | 1 | 1 | 1 |
| ---------------------------------- | ------------- | ------------- | ------------- | ------------- | ------------- | ------------- |
| Digit 1 | 33.9% | 30.9% | 30.4% | 28.3% | 34.7% | 32.0% |
| Digit 2 | 19.6% | 19.0% | 18.9% | 20.2% | 19.7% | 16.8% |
| Digit 3 | 12.6% | 13.2% | 13.1% | 14.9% | 13.0% | 11.1% |
| Digit 4 | 8.9% | 9.8% | 9.8% | 10.9% | 9.3% | 8.9% |
| Digit 5 | 6.7% | 7.6% | 7.6% | 7.9% | 6.9% | 7.8% |
| Digit 6 | 5.4% | 6.1% | 6.2% | 5.9% | 5.4% | 6.9% |
| Digit 7 | 4.6% | 5.1% | 5.2% | 4.6% | 4.3% | 6.1% |
| Digit 8 | 4.2% | 4.4% | 4.6% | 3.9% | 3.6% | 5.5% |
| Digit 9 | 4.0% | 4.0% | 4.1% | 3.4% | 3.1% | 4.9% |

FIG 4: Average of Averages Schemes **(country)**

The resultant LD distribution here is much superior to the simple averaging scheme (Stigler's Law), as it is much closer to the logarithmic. Yet, it's not quite close enough to the logarithmic. Also of note is the last column on the right, where the arbitrary numbers of 4500, 800, and 30 did NOT disrupt results really. This is quite surprising since we have strongly argued earlier that all ranges should start and end on IPOT numbers to avoid any bias towards low or high digits!

We have successfully aggregated and measured Continental USA, and one can claim that this data represent something 'more random' in a sense than data for just one single city. Now, what about other countries?! What we now seek is a global LD distribution of address data! And so naturally, the next higher order scheme is the averaging the averages of the averages, and few such calculations are enclosed here:

| | | | | | | |
|---|---|---|---|---|---|---|
| max UB Highest ending | 250 | 400 | 500 | 333 | 500 | 1,005 |
| max UB Highest starting | 100 | 55 | 250 | 55 | 50 | 995 |
| max UB Lowest | 30 | 50 | 60 | 17 | 50 | 30 |
| min UB | 1 | 1 | 30 | 1 | 1 | 15 |
| LB | 1 | 1 | 1 | 1 | 1 | 1 |
| ---------------------------------- | ------------- | ------------- | ------------- | ------------- | ------------- | ------------- |
| Digit 1 | 29.6% | 29.9% | 29.3% | 31.5% | 30.6% | 30.6% |
| Digit 2 | 18.4% | 17.7% | 16.4% | 17.8% | 17.6% | 19.2% |
| Digit 3 | 13.4% | 13.1% | 12.6% | 12.6% | 12.7% | 13.3% |
| Digit 4 | 10.2% | 10.3% | 10.1% | 9.7% | 9.9% | 9.8% |
| Digit 5 | 8.0% | 8.1% | 8.4% | 7.8% | 7.9% | 7.5% |
| Digit 6 | 6.5% | 6.6% | 7.1% | 6.4% | 6.6% | 6.1% |
| Digit 7 | 5.4% | 5.5% | 6.1% | 5.4% | 5.6% | 5.1% |
| Digit 8 | 4.6% | 4.7% | 5.3% | 4.7% | 4.8% | 4.4% |
| Digit 9 | 3.9% | 4.1% | 4.7% | 4.1% | 4.2% | 4.0% |

FIG 5: Average of Averages of Averages Schemes **(global)**

Results are even more logarithmic now – even when we strayed far away from the safety of unbiased IPOT borders! Again, one can claim that this global data represent something 'more random' in a sense than data for just one single country. Now, it is by the same logic that we went from city, to country and then globally that we find it necessary to continue to higher and higher order averaging schemes to be able to capture and model everyday numbers perfectly well, especially if we wish to incorporate other types of more complex data as oppose to merely simple count data of street addresses. As expected, the higher the order of the scheme the closer computerized results are to the logarithmic, and convergence is rapidly achieved!
Interestingly, as order increases, results become completely independent of the values for the limits chosen, even non-IPOT numbers when chosen do not result in any significant deviation from the logarithmic, namely: the scheme is totally independent of parameters!

B.J. Flehinger presented a rigorous proof in her 1966 article [F] that an iterated averaging scheme (exactly as done here) in the limit approaches the logarithmic. The above discussion may hopefully serve perhaps as the conceptual framework for going in such round about way in calculations and having the logarithmic then emerged, and that such a scheme could refer also to usage of numebrs and not merely to the number line itself.

One serious dichotomy though (overlooked earlier) exists between the simple averaging scheme and (city's) street address data. By taking the simple average of LD distributions of all the streets we allot each street equal importance or weight. But how could that be so if some are long and some are short?! There are many more houses on long streets than on short ones, and the postman surely must be using house numbers pertaining to long streets much more often than those pertaining to short ones?! Taking this fact into account would result in a milder LD distribution that is by far less skewed in favor of digits 1. For example, if a Greek-like Parable would represent street addresses for which all streets start at 1, the shortest is just one house long standing alone, and the longest street is of 9 houses, yielding 45 houses in total on 9 streets, (1+2+3+4+5+6+7+8+9 = 45), then direct calculations of LD on those 45 houses shows a linear fall in LD distribution as follow: (20%, 17.8%, 15.5%, 13.3%, 11.1%, 8.9%, 6.7%, 4.4%, 2.2%). There is one exception though, in the case where a selected group of people living in those houses are chosen (say famous people) and where one can explicitly assume, for one reason or another, that they all live on different streets, in which case equal weight for all streets is a valid assumption. For example, if we select for each of 200 cities say the most famous person there, one 'representative' famous person per city, then their street addresses are as in the averaging (of the averages) schemes above, as the assumption of equal weight of each street within the entire scheme is valid.
A straightforward way to avoid this dichotomy is to assume that short streets occur with much higher frequency than long ones, a reasonable assumption in a way. Certainly there are many more streets with say 30 or 40 houses than ones with 35,000 houses like 5th or 2nd Avenues in New York City. Yet it is very hard to argue that streets with a single house in them are the most numerous and come with the highest probability!? While the fact that the longest street, say Broadway in NYC, having perhaps 618,000 houses comes with the lowest probability is certainly acceptable - being uniquely the longest in the whole city, is surely the rarest. To eliminate this dichotomy for a Greek-like Parable representing street addresses, it would be necessary to assume that there are 9 as many streets with 1 house than those with 9 houses. That is:
Probability(street with 1 house) = 9 * Probability(street with 9 houses) , and also that:
Probability(street with 3 houses) = 9/3 * Probability(street with 9 houses) , and that:
Probability(street with 7 houses) = 9/7 * Probability(street with 9 houses) , and so forth.

The essential aspect of the simple averaging scheme is incorporating multiple intervals that are growing in focus and thus having **different** length each, yet all of them are fastened to the **same** much lower LB, typically 1. Calculations of schemes that are designed totally differently, and where intervals come with fixed length (so that $UB_i - LD_i$ = constant, for all of them) but which differ in their position (i.e. shifting each interval one integer forward to a higher location on the x-axis) lead to LD distribution that do not resemble the logarithmic distribution at all, no matter how we vary (in unison) those upper and lower bounds of fixed distances. Fundamental to Benford's Law is the fact that usage of numbers has one unique direction - up, thus any single set of numbers regarding a particular topic containing numbers on say ($6*10^{13}$, $7*10^{13}$) where digit 6 leads, would most likely also contain numbers on say ($5*10^{13}$, $6*10^{13}$) as well where digit 5 leads (assuming lower bound is somewhere near the origin) whereas the converse is not true. This example shows how the higher digit 6 frequently must share the spoils with lower digit 5, while the converse is not true. This argument can be generalized to any two competing adjacent digits, hence explaining in general why we do get this monotonically decreasing LD distribution.

*In other words: To get to high digits one must 'pass' through low ones, and this in essence is what deprives high digits of equality in leadership.*

## [6] Chains of Distributions

An alternate point of view of the simple averaging scheme (and ultimately that of Flehinger's iteration as well) is to consider those varying upper bounds as originating from a random process. Furthermore, each interval in itself is to be viewed as a random variable, not merely a set of integers. And so, instead of attempting to record digital leadership for a variety of intervals containing integers and having different upper bounds and then average them all, the corresponding view here is to take the LD pulse of random numbers from the continuous uniform on (0, B) standing for the generic single interval, and where B, the upper bound, is itself a random number drawn from another continuous uniform distribution on (0, C). We may call the latter 0 'the min upper bound' and the constant C 'the max upper bound', and view the latter uniform distribution as existing merely to provide us with a random parameter. Schematically this is written as U(0, U(0, C)). Here the continuous uniform is used, as oppose to the discrete uniform as in the averaging schemes where only integers were considered. Also conveniently, the ranges for the chains start from 0 as oppose to the usage of 1 being the LB in all the averaging schemes, yet we disregard and overlook this difference which is quite inconsequential when C is large enough, and thus we argue that these types of chains of distributions are essentially - at least in some sense - mirror images of the averaging schemes! [There is actually a **fundamental** difference between U(0, U(0, C)) and U(1, U(1, C)) regarding their LD whenever C is not large enough as compared with 0, 1 or any other LB, this is so because the range (0, 1) is quite apart from all other x-axis ranges regarding LD, and also in general, as it contains infinitely many IPOT numbers as well as an infinite log spread, but more on that in the section at the end detailing many more results regarding the chains.] Note that on [6, 7) for example digit 6 leads by as much as digit 8 leads on [8, 9), so the switch from averaging schemes of integers to simulations of chains of continuous variables does not alter LD relative distribution. Why do we choose here (twice) the uniform as oppose to other distributions? The answer lies in the need to correspond to the assumption we made earlier that integers come with equal weights within intervals, and that all intervals are equally considered in the entire scheme. Attempts here to perform these simulations resulted in various LD distributions depending on the particular value of C as might be expected. Note that here a choice of 200 say for C would yield the same result as a choice of 200,000. Also note that for C=100 (an IPOT number) results were very close to Stigler's law! Here are 9 such simulation results:

| Digit | C=100 | C=200 | C=300 | C=400 | C=500 | C=600 | C=700 | C=800 | C=900 | Average |
|---|---|---|---|---|---|---|---|---|---|---|
| 1 | 24.2% | 31.7% | 36.7% | 34.9% | 33.6% | 31.0% | 28.6% | 27.7% | 25.7% | **30.5%** |
| 2 | 18.2% | 12.1% | 16.1% | 20.4% | 20.6% | 21.0% | 20.8% | 20.2% | 19.1% | **18.7%** |
| 3 | 14.5% | 11.1% | 9.1% | 10.8% | 13.6% | 14.6% | 15.6% | 15.1% | 15.3% | **13.3%** |
| 4 | 11.9% | 9.6% | 7.6% | 7.2% | 7.7% | 10.4% | 10.6% | 11.2% | 11.6% | **9.8%** |
| 5 | 9.5% | 8.9% | 6.8% | 5.9% | 5.8% | 6.3% | 7.7% | 8.6% | 8.8% | **7.6%** |
| 6 | 7.7% | 7.8% | 6.7% | 5.8% | 4.9% | 4.5% | 5.2% | 6.5% | 7.5% | **6.3%** |
| 7 | 6.0% | 7.0% | 6.3% | 5.4% | 4.6% | 4.5% | 4.0% | 4.1% | 5.3% | **5.2%** |
| 8 | 4.6% | 6.3% | 5.6% | 4.7% | 4.7% | 3.9% | 4.0% | 3.4% | 3.6% | **4.5%** |
| 9 | 3.4% | 5.5% | 5.2% | 4.9% | 4.6% | 3.8% | 3.5% | 3.1% | 3.1% | **4.1%** |

FIG 6: Picking Numbers Randomly from a Chain of Two Uniform Distributions, U(0, U(0,C))

The next natural step is to pick up numbers randomly from the uniform on (0, B) where B itself is being picked randomly from the uniform on (0, C) and where C is randomly chosen from the uniform on (0,K). Schematically this is written as U(0, U(0, U(0, K))). LD of this simulation would mirror the average of averages scheme (country). Performing such simulation yielded better conformity with the logarithmic and the choice of the last constant K was less crucial to the resultant digital distribution than was the choice of C in the previous simulation where changes in C there set off larger swings in LD distribution. But clearly the value of K is still playing an important role here.

Results for the next logical step, namely U(0, U(0, U(0, U(0, M)))), a chain of 4 uniform distributions standing for the average of averages of averages scheme (global), came out very close to the logarithmic as can be seen in the table below, and the choice now of the constant M matters much less! Here are 7 such simulation results:

| Digit | M=100 | M=200 | M=300 | M=400 | M=500 | M=600 | M=700 | M=800 | M=900 | AVG | Benford |
|---|---|---|---|---|---|---|---|---|---|---|---|
| 1 | 31.3% | 29.9% | 29.5% | 29.9% | 29.6% | 29.8% | 30.9% | 31.0% | 31.1% | 30.3% | 30.1% |
| 2 | 17.5% | 18.1% | 17.5% | 16.9% | 17.7% | 17.1% | 17.7% | 17.8% | 18.0% | 17.6% | 17.6% |
| 3 | 12.0% | 13.0% | 13.2% | 12.3% | 11.7% | 12.6% | 11.8% | 12.3% | 12.4% | 12.4% | 12.5% |
| 4 | 9.5% | 9.9% | 9.7% | 10.2% | 9.6% | 9.9% | 9.4% | 9.1% | 9.0% | 9.6% | 9.7% |
| 5 | 7.8% | 7.9% | 8.0% | 8.1% | 7.8% | 8.2% | 7.2% | 7.5% | 7.6% | 7.8% | 7.9% |
| 6 | 6.2% | 6.3% | 6.8% | 6.9% | 6.8% | 6.8% | 6.8% | 6.7% | 6.7% | 6.7% | 6.7% |
| 7 | 5.9% | 5.8% | 5.7% | 6.0% | 6.5% | 5.9% | 5.9% | 5.7% | 5.9% | 5.9% | 5.8% |
| 8 | 5.1% | 4.8% | 5.0% | 5.4% | 5.4% | 5.2% | 5.6% | 5.4% | 5.1% | 5.2% | 5.1% |
| 9 | 4.7% | 4.3% | 4.6% | 4.3% | 4.9% | 4.5% | 4.7% | 4.6% | 4.4% | 4.5% | 4.6% |

FIG 7: Picking Numbers Randomly from a Chain of Four Uniform Distributions, U(0, U(0, U(0, U(0, M))))

Clearly, what we need here is an infinite chain to obtain the logarithmic, and that is exactly what Flehinger's scheme is really all about! Her Scheme is the mirror image of an infinite chain of the continuous uniform distributions, all anchoring at 0 **(but not at 1 though!)**, each depending sequentially on the next distribution to obtain the parameter. Note that with each higher order chain the actual value of the last constant becomes less important, and LD distributions become more similar regardless of the particular choice of that constant.

This result shows that the process of picking a number randomly from a truly random interval, one with an infinite successive arrangement of random upper bounds yields exactly the logarithmic. This very random nature we pick numbers here, not even specifying the range clearly, but rather keep referring its upper bound to yet another random process whose range is again being referred to another one and so forth, suggests the vague and mathematically undefined concept one may call 'super random number'.

The uniform U(a, b) has two parameters, a, the lower bound which we insist on fixing at 0 here, and b, the upper bound which we refer to another uniform distribution of the same type. Now, instead of employing just the uniform distribution, a more general viewpoint would be to select numbers randomly from <u>any</u> distribution with all its parameters being random numbers of other distributions (of the same/different type), whose parameters in turn are also derived from yet other distributions, and so forth. The conjecture here then is that this algorithm would yield the logarithmic distribution so long as there are enough sequences in the chain. The conjecture is stated for most distributions and parameter types, but there are exceptions. Furthermore, after numerous experimentations with several such simulations it appears that convergence is relatively quite rapid and after 3 or 4 sequences of distributions we come pretty close to the logarithmic. This more general process should be the one earning the term ‘super-random number’, while Flehinger’s uniform scheme is just one very special case of this more general algorithm. In a sense, we may think of the primary distribution (the first one, whose existence does not serve as parameter for another) as a distribution with some deep roots of uncertainly - one that is most thoroughly random, in a sense.

A notable demonstration of the very general nature of the chains is a simulation with a 4-sequence distribution chain performed as follows: a normal distribution is simulated where: **(A)** the mean itself is simulated from the uniform with a=0 and b=chi-sqr whose d.o.f. is derived by throwing a die, **(B)** the s.d. of the normal is also not fixed but rather simulated from the uniform(0,2). Note that here both parameters are being chained. Schematically this is written: **Normal(Uniform(0, chi-sqr(a die)), Uniform (0,2))**.

The table below shows the LD distribution results of six different simulation batches. Each batch of simulation represents the LD distribution of 2000 numbers being simulated.

| Digit | 1st Batch | 2nd Batch | 3rd Batch | 4th Batch | 5th Batch | 6th Batch | Average | Benford |
|---|---|---|---|---|---|---|---|---|
| 1 | 29.7% | 31.1% | 30.1% | 28.9% | 29.3% | 30.5% | 29.9% | 30.1% |
| 2 | 17.8% | 17.9% | 18.0% | 18.2% | 18.9% | 17.4% | 18.0% | 17.6% |
| 3 | 14.7% | 13.6% | 12.9% | 14.3% | 13.8% | 13.1% | 13.7% | 12.5% |
| 4 | 10.2% | 10.8% | 10.2% | 10.8% | 10.0% | 10.9% | 10.5% | 9.7% |
| 5 | 8.1% | 8.1% | 7.8% | 7.9% | 7.8% | 8.1% | 8.0% | 7.9% |
| 6 | 6.6% | 5.9% | 6.4% | 5.6% | 6.3% | 6.3% | 6.2% | 6.7% |
| 7 | 4.9% | 5.8% | 5.2% | 6.0% | 5.6% | 5.4% | 5.5% | 5.8% |
| 8 | 4.0% | 3.9% | 4.7% | 4.5% | 4.6% | 4.9% | 4.4% | 5.1% |
| 9 | 4.3% | 3.1% | 4.9% | 4.0% | 3.9% | 3.7% | 4.0% | 4.6% |

FIG 8: LD of Chain of Distributions, a Die, chi-sqr, Uniform and the Normal, 6 Runs

The agreement with the logarithmic distribution is quite good, even though there were only four sequences in the chain here.

**Initially as the idea of the chains was explored, the conjecture required ALL parameters be chained, leaving none stuck as a constant.** In a sense, we may think of the primary distribution (the first one, whose existence does not serve as a parameter for another) as **a distribution with some deep roots of uncertainly - one that is most thoroughly random, having such uncertain and shifting parameters.** For distributions with more than one parameter, the general image of the chain is of a pyramid-like arrangement of densities, with the primary sitting on top.

**Yet, this simplistic stipulation turned out to be wrong** as it depends on the exact nature of parameters/distributions in question. In some cases this stipulation is true, yet in others this is not the case as some parameters are simply totally or partially indifferent to chaining and do not contribute to any logarithmic convergence. The vast majority of all typical distributions though converge to the logarithmic if ALL parameters are blindly and mindlessly (regardless of types) fasten to others in chains. Only a few rare types/cases are totally immune to chaining. For a 2-parameter distribution say, there are typically 3 scenarios, (A) both need to be chained in order for convergence, (B) only one needs to be chained while the other one is totally irrelevant, (C) nothing works here and even if both are chained we do not see any convergence whatsoever. Hence, **the chain conjecture is not a universal rule but rather restricted to particular types of distributions/parameters. Very roughly, scale and location parameters are chain-able, while shape is not.** But this is way too simplistic, in reality the full account is quite complex. Fortunately, a rough conjecture regarding some general principles here was found. For a better exposition and flow here, the full account is given in a separate sub-section at the end.

A **second conjecture** regarding chains of distributions is made here (see extra section at the end) pertaining to the shortest chain possible having only two sequences as oppose to infinitely many (i.e. the other extreme in length). The 2nd conjecture claims that any 2-sequence chain in which ALL parameters of the primary density are derived from logarithmic distributions is logarithmic there and then without any need to expand infinitely. In other words, if distributions serving as parameters are logarithmic in their own right then the chain is also logarithmic. Numerous simulations confirmed the 2nd conjecture, and which appear true exactly for the very same types of parameters and distributions that the first (infinite) conjecture refers to.

Does the chain of distribution have anything to do with typical everyday data? In other words, could it represent some sort of an explanation for the logarithmic phenomena for everyday data? No! but in some limited sense it might be so partially for some very particular class of everyday data if we ponder the preponderance of causality in life, the multiple interconnectedness and dependencies of entities that we normally measure and record, and that so many measures are themselves parameters for other measures. For example, lengths and widths of rivers depend on average rainfall (being the parameter) and rainfall in turns depends on sunspots, prevailing winds, and geographical location, all serving as parameters of rainfall. Weights of people may depend on overall childhood nutrition, while nutrition in turns may depend on overall economic activity and so forth. The chain can indeed explain logarithmic behavior of some pieces of everyday data arising from sources that mimic the chain structure (and only approximately so when just a few-sequence-deep chain could serve as a model), yet only a small portion of the entirety of our everyday data could be covered here. High-order averaging schemes mentioned earlier though (for whom the chains were corresponding to) are a bit more readily viewed as a (failed) explanation of sorts as it attempts directly to measure (to average out) aggregate LD of some approximating model. Yet, the only concrete and full explanation for the phenomena is Ted Hill's proof regarding the LD of a uniquely defined distribution consisting of multiple distributions, a topic that would constitute our next section.

## [7] Hill's Super Distribution

After many decades since the re-discovery of the phenomena in 1938 by Frank Benford, and after numerous failed attempts at proofs, a rigorous mathematical explanation was given in 1995 by Theodore P. Hill, demonstrating that a distribution consisting of infinitely many distributions is logarithmic. Hill's argument then is that the totality of everyday data is simply a composition of a variety of variables and distributions, and that picking numbers from such real life data is schematically equivalent to picking from his abstract random number arising from a distribution of distributions. Since his abstract mathematical model allows such a wide variety of classes and forms of distributions (restricted to those defined on the positive x-axis), the perfect applicability to everyday real data is obvious.

One could view Hill's random number from his distribution of distributions as some kind of a super random number, just as was said about the infinite chain of distribution, yet there are two important differences. The first difference is that the statistician performing Hill scheme is thought of as a person standing at a center and facing 'horizontal' distributions, each independent of all the others, while for the chain of distributions scheme the statistician is thought of as standing at the very end (top) of a long 'vertical' chain of distributions, all aligned sequentially, all depending on one another. The second difference is that for the chain of distributions scheme, the final outcome is simply a collection of distributions of the same form (that of the top disrtibution not serving as a parameter for any other distribution), and differ only with respect to their parameters, while Hill's scheme allows a diversity of forms as well. It is worth noting that Hill's super distribution is parameter-independent, or parameter-less, just as the infinite order averaging schemes and the chains of distirbutions are.

Putting Hill' super distribution under the microscope, one sees numerous 2-sequence chains hiding there, stubbornly avoiding the limelight, moreover, actually that is **all** one could see there, chains exclusively! why? Let's take the exponential distribution for example, one minor item on his overall super distribution agenda; exp(56.7) is just as likely to be part of Hills's scheme as exp(0.0139) say, and so forth. And exponentials on [0, 1] say should are as likely to be included as exponentials on [1, 2], so what we really have here on hands is simply the 2-sequence chain **exponential(unifrom(0, +oo)).**

## [8] The Scale Invariance Principle

A novel approach to leading digits is one based on the scale invariance principle and on the base invariance principle. Scale and base invariance arguments claim that if there is indeed any universal law for significant digits, then it should be independent of the units and scales as well as the base system employed by society, because base and scales are cultural, arbitrary, and do not represent any fundamental properties of numbers. In other words, if Benford's Law were to be dependent on humanity using the meter as oppose to the foot, then the law would not be universal. The difficulty with such an argument is that Benford's Law was indeed about every day and scientific data as we have it in our society, not about our number system or pure numbers, and that the law could in principle depend on scale and base. Such an argument relies on the unproven assumption that some invariant first digit law exists in the first place.

In any case, the scale invariance principle, **if assumed** implies the logarithmic distribution. Roger Pinkham [P] have demonstrated that the logarithmic distribution is scale-invariant, and that it's the **only** distribution with such a property; therefore any first digit law that is independent of choices of scale must be the logarithmic!
Base invariance principle has also been used to prove the logarithmic distribution.

Derivation from scale invariance principle has been criticized for making unfounded assumptions; nonetheless let us justify the scale invariance principle based on some very general consideration about the nature of typical usage of numbers in everyday data and the implied leading digits distribution.

Units or scales are merely arbitrary measuring physical rods against which we compare other entities. We define here the term **'measured category'** as some very specific physical entity of a particular category that needs measuring, such as heights of people, heights of horses, heights of buildings, lengths of rivers, distances between global cities, ages of people or lifespan of a single bacteria (these examples constitute 7 different measured categories in all). Units or scales are used in quoting measured categories.
Let us assume for a moment that society has standardized all measurements and does not tolerate using different scales/units for different measured categories of the same dimensions, that only one single universal scale/unit is used to measure all measured categories regarding any one particular dimension.

On the face of it, the scale invariance principle appears totally unfounded. For example, people's height measured in feet yields mostly 5 and 6 as leading digits while on the meter scale digit 1 takes a strong lead. But the rationale for the principle rests on the assumption that our measured categories are not just numerous, but enormous. It also rests on the notion that changes in scale will have such varying and independent effects on measured categories in terms of leading digits that it will all sum up to nothing. In other words, even though a change in scale would revolutionize digital leadership for almost each and every measured category (individually) yet, the changes for almost all measured categories do take totally different turns, canceling each other's effects and leaving the net results on digital leadership unaffected. Had we used only 4 measured categories say for the length dimension, then a change in scale would certainly not thoroughly shuffle all numbers in such a way as to leave digital leadership unaltered. Intuitively, any single shuffle of scale that would result in digit 7 for example taking a stronger lead all of a sudden in xyz measured category would be offset somewhere else (another measured category) where digit 7 would be all at once deprived of some existing leadership, balancing things out. All this is relevant because Benford's Law applies to the totality of everyday data, not to any subsets thereof or merely one measured category. For a more specific example, consider two measured categories, heights of people in feet where most quotes are assumed to be on (4.5 ft, 6.7 ft) and heights of building in a certain region where most quotes are on (65 ft, 193 ft). Now, changing our scale to meters using the conversion formula (FEET measurements)*0.30919 = (METER measurements) would give roughly the new intervals (1.39m, 2.07m) and (20.1m, 59.7m). Clearly digit 1 gained a lot of leadership with people's height but lost badly with the buildings! This trade-off is due to measured categories being something real/physical, and that the two intervals in the example above relate to each other in a fixed manner independent of the scale chosen. The fact that LD distribution of the totality of all our measured categories remains unaffected under a unified scale change is fairly intuitive, yet a rigorous proof of this grand trade-off should be given.

It is quite intuitive that under the assumption of uniform scale for all measured categories, any possible 'digital leadership bias' about the choice of scale during earlier epochs when it was decided upon, would be totally ineffective due to the enormous size of the number of measured categories now in existence. In other words, under the assumption, society could never find one magical scale that would give leadership preferences to some digits at the expense of others.

We could relax the above assumption of uniform scale for all measured categories, and even go to the other extreme by imagining a society that is extremely fond of multiple units, and in addition tends to assign units in haphazard and random fashion without any particular calibration with regard to leading digits or any other sensible manner, a society that liberally invents different scales for each and every measured category. Given the assumption of the enormity of number of measured categories, even for such idiosyncratic and inefficient society, the invariance principle still holds. That reshuffling of all their units in one single epochal instant (either uniformly by a fixed factor for all measured categories, or by different factors assigned randomly for each measured category) would not result in any change in overall leading digits distribution. Such a society also passes the scale invariance test.

But let us considering a third case, of another society, one that is also very fond of multiple units, and liberally invents different scales for each and every measured category, but doesn't assign units in haphazard and random fashion; rather it tends to carefully calibrate things in such a way leaving them to quote numbers with mostly digit 1 and 2 leading. For this society, the scale invariance principle has no validity, because in such a society, reshuffling all scales simultaneously would drastically alter leading digits distribution (in favor of the logarithmic if done randomly and differently for all measured categories).

And what about our society during this particular epoch? Surely our society does limit the diversity of scales a great deal. We don't liberally invent different scales for each and every measured category, yet we do not restrict ourselves completely. It is also true that we do not usually assign units completely haphazardly and randomly when inventing new ones. There are numerous examples of convenient yet peculiar and drastic adjustments of scale, such as astronomical units of length being light years, Avogadro's number in chemistry, and many others. A typical example is the scale for some measured categories regarding weights: tons are used at international freight companies, milligrams in bio labs, and kilos at butcher shops, with the overall result that many quotes of weight are with low-digit-led numbers (1.34 ton of crude oil delivered to HK, 3.5 milligram of $C_6H_3O_5$ solution, 1.5 kilo beefsteak.) It might even be argued that people unknowingly like dealing with low-digit-led numbers, and that they at times invent scales that promote such outcome. In any case, even if that is true in some instances, the fact that we severely limit multiple usages of units as compared with our numerous number of measured categories implies that no matter how we calibrate scale for this or that particular measured category, inevitably numerous other measured categories (all united in paying homage to one single scale) will have the final word in overall leading digits distribution and will swamp them by their sheer size.

Clearly, a change in scale for one dimension implies the transformation of numerous measured categories by different multiplication factors. For example, a change in the length dimension that reduces standard unit to half its original size would require the factor of 2 for simple lengths, 4 for areas, 8 for volume, 1/2 for kilograms per length and 1/8 for minutes per volume, and so forth. Similar changes in the units of mass and time dimensions would result in different factors. Nonetheless, the scale invariance principle can also be invoked just the same to state that one single unified transformation of all measured categories (of the totality of our data) by ANY SINGLE multiplicative factor would not alter existing logarithmic distribution whatsoever! **As a corollary, transforming logarithmically well-behaved and large enough pieces of data by multiplying them by ANY single factor would barely nudge their LD distributions!** Readers are encouraged to attempt to perform this rather striking demonstration of the scale invariance principle!

It is noted that certain data types have got nothing to do with scales. Example of such entities are periods, houses, people, and so forth, with a particular, non-optional unit of measurement, where the scale invariance principle is totally irrelevant. Yet, it could be argued that such count data, as large as it is, is being overwhelmed by the sheer volume of other scale-based data preserving and guaranteeing the overall correctness of the scale invariance principle.

## [9] The Logarithmic as Repeated Multiplications

Completely different aspect of Benford's law was given by Hamming [H-b], Raimi [R-76], and others regarding sequences of multiplications and divisions. Repeated multiplication processes effectively drive numbers toward the logarithmic distribution in the limit. Note the deterministic nature of this process, as oppose to the random.

Consider the arbitrarily chosen value of 8 being repeatedly multiplied by another arbitrary number 3. This yields the follwing table:

| $8*3^i$ | $8*3^i$ series | 1st digit |
|---|---|---|
| $8*3^0$ | 8 | 8 |
| $8*3^1$ | 24 | 2 |
| $8*3^2$ | 72 | 7 |
| $8*3^3$ | 216 | 2 |
| $8*3^4$ | 648 | 6 |
| $8*3^5$ | 1944 | 1 |
| $8*3^6$ | 5832 | 5 |
| $8*3^7$ | 17496 | 1 |
| $8*3^8$ | 52488 | 5 |
| $8*3^9$ | 157464 | 1 |
| $8*3^{10}$ | 472392 | 4 |
| $8*3^{11}$ | 1417176 | 1 |
| $8*3^{12}$ | 4251528 | 4 |
| $8*3^{13}$ | 12754584 | 1 |

| Digit | Proportion |
|---|---|
| 1 | 35.7% |
| 2 | 14.3% |
| 3 | 0.0% |
| 4 | 14.3% |
| 5 | 14.3% |
| 6 | 7.1% |
| 7 | 7.1% |
| 8 | 7.1% |
| 9 | 0.0% |

FIG 9 : Simple Mutliplication Process, $8*3^i$ and its LD Distribution

Other arbitrarily chosen values give similar results especially if extended much farther than the mere 14 sequences in the above table, so that the above result is quite general and representative, and in the limit the logarithmic emerges, except in some very rare anomalous cases discussed in a separate section towards the end of this article.

In any case, geometric series and exponential growth exhibit nearly perfect logarithmic behavior right from the start, and that can be said for almost any growth, any factor, and any starting value, all yielding the logarithmic almost exactly given that enough elements are being considered. It is important to note that geometric series and exponential growth are nothing but repeated multiplication processes albeit ones with a single (constant) factor applied.

Such logarithmic behavior is striking at first. Here geometric series and exponential growths keep pacing forward, with no end in sight, leaving a long trail of logarithmically-well-behaved numbers from the very start and exhibit such behavior for any interval cut from the whole series at any two arbitrary points (assuming they are sufficiently apart). Somehow as if by magic the series jump forward using carefully selected choreographic steps that not only favor overall low-digit-led numbers, but also does it in a way as to yield the logarithmic exactly!

Examination of related LOG and Mantissa distributions of multiplication processes leads to a much better understanding of their logarithmic behavior, and will be thoroughly discussed in later sections.

## [10] The Case of k/x Distribution

The probability density function f(x) = k/x occupies a central position and importance in the study of Leading Digits and Benford's Law, and studying all aspects of this distribution is essential for a complete understanding of the phenomena.

For the probability density function of the form f(x) = k/x over the interval $[10^S, 10^{S+G}]$ where S is **any** real number and G is **any** positive integer, k is a constant depending on the values of S and G, the sum of all the areas under the curve where digit d is leading is indeed $LOG_{10}(1+1/d)$, namely perfectly logarithmic. If G is fairly large (say over 25 perhaps, but depending on precision required), then the requirement that G must be an integer can be relaxed. Note that the distribution considered here has its long tail to the right abruptly cut at $10^{S+G}$, its head abruptly launched at $10^S$, and that exponents of 10 representing the two boundaries of the interval differ exactly by an integer, namely G.

[We shall use the notation ln(q) or ln q to mean the natural logarithm base e of the number q.]

Let us prove this assertion. We first note that the entire area should sum to one, that is $\int k/x\, dx = 1$ over $[10^S, 10^{S+G}]$, therefore $k[\ln(10^{S+G}) - \ln(10^S)] = 1$, or $k[(S+G)\ln 10 - (S)\ln 10] = 1$, so that $k*\ln 10*[(S+G) - (S)]=1$ and $k*\ln 10*G = 1$, so $k = 1/[G*\ln 10]$. Notice, that this determination of k was in total generality, where G can assume any value, not necessarily only an integer, and that G represents the difference in the base 10 exponents of the two boundaries spanning the entire interval length of x in question.

Secondly, given a particular p.d.f. of the form k/x on (a, b), we note that for any two **sub**intervals of (a, b) having the same exponent difference, their areas under the curve are identical. Given $[10^P, 10^{P+I}]$ and $[10^Q, 10^{Q+I}]$, both contained inside (a, b), P and Q being <u>any</u> set of numbers, not necessarily integers, we know that the values of their related constant k are identical since they belong to the same distribution defined on (a, b). The areas under the curve are $k[\ln(10^{P+I}) - \ln(10^{P})]$ and $k[\ln(10^{Q+I}) - \ln(10^{Q})]$ respectively, or simply k[(P+I)ln10 – (P)ln10] and k[(Q+I)ln10 – (Q)ln10], simplifying we get: k*ln(10)*[(P+I) – (P)] and k*ln(10)*[(Q+I) – (Q)], which yields k*ln10*I as the **same** area for each subinterval! If the whole interval (a, b) is expressed as $[10^{S},10^{S+G}]$ so that k = 1/[G*ln10] then area for each is simply I/G, namely the ratios of exponent difference of the subinterval to the exponent difference of the entire range. For example, for k/x defined on (1, 10000), [1, 10] and [100, 1000] have equal areas, [1, 10] is narrower on the x-axis but with a very high p.d.f. value, while [100, 1000] is extremely long on the x-axis in comparison, but it comes with a much shorter p.d.f. value. This trade-off cancels out each effect exactly so that areas end up the same.

To set the stage for the proof, we represent the interval in question here as $\mathbf{[10^{N+f}, 10^{N+f+G}]}$, where N is zero or some positive integer, f is some possible fractional part, and G is a positive integer representing the integral part of the difference in exponents for the interval.

We initially let f = 0, a restriction that would be relaxed later. Considering any digit D, the probability that D leads is given by the area under f(x)=k/x on the following intervals: $\{[D10^{N}, (D+1)10^{N}],\ [D10^{N+1}, (D+1)10^{N+1}],$ ...**G times**...$[D*10^{(N+G-1)},(D+1)*10^{(N+G-1)}]\}$. This is so because it is on these intervals and these intervals alone that D leads. Calculating the various definite integrals we obtain:

$k[\ln((D+1)10^{N}) - \ln(D10^{N})] + k[\ln((D+1)10^{N+1}) - \ln(D10^{N+1})]$
+....**G times** ....+ $k[\ln((D+1)10^{N+G-1}) - \ln(D10^{N+G-1})]$ or:

k[ln(D+1)+N*ln(10)-ln(D)-N*ln(10) ] + k[ln(D+1)+(N+1)*ln(10)-ln(D)-(N+1)*ln(10)]
+ ....**G times**.... + k[ln(D+1)+(N+G-1)*ln(10)-ln(D)-(N+G-1)*ln(10)]

Canceling out like terms (not involving D) we are left with:

k[ln(D+1)-ln(D)] + k[ln(D+1)-ln(D)]+ ...**(G times)** ... + k[ln(D+1)-ln(D)]
k[ln[(D+1)/(D)]] + k[ln[(D+1)/(D)]]+ ...**(G times)** ... + k[ln[(D+1)/(D)]]

that is: G*k*ln[(D+1)/(D)]. Substituting here the expression for k above we obtain: G*(1/[G*ln10])*ln[(D+1)/(D)] = ln[(D+1)/D]/ln10. This expression uses the natural logarithm base e, and applying the logarithmic identity $LOG_{A}X=LOG_{B}X / LOG_{B}A$ to convert this ratio to the common base 10 we get:

$LOG_{10}[(D+1)/D]/LOG_{10}[e]$ / $LOG_{10}[10]/ LOG_{10}[e]$
$LOG_{10}[(D+1)/D]/LOG_{10}[e]$ / $1/ LOG_{10}[e]$
$LOG_{10}[(D+1)/D]$
$LOG_{10}[1+1/D]$!

Let us now see why $f \neq 0$ should not yield any different result. As shown above k depends solely on the difference between the exponents being that its expression is $k = 1/[G*\ln 10]$, hence since $f \neq 0$ ( a slight shift to the right from the $f=0$ situation) does not alter that difference here, and k is still of the same value (whether f equals to zero or not), and thus function x/k is not altered either. Now, since areas under the curve are identical for any two subintervals of the same exponent difference it follows that any nonzero f that requires an additional area to the **right** of the intervals
$\{[D10^{N}, (D+1)10^{N}],\ [D10^{N+1}, (D+1)10^{N+1}],\ \ldots,\ [D*10^{(N+G-1)}, (D+1)*10^{(N+G-1)}]\}$.
would then also require an **<u>identical</u>** subtraction on the **left** side of that intervals! This completes the proof.

An integral difference, G, in the exponents of the interval $[10^{S}, 10^{S+G}]$ was required for logarithmic behavior, but if G is fairly large (depending on precision required), then the requirement can be relaxed without effecting much the above proof. This is so because the value of the constant k is $1/[G*\ln 10]$ and as such it is inversely proportional to G, hence each term in the earlier expression
$k[\ln[(D+1)/(D)]] + k[\ln[(D+1)/(D)]] + \ldots$ **(G times) ...** $+\ k[\ln[(D+1)/(D)]]$
becomes quite marginal in the overall scheme of things for large values of G (and the implied small values of k).

The conceptual assertion that k/x on (a, b) is always or continuously logarithmic, can be justified or illustrated, if we slide to the right or to the left a flexible imaginary lens, focusing solely on the subinterval $\mathbf{[10^{S}, 10^{S+1.00}]}$ contained anywhere inside the entire (a, b) interval, by letting S vary. Here exponent difference is always 1 no matter what value S assumes. The width of the lens (on the x-axis) widens as we slide to the right (but height falls) and shrinks as we slide to the left (but height rises), while data generated from it in accordance with local density is always exactly and perfectly logarithmic!

## [11] Uniform Mantissa, Varied Significand, and the General Law

In the first paper ever acknowledging the phenomena, Simon Newcomb asserted that “the law of probability of the occurrence of numbers is such that all mantissae of their logarithms are equally probable” [N]. The definition of the mantissa of any number X is the unique solution to $\mathbf{|X| = [10^{W}]*[10^{mantissa}]}$, where W - called the ‘characteristic’ - is the result of rounding **down** log(|X|) to the nearest integer, namely the largest integer less than or equal to log(|X|), or equivalently the first integer to the left on the real number line $(-\infty, +\infty)$ of the log-axis regardless whether log(|X|) is positive or negative.

Similarly, significand of X is defined as the unique solution to $\mathbf{|X| = [10^{W}]*[S]}$, where W is the result of rounding down log(|X|) to the nearest integer (the ‘characteristic’). In other words, significand is simply the way the number is written and expressed (using digits) disregarding the decimal point. For example, if X = 175.831, its significand is 1.75834. Hence significand is that first part of scientific notation with the decimal power totally ignored. Note that significand □ [1, 10).

From the two definitions above it follows that the relationships between mantissa and significand are:

**significand = $10^{\text{mantissa}}$** & **mantissa = $\log_{10}$( significand )**

Following the curve of the mantissa and significand (as a function of x) along their journeys throughout the x-axis, we notice a perfect repetition after each IPOT number. The two equations below (which follows directly from the definitions) as well as the following graph for the significand defined over (10, 1000) illustrate this principle:

**mantissa of ($10^N$*X) = mantissa of (X) [N being an integer]**
**significand of ($10^N$*X) = significand of (X) [N being an integer]**

We have now the tools to properly define Benford'a Law in its general form, detailing exactly how digits appear in total generality (!), that is: how higher order digits are distributed and the relationship and dependencies between the orders themselves.

**The general form of Benford's Law:**
**Probability( significand ≤ $S_0$ ) = $LOG_{10}(S_0)$ where $S_0$ is an element of [1,10).**

It is truly exhilarating to see the beauty and elegance in having such a concise and compact law covering the entire phenomena!

**The general form of Benford's Law implies $LOG_{10}(1+1/d)$.**
**The general form of Benford's Law implies uniformity of mantissa.**

## [12] Related Log and Uniqueness of k/x Distribution

Let us now turn our attention to the most obvious way of obtaining uniformity of mantissa, namely the special case of uniformity of log of data itself with its range on the log(x) axis measuring exactly some integral value.

Assume X has p.d.f. f(x)=k/x over the interval [$10^S$, $10^G$], with both S and G being **any** number, large, small, integer or non-integer. The Distribution Function Technique can be employed to find the p.d.f. of the transformation Y=$LOG_{10}$X. The Distribution function of Y, F(y), is then P[Y<y, for S≤y≤G] = P[$LOG_{10}$X<y, for S≤y≤G]. Taking 10 to the power of both sides of the inequality we get: F(y) = P[$10^{LOG\,10\,X}$<$10^y$, for S≤y≤G], or F(y) = P[X<$10^y$, for S≤y≤G] = ∫k/x dx over the interval ($10^S$,$10^y$), so that F(y) = k[$\log_e(10^y)$-$\log_e(10^S)$] = k[(y-S)*$\log_e$10]. Since k = 1/(I* $\log_e$10) where I represents the difference in the exponents of the two boundaries spanning the entire interval length of x, as was demonstrated earlier regarding k/x, we finally obtain: F(y) = [1/([G-S]* $\log_e$10)]*[(y-S)*$\log_e$10] = (y-S)/(G-S). Differentiating with respect to y and employing the relationship f(y)=dF(y)/dy, we finally get: p.d.f.(Y) = 1/(G-S), namely a constant probability density function of the form 1/(exponent difference). Hence:

**<u>Proposition I</u>: If X has p.d.f. f(x)=k/x over the interval [$10^S$, $10^G$], S and G being any number, then its related log distribution, namely Y=$LOG_{10}$X over [S, G] is uniformly distributed with p.d.f.= 1/(G-S). (No consideration is given here to Leading Digits).**

**<u>Corollary I</u>: If X has p.d.f. f(x)=k/x over the interval [$10^S$, $10^{S+N}$] where N – the difference in exponents – is exactly an integer, then X is logarithmic according to the General Law.**

This follows directly from Proposition I implying that related log of X is uniformly distributed while difference in exponents is assumed to be of an integral value, together rendering mantissa of X uniformly distributed. As an example, imagine log of X being uniformly distributed on [3.88, 6.88], with exponent difference of 3, then clearly mantissa is uniformly distributed as well [just exchange the margin on the left before the integral value of 4 on (3.88, 4) with right margin filling in (6.88, 7) space to obtain the wholesome interval (4, 7) with its uniformity intact.]

In the earlier section regarding the case of k/x distribution, it was demonstrated that if the interval over which k/x is defined is of an integral exponent difference than its first leading digits are logarithmically distributed as in $LOG_{10}$[1+1/d]. Here we go a step farther asserting that an integral exponent difference in its interval implies the logarithmic according to the General Law with all higher orders considered.

Assume Y is uniformly distributed over [R, S], with the difference S-R being **any** number, having p.d.f. 1/(S-R). The Distribution Function Technique can be employed to find p.d.f. of the transformation X=$10^Y$ over the interval [$\mathbf{10^R}$, $\mathbf{10^S}$]. The Distribution function of X, F(x), is then P[X<x, for $10^R \le x \le 10^S$] = P[$10^Y$<x, for $10^R \le x \le 10^S$]. Taking log of 10 to both sides of the inequality we get: F(x) =
= P[$LOG_{10}(10^Y) < LOG_{10}\,x$, for $10^R \le x \le 10^S$], or F(x) =
P[Y< $LOG_{10}\,x$, for $10^R \le x \le 10^S$] = ∫1/(S-R) dy over the interval [R, $LOG_{10}\,x$],
so that F(x) = [($LOG_{10}\,x$) - R]*[1/(S-R)]. Differentiating with respect to x and employing the relationship f(x) = dF(x)/dx, we finally get:
p.d.f.(X) = [1/$LOG_e$10]*[1/(S-R)]*[1/x] distributed over [$\mathbf{10^R}$, $\mathbf{10^S}$]. Hence:

**<u>Proposition II</u>: If Y is uniformly distributed over [R, S], with the length S-R being any number, then X=$10^Y$ over [$10^R$, $10^S$] is distributed with p.d.f. of the form f(x)=(1/[(S-R)* $LOG_e$10])*(1/x). (No consideration is given here to Leading Digits).**

**<u>Corollary II</u>: If related log of X (called Y) is uniformly distributed over [R, S] with the length S-R (exponent difference) being an integer, then X itself (the transformation X=$10^Y$) is logarithmic according to the General Law if considered over its entire range of [$10^R$, $10^S$]. Moreover, the manner in which X is logarithmic is of the type k/x, where logarithmic behavior is consistent and continuous, being observed also on ANY subinterval having an integral exponent difference anywhere in [$10^R$, $10^S$].**

This follows straight from Proposition II and the association of uniformity of mantissa with the General Law itself, as well as from what was observed earlier about k/x distribution [and employing the fact that uniformity of log over an integral interval implies uniformity of mantissa as well].

**It is very important here to note that the integral restriction on exponent difference of Corollary I and Corollary II becomes superfluous and can be waived if that exponent difference is quite larger, depending on precision desired. Very roughly, a value over 25 or so yields LD distribution that is quite close to the logarithmic even though it is not an integral value.** The reason for this waiver is straightforward: as the entire related log interval is enlarged, that extra piece of data stemming from the fractional part of log over and above an integral interval becomes exceedingly just a tiny fraction of the overall data (in comparison with that much larger part of data stemming from the main integral interval - a large proportion of data which in itself is perfectly logarithmic.) Yes, this extra piece ruins the chance for any perfect logarithmic behavior to be observed, although only slightly so.

If data or distribution is strictly between any two adjacent/consecutive integral powers of ten (such as 10&100, 1&10, 100&1000, etc.) then the ONLY distribution consistent with The General Benford's Law (all higher orders considered) is k/x! No other distribution could satisfy the general logarithmic condition, hence its uniqueness in LD. In fact over such an interval it can be said that k/x is the definition of being logarithmic in the general sense. to prove this, note first that if data is between TAIPOT then its related log is then stuck narrowly between two consecutive integers. Now, for logarithmic behavior to exist uniformity of mantissa is required, thus it follows that related log itself is also uniformly distributed here because log only varies on the short length of unity on log(x)-axis (within such a 'narrow range' of unity length only uniformity of log could ever let mantissa inherent uniformity from the log). But uniformity of related log density in turns implies data is distributed in the form of k/x as per Proposition II. Moreover, same uniqueness of k/x applies to any data or distribution between two values whose exponent difference is unity, that is whenever interval is on $[10^S, 10^{S+1}]$ S being any number. For example, any distribution on (3.47, 34.7), namely on $(10^{0.54033}, 10^{1.54033})$, is logarithmic if and only if it's distributed as k/x. and what about distributions between two non-adjacent integral powers of ten such as (1, 1000) for example? There, no distribution can claim uniqueness, as would be shown in the later section on the lognormal, exponential and k/x. There is an infinite number of distributions in theory all perfectly logarithmic in the general sense on such 'wider interval', and where k/x is but one such manifestation of a density behaving logarithmically.

Even though k/x density is always diminishing along the positive x-axis of course, having a conspicuous tail to the right, yet it goes on and on, not terminating until its upper limit of defined range is reached, and there, an abrupt and decisive end to concentration occurs. In fact, concentration of k/x is constant between **any** TAIPOT points! That means that there are as many numbers on (10, 100) as on (100, 1000)! [given that k/x is defined on say (1, 10000) which incorporates all of these sub-intervals between TAIPOT.]

The same constant concentration is found between any two points whose exponent difference is unity, such as for example (3.47, 34.7), namely ($10^{\mathbf{0.54033}}$, $10^{\mathbf{1.54033}}$). In fact related log density of k/x is constant, never diminishing, as it is uniformly distributed! By contrast, real life data always terminate no matter what, and do so gradually, not abruptly! Weights of people terminate gradually somewhere around 1.9 to 2.5 meters, and so forth. Even more so for the aggregated whole of the variety of real life usages and measurements, pointing to some smooth and gradual decline in density towards the truly upper limits of everyday data quoted in trillions or more.

## [13] Related Log Conjecture

Little reflection is needed to realize that uniformity of log of data over an integral interval (such as in the k/x case) is not the only circumstances yielding uniformity of mantissa. In situations where density curve of log of data is continuous and has wide enough range on the log(x) axis (roughly over 3 or 4 say), and when it rises from the log(x)-axis itself until it reaches a certain plateau or zenith, followed by a fall all the way back to it, it is conjectured that uniformity of mantissa may also be achieved, especially when log(x) is symmetric and where the rise and fall are not too steep but rather gradual.

In such cases, the log function monotonically increases up to the central point (meaning that **high** digits are benefiting somewhat as compared with the logarithmic situation) and then steadily falls off from there (meaning that now **low** digits get an advantage over and above the logarithmic situation), and it is easy to envision overall mantissa ending up uniform assuming that log traverses plenty of log(x) distance. It is plausible to argue that whatever low digits lose (in relation to the logarithmic) on the left of the highest point (rising) is exactly what they gain to the right of it (falling), and regardless of location of the center! On the other hand, for symmetrical distributions representing log of data that are too narrowly focused on the log(x) axis, there is no such meaningful trade-off. For example, for related log strictly bordered by 7.904 and 8.000 (hence generating mantissa on (0.904, 1.000) exclusively), the resultant LD distribution of the data itself necessarily consists of only 8 and 9 digits, and regardless of the shape of its related log curve!

Next to be depicted here is a chart showing hypothetical density of log of some data or distribution. Density curve starts from the very bottom, touching it, and ends all the way down there as well. On the extreme left **high** digits in a rare show of force are winning slightly, then digital equality in achieved exactly only at a single point if curve is truly smooth, and it is approximated here to be as such on a wide interval. Further on, low digits win slightly, followed by exact logarithmic behavior, and then low digits win even more than that, and so forth. Gradually we end up on the extreme right side where low digits almost completely dominate all other digits, much more than their normal logarithmic advantage. Intuitively all this suggests some grand trade-off between the left region and the right region, leaving the logarithmic center as a good representative of the <u>entire</u> range.

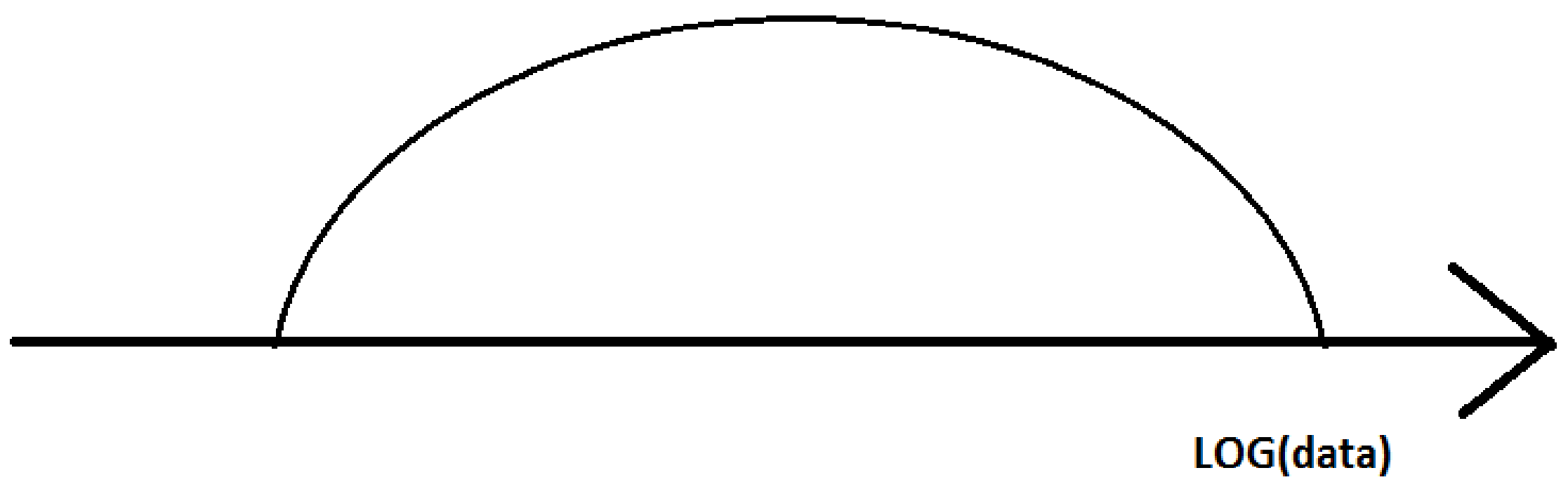

FIG 10: Symmetrical Related Log of Data or Distribution

It should be noted that overall mantissa here is the result of numerous competitions between the digits on each of those sub-intervals between integral values on the log axis, and thus the conjecture here then argues that by symmetry aggregating individual densities of all mantissa subintervals between the integers (in the chart above) points to a uniform overall distribution, namely: uniform mantissa, and thus logarithmic.

Let us examine carefully via simulations the conjecture above.

Three obvious examples for such symmetrical log distributions come to mind: (1) the normal distribution, (2) the triangular, and (3) the semi-circular.

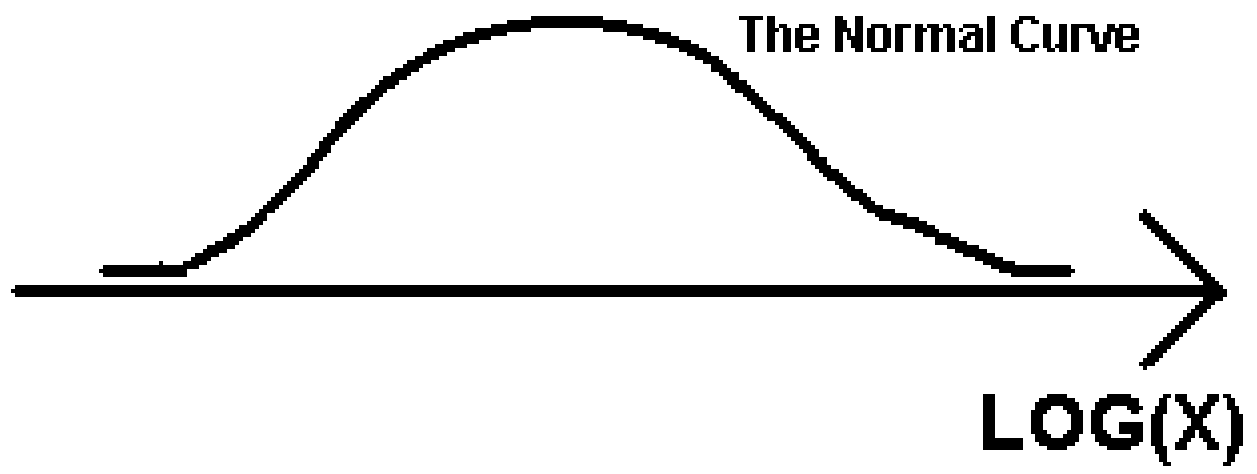

FIG 11: The Normal Curve as the Distribution of Related Log of Data

Simulations show that if related log (base 10) of data is normally distributed then data itself (that is $10^{\text{Normal}}$) is nicely logarithmic as long as s.d. of the normal is (roughly) 0.4 or larger, and regardless of what value the mean takes! Why? Because a large value for the s.d. here guarantees wide enough spread over the log(x)-axis, while too small a value leaves distribution too narrowly focused on some very small interval. By Chebyshev's theorem, 8/9 of distribution is within 3 s.d. Applying this here to the minimum s.d. value of 0.4 above encountered in simulations, we then have roughly 89% of distribution over a range of at least 2.4 units [0.4*(3+3)], so perhaps we might wish to extrapolate this and propose that <u>in general</u> what we need here is a minimum range of 2.5 to 3 log(x) units for sufficient agreement with the logarithmic. This particular arrangement of distribution outlined above is actually something quite familiar in statistics, since if $X = 10^{\text{Normal}}$ then log(X) is normally distributed, and hence log of any other base would also yield the normal, therefore $e^{\text{Normal}}$ would also be normal, But $e^{\text{Normal}}$ is that well-known lognormal distribution, which is logarithmic given parameters are of certain values, a topic that will be further explored in a later section.

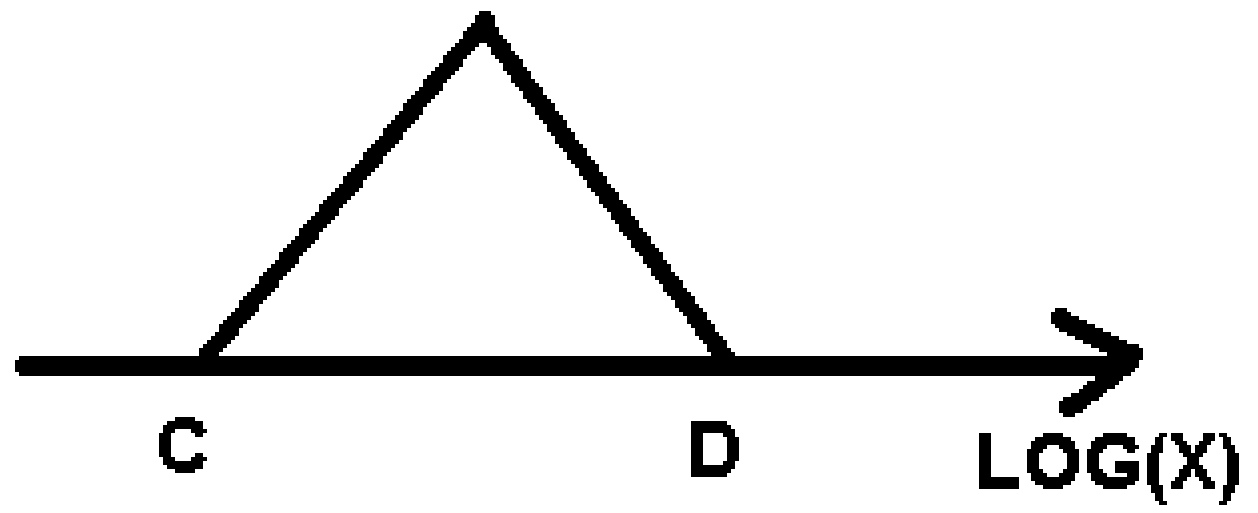

FIG 12: Triangular Distribution of Related Log of Data

Simulations of the triangular distribution (as in the figure above) use the following programming command to obtain a random LOG(X) value:
IF rand<0.5 THEN: C+SQRT((rand/2)*((D-C)$^2$)) ,
ELSE: D-SQRT(((1-rand)/2)*((D-C)$^2$)) ),
where D and C are the upper and lower bounds respectively, rand is a simulated number from the uniform on (0,1) representing a particular value for the cumulative distribution function, and SQRT is the square root. Note that results are independent of the height h and so it is omitted. These simulations show a near perfect logarithmic behavior for $\mathbf{10}^{\mathbf{Triangular}}$ whenever D-C > 2. In other words, for those log distributions whose ranges are of length 2 or higher. It is interesting to note that logarithmic behavior here for D-C > 2 holds regardless of the location of the center almost, that is, results are almost independent of the value of (D+C)/2, and depends almost solely on the value of (D-C), that is, the width. Although for rangers narrower than 2, fractional translation of location does matter. If center is exactly at an integer and if entire range also spans exactly integral units to the right and to the left then it could be shown geometrically that overall mantissa is exactly uniform, that is, that the rise on the left and the fall on the right cancel each other exactly.

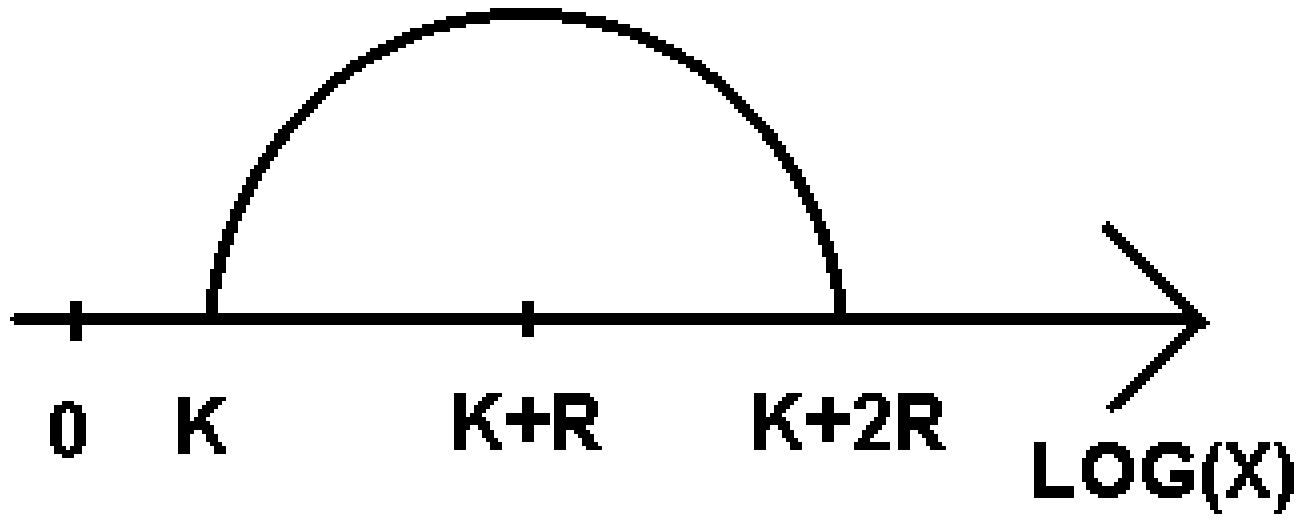

FIG 13: Semi-Circular-like Distribution of Related Log of Data

Simulations of the semi-circular-like distribution (with diameter laying on the log(x)-axis) are quite difficult to perform due to the form of the density function here, hence they were abandoned in favor of direct calculations. The actual shape of the distribution here is not truly circular, since area under the curve must sum to 1, whereas area of the semi-circle is $\frac{1}{2}\pi R^2$ while we seek total flexibility in choosing any value for the radius R, not only the solution to $\frac{1}{2}\pi R^2 = 1$. There is a need to introduce the adjusting factor $(\frac{1}{2}\pi R^2)^{-1}$ leading to $(\frac{1}{2}\pi R^2)^{-1}(R^2 - (\log(x) - (K+R))^2)^{(1/2)}$ as the expression for the density function of log(x). Hence the adjusting factor bends and contorts the circular area downward.

Direct calculations of the semi-circular-like distribution of log of data shows satisfactory logarithmic behavior for $\mathbf{10}^{\textbf{Semi-Circular}}$ whenever R is roughly larger than 1.1 (i.e. having a total range – its diameter – over 2.2). It is very important to note the near perfect agreement with the logarithmic for large enough values of R, say over 4 or 5! Note that almost the same results are obtained (for large enough R) whenever center of distribution is displaced to the left or to the right by **any** amount; hence the idea presented here is quite general and independent of the location on the log(x)-axis.

Taking stock now, we note that in general, a range of roughly 2.5 was quite sufficient for the normal, triangular and semi-circular-like related log distributions to bequeath a near logarithmic behavior to the data. It is conjectured here that this 2.5-value is quite general and should also be sufficient for other such symmetrical log densities to yield a near logarithmic behavior for the data, and much closer to the logarithmic on longer ranges say over 4 or 5.

A crude simulation of a mini Hill's super distribution was performed. It employed just 6 different distributions and a throw of a normal unbiased die (having 6 sides) 50,951 times to decide which to draw from. The distributions were made quite rare or unusual on purpose so as to at least attempt to capture the spirit of Hills's idea. Care was taken to avoid distributions with negative values. Below is the list of these 6 distributions.
U(0, 1) signifies the uniform on (0, 1). Note that some are short chains of distributions themselves:

$DIST_1 = 5.4*U(0, 1) + 6.034*U(0, 1) + 0.054*U(0, 1)$

$DIST_2 = |( - 0.0042312*LOG_e(1- U(0, 1)) -5) + 0.0042312 * U(0, 1)|$

$DIST_3 = |Normal( 5 , 3)|$ the absolute value of the Normal with mean=5 s.d.=3

$DIST_4 = |\,|\,|( - 0.345* LOG_e (7 + U(0, 1)) - 3| - 0.345*U(0, 1))| - 1.37*U(0, 1)|$

$DIST_5 = |Normal( 0.002442281, 0.256533505)|$

$DIST_6 = |(7*U(0, 1)) - (3*U(0, 1))|$

LD distribution of this simulation did not come out exactly logarithmic of course, but it was fairly close to it. The histogram of its log (base 10) resembles slightly the normal distribution perhaps! It's quite symmetrical! Also the range of the log is wide enough as compared with unit length, roughly from –1 to 2, spanning about 3 units of length. This simulation result hints about **a possible connection between Hill's super distribution and the lognormal, or at least a related log curve resembling the Normal, and hence between the totality of our everyday data (or any data derived from multiple sources) and the lognormal.**

(Even though the lognormal involves taking base e to the normal distribution while the graph here is of base 10, it doesn't matter much; natural log base e to the power of Hill's mini data gave almost the same shape.)

The result may be due to the particular choices of the 6 distributions above of course, and not general and so more representative simulations of Hill's distributions should be performed. In any case, a good simulation however finite, should certainly involve many more distributions than merely six. Here is the histogram of the log base 10 of this mini Hill's super distribution:

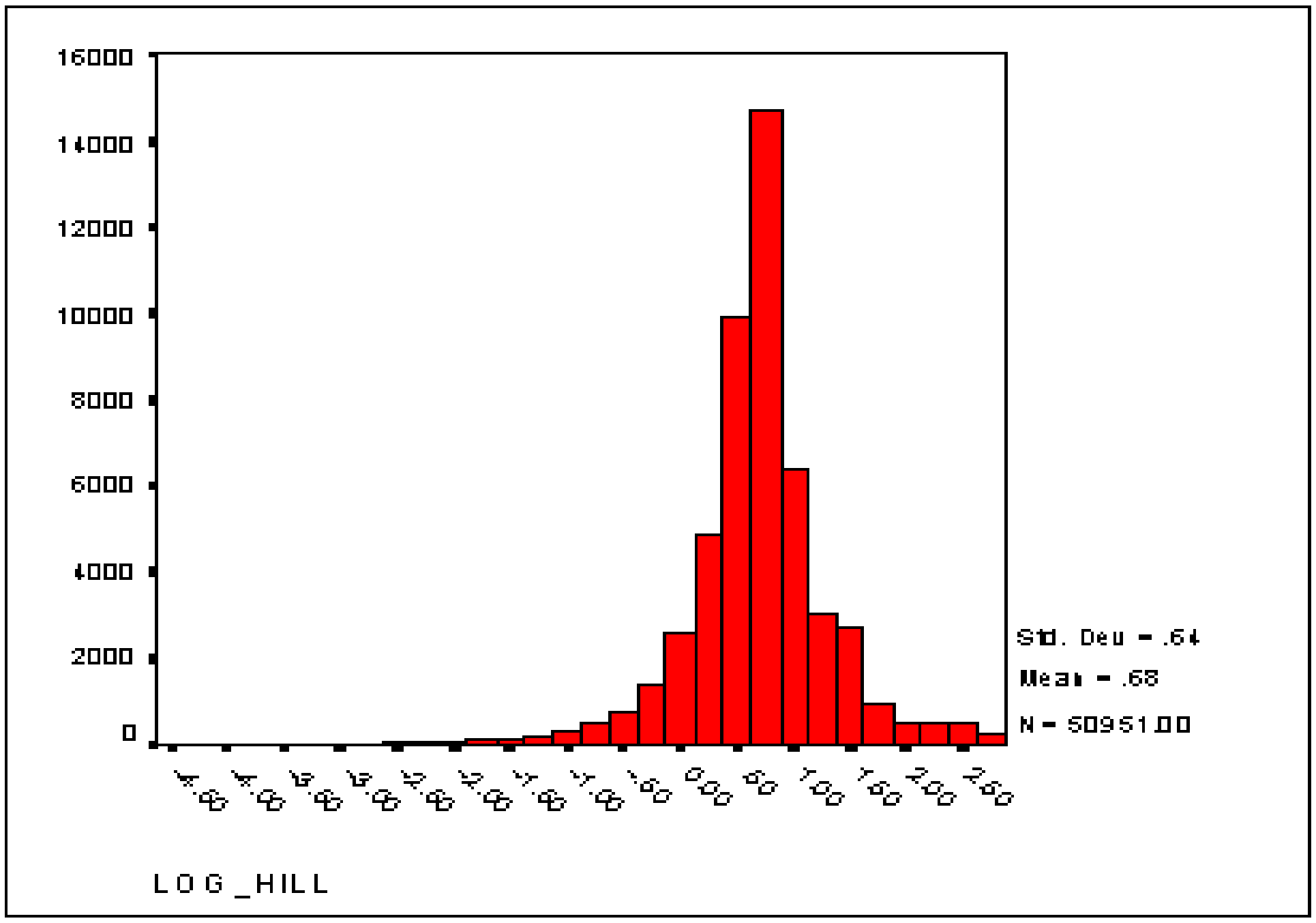

FIG 14: Histogram of the Related Log of a Mini Hill's Super Distribution Simulation

A chain of 5 uniform distributions was simulated, 19,975 simulations in total, with the last one as uniform on (0, 10000). Schematically: U(0, U(0, U(0, U(0, U(0, 10000))))). The resultant distribution appears to come with a very long one-sided tail to the right and its LD behavior is nearly perfectly logarithmic. <u>Its related log distribution is not symmetrical, and its tail to the left is longer than the one to the right.</u> In any case, its range is wide enough as compared with unit length **(at the very comfortable 5-unit length, compensating perhaps for its asymmetry!)** and it starts and ends on the log(x)-axis itself without too steep a rise or a fall. Here is the unadjusted p.d.f. (histogram) of related log of this particular simulation scheme:

[Note: This may appear a mirror image of Flehinger's scheme albeit with a much shorter 5-squence length, yet it is **not** really equivalent, since for the chain of uniform distributions LB starts at 0, while for Flehinger's scheme LB starts at 1, and there are infinitely many IPOT numbers in that LD-special interval of (0, 1)!]

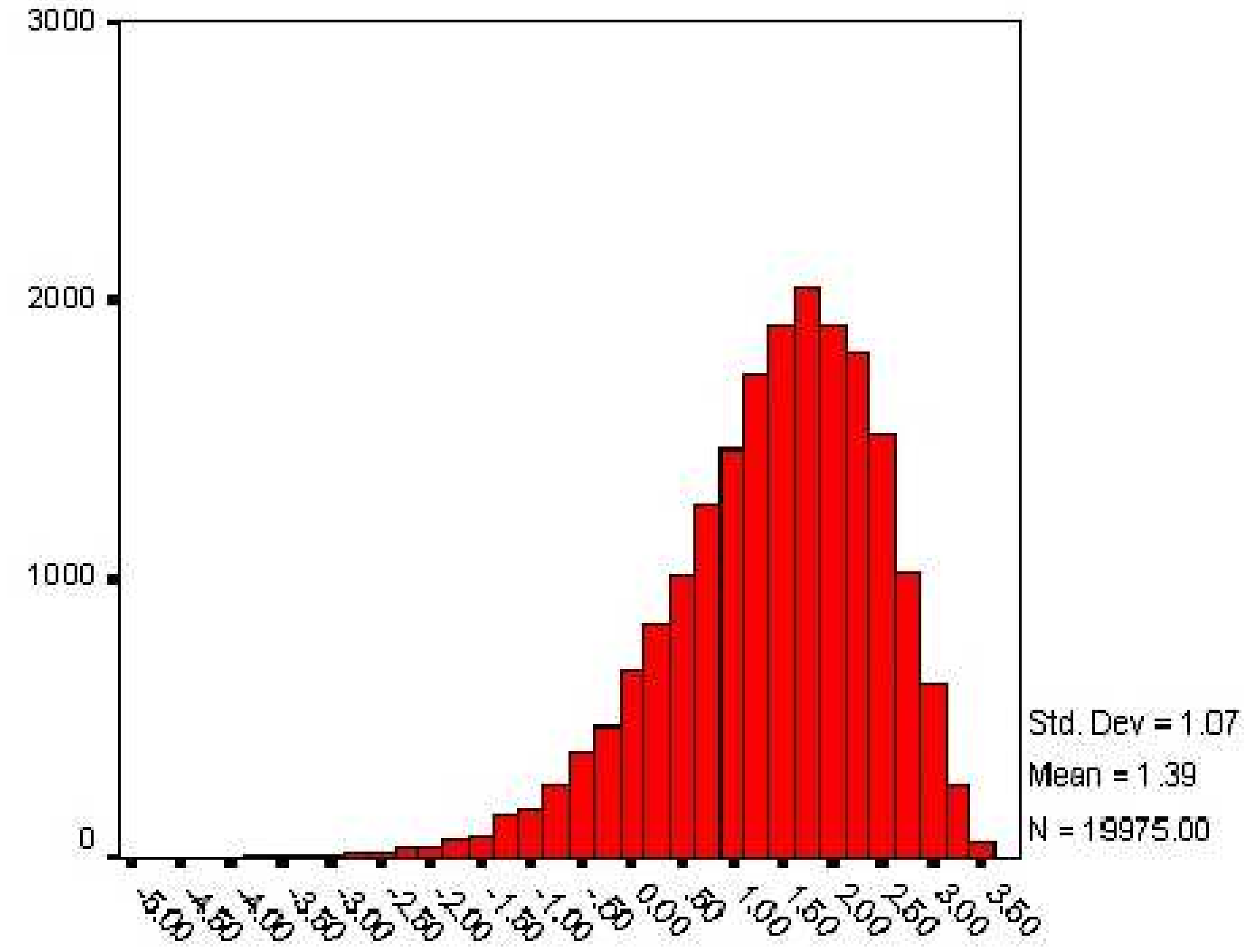

FIG 15: Histogram of Related Log of a Simulation of a **Chain** of 5 Uniform Distributions

Inspired by such non-symmetrical related log distributions that nonetheless bequeath nearly logarithmic behavior to its data, we are led to explore simulations of an asymmetrical triangular distribution, as in the following figure:

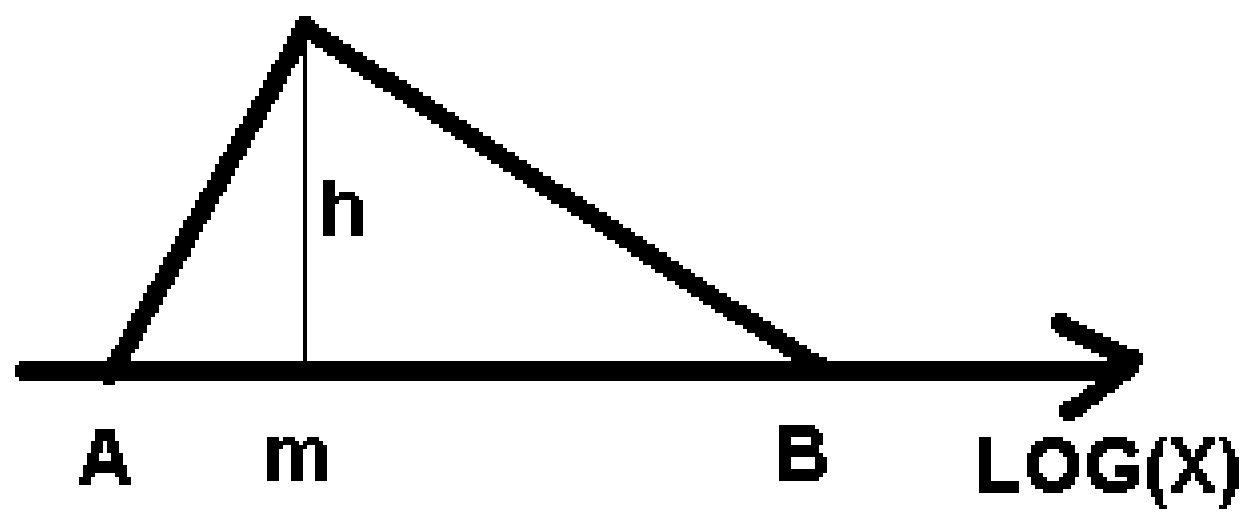

FIG 16 : Asymmetrical Triangular Distribution of Related Log of Data

This is done using the following programming command to obtain a random log(x) value:

```
IF      rand<(m-A)/(B-A)       THEN A+SQRT(rand*(m-A)*(B-A))
ELSE      B-SQRT((1-rand)*(B-m)*(B-A))
```

where B and A are the upper and lower bounds respectively, m yielding the highest density, h the height at m, and rand a simulated value from U(0, 1). These simulations point to logarithmic behavior for $\mathbf{10}^{\textbf{Asymmetrical Triangular}}$ whenever (B – A) is roughly larger than 3 or 4 (hence: range is wide enough) and m not too close to either A or B, that is, m approximately not within 0.4 units of distance from either A or B (hence: curve is not too steep, and not too asymmetric). It might be conjectured that here again, as for the chain of uniform distributions scheme above, **the system needs to compensate for its asymmetry by using a slightly wider range!**

But how does it all work out to uniformity of mantissa in spite of its asymmetry? The conceptual explanation is that to the left of the **m** point, the curve is very steep, so low digits suffer greatly (in comparison with the logarithmic situation) but very briefly, while on the right side density incline quite gently helping low digits only mildly there, yet doing so for a much longer duration (longer stretch) resulting in that same perfect trade-off just as symmetrical densities do.

One consistent result that was observed in all of these simulations is that $\mathbf{e}^{\textbf{Related Log}}$ demanded much larger range on the (common/decimal) log(x)-axis in order to bequeath logarithmic behavior to the data than $\mathbf{10}^{\textbf{Related Log}}$ had demanded.
For example, simulations for related log(x) in the form of the symmetrical triangular distribution required a range of about 2 units on the log(x)-axis for a reasonable logarithmic behavior for $\mathbf{10}^{\textbf{Symmetrical Triangular}}$, while requiring 4 such units (double) for $\mathbf{e}^{\textbf{Symmetrical Triangular}}$ to behave about logarithmically. Surely, one should not generally mix differently based logs, nor confuse one base (10) with another (e), and by "related log" we meant all along common (decimal) log base 10, but there are some very important implications here for the lognormal distribution which is defined as $\mathbf{e}^{\textbf{Normal}}$ !

Also of note is that ANY transformation of related log distributions by a translation by an integral value didn't change LD distribution of x in any way. In other words, moving density of related log to the right or to the left by an integer had no impact whatsoever on LD of x! Of course there is no surprise here at all, as this integral log transformation translates into a common factor of an IPOT number (such as 10, 100, and so forth) applied on all values of x equally, and such a transformation obviously is LD-neutral!

Moreover, in general, transformations of related log distributions by a translation by a NON- integral values, i.e. by **fractions**, doesn't change LD distribution much –given that range on the log(x) is wide enough, as mentioned earlier, and thus lending a much more general importance to what was observed in this section.

An extremely important feature in ALL of these simulations and calculations is that while those related log densities (be it the normal, semi-circular, or the triangular) came in forms either nicely symmetrical or just slightly skewed, x itself for all of them without an exception always came out with an asymmetrical density having a decisive and sharp fall as we progress to the right on the x-axis, namely: having a clear and prominent tail to the right, an issued further explored in our next section.

## [14] Density is Falling - Having One-Sided Tail to the Right

As noted earlier, in all simulations of symmetrical and non-symmetrical related log distributions, density of data itself always diminishes along the x-axis. The same observation about a falling density having one-sided tail to the right was made in the simulations of Hill's super distribution, in those of the chains and the averaging schemes; hence this appears to be an intrinsic and necessary property of any logarithmic behavior!

## [15] Fall in Density Must Be Well-Coordinated Between IPOT Values

A density falling off as sharply as in the case of k/x is NOT a sufficient condition for logarithmic behavior AT ALL! An essential feature or criteria of logarithmic behavior is that such a fall must be properly aligned and coordinated with those relevant intervals between integral powers of 10 such as 1, 10, 100, and acknowledging this fact is a crucial factor in the understanding and prediction the LD phenomena.

The following two illustrative examples help to shed some light on this crucial issue. The chart below depicts one possible serious pitfall wherever density of data is falling over the 'wrong' interval; a distribution having the same shape as k/x but defined over (-3, 6) instead of the more appropriate one of say (1, 10), as it is shifted to the left by 4 units. Its density is p.d.f. = k/(x+4) defined over (-3, 6). In this case LD is decidedly non-logarithmic, and furthermore, it's not even monotonic, as digit 2 leads the most, gathering 39% of numbers. The exact LD distribution here is {28%, **39%**, 8%, 8%, 7%, 2%, 2%, 3%, 3%}, where digit 2 by far takes the most leadership. A similar yet different pitfall is the second example of shifting k/x to the right by 4 units as in p.d.f.= k/(x - 4) over (5, 14) instead of the more appropriate one of say (1, 10). This shift does not involve any negative values whatsoever, yet its LD distribution does not resemble the logarithmic in any way, it's not even monotonic, LD here is: {22%, 0%, 0%, 0%, **30%**, 18%, 12%, 10%, 8%}, with digit 5 drawing the most leadership of course. Why is digit 5 the winner here? Because the 'wall' of numbers starts abruptly at 5.00 and falls off from there steadily throughout, hence area on (5, 6) necessarily is the largest by far and digit 5 wins strongly.

In all these examples, we encounter **distributions having the appropriate density shape for logarithmic behavior, falling at the correct logarithmic rate, yet badly-oriented on the x-axis and falling over the wrong intervals, in other words, the densities are being out of phase with those intervals bordered by adjacent integral powers of ten and hence not logarithmic.** All of this is an extremely important LD-lesson in general that the range over which distribution is defined is as important in LD as its shape of course, and that only the right combination of both - shape over appropriate interval - is what determines resultant LD distribution.

A distinction should be made between the following two types of logarithmic densities:

**(I) A density that starts out abruptly very high at a certain point on the x-axis and then keeps falling consistently having a steady tail to the right.**
**(II) A density that first rises up gradually, then hangs up there about flat, and then reverses direction and consistently falls.**

For type (I) most of the data (area under the curve) is near the starting point, and ***it very hard to envision it logarithmic unless starting point is some IPOT value (perhaps 0 or 1!)***, otherwise how could digit 1 obtain its large 30.1% leadership!? The only way this could be done logarithmically without starting at an IPOT point is to begin abruptly on some non-IPOT value say $\mathbf{10^{Z}}$ where Z is **not** an integer but end just as abruptly on $\mathbf{10^{Z+integer}}$, in other words having an integral exponent difference! Note that what we called here the 'starting point' is what was called Lower Bound (LB) in all our averaging schemes earlier. For type (II) not much data is near the starting point, in fact most of it is much farther to the right, and for such densities it does not matter much as far as LD is concerned where the starting point is at, as it could definitely be at some NON-IPOT value and yet quite logarithmic. ***Typical pieces of everyday data are for the most part of type (II), and thus starting point is of no consequence mostly. Yet, for that small portion of real life data of type (I) that is still logarithmic in its own right, it may be argued that at times the natural pull that the numbers 0 and 1 exert plays a crucial rule here utilizing them as anchors serving to coordinate and align the fall in the density with those intervals bordered by adjacent integral powers of ten (but this is true only for some particular types of data, e.g. address data). Certainly 0 and 1 are the starting values of so much of our data, 0 for measurements, and 1 for count data, and fortunately for a few types of data in LD they just happened to be IPOT numbers!! Yet it also should be noted with caution that on the other hand often real data congregate well above the starting points of 0 and 1.***
For a counter example where count data mostly aggregate well above its potential starting point of 1, imagine the count of annual accidents in a given city, as but one example.

For a real-life illustration of the enormous importance of 1 or 0 as anchors in LD for densities of type (I), let us explore the adverse implication to LD of their absence. Imagine a country where the king arbitrarily decrees that all streets must start with the number 5 being a holy or lucky number in their religion or folklore. We perform a Greek-Parable-like scheme to model this hypothetical street address data, similar to what was done in the simple averaging schemes, but with LB strictly at 5. All streets have the sign at the door of the first house showing #5, and then for consequent houses the numbers 6, 7 ,8 , and so forth. Street length varies from 1 to 9 houses long. For simplicity not all streets get equal importance in terms of having the same # of houses, hence the shortest street represents only 1 house (showing the strange sign #5) while the longest street represents 9 houses, but LD results of schemes giving equal importance to all streets point to the same conclusion almost. LD distribution for those 45 houses is:
{10/45, 0/45, 0/45, 0/45, 9/45, 8/45, 7/45, 6/45, 5/45} or
{**22%**, 0%, 0%, 0%, 20%, 18%, 16%, 13%, 11%}
Hence digit 1 still squeezed a narrow and unexpected victory, but a lot of focus is given to digit 5, standing at the 2nd place, because 'density' starts with this digit!
Clearly there is not even a resemblance to the logarithm! Not even to Stigler's Law! Absent here even that monotonically declining probabilities pattern. And this is the price we pay obeying the king and losing those crucial natural anchors of 1 and 0!

## [16] Some Synthesis Between The Deterministic and The Random

It may first appear far-fetched that two totally different processes; deterministic multiplication processes on one hand (including exponential series, the Fibonacci series, and others) and probabilistic large and varied collection of data or distributions on the other hand, both are stamped with that identical logarithmic LD signature for all digits! There must be not only some minor LD commonality here but rather another much deeper and intrinsic correspondence between the two! What is it?

Let us then select two representatives, one from each camp, and search for some fundamental common denominator: **(1)** The simplest of the simple averaging schemes, namely the Greek parable, representing collections of distributions. **(2)** The first 21 elements of 12% exponential growth starting at 1, representing multiplication processes.

The idea here is to **UNIFY** both processes by viewing them as datasets, hence any deterministic process is to avail itself also as probabilistic random variable in a sense of being viewed as a collection of numbers all with equal discrete probabilities. For example: prime numbers up to 50 could be viewed as in picking randomly from the set of elements {2, 3, 5, 7, 11, 13, 17, 19, 23, 29, 31, 37, 41, 43, 47} all having equal probabilities. Here Probability(P<11) = 4/15, and so forth. We would then be able to contemplate and construct density curves for <u>both</u> processes and make comparisons.

Let us reconsider the simple averaging schemes and instead of averaging the various LD distributions of those many individual intervals (all accorded equal importance), let us first aggregate all the various intervals into one single dataset (in some reasonable and consistent manner) and then measure LD pulse (only once) of that grand dataset. Since each interval is accorded equal weight and since intervals come with difference lengths, it is necessary to 'stretch' them out by using repetition until all are of equal length (namely the length of the longest interval). In other words, treating the intervals as datasets and combining them all in one grand collection of numbers. The motivation here is to attempt to glance at the actual distribution form of the imaginary data itself (and its other properties perhaps) instead of just focusing on resultant LD distribution. For the Greek parable we assume that 1 out of 9 topics will be about a spouse, 1 out of 9 topics will be about slaves, and so forth. To accord equal importance, we stretch topics to 9 by assigning 9 conversations per topic, resulting in a total of 81 typical numbers in typical conversations (9 topics * 9 conversations). In other words, let us <u>extrapolate</u> the original setup to achieve uniform length of 9 of each topic.

Here is a table representing one such attempt:

| | | | | | | | | | |
|---|---|---|---|---|---|---|---|---|---|
| **spouse** | 1 | **1** | **1** | **1** | **1** | **1** | **1** | **1** | **1** |
| **houses** | 1 | 2 | **1** | **2** | **1** | **2** | **1** | **2** | **1** |
| **slaves** | 1 | 2 | 3 | **1** | **2** | **3** | **1** | **2** | **3** |
| **oranges** | 1 | 2 | 3 | 4 | **1** | **2** | **3** | **4** | **3** |
| **sheep** | 1 | 2 | 3 | 4 | 5 | **1** | **2** | **3** | **4** |
| **chicken** | 1 | 2 | 3 | 4 | 5 | 6 | **2** | **4** | **6** |
| **dogs** | 1 | 2 | 3 | 4 | 5 | 6 | 7 | **3** | **6** |
| **olives** | 1 | 2 | 3 | 4 | 5 | 6 | 7 | 8 | **5** |
| **Gods** | 1 | 2 | 3 | 4 | 5 | 6 | 7 | 8 | 9 |

FIG17: The Greek LD Parable as a Dataset

Care should be taken in filling in the blanks for the numbers on the upper-right corner. Surely we ought to imitate whatever is on the left side of each row, row by row. We first simply repeat them in order, and then - when there is no more room for another repetition of the entire set – we pick at random from the numbers on the left. Hence we note the arbitrary manner in which the few remaining numbers on the extreme right side are chosen, but variations have very minor overall effect on results.

Here is the grand dataset for this particular attempt presented in the table above: **{1,1,1,1,1,1,1,1,1,1,2,1,2,1,2,1,2,1,1,2,3,1,2,3,1,2,3,1,2,3,4,1,2,3,4,3,1,2,3,4,5,1,2,3,4,1,2,3,4,5,6,2,4,6,1,2,3,4,5,6,7,3,6,1,2,3,4,5,6,7,8,5,1,2,3,4,5,6,7,8,9}.** The LD distribution of this grand dataset is of course extremely close to the result of the original averaging scheme where LD was obtained by averaging out the 9 different mini LD distributions of each topic.

Treating as a dataset the collection of the first 21 elements of 12% exponential growth series that starts at 1 we obtain: **{1.00, 1.12, 1.25, 1.40, 1.57, 1.76, 1.97, 2.21, 2.48, 2.77, 3.11, 3.48, 3.90, 4.36, 4.89, 5.47, 6.13, 6.87, 7.69, 8.61, 9.65}.**

Now let us glance at these two very different datasets using a variation on stem and leaf plot, a histogram of sorts, as in the two following figures:

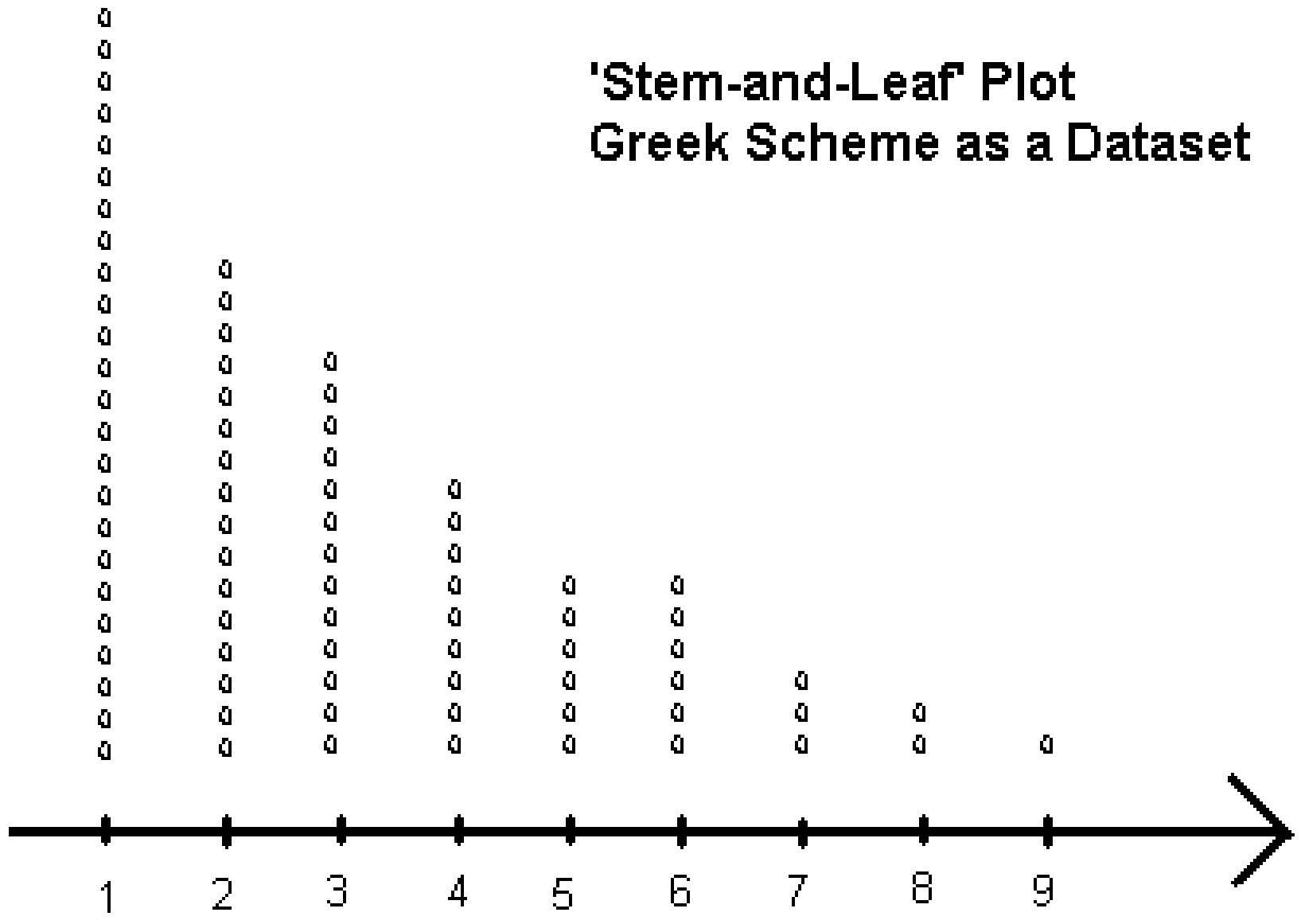

FIG 18: **THE RANDOM PROCESS,** A Variation on Stem-and-Leaf Plot, Greek Scheme

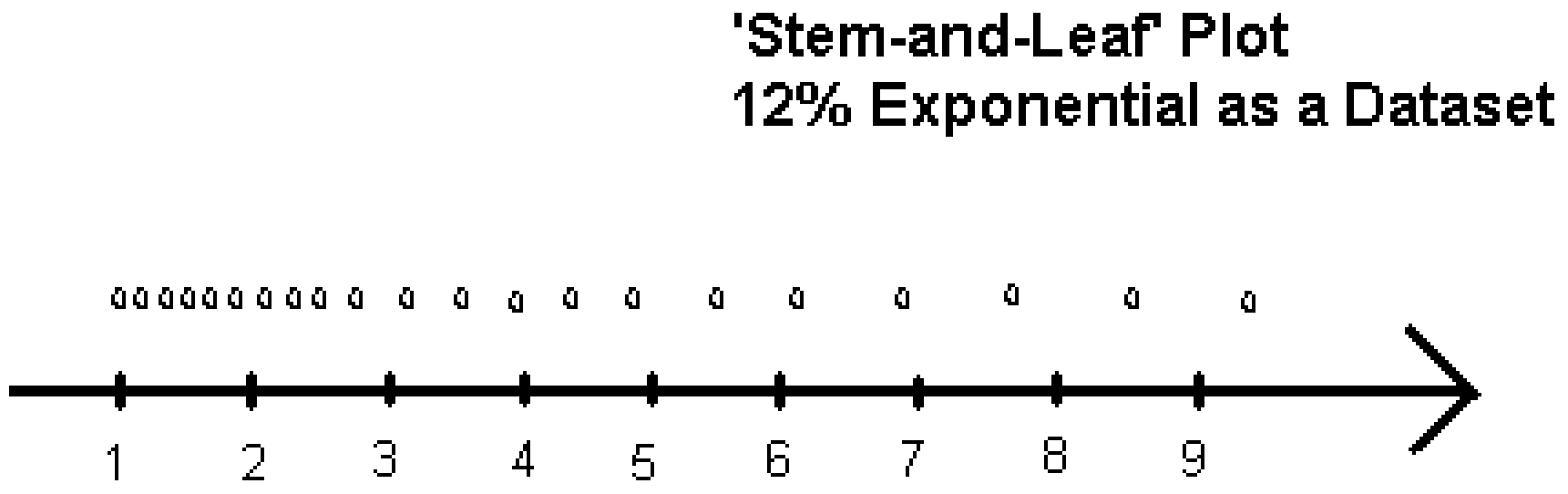

FIG 19: **THE DETERMINISTIC PROCESS,** A Variation on Stem-and-Leaf Plot, 12% Growth

The common denominator here can be easily visualized; it's simply the frequency per unit length, or density, which is diminishing for both datasets. Both appear as having lop-sided densities of sorts with tails to the right. Left unexplained are two points: (A) The fall in the density <u>for both of them</u> is nicely aligned and synchronized between IPOT values, and (B) the sharpness in the fall of the density is <u>the same for both processes</u>.

The idea imbedded in the Greek parable ultimately points to - and converges to - Stigler's law which is quite close to the logarithmic and from which further iterations [Flehinger's] lead directly to the correct LD distribution. Also the model of the averaging scheme is similar to Hill's super distribution, in that it is basically a collection of intervals (thought of as distributions as well), hence it should be considered a fair representative here.

## [17] LD-Dichotomy Between Deterministic & Random Processes

Turning our attention again to exponential growth series in general, we note that distances between elements of exponential growth are expanding of course, hence we observe that their 'discrete density' curve is falling, as seen in the previous section:

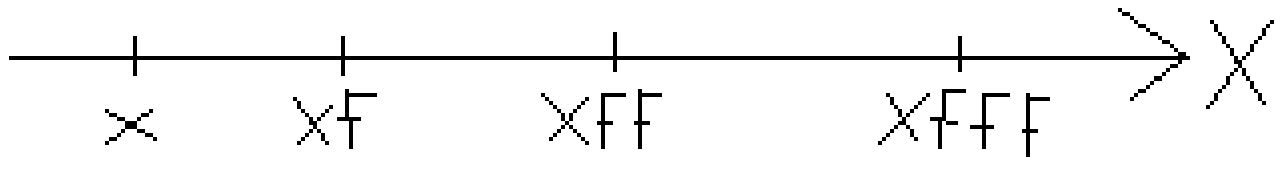

FIG 20: Exponential Growth Series on the X Axis

$$xf - x = x*(f-1)$$
$$xff - xf = x*(f-1)*f$$
$$xfff - xff = x*(f-1)*f*f$$

This is so because $f = (1 + percent/100) > 1$. For example, 6% growth has factor $f = 1.06$. Hence, occurrences on the x-axis become progressively less and less frequent (less dense) as per any fixed interval of some arbitrary length. Furthermore, let us examine the related log series of an exponential series to gain a better insight here. For exponential growth such as $\{B, Bf, Bf^2, Bf^3, \ldots, Bf^N\}$ the related log series is simply: $\{ LOG_{10}(B),\ LOG_{10}(B) + LOG_{10}(f),\ LOG_{10}(B) + 2{*}\ LOG_{10}(f),\ \ldots\ ,\ LOG_{10}(B)+N{*}\ LOG_{10}(f) \}$. Now, since this is nothing but the steady additions of the same constant, namely $LOG_{10}(f)$, from a fixed point, namely $LOG_{10}(B)$, we conclude that: For exponential growth in particular and multiplication processes in general, their related log is uniformly distributed (albeit discretely, not continuously). But by Proposition II, the approximate or most appropriate 'density function' here is then k/x. Hence the following important conclusion: Exponential growth in particular and multiplication processes in general are related to the k/x distribution, having consistent, steady and continuous logarithmic behavior all along their entire interval. Their approximate 'density' is seen as falling steadily at the constant logarithmic rate, and having a tail to the right all throughout their entire range. This feature in multiplication processes - having a steady LD behavior throughout cannot be found in random processes. Note that since exponential series are almost always distributed (i.e. considered) over very large intervals as compared with unit length of its related log function, LD requirement of Corollary I for an integral value of exponent difference can be easily waived, implying a near logarithmic behavior in most such typical cases. Exponential series, with their constant multiplicative factor, are good representation of multiplication processes in general, because the other extreme, that of random multiplicative factor, yields overall uniform related log distribution just the same, by the same line of reasoning above, because its related log distribution can be thought of as arising from slight disturbances of an original equidistance arrangements of values as in the exponential growth series, where some values have been randomly displaced to the left and others to the right, and all by some random distance of displacement, thus leaving overall uniformity of log density intact.

Related log densities has a way of revealing much about data that neither density itself nor LD distribution can show, and therein lies the fundamental difference between the Random and the Deterministic. The Random process referred by Benford's Law is a collection of distributions, data composed and aggregated from a variety of sources as in Hill's proof, and there related log starts out from the very bottom of the lox(x)-axis, rises up gradually, then falls gradually as well all the way back to the bottom, ending there, imitating slightly the shape of the Normal or the contorted semi-circular in some way. It has a perfect logarithmic behavior if considered in the aggregate over its entire range, but it's local LD behavior over small segments within the entire range actually evolves from favoring high digits, to favoring low ones slightly, to the logarithmic, and then on to extreme bias against high digits. Deterministic multiplication processes on the other hand see their related log densities uniformly distributed, always hanging above log(x)-axis horizontally and steadily, with LD behavior being completely steady and consistent, true for all local subintervals anywhere. The Random density itself appears a bit similar to the lognormal, while density of the Deterministic itself is clearly of the k/x type. A note: no matter how uniquely and perfectly logarithmic the distribution k/x is considered to be, it is NOT at all the one that is appropriate for typical everyday real data for which Benford's Law refers to! The k/x distribution is exclusively connected with deterministic multiplication processes!

It is worth noting that any subset cut out from the whole of some logarithmic data is still logarithmic given that extraction from the whole is random, or equivalently, whenever the whole has been well mixed before an ordered extraction on the newly shuffled data took place. On the other hand, when dataset is ordered (i.e. sorted low to high say), cutting off pieces from it in an orderly fashion will not result in logarithmic behavior unless logarithmic data is deterministically-driven and enough elements are included. For this reason, segments of exponential growths or decay are always more readily logarithmic even though the infinite series is being sampled with its order fully intact, but this would be so only when LD cycle is much shorter relative to the length of the entire sequence being examined. In the same vein, subsets from the exponential and lognormal distributions should be taken randomly from the whole range of $(0, +\infty)$ without any focus on specific sub-intervals.

Now, an alternative point of view about the remarkable correspondence between these two disparate processes can be given via the perspective of mantissa, namely that these processes both come with uniformity of mantissa because of their particular related log densities; (I) the deterministic process achieving uniform mantissa by the direct way of related log itself being flat and uniform over an integral interval; (II) and the random achieving uniform mantissa by having its related log arising and falling gradually from the log(x) axis itself over large enough a range.

## [18] LD of Lognormal, Exponential, and k/x Distributions

The lognormal distribution represents a variable obtained from a product of many small independent and identical factors and is modeled on the *multiplicative central limit theorem* in mathematical statistics hence logarithmic behavior is expected when considered over its entire range of $(0, +\infty)$ and given the right parameters. On the other hand, the normal distribution (whose curve is symmetrical), being a model of a variable obtained from numerous additions of small independent and identical factors is not expected to show any logarithmic behavior and indeed does not show any such behavior no matter what parameters are chosen,.

Note that if we represent the normal as $N = x_1+x_2+x_3+\ldots+x_n$ over $(-\infty, +\infty)$ where $x_i$ are all independent and identical variables, then

$$e^{\text{normal}} = e^{(x1+x2+x3+\ldots+xn)} = e^{(x1)} * e^{(x2)} * e^{(x3)} * \ldots * e^{(xn)} = \text{LogNormal}$$

defined over $(0, +\infty)$ namely a process of repeated multiplications, and thus under the right parameters a logarithmic LD distribution should follows.

*In brief*: (1) lognormal = $e^{\text{normal}}$ & (2) normal = ln(lognormal).

Often, everyday data is modeled on and compatible with the lognormal distribution rather than the normal, and this fact is employed by those who assert that a lot of our everyday data is derived from repeated multiplications (hence the lognormal). Yet that assertion is debatable, and disputed by many.

Simulations of the lognormal distribution show some good logarithmic behavior whenever the shape parameter is over 0.7 (roughly), and regardless of the value assigned to the location parameter, while a near perfect logarithmic behavior is observed whenever parameter is over 1.2 roughly. For shape parameter in the range between 0 and 0.5 there are wild fluctuations in its related LD distribution that also depend to a large extend on the location parameter, and for the most part doesn't resembled the logarithmic at all. It is worth noting that for very low values of the shape parameter the distribution is almost symmetrical and becomes progressively more skewed to the right as the shape parameter increases in value; coinciding with progressively better conformity to the logarithmic! This is nicely consistent with what was noted earlier, namely that asymmetrical one-sided tail to the right is essential for logarithmic behavior. It is also consistent with all else that was observed in earlier sections, namely that for any base B, $B^{Symmetric}$ is logarithmic given that the range of the symmetric distribution is wide enough. The shape and location parameters of the 'resultant' lognormal are nothing but the s.d. and the mean respectively of the 'generating' normal distribution, hence low values of the shape parameters are situations where the symmetrical log(X) distribution (the normal in this case) is focused on some too narrow a range, which leads to non-logarithmic distribution for X!

The exponential distribution with parameter p is defined as $f(x) = (p)*(1/e^{px})$ for x on the interval $(0, +\infty)$. It has an asymmetrical one-sided long tail to the right that falls off proportionately to x (or to $e^{px}$ rather). Its overall LD distribution (no matter what value the parameter takes) is quite close to the logarithmic but never exactly so when considered over its entire range of $(0, +\infty)$. A more careful investigation shows that individual digits nicely oscillate closely around their logarithmic values indefinitely in an ever-widening cycle as parameter p increases, but are always out of phase with each other, hence never truly logarithmic. In order to avoid a loss of continuity, a separate section is added at the end of the article to give a more detailed account of those interesting and unique patterns in LD behavior of the exponential. The special case of the perfectly logarithmic k/x defined over integral exponent difference was discussed earlier. Let us now put the <u>exponential</u>, <u>lognormal</u>, and <u>k/x</u> distributions under a closer scrutiny to see what happens locally to their mini-LD distributions on those sub-intervals bordered by adjacent integral powers of ten, as this leads to significantly better understanding here. Here are some LD distributions, interval by interval, of 13,000 simulations from the exponential with parameter 0.02547: [Note: 4 values falling below 0.01 were discarded.] Between 0.1 and 1: {0.09, 0.11, 0.17, 0.10, 0.10, 0.09, 0.13, 0.09, 0.11}. Between 10 and 100 : {0.24, 0.20, 0.15, 0.12, 0.09, 0.07, 0.06, 0.04, 0.03}, while between 100 and 100 we get:{0.93, 0.07, 0,0,0,0,0,0,0}. Except for intervals very near the origin where data is quite thin with less than 3% or 4% of overall values and a bit unreliable, the pattern seen above is totally consistent across all exponential distributions having other parameters, so all this is quite general. Low digits win on each of these intervals, if given close scrutiny and large enough simulated values are generated, but their advantage over high digits is weak near origin and gets progressively stronger later on, finally achieving at the end a complete dominance of digit 1 over all others.

Here are some LD distributions, interval by interval, of 13,000 simulations from the lognormal with shape parameter 2.3 and location parameter 1. [Note: 4 extreme values were discarded.] Between 0.01 and 0.1: {0.14, 0.12, 0.13, 0.12, 0.09, 0.11, 0.08, 0.11 0.09}. Between 1 and 10: {0.31, 0.18, 0.13, 0.10, 0.07, 0.06, 0.05, 0.04, 0.05} while between 100 and 1000:{0.52, 0.17, 0.11, 0.05, 0.05, 0.05, 0.03, 0.01, 0.02}. Again, this general pattern observed in the shift of LD of the lognormal from interval to interval is general and is true for all other relevant parameters. Here though on intervals near the origin, low digits lose leadership, curve is actually ascending, but they very soon achieve equality, and immediately afterward start to win strongly until a near dominance of low ones over high ones is observed. [The discussion here regarding the lognormal distribution focuses solely on ones with shape parameter larger than 0.7 which are nearly logarithmic.] And what about LD distributions of the k/x distribution interval by interval? Clearly, there LD distributions are all the same, namely the logarithmic itself exactly, with perfect repetition on each such sub-interval, as seen earlier.

Let us now state the following in general about **any** distribution. For each subinterval bordered by TAIPOT numbers, complete uniformity of LD distribution there locally (non-logarithmic) implies uniformity of the density function itself namely a flat curve, and visa versa, assuming smoothness. A logarithmic LD distribution on such subinterval implies a tail to the right with density falling just as steeply as that in the k/x case. And cases where digit 1 leads much more than the usual 30% and high digits are almost totally excluded from leadership, are areas where the density falls even more rapidly or sharply than k/x. Clearly, overall LD distribution is obtained after averaging out (the weighted) different shapes/slopes on the entire curve incorporating all subintervals. It should be emphasized that our old familiar concept of slope from calculus does not translated here literally in the context of LD. It is certainly true that for any two competing density curves of equal average height on any given fixed interval between TAIPOT points, the one with the steeper negative slope endows more leadership to low digits. Yet in general, slope can't be divorced from interval's location on the x-axis in the context of LD. A good conceptual example is k/x distribution, where derivative is $-k/x^2$, hence slope is constantly negative but getting flatter even though local LD distributions on those subintervals between TAIPOT points remains constant throughout!

The chart below depicts all 3 distributions above superimposed. The focus of all these 3 distributions is on **(1, 100)**, and very little else is missed by neglecting to show the rest over 100. The choice of parameters was deliberate so as to have all three distributions having the same median, namely 10. Distribution k/x is defined over $(1, 10^2)$ or (1, 100), so G=2, k = 1/[G*ln10] = 0.21715. Median there is always $\mathbf{10^{G/2}}$ whenever lower bound is 1, hence here the median is $\mathbf{10^{2/2}}$, namely 10. The lognormal distribution over (0, + ∞) has here a location parameter 2.303 and shape parameter 1.11. The median there is always $\mathbf{e^{location\text{-}parameter}}$, or $\mathbf{e^{2.303}}$ = 10 in our case. The exponential distribution over (0, + ∞) has parameter 0.069314718. The median there is always **(1/p)*ln(2)** hence median here is 10.
LD distribution of this lognormal over the area of (1, 100) is close to being logarithmic, that of the k/x is perfectly such over (1, 100), while that of our exponential is a bit off and comes with the following LD distribution:
{31.8%, 18.9%, 12.7%, 8.8%, 6.9%, 6.2%, 5.1%, 4.7%, 4.6%}.

The chart attached here terminates abruptly at X=70, neglecting to show a only 7.7% of the area of the rest of k/x distribution to the right, just 4% of the rest of the lognormal, and a mere 1% of the rest of the exponential, so that the overall picture can still be visualized clearly. Also, the abrupt termination of the roof at 0.08 masks a slight rise of k/x up to 0.21 when x is 1.

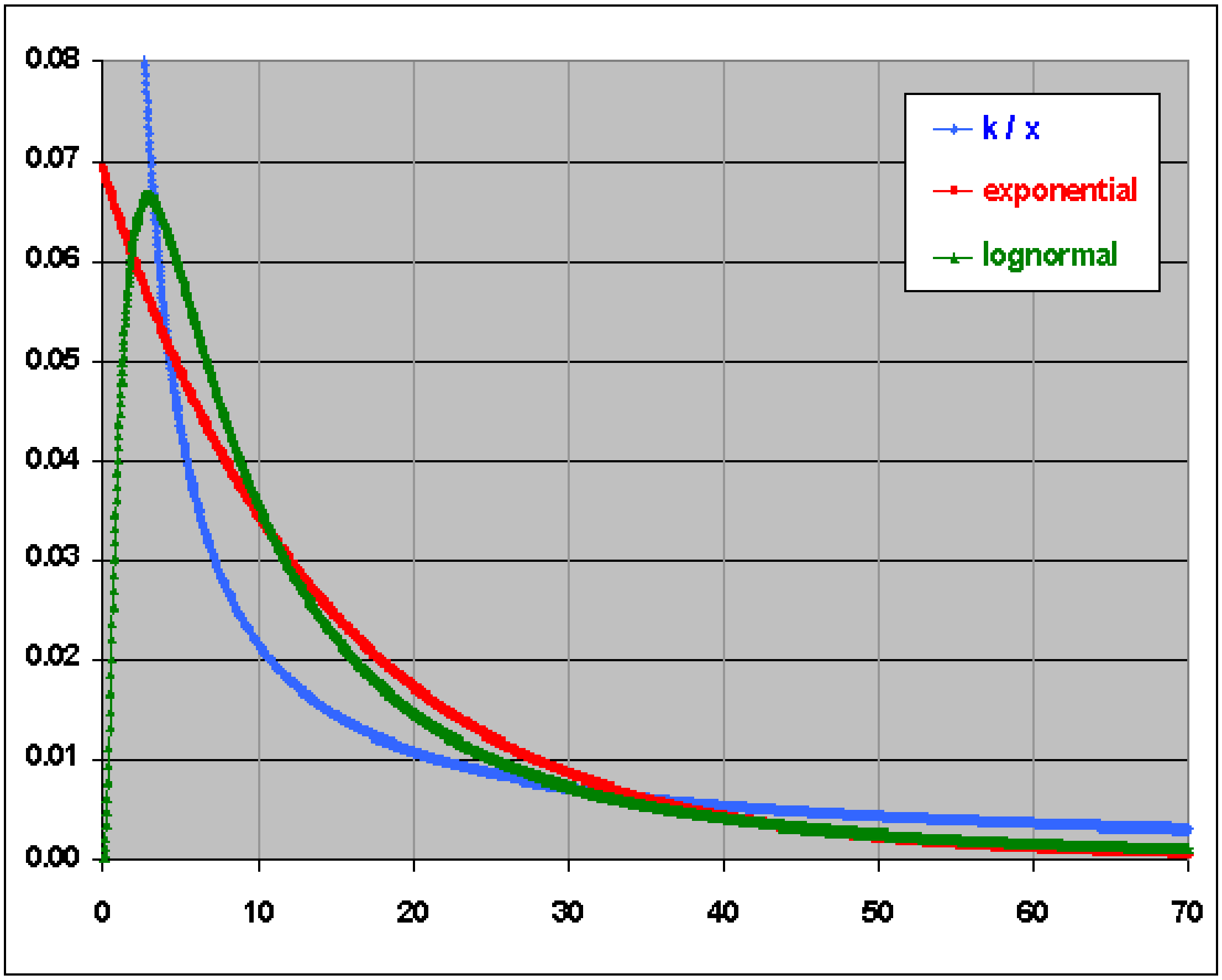

FIG 21: The Exponential, the Lognormal and k/x Distributions all centered on 10 (via the median)

All 3 distributions need to transverse these 2 sub-intervals between TAIPOT points namely (1, 10) & (10, 100) in a way that would result in the logarithmic or nearly so. While **k/x** has a steady (LD-wise) fall over these 2 sub-intervals, the **exponential** starts out on around (1, 10) not as steep as k/x, only to reverse course later and to descend even more than k/x does on roughly around (10, 100). The **lognormal** which starts actually ascending briefly around (1, 3), quickly reverses course and bow to the inevitable descend, then finally join the exponential on around (10, 100) in descending even more than k/x does.

## [19] LD-Inflection Point

The common (decimal) log density of k/x distribution is uniform, as seen in Proposition I. The common log function of the lognormal is schematically as in: **LOG $_{10}$ (e$^{NORMAL}$),** a mix of the natural and common logs and bases, and is simply **[LOG $_{10}$ e]*Normal**, which is a normal distribution also, hence the common log function of the lognormal is normally distributed. The common log of the exponential with parameter p is distributed as in the following expression: $[p]*[\ln 10]*[10^y]*[e^{-p*(10 \text{ to the } y)}]$. This is derived by employing the Distribution Function Technique for transformations. Clearly, this related log function is not symmetrical, rather, it appears (for ALL values of parameter p) as a contorted bell-shape-like distribution skewed to the left, and its entire range is always roughly 3, regardless of parameter p. Here is the chart of one typical example where p = 0.271:

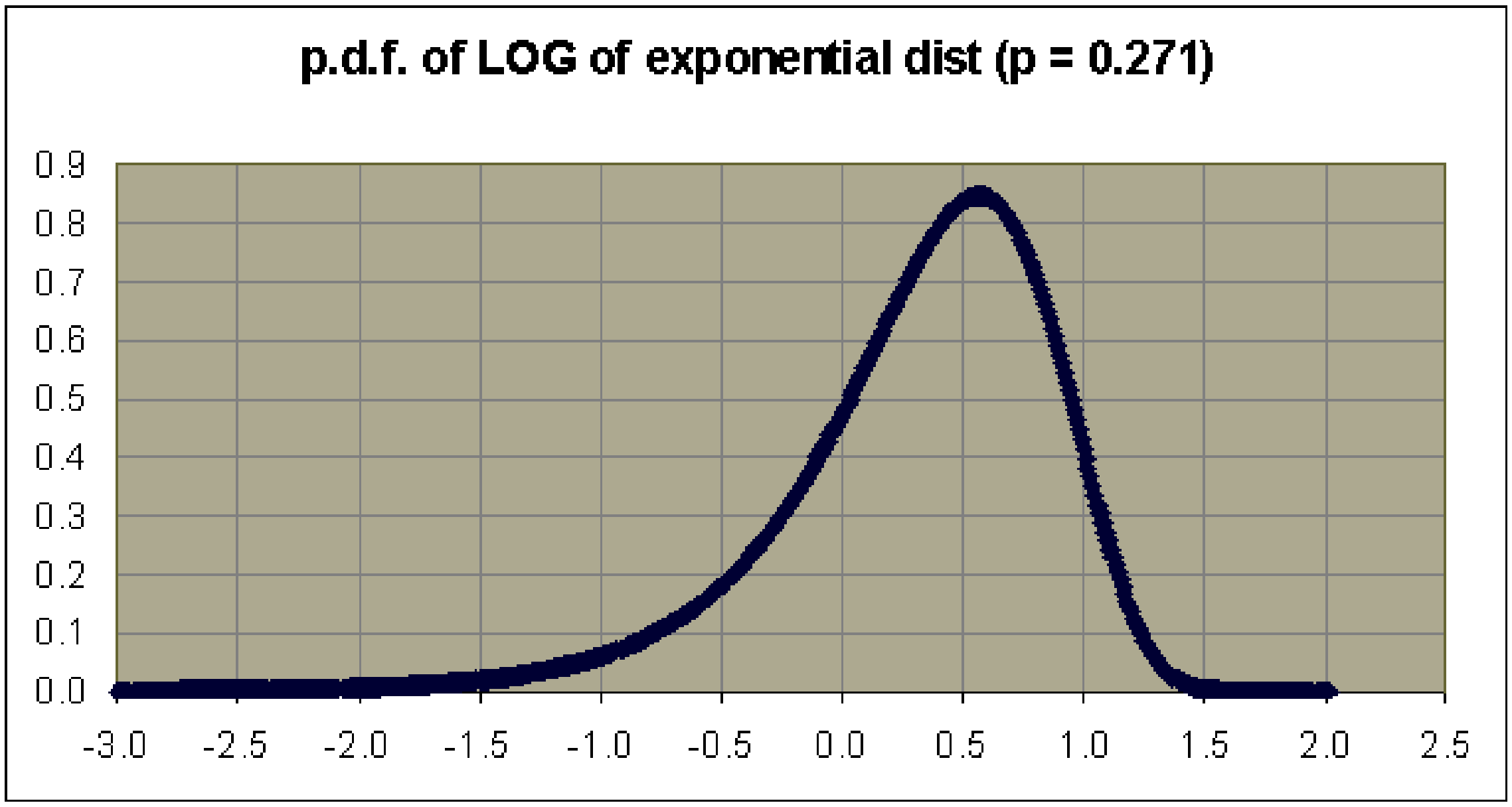

FIG 22 : Typical Distribution of the Related Log of the Exponential Distribution

As conjectured above, and encountered again here in this case of the related log of the exponential, the system can't deliver strong logarithmic behavior to x with such an asymmetrical shape unless a bit of a wider spread/range on the log(x) is available. The entire log(x) range here is always roughly no more than 3-unit length, regardless of the value of parameter p, and this is not sufficient, hence x is not fully logarithmic.

The point on the log(x)-axis of the exponential yielding maximum height (**maximum density of related log that is - not its maximum value on the extreme right**) corresponds to an important turning point for LD on the x exponential curve itself, a point we may call **LD-inflection point**. Such a point is roughly log(x)=0.6 on chart here above. Setting the derivative of log of the exponential equal zero, that is, setting $[p] * [\ln 10] * [10^y] * [e^{-p*(10 \text{ to the } y)}] = 0$ yields: $\log_{10}(\text{exponential})_{max} = \log_{10}(1/p)$ which corresponds to $1/p$ on the exponential density curve itself. In the chart above 1/p is 1/**0.271** = 0.57, quite near our 0.6 estimate earlier. From what was shown earlier we note that to the left of 1/p, although low digits overtake high ones, LD are bit more equitably distributed as compared with the logarithmic, and to the right of 1/p dichotomy between high and low digits is even more acute than that of the logarithmic inequality. Several simulations confirm this fact regardless of value of parameter, and it's quite clear that exactly at that infinitesimal flat point on the log(x) curve – namely at the max – the corresponding point on the x-axis itself is in a logarithmic-like mode borrowing this property from the k/x distribution (Proposition II). The argument here is that at max, log(x) curve has zero slope and so it's perfectly horizontal and suggestive of course of the uniform distribution (for the related log) at this infinitesimal location. It should be noted that the slope of the exponential itself is negative everywhere as seen from its derivative $f'(x) = -p^2e^{-px}$, and thus curve is always falling, with the obvious result that high digits can never hope to win the leadership game even just locally anywhere. Moreover, even though slope is constantly getting flatter (since 2nd derivative is $+ p^3e^{-px}$), nonetheless low digits are continuously taking more leadership from high digits as we move to the right.

Let us recap, the exponential distribution, regardless of parameter, has definite LD patterns over those subintervals bordered by TAIPOT. The exponential starts out giving low digits some very slight advantage, as there is almost LD-equality near the origin; then progressively low digits take stronger lead on subsequent such subintervals; once LD-inflection point value of 1/p is reached within a certain subinterval, an approximate logarithmic-like mini distribution is observed there locally; and finally on subsequent subintervals low digits take a commanding lead leaving none for high digits farther on. Yet, the important result here is that overall LD distribution, aggregating over all such subintervals, ends up quite close to the logarithmic. Note that the largest weight or portion of overall data is given to that subinterval containing LD-inflection point!

And what is LD-inflection point for the lognormal? What we are seeking here is **Max { density of $LOG_{10}$lognormal }** or **Max { density of $LOG_{10}e^{NORMAL}$** } that is: **Max { density of ($LOG_{10}e$)*(Normal)**} which is achieved whenever the generating normal itself is at max, namely whenever it's at its mean parameter [since mode - the highest point or max - median and average are all the same in any normal distribution], therefore this maximum value of the related log of the lognormal is simply its location parameter times **$LOG_{10}\, e$**. Finally, on the curve of the lognormal itself LD-inflection point is at $\mathbf{10}^{[LOG_{10}e \,*\, \text{location parameter}]}$, or simply: $\mathbf{e}^{\text{location parameter}}$.

Let is recap, the lognormal starts out actually ascending at the very beginning, hence low digits are either losing or winning very slightly, but this weakness of low digits is fleeting and only over a very small portion of overall data; then curve begins to fall off a bit and low digits get a bigger portion ; once LD-inflection point value of $\mathbf{e}^{\textbf{location parameter}}$ is reached the logarithmic is exactly observed momentarily (but over a large portion of overall data!) and then further on the fall is much steeper and low digits lead even more strongly there. Miraculously, overall LD for the lognormal in the aggregate is always **almost** perfectly logarithmic in spite of that LD emotional roller coaster.

Distribution k/x on the other hand does **not** have any LD-inflection point as it is smoothly and continuously logarithmic.

Note that while k/x distribution is exactly logarithmic over proper range, the lognormal is only very nearly logarithmic for high values of the shape parameter but can never attain a truly exact logarithmic LD behavior if one is being too strict in measuring it.

## [20] LD of the Wald, Weibull, chi-sqr, and Gamma distributions

**The Wald distribution**: the higher the mu (the location parameter) the more logarithmic the distribution is for a given lanbda. The lower the lambda the more logarithmic it is for a given mu. Hence for a logarithmic behavior here we need to have a low ratio of lambda/mu. Roughly, whenever lambda/mu is lower than 0.85 we find a near perfect logarithmic behavior there approximately.

**The Weibull distribution** is nearly perfectly logarithmic whenever the shape parameter k is very low, roughly whenever it’s 0.7 or lower, and regardless of what happens to the lambda parameter (the scale).

**The chi-sqr distribution** is nearly logarithmic whenever value of degrees of freedom is 1, approximately logarithmic when it’s 2, and becomes progressively less and less so for higher parametrical values.

**The gamma distribution** is nearly logarithmic whenever shape parameter alpha is lower than 0.7 (roughly).

## [21] Singularities in Exponential Growth Series

Let us reiterate; for exponential growth such as { **B, Bf, $Bf^2$, $Bf^3$, … , $Bf^N$** } where **B** is the base value, **P** is the % growth, and **f** is the constant factor multiplying, the relationship **f = (1+P/100)** holds. The related (common) log series are simply: **{ $LOG_{10}(B)$, $LOG_{10}(B) + LOG_{10}(f)$, $LOG_{10}(B) +2^* LOG_{10}(f)$, … , $LOG_{10}(B)+N^* LOG_{10}(f)$ }**. Due to their discrete nature, (finite) exponential series are never exactly logarithmic, only in the limit - possibly. Discreteness in Leading Digits at times causes deviations from the logarithmic, and as often could be remedied by considering infinitely many numbers and the limiting case. The understanding gained in the previous section about the nature of k/x distribution and its connection to multiplication processes hints at the necessity to stand between TAIPOT for exponential series as well – its discrete nature notwithstanding - and that is exactly what is found! Two factors facilitate and determine logarithmic behavior here: **(I)** The length of the series in focus, and normally enough sequences must be considered here for a close logarithmic fit, not just a few. **(II)** The value of exponent difference between the first and the last elements considered (*or equivalently - the difference in the two extreme values of related log*) should be ideally as close to an integer as possible for a better logarithmic fit. In general, low growth series that adheres closely enough to just one of the above requirements can be generally excused for neglecting the other and still be close to the logarithmic. For very high growth series, where LD cycle is very short, being between TAIPOT is not sufficient at all; rather it is necessary is to have truly a lots of sequences, and only for such long series the logarithmic is nearly observed.

It is certainly proper to think of **$LOG_{10}$( f )** namely **$LOG_{10}$(1+P/100)** as being a fraction, as it is so in any case up to 900% growth of P, else the integral part can be ignored as far as LD is concerned. The related log series of exponential series then represents constant **additions** (accumulation) of log(f) from an initial base of log(B). While this related log grows ever larger, mantissa on the other hand is oscillating between 0 and 1 of course, as it takes many small steps forward, then suddenly one large leap backwards whenever it overflows 1, and so forth. This is so since mantissa is obtained by constantly removing the whole part each time log overflows 1 *(whenever series ≥ 1, that is log(series) is positive or zero - an assumption which we can take for granted)* and it is easily seen as being uniform on (0,1) as more and more points of newly created mantissa keep falling there on a variety of points, covering an ever increasing portion of the entire (0,1 ) range.

An example of a normal (typical) series with base 3 and 30% growth rate having factor f=1.30 and the implied log(f) = 0.1139433. Notice how we always re-enter into the (0, 1) interval at different or newer locations, namely at 0.047, 0.072, and 0.098, and so on. This is a necessary condition for logarithmic behavior! Interestingly, such an indirect way of looking at exponential growth series leads to the detection of some peculiar LD singularities there. The argument above falls apart when the fraction log(f) happens to be such that exactly M whole multiples of it add up to unity, as in the fractions 0.50, 0.25, 0.10, 0.125, 0.05, etc., in which case that constant additions here leads to re-entering (0, 1) always at the same point and taking the same type of steps over and over again focusing only on a few selected lucky points having some very strong concentration or density, and thus ignoring all other points or sections on the interval (0, 1). All this results in some very uneven distribution on (0, 1) and thus yielding a non-logarithmic LD distribution for the exponential series itself.

An example of such an anomalous series with base **8** and **77**.8279% growth rate having factor f = 1.778279 and the implied problematic log(f) = 0.25 where exactly 4 (called M) multiples of it add up to 1. Notice how we always re-enter into the (0, 1) interval at the same location, namely at 0.153, and then always take the same subsequent <u>long</u> steps skipping so many points in between. Such state of affairs can not result in any logarithmic behavior.

Yet, even for those rebellious rates, some comfort can be found whenever **$LOG_{10}$(f)** - *the length of the steps that mantissa of the series advances along the log axis* - is quite small as compared with unit length of the log-axis, because then no matter how repetitive, peculiar, and picky the mantissa chooses the points upon which to stamp on in (0, 1), each step is still so tiny that it has no choice but to walk almost all over the interval covering most corners and locations of (0, 1)! Put simply: it has legs too short to be really picky and so it can not jump and skip much, thus it is reduced to walking almost all over the interval, willingly or unwillingly!

Therefore, LD of <u>even</u> these anomalous and picky series begin to resemble monotonically decreasing LD-distributions favoring low digits whenever log(f) values are below 0.1428 (1/7) approximately, and slowly approach the logarithmic as these values decreases even further. When log(f) fraction is very small, say 0.01 (a rational 1/100, called 1/M), its spread over (0, 1) is good enough, so that its LD distribution is quite near the logarithmic.

In order to generate a table of **all** these anomalous series up to log(f) = 0.01, we use the relationship f = (1+P/100) so that the argument above translates into:

**$LOG_{10}$(1+P/100) = 1/M** *(M an integer).* Now solving for P (the % growth) by first taking 10 to the power of both sides of equation we get:

$$(1+P/100) = 10^{1/M}$$
$$P/100 = 10^{1/M} - 1$$
$$P\% = 100*(10^{1/M} - 1) \quad M = 1, 2, 3, \text{etc.}$$

The table below uses the above relationship to generate those anomalous rates and comments on their LD distributions by varying integer M from 1 to 25 as well as evaluating it at M equals 35, 40, 50 and 100. Calculations using only the first 1000 elements were used, therefore it is important to note the slight dependency of exact resultant LD distributions on initial value of the series (the base).

| Integer M | The implied mantissa fraction log(1+P/100) | The implied P % growth | LD behavior of the related Exponential Growth | chi-sqr test (initial value=3) |
|---|---|---|---|---|
| 1 | 1 | 900.0000 | LD Neutral, a single digit leading all. | 7003.9 |
| 2 | 0.5 | 216.2278 | Only 2 digits leading, the rest never lead. | 6464.6 |
| 3 | 0.3333 | 115.4435 | Only 3 digits leading, the rest never lead. | 1918.1 |
| 4 | 0.25 | 77.8279 | Only 4 -5 digits leading, the rest never lead. | 1659.7 |
| 5 | 0.2 | 58.4893 | Only 4 -5 digits leading, the rest never lead. | 954.2 |
| 6 | 0.1667 | 46.7799 | Only 5 -6 digits leading, the rest never lead. | 785.1 |
| 7 | 0.1429 | 38.9495 | Only 5 -6 digits leading, the rest never lead. | 418.9 |
| 8 | 0.1250 | 33.3521 | Only around 6 -7 digits leading. | 518.2 |
| 9 | 0.1111 | 29.1550 | Only 7 digits leading. Fairly monotonic. | 416.9 |
| 10 | 0.1 | 25.8925 | Only 7 digits leading. Fairly monotonic. | 296.4 |
| 11 | 0.0909 | 23.2847 | 7 -8 digits leading, begins to resemble logarithmic. | 166.0 |
| 12 | 0.0833 | 21.1528 | resembling the logarithmic better. Definitely monotonic. | 122.5 |
| 13 | 0.0769 | 19.3777 | resembling the logarithmic better. Definitely monotonic. | 103.1 |
| 14 | 0.0714 | 17.8769 | resembling the logarithmic better. Definitely monotonic. | 129.0 |
| 15 | 0.0667 | 16.5914 | resembling the logarithmic better. Definitely monotonic. | 95.6 |
| 16 | 0.0625 | 15.4782 | More logarithmic in appearance. | 36.9 |
| 17 | 0.0588 | 14.5048 | More logarithmic in appearance. | 133.0 |
| 18 | 0.0556 | 13.6464 | More logarithmic in appearance. | 53.3 |
| 19 | 0.0526 | 12.8838 | More logarithmic in appearance. | 110.1 |
| 20 | 0.05 | 12.2018 | More logarithmic in appearance. | 20.3 |
| 21 | 0.0476 | 11.5884 | All digits are represented, closer to logarithmic. | 105.2 |
| 22 | 0.0455 | 11.0336 | All digits are represented, closer to logarithmic. | 26.5 |
| 23 | 0.0435 | 10.5295 | All digits are represented, closer to logarithmic. | 81.4 |
| 24 | 0.0417 | 10.0694 | Converging. | 33.2 |
| 25 | 0.04 | 9.6478 | Converging. | 51.9 |
| 35 | 0.0286 | 6.8000 | Converging. | 12.1 |
| 40 | 0.025 | 5.9254 | very close to the logarithmic for all starting points. | 24.3 |
| 50 | 0.02 | 4.7129 | very close to the logarithmic for all starting points. | 11.1 |
| 100 | 0.01 | 2.3293 | Nearly logarithmic. | 1.5 |

FIG 23: Some Singularities in LD Distributions of Exponential Growth Series (Basic Type)

Even slight deviations from those anomalous rates above (*moving back to normal rates so to speak*) resulted again in a near perfect logarithmic behavior. For example, for the rebellious rate **$100*(10^{(1/12)} - 1)\%$** series (approximately 21.152765862859%), an addition or subtraction of about 0.02% (i.e. about 21.1727% and 21.1327%) was enough to better its LD behavior and to arrive very close to the logarithmic. For rates that are extremely close to a given rebellious one (yet not exactly equal), logarithmic behavior will not be obvious unless thousands or millions of elements are considered in the series – depending on how close it is and the magnitude of the rebelliousness of the rate itself. In any case, at M=100, that is, for the fraction 0.01, or roughly 2.3292992% growth rate, LD becomes very nearly logarithmic and very little improvement is obtained by raising or lowering the rate slightly.

In contrast with this serous pitfall, multiplication processes of **random numbers** on the other hand (which can be thought of as non-constant exponential growth/decay) are nicely logarithmic. Here, we never stumble upon the peril encountered above, namely the problem of fractions whose multiples add exactly to unity on the mantissa axis, since it is very clear that repeated additions of totally different fractions result in covering the entire (0, 1) interval nicely, evenly and fully whenever there are plenty of such fractional accumulations.

Such an esoteric or indirect way of finding singularities in exponential series via their mantissa begs another more direct or straightforward explanation. The natural suspicion here is to lay blame on the LD neutrality of integral powers of ten in multiplication processes, and this is exactly what is found! For **all** the anomalous rates above, and **only** for those rates, there exists an **integer N** such that **$(1+P/100)^N = 10$**. This equality implies that applying the factor (1+P/100) N times yields the LD-neutral number of 10. For all such rebellious rates, the cumulative effect of N such multiplications is the LD-neutral factor of 10, and thus a particular digit continues to lead once again at the Nth sequence no matter. The emergence of these particular cumulative factors whose values are powers of ten is cyclical, with each particular rate having a particular period, namely N. Clearly, the larger the value of N the more diluted is the net effect on overall LD distribution of the series since a large value of N implies that this LD repetition happens less often, and this is nicely consistent with what was established before, namely the lessening in severity of the deviations from the logarithmic for lower growth rates of those rebellious series (whose values of N are much larger, implying much longer periods and diluted effects)!

Little reflection is needed to realize that N above (the Nth cumulative factor) is simply M, the number of whole multiples of the faction log(f) adding exactly to unity. That is, if exactly M whole multiples of the fraction log(f) is unity, then the Mth (Nth) cumulative factor is then 10, since returning to the very same mantissa means that LD situation hasn't changed and therefore it must be that the cumulative factor is 10 or perhaps another higher power of ten.

**General types** of anomalous exponential rates are found when exactly T whole multiples of the fraction log(f) add up exactly to any integer L, and not just to unity. For example, 5 whole multiples of 0.4 yield exactly the value of 2 units of distance spanning log/mantissa interval, hence its related 151.1886% growth series has a non-logarithmic LD distribution, its $5^{th}$ cumulative factor is 100, its $10^{th}$ cumulative factor is 10,000, and so forth. The general rule for any fraction log(f), a particular log/mantissa interval of length L, and T (thought of as being larger than L) as an integral multiple of some small interval, is as follow: whenever the fraction log(f) equals the **rational number L/T**, a non-logarithmic LD behavior results. That is, whenever $\mathbf{LOG_{10}(1+P/100) = L/T}$ and L & T being positive integers. Solving for P we get: $\mathbf{P\% = 100*(10^{L/T} - 1)}$. For this general type of anomalous rates, the direct explanation blaming multiple integral powers of ten for its rebellious LD behavior is found in $\mathbf{(1+P/100)^T=10^L}$. The larger the value of T the less severe LD deviation is observed. For T values over 50 (roughly), the behavior is quite nearly logarithmic regardless of what value L takes. For a given value of T, deviations are almost of the same magnitude for all values of L. Cumulative factors of integral powers of tens emerge cyclically with T as the period, and are of the forms $\mathbf{10^L, 10^{2L}, 10^{3L}}$ and so forth. The tables below lists some such anomalous rates:

[Note: chi-sqr tests were obtained with initial values of 3. Also,"***1st C.F. of P. 10***" means "The first cumulative factor of integral powers of ten", namely $10^L$.]

| Interval L | # of parts T | the fraction L/T | % Growth $100*(10^{L/T} - 1)$ | chi-sqr Test | 1st C.F.of P.10 $10^L$ |
|---|---|---|---|---|---|
| 2 | 5 | 0.4 | 151.1886 | 7003.9 | 100 |
| 2 | 7 | 0.2857 | 93.0698 | 584.4 | 100 |
| 2 | 9 | 0.2222 | 66.8101 | 416.9 | 100 |
| 2 | 13 | 0.1538 | 42.5103 | 89.1 | 100 |
| 2 | 17 | 0.1176 | 31.1134 | 92.2 | 100 |
| 2 | 23 | 0.0870 | 22.1677 | 80.7 | 100 |
| 2 | 25 | 0.08 | 20.2264 | 51.9 | 100 |
| 2 | 47 | 0.0426 | 10.2943 | 5.0 | 100 |
| 2 | 67 | 0.0299 | 7.1151 | 8.1 | 100 |
| 2 | 71 | 0.0282 | 6.7011 | 5.3 | 100 |
| 2 | 120 | 0.0167 | 3.9122 | 2.4 | 100 |
| 2 | 214 | 0.0093 | 2.1753 | 7.1 | 100 |
| 2 | 344 | 0.0058 | 1.3477 | 5.7 | 100 |
| 2 | 657 | 0.0030 | 0.7034 | 1.8 | 100 |

FIG 24: Other Anomalous Series (General Type) with L as a constant 2 and increasing T

| Interval L | # of parts T | the fraction L/T | % Growth $100*(10^{L/T} - 1)$ | chi-sqr Test | 1st C.F.of P.10 $10^L$ |
|---|---|---|---|---|---|
| 3 | 4 | 0.75 | 462.3413 | 7003.9 | 1000 |
| 3 | 5 | 0.6 | 298.1072 | 950.8 | 1000 |
| 3 | 7 | 0.4286 | 168.2696 | 418.9 | 1000 |
| 4 | 7 | 0.5714 | 272.7594 | 417.7 | 10000 |
| 6 | 7 | 0.8571 | 619.6857 | 410.3 | 1000000 |
| 4 | 9 | 0.4444 | 178.2559 | 324.6 | 10000 |
| 3 | 11 | 0.2727 | 87.3817 | 225.5 | 1000 |
| 4 | 11 | 0.3636 | 131.0130 | 225.5 | 10000 |
| 5 | 12 | 0.4167 | 161.0157 | 124.8 | 100000 |
| 4 | 13 | 0.3077 | 103.0918 | 80.8 | 10000 |
| 4 | 17 | 0.2353 | 71.9072 | 130.3 | 10000 |
| 3 | 25 | 0.12 | 31.8257 | 51.9 | 1000 |
| 4 | 25 | 0.16 | 44.5440 | 43.4 | 10000 |
| 6 | 25 | 0.24 | 73.7801 | 51.9 | 1000000 |
| 7 | 25 | 0.28 | 90.5461 | 51.9 | 10000000 |
| 8 | 25 | 0.32 | 108.9296 | 51.9 | 100000000 |
| 11 | 25 | 0.44 | 175.4229 | 43.2 | 1E+11 |
| 13 | 25 | 0.52 | 231.1311 | 53.5 | 1E+13 |
| 17 | 25 | 0.68 | 378.6301 | 53.1 | 1E+17 |
| 19 | 25 | 0.76 | 475.4399 | 57.0 | 1E+19 |
| 23 | 25 | 0.92 | 731.7638 | 52.7 | 1E+23 |
| 24 | 25 | 0.96 | 812.0108 | 49.7 | 1E+24 |
| 3 | 67 | 0.0448 | 10.8603 | 9.9 | 1000 |
| 6 | 67 | 0.0896 | 22.9001 | 8.2 | 1000000 |
| 20 | 67 | 0.2985 | 98.8417 | 7.6 | 1E+20 |
| 50 | 67 | 0.7463 | 457.5305 | 12.5 | 1E+50 |
| 20 | 123 | 0.1626 | 45.4125 | 0.6 | 1E+20 |
| 20 | 133 | 0.1504 | 41.3761 | 1.1 | 1E+20 |
| 3 | 344 | 0.0087 | 2.0284 | 1.8 | 1000 |
| 4 | 344 | 0.0116 | 2.7136 | 6.2 | 10000 |
| 6 | 344 | 0.0174 | 4.0979 | 1.6 | 1000000 |
| 20 | 344 | 0.0581 | 14.3246 | 5.4 | 1E+20 |
| 50 | 344 | 0.1453 | 39.7490 | 1.5 | 1E+50 |
| 20 | 345 | 0.0580 | 14.2802 | 6.1 | 1E+20 |
| 277 | 600 | 0.4617 | 189.5121 | 1.3 | 1E+277 |
| 20 | 1223 | 0.0164 | 3.8373 | 0.8 | 1E+20 |
| 277 | 3000 | 0.0923 | 23.6896 | 0.4 | 1E+277 |
| 500 | 11200 | 0.0446 | 10.8263 | 1.3 | == |
| 747 | 13577 | 0.0550 | 13.5062 | 0.1 | == |

FIG 25: Other Anomalous Series (General Type) with increasing T

We can confidently expect not merely most, but rather almost all exponential series to behave logarithmically. This is so for two reasons: **(I)** For a rate to be rebellious, the fraction log(1+P/100) must exactly equal some rational number L/T, and this is very rare, in fact the probability of obtaining a rational number when one is picked at random from any continuous set of real numbers is zero, although often financial, economics and other rates are quoted rationaly as fractions. Nevertheless, merely being (very) close to a rebellious rate is problematic and logarithmic behavior is dirupted unless some very large numbers of elements are to be considered - depending on the intensity of the rebellionness of the anomalous rate and on the closeness to it. **(II)** Even if a given rate is rebellious or close to one, there is a cap on deviation given by the value of T in the equation log(1+P/100) = L/T. As argued and seen above, whenever T is roughly over 100 deviation from the logarithmic is fairly small.

As an additional general check on LD behavior of exponential series seeking confirmation, a computer program was run. It started simulating exponential growth series from 1% all the way up to 600% in increments of 0.01% checking all the series for their LD behavior. The result was just as expected, even though it may sound puzzling or parodoxical to the social scientist: Any **deviant** DL-behavior differing from the logarithmic was always exclusively associated and correlated with the **rationality** of its related log(f) fraction. There seen to be nothing else adversely affecting LD behavior at all, except for the argument given in this section.

In the section on multiplication processes the reader was reminded that exponential growth is not free of units or scale, but rather depends on the length of the time used in defining the period, and that any high growth that is quoted using some fixed period of time (say a year) could in principle be quoted as lower growth if defined over a shorter period of time (say a month). Let us see how all this may be related to the anomalous rates above without any contradictions: **An annual rebellious rate implies also some rebelliousness on the part of its related (equivalent) monthly rate as well, albeit with reduced intensity.** Equivalency implies that $(1+YR/100)=(1+MT/100)^{12}$ hence $LOG_{10}(1+YR/100)=LOG_{10}(1+MT/100)^{12}$ or $LOG_{10}(1+YR/100)=12*LOG_{10}(1+MT/100)$. Now, if $L/T =LOG_{10}(1+YR/100)$ where L/T is rational, then $L/T =12*LOG_{10}(1+MT/100)$ or that $L/(12*T)=LOG_{10}(1+MT/100)$ which is also rational (12 is an integer!) hence rebellious as well monthly-wise. But deviation from the logarithmic for the monthly rate is less severe as it comes with a longer period. Put another way: the related cumulative factors of integral powers of ten for the monthly rate occurs with much less frequency (snapshots are taken 12 times more often than the annual series) and thus are much closer to the logarithmic than its related annual rate.

In general, any basic type series having $(1+P/100)^N = 10$ points to another basic type of a lower rate, simply by writing it as $(1+P/100)^{(1/R)*R*N}=10$ or $((1+P/100)^{1/R})^{*R*N}=10$, so that the Rth root of (1+P/100) points another basic type anomalous rate. Also, any basic type series having $(1+P/100)^N = 10$ points to another higher rate (general or basic) type series. We simply raise both sides of the above equation to the Lth (integral) power yielding $(1+P/100)^{NL} = 10^L$ or $((1+P/100)^L)^N = 10^L$ so that the Lth power of (1+P/100), namely $(1+P/100)^L$ points to another anomalous series, provided L is an integer. When L is larger than N it points to a general type. When L is smaller than N, the type depends of whether or not L is a proper factor of N (basic) or not (general). If L is a proper fraction of N, say QL=N, Q is an integer, $1 < Q < N$, then $((1+P/100)^L)^N = 10^L$ can be written as $((1+P/100)^L)^{QL} = 10^L$ and upon taking the Lth root of both sides we obtain $((1+P/100)^L)^Q = 10$ pointing to a basic anomalous series. Examples taken from the tables in this section are: $1.066050 = \sqrt[9]{(1.778279)}$, and $1.136464 = \sqrt{(1.291550)}$, Also, the basic type rate of 1.110336 points to $1.110336^2 = 1.23284$ (basic type, 2 is a factor of 22) and to $1.110336^6 = 1.87381$ (general type, 6 is not a factor of 22).

It is worth noting that exponential growth and exponential decay clothed in probabilistic curves of 'density' both point to the exact same curve, a curve with a one-sided tail to the right. It's only the time dimension that differentiates between the two, as it transverses different directions on the same curve. The same analysis and same results derived in this section applies to exponential decay. In summery, exponential growth and exponential decay, together with the purely random multiplication process should be considered as covering the entirlety of the mutlipliction processes phenemena and thus unanimous conclusions and results from these topics should be considered as good representative on the whole of multiplication processes study in Leading Digits.

## [22] Higher Order Significant Digits

The general law implies precise significant high orders digits distributions, now including the 0 digit as well. That also implies a dependency between the orders. Here is a table of the **unconditional** first 3 orders of distributions:

| DIGIT | 0 | 1 | 2 | 3 | 4 | 5 | 6 | 7 | 8 | 9 |
|---|---|---|---|---|---|---|---|---|---|---|
| 1st LD | 0.0% | 30.1% | 17.6% | 12.5% | 9.7% | 7.9% | 6.7% | 5.8% | 5.1% | 4.6% |
| 2nd LD | 12.0% | 11.4% | 10.9% | 10.4% | 10.0% | 9.7% | 9.3% | 9.0% | 8.8% | 8.5% |
| 3rd LD | 10.2% | 10.1% | 10.1% | 10.1% | 10.0% | 10.0% | 9.9% | 9.9% | 9.9% | 9.8% |

FIG 26: 1st, 2nd, and 3rd Order Leading Digits, **UN**conditional Probabilities

Third order LD distribution is nearly uniform, and still higher orders converge to the uniform fairly quickly.

If we assume that all our numbers are made of 4 digits, and if we approximate all fourth order significant digits probabilities as 10%, equal for all ten digits, then the overall frequency in the usage of the digits, regardless of order, can be easily calculated as the simple average of all four distributions given in the following table:

| DIGIT | 0 | 1 | 2 | 3 | 4 | 5 | 6 | 7 | 8 | 9 |
|---|---|---|---|---|---|---|---|---|---|---|
| Overall usage | 8.0% | 15.4% | 12.2% | 10.7% | 9.9% | 9.4% | 9.0% | 8.7% | 8.4% | 8.2% |

FIG 27: Digital Usage for 4-digit-long Numbers

If we assume that all our numbers are made of 7 digits, and if we approximate all fourth order and higher significant digits probabilities as 10%, equal for all ten digits, then the overall frequency in the usage of the digits, regardless of order, can be easily calculated as the simple average of all seven distributions given in the following table:

| DIGIT | 0 | 1 | 2 | 3 | 4 | 5 | 6 | 7 | 8 | 9 |
|---|---|---|---|---|---|---|---|---|---|---|
| Overall usage | 8.9% | 13.1% | 11.2% | 10.4% | 10.0% | 9.7% | 9.4% | 9.2% | 9.1% | 9.0% |

FIG 28: Digital Usage for 7-digit-long Numbers

Therefore, the effect of Benford's Law on overall keyboard usage or digital storage in computers say is not at all as dramatically skewed in favor of the low digits as is the case in the occurrences of the first digits where digit 1 is almost 7 times more likely to occur than digit 9.

## [23] Distribution Assuming Other Bases

Another generalization of Benford's law is LD distribution of other number systems with different bases. The probability of the first leading digit d (in the range 1 to B-1) for the base B is: P[d is first] = $LOG_{10}(1 + 1/d) / LOG_{10}(B)$, while other formulae such as those for higher order significant digits are carried over to other bases with the single adjustment factor of $1/LOG_{10}(B)$.

For binary numbers the probability of the first leading digit being 1 is 1.000 since the leading non-zero digit of a binary number is necessarily 1, and in which case Benford's law becomes a tautology.

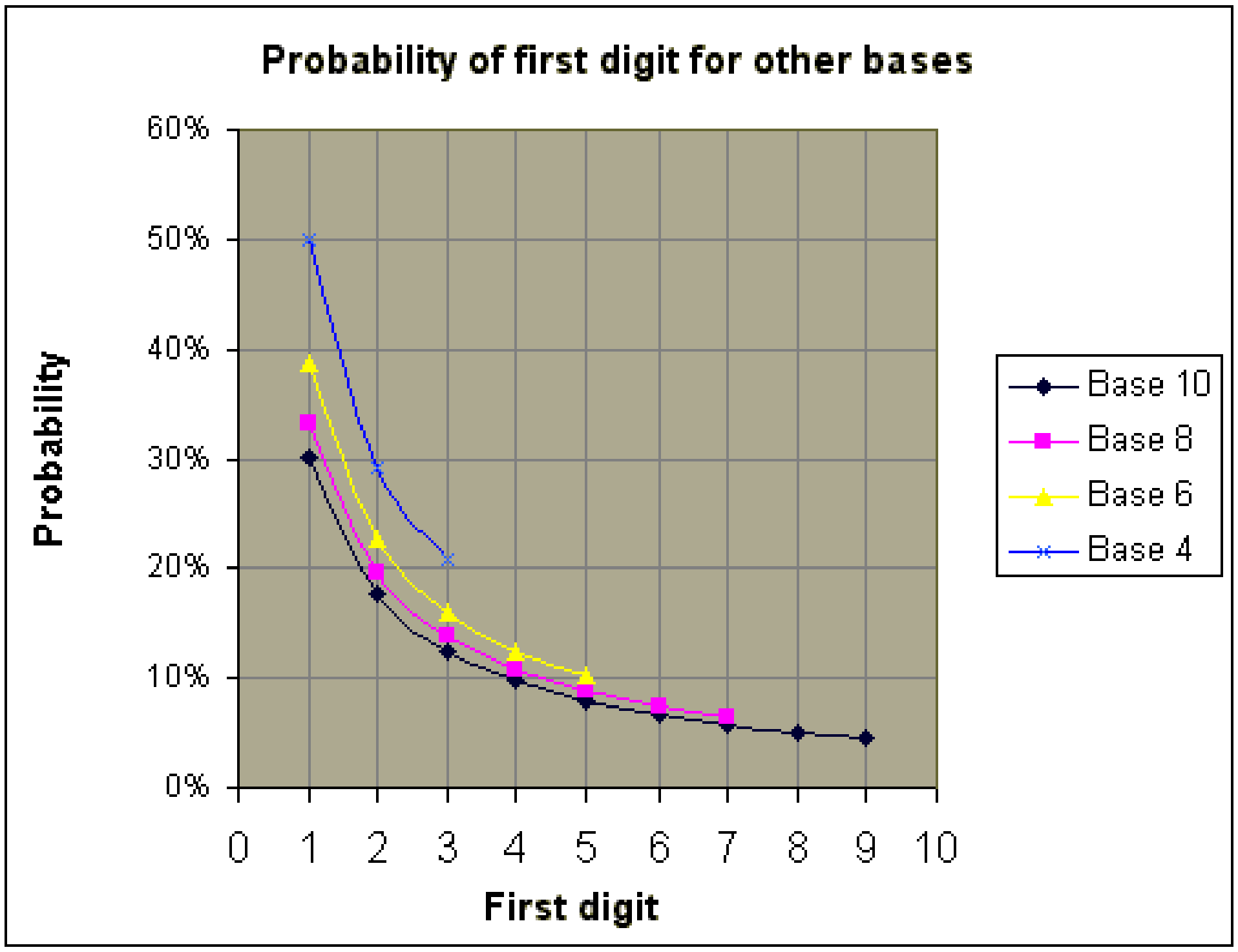

FIG 29: LD distributions for various number systems

## [24] Applications

**[I]** Accounting data is known to conform very closely to Benford's law, therefore when accountants or individuals fake such data it is possible at times to detect it. Such a method of detection is widely used at the Internal Revenue Service. Often, faked accounting data has equal LD distribution.

**[II]** Data analysts looking for data entered in error use Benford's law to spot its occurrences. There are other applications in data analysis.

**[III]** Theodore Hills writes [H-d]: "It is possible to build computers whose designs capitalize on knowing the distribution of numbers they will be manipulating. If 9's are much less frequent than 1's then it should be possible to construct computers which use that information to minimize storage space, or to maximize rate of output printed. A good analogy is a cash register drawer. If the frequency of transactions involving the various denominations of bills is known, then the drawer may be specially designed to take advantage of that fact by using bins of different sizes. In fact, German mathematician Peter Schatte has determined that based on an assumption of Benford input, the best computer design that minimizes expected storage space (among all computers with binary-power base) is base 8, and other researchers are currently exploring use of logarithmic computers to speed calculations."

**[IV]** When a mathematical model is proposed regarding a set of data that is known to follow Benford's law and when its task is to predict future values (such as stock prices or number of people infected in a pandemic) then a possible test disproving its validity outright can be performed by checking whether or not predicted values are still logarithmic. When model is theorizing about input-output scenarios, the test would involve what is called the "Benford-in-Benford-out" principle.

## [25] More on Chains of Distributions

Here is a brief summary of several other simulations of <u>short</u> chains:

Uniform**(***Normal***(**chi-sqr(die), Uniform(0,3)**),** *Normal***(**Uniform(77,518), Uniform(0,2)**)).**
10,000 simulations gave: {0.331, 0.213, 0.140, 0.081, 0.051, 0.049, 0.050, 0.044, 0.043}.

Normal(Normal(Normal(Normal(Normal(Normal(Uniform(0, 7), 13), 13), 13),13),13),13).
12,000 simulations gave: {0.250, 0.215, 0.161, 0.119, 0.089, 0.057, 0.043, 0.036, 0.026}.

Normal(*Uniform*(50,100), *Unifrom*(0,133)). 14,000 simulations gave overall LD distribution of: {0.270, 0.158, 0.127, 0.108, 0.085, 0.077, 0.064, 0.060, 0.048}.

Normal**(***Uniform***(**0,Unif(0,Unif(0,Unif(0,55)))**),** *Uniform***(**0,Unif(0,Unif(0,Unif(0,Unif(0,3))))**)** ).
12,000 simulations gave: {0.300, 0.175, 0.125, 0.100, 0.079, 0.066, 0.061, 0.048, 0.046}.

Normal(0, Uniform(0,3)). This is an amazingly short yet rapidly converging chain!
10,000 simulations gave:{0.320, 0.194, 0.121, 0.095, 0.071, 0.053, 0.056, 0.047, 0.041}.

Rayleigh(Uniform(0, Exponential(Rayleigh(Uniform(0, Exponential(Rayleigh… *etc.*
*9 full such cycles of 3 sequences, totaling 27-sequence chain.*
10,000 simulations gave:{0.301, 0.177, 0.127, 0.097, 0.080, 0.070, 0.054, 0.055, 0.039}.

Weibul**(** *Normal***(** Wald(2, 25), Weibull(0.1, 1)**)**
**,** *Uniform***(**0, Rayleigh(Rayleigh(Weibull(Uniform(0, 65), Normal(87, 5)))))**)** ).
10,000 simulations gave: {0.301, 0.172, 0.122, 0.095, 0.081, 0.073, 0.054, 0.055, 0.047}.

Uniform(0,U(0,U(0,U(0,U(0,U(0,U(0,U(0,U(0,17))))))))). **Flehinger's iterated scheme-like** in disguise! Except that LB is at 0 here while that of Flehinger was at 1, a crucial difference yielding very different results than her infinite schemes. Recall that this short yet potent distance from 0 to 1 contains an infinite number of intervals between IPOT (integral powers of ten) numbers, as well as an infinite numbers of related log negative integers! This chain has only finite 9 sequences, yet this is considered more than plenty here for a very good convergence. 14,000 simulations gave:
{0.303, 0.178, 0.129, 0.095, 0.080, 0.065, 0.057, 0.049, 0.045}.

Uniform(0, ***1/(x*ln10) dist over (10, 100)*** ). This is nothing but **Frank Benford's own attempt at a proof** of the logarithmic distribution in disguise! 12,000 simulations gave the logarithmic almost exactly. Hence, Benford's curious and unsuccessful attempt can now be viewed as simply one short yet powerful chain of distributions! Unfortunately this doesn't explain the phenomena as it would be very hard to argue that everyday data mimic this unique statistical process. More on F. Benford's attempt at an explanation - and the correct mathematical reasoning of why this short chain is actually logarithmic as it stands - in a later section below.

limit $N \rightarrow \infty$ Uniform(0, Uniform($10^N$, $10^{N+1}$)), N is an integer. This 'chain limit' is nothing but **Stigler's law** in disguise! The proper view here perhaps is the other way around, namely that Benford's attempt, Flehinger scheme, and Stigler's law are nothing but chains of distributions in disguise! The chain idea can be thought of as the common thread going through all these schemes.

Conceptually it is easy to envision or intuit how we get approximately a k/x-like-shaped distribution, or at a minimum something with a definite tail to the right. Consider the 2-sequence short chain **exponential( *Uniform(0, 1)* )**, pointing to a set of exponentials all starting at 0 (by default) and having progressively different focus on the (0, ∞) range. This chain is not long enough to converge yet it is quiet close to the logarithmic as it stands. Let us recall that exponential is defined over (0, +∞) and that the mean is 1/parameter. We can attempt to envision the overall shape of the chain by aggregating 10 typical realizations of exponential distributions of the chain having parameters nicely spread over (0, 1). We choose parametrical values of {0.1, 0.2, 0.3, 0.4, 0.5, 0.6, 0.7, 0.8, 0.9, 1.0}.

Since all individual distributions have a common lower bound (namely zero, as for all exponentials) while upper bounds vary (depending on parameter), the result is a definite fall in the density of the aggregate, with a tail to the right, and not just because of the fact that each exponential has similar such shape to begin with, because very similar results are obtained with **U(0, *U(0, 1)* )** for example, which comes with a flat curve for each uniform! This certainly reminds us of the dichotomy seen earlier between LB and UB concerning the averaging schemes. The chart of the combined or aggregate density curve for all these 10 exponentials illustrates the resultant near-logarithmic curve.

**The infinite chain conjecture assumes that the primary density and all intermediary distributions are defined on (0, + ∞)**, **(- ∞, 0)**, or **(- ∞, + ∞)**. Conceptually, it might be difficult to see why distributions drawing from both sides of the origin such as the Normal converge in simulations, but it DOES converge! A possible explanation perhaps is to think of a Normal sitting firmly on some positive territory mostly, and considering Chebyshev with 10 s.d. on either side we may forgo the rest of the tail. This easy explanation though is superfluous; as convergence to Benford is seem also when values are well-mixed with plenty of positive and negative numbers coming out simultaneously in simulation results.

The infinite chain conjecture can then be stated formally as follow: Consider the following set of distributions; all having any number of parameters; **all parameters are real, not discrete**; all parameters and all distributions are defined on such ranges as: (0, + ∞), (- ∞, 0), or (- ∞, + ∞), and where **all moments/means have finite values**, then, for chains made of such distributions (limited to those types of parameters discussed below) where ALL parameters are tied up to yet other distributions, and so forth, the primary distribution (the one sitting on top of the pyramid not serving as a parameter for others) is Benford in the limit as the number of sequences goes to infinity, and regardless of course of the exact values of the parameters at the very (infinitely deep) bottom of the chain (if anybody can see that far!)

Surprisingly, chaining distributions defined on **(a, + ∞)** or on **(- ∞, -a)** yield convergence to the logarithmic just as well, and regardless of whether **a** is an integral power of ten or not. Conceptually though one finds it hard to explain or intuit such convergence except in cases of distributions defined on **(+10 $^{integer}$, + ∞) or (- ∞, -10 $^{integer}$)**. But this behavior can be understood when distribution is not abruptly launched at **a** and then falls off steadily (the way the exponential is from 0 and k/x distribution from the left-most point) but rather gradually rises from **a,** and then falls, meaning that not a big portion of data is near **a** at all.

Chains inadvertently possess some sort of 'proper' direction in order for them to converge. For example, U(0, U(0, U(0,M))) taken further to infinity is Flehinger-like successful scheme (with lower bound stuck at 0 as oppose to 1), while U(U(U(A, B), C),D) taken to infinity as well, where D>C>B>A, is not logarithmic at all! The former fastens LB of the chain to 0 and expand UB outward, while the latter fastens UB of the chain to the constant D and contract LB gradually up to A. Certainly it is plausible to argue that for everyday real life datasets such as length of rivers, population of towns, quantities bought and sold and so forth, an extremely natural lower bound would be zero or one, and that distributions expand outward in the positive direction, much more like the former scheme and unlike the latter.

It is important to note that the 2-sequence uniform chain **U(L, U(L, M))** converges nicely to the logarithmic NOT whenever difference **M-L** is large, but rather whenever log difference is large! That is whenever **$LOG_{10}M$ - $LOG_{10}L$** is larger than roughly 8 or 9 depending on logarithmic accuracy desired. Hence **U(0, U(0, 1))** is by far superior to **U(1, U(1, 999))** say, as there are INFINITELY many IPOT values on the "short" interval between 0 and 1 !!! That is, 'log distance' between 0 and 1 is "infinitely log-long"!The chain **U(1, U(1, M))** needs to have its M roughy with the value of 10,000,000 to see a good agreement with the logarithmic!

Another 'chain paradox' here is the observation that **U(1, U(1, 10))** gave LD of {36.0%, 19.6%, 14.7%, 10.1%, 7.7%, 5.3%, 3.3%, 2.6%, 0.7%} with good spread from 1 to about 5 or almost 7, while **U(1, U(1, U(1, U(1, 10))))** gave LD as {81.0%, 11.6%, 4.3%, 1.7%, 1.0%, 0.3%, 0.1%, 0.0%, 0.0%} and most values were congregating just over 1 with the vast majority below 2, and without any meaningful spread. It is as if the longer the chain is the more values are being 'pushed over' to regress to 1! These two last examples show that we could get a worsening logarithmic result the longer the chain is being chained (the shorter 2-sequence chain had a better logarithmic result than the longer 4-sequence one!?)

But when fastened to zero the chain with more sequences (longer) is superior! That is, the longer chain of **U(0, U(0, U(0, U(0, 1))))** is closer to the logarithmic than the shorter chain of **U(0, U(0, 1))**.

Finally, the confluence of two observations here points to a second chain conjecture. First, some particular infinite chain is imagined and is being called CHN. Since at the end of the day CHN is nothing but some particular distribution (and logarithmic as per our infinite conjecture), it is now being used to represent the parameter for some single-parameter distribution called Z. Clearly, this new chain Z(CHN), having infinite +1 sequences, is logarithmic just the same. Yet, a different perspective here is to view Z(CHN) merely as a 2-sequence chain, and to proclaim that its LD logarithmic behavior springs from the fact that its parameter is being chained to a density that is logarithmic in its own right. A second similar insight can be inferred from Benford's own attempt in 1938 to establish the logarithmic distribution by letting UB vary exponentially. His semi-logarithmic scale of P vs. UB and his calculations of the area under this curve to arrive at its average height are mathematically equivalent to the use of exponential growth series as upper bounds. His scheme could then be interpreted as the chain Uniform(0, ***exponential series*** ). This short chain is noted for its special feature of inserting a distribution that is logarithmic in its own right, namely exponential series, to serve parameter **b** of the uniform distribution there. Admittedly, parameter **a** is being made fixed at 0, as opposed to being chained as well, but in this case it is useful **not** to chain **a**, only **b.**

**The 2nd conjecture then claims that any 2-sequence chain having ALL parameters of the primary density derived from logarithmic distributions is logarithmic there and then without any need to expand infinitely**.

Three chains that are closely related to Benford's own scheme have been simulated, and the results lend strong support to this 2nd conjecture.
**(A):** ***Unifrom(0, lognormal)*** is progressively more logarithmic as shape parameter is gradually increases and carefully observed all the way from 0.1 to 0.7, culminating in a near perfect such behavior of the chained uniform for shapes over 1.25.
**(B):** ***Unifrom(0, exponential)*** is roughly logarithmic, since any exponential is only approximately so.
**(C):** ***Unifrom(0, a chain of 15 Rayleighs)*** is nearly logarithmic. Yet, strictly speaking, none of these three examples support the conjecture completely, since none of these parametrical densities (lognormal, exponential, 15-sequence-Rayleigh chain) is perfectly Benford no matter what parameter is chosen, rather they are just very close to it. Yet, they certainly do suggest the principle. The fact that for all these three examples, parameter **a** is stuck at 0 does not imply any lack of support for the conjecture, because had we actually chained parameter **a** to a distribution it would still be necessary to have it bounded by LB<min[chain of b-parameter], since **a** must be less than **b**, and this would most likely point to 0 anyhow as the value for **a**.

Other simulations suggestive of the conjecture are:
**chi-sqr(integer of exponential growth series)** is nearly logarithmic regardless of rate and base, but provided that the exponential series is long enough;
**chi-sqr(integers of a lognormal)** is progressively more logarithmic as shape parameter increases from 0.1 to 0.7, culminating in near perfect logarithmic behavior for shape over 1.25;
**chi-sqr(integers of a chain of 13 uniforms - *all with parameter* a *stuck at 0*)** is nearly logarithmic. [Note that this is so in spite of our earlier prohibition against using discrete parameters such as the d.o.f. of the chi-sqr!!! Indicating that the conjectures are too timidly and too cautiously stated, that the chain is applicable in much more situations. Also, this chain above is logarithmic in spite of d.o.f. being essentially a non-chainable shape parameter of the gamma distribution as will be explained later in the other section.]

Here is a different set of simulations that are also suggestive of the second conjecture; the following 2-sequence chains here all came out quite nearly logarithmic, while their parameters are being derived from densities that are nearly or exactly Benford:

**Weibull(** lognormal shape>1.25, lognormal shape>1.25 **)** ;
**Weibull(** chain of uniforms, chain of uniforms **)** ;
**Weibull(** 1/(x*ln10*M) on $(10^{H}, 10^{H + \text{integer_}\mathbf{M}})$**,**
1/(x*ln10*N) on $(10^{G}, 10^{G + \text{integer_}\mathbf{N}})$ **)** ;

**Rayleigh(**chain of uniforms**)** ;
**Rayleigh(** lognormal with shape>1.25 **)** ;
**Rayleigh(**1/(x*ln10*N) over $(10^{F}, 10^{F + \text{integer_}\mathbf{N}})$ **)** ;

**Wald(**lognormal shape>1.25, lognormal shape>1.25**)** ;
**Wald(**chain of uniforms, chain of uniforms**)** ;
**Wald(**chain of Rayleighs, chain of Rayleighs**)** ;
**Wald(**1/(x*ln10*M) on $(10^{H}, 10^{H + \text{integer_}\mathbf{M}})$**,** 1/(x*ln10*N) on $(10^{G}, 10^{G + \text{integer_}\mathbf{N}})$**).**

The general form of the 2nd conjecture could be succinctly expressed as follows: **AnyDensity(AnyBenford) is Benford**, although this is true for a limited class of distributions/parameters to be discussed later on in the next section.

One could hardly argue that these are special cases, and that what drives the nearly perfect logarithmic behavior for these 2-sequence chains is not the fact that they have chained their parameter to a logarithmic distribution, but rather the particular density form used or the specific chain arrangement. Such claims are almost thoroughly refuted by using quite convincing demonstrations such as the case of the chain **Rayleigh(lognormal)** as but one example. Here is a table of such simulations with 5000 values per batch, having location parameter stuck at 3, while shape (the one parameter exclusively controlling logarithmic behavior) is allowed to vary:

| shape of lognormal | 0.1 | 0.2 | 0.3 | 0.4 | 0.5 | 0.6 | 0.7 | 0.8 | 0.9 | 1.0 |
|---|---|---|---|---|---|---|---|---|---|---|
| chi-sqr - lognormal | 6380.9 | 5810.5 | 4095.6 | 2467.4 | 1281.1 | 645.1 | 264.6 | 83.3 | 35.3 | 11.4 |
| chi-sqr - chain | 1027.9 | 716.7 | 517.6 | 216.7 | 141.9 | 70.8 | 31.3 | 15.7 | 14.1 | 10.1 |

FIG 30: Shape Parameter Manipulating Logarithmic Behavior of the Chain **Rayleigh(lognormal)**

For each value of the shape parameter of the lognormal serving as the parameter of the Rayleigh, we will give a measure of the degree of logarithmic behavior of the chained Rayleigh by a chi-sqr value, testing agreement with the logarithmic distribution, where a very low chi-sqr value indicates a good agreement. The result here shows a total dependency of the chain in being logarithmic on how far the lognormal itself is so! A perfect correlation! Here, the form of the distribution serving the parameter – namely the lognormal - and the whole setup of the chain is a **constant**, while the **only** factor that varies here is Benfordness of the lognormal, which clearly induces and controls Benfordness of the chained Rayleigh!

This strongly suggests in general that for the chain P(Q), it's the Benfordness of chainer Q that determines Benfordness of chainee P, or at a minimum, that Benfordness of chainee P is not indifferent to Benfordness of chainer Q. And this also strongly suggests a daring **extrapolation of the 2nd conjecture** by applying it (in total generality) also to **any** Benford or non-Benford densities serving a parameter, and asserting that Benfordness of any chain P(Q) is at a minimum equal to that of Q, but probably slightly higher, never less, and thus Benfordness is always conserved. In other words, that [chi-sqr of Q] ≥ [chi-sqr of P(Q)], so that there is never a retreat or backtracking in Benfordness in the act of chaining any two distributions.
We can mathematically define relative "**Benfordness**" not only using the chi-sqr test, but directly, in terms of how severely the digits are skewed, favoring low digits over high ones, and there are many ways of going about it. One dilemma would be what to name 'low digits', would they be 1 & 2, or perhaps 1 & 2 & 3, and so forth.

Another simulation in support of this last statement is displayed here with sequential chi-sqr values:

Uniform(0, Rayleigh(Rayleigh(Weibull(Uniform(0, 65), Normal(87, 5) ))))
**6.9** **91.1** **928.5** **28924.2** 1739.7 **42593.8**

What this 5-sequence chain shows us is that Benfordness is continuously being improved at each sequence in the chain, that there is never a pause or a retreat! [chi-sqr value of 1739.7 of the Uniform(0, 65) can be overlooked, it pertains to the shape parameter of the Weibull which does **not** play any rule in chaining, as will be shown later in the next section.]

Consider the logarithmic chain:

**gamma( infinite chain of uniforms with a at 0, infinite chain of exponentials) ;**
This could be viewed either as:
**(I)** simply an infinite chain, and thus logarithmic as per our 1st infinite conjecture.
**(II)** gamma(logarithmically distributed param, logarithmically distributed param), and thus logarithmic as per our 2nd conjecture. This flexibility in our point of view demonstrates how the second conjecture synthesizes and harmonizes nicely with the original infinite chain conjecture. But actually this is necessary indeed in order to hold the whole edifice of the chain idea intact, since contradictions would arise otherwise.

In general, assuming the **preposterous** claim that every logarithmic density can be expressed as or represented by an infinite chain, then the second 2-sequence conjecture can be directly derived from the infinite-sequence first conjecture. To formally prove this, simply write in the abstract **AnyDensity(logarithmic distribution #1, logarithmic distribution #2)** as **AnyDensity(infinite chain #1, infinite chain #2**) and hence itself an infinite (+1) chain as well, which then completes the assertion.

An even more harmonious relationship and mutual dependency between both the 1st and the 2nd conjectures can be found employing the following reasonable argument showing how the <u>extrapolated</u> 2nd conjecture mentioned earlier explains or contributes to the infinite-sequence conjecture by employing two assumed factors:
**[I]** We have widen the scope of the 2nd conjecture and asserted that Benfordness is always conserved in the act of chaining, that there is never a retreat or backtracking on Benfordness.
**[II]** Moreover, focusing on the micro level of an infinite-sequence chain, we postulate that there is an active mechanism in most sequential jumps resulting in a dynamic tendency for the outcome to end up with some 'improvement' in Benfordness. The confluence of these two factors above on an infinite chain is such as to insure a slow but certain drift towards complete and final Benfordness!

An immediate corollary of the 2-sequence conjecture is that any <u>finite</u> chain, however long, but having ALL its last (lower) parameters derived from logarithmic distributions is logarithmic. This is so via repeated applications of the second conjecture from the bottom up. [The symbol **B** to be used here would mean **B**enfordness, that is, a logarithmic distribution]. For example:
K**(**R(S(N(M(P(**B**))))),V(H(**B**, **B**))**)** is reduced to K**(**R(S(N(M(**B**)))),V(**B**)**)**, then to K**(**R(S(N(**B**))),**B)**, and so forth, to K(**B**, **B**), and finally to **B** . The extrapolated 2nd conjecture guarantees a quicker near-convergence of the original infinite conjecture for finite chains whenever all lowest (parametrical) densities are themselves closer to Benford to begin with, or at least having monotonically decreasing LD distributions resembling it. In other words:
If one's starting point is at a higher plateau then one reaches the summit sooner.

## [26] Chain-able Distributions/Parameters

An attempt to find a general rule describing what kind of parameters (or combination thereof) respond favorably to chaining was partially successful. Several general rules were found covering a wide range of distributions.

Our short 2-sequence chain conjecture offers an easy way of checking compliance. The setup for computer simulation is straightforward as parameters are tied either to **$1/(x*\ln10*N)$** over **$(10^F, 10^{F + \text{an integer } N})$**, or to the lognormal, turning the shape parameter there high and low as a knob and observing effects (the latter being a much more convincing demonstration). In contrast, approximating an infinite chain by using a finite number of chained distributions is cumbersome and leads to the constant dilemma of figuring out the correct number of necessary sequences for sufficient accuracy. Since both conjectures are intimately connected and in a way mirror each other, we shall assume for convenience that result for one implies the same for the other! In other words, that the rules governing convergence for the 1st conjecture are identical for the 2nd conjecture!

General principles:

**Scale** parameters such as **$\lambda X$** or **$X/\lambda$** (**divisions & multiplications**) in the absence of any location parameter always responds vigorously to chaining and yield the logarithmic. The requirement here is that each X is always consistently accompanied by such $\lambda$, and that there isn't any other $\lambda$ combining with X using any other arithmetic operation anywhere else in the PDF expression (nor reappearing alone as addition or subtraction, such as in $\lambda X$-$\lambda$, and so forth.)

**location** parameters such as **$X$-$\mu$** (**subtractions**) in the absence of any scale parameter responds vigorously to chaining and yield the logarithmic only where the range of the chained $\mu$ is high and hovers above a certain level. Low distributed values of the chained $\mu$ do not yield the logarithmic, yet it responds to chaining a great deal, approaching the logarithmic but not quite reaching it. This curious state of affairs is hard to explain!

For the case of a simultaneous presence of **location and scale** parameters everywhere in the PDF expression of the form **$(X-\mu)/\lambda$**, both parameters must be chained in order to obtain the logarithmic. Again, the requirement here is that each X is always accompanied by $\lambda$ and $\mu$ in the same arithmetic way everywhere in the PDF expression. Chaining only one parameter while leaving the other as a constant does not yield the logarithmic unless done in a way that yields much higher values (in comparison) for the chained parameter than the values of the other parameter being left a constant. Yet there is a sharp dichotomy here between $\lambda$ and $\mu$ when only one parameter is being chained; chaining only scale does not yield anything resembling the logarithmic unless done in a way that results in much higher values for scale than location. On the other hand, chaining only location always yield something quite close to the logarithmic, and regardless of the relative level in values between location and scale.

When the parameters above are of the form **$\lambda(X-\mu)$**, <u>both</u> parameters must be chained in order to obtain the logarithmic. Chaining only one parameter while leaving the other as a constant does not yield the logarithmic unless done in a way that yields much higher values (in comparison) for the location parameter than the values of the scale parameter. This draws a distinction between this case (where location has to be large no matter what is being chained) and the previous case where $\lambda$ appears as a denominator (where the parameter being chained has to be large). And just as in the previous case, chaining only location always yield something close the logarithmic, while chaining only scale does not have any effect unless location is much larger than scale.

**<u>Shape</u>** parameters such as $X^k$ (**powers**) do NOT yield the logarithmic when chained. Interestingly, chaining only the shape of the Weibull distribution does not lead to <u>any</u> kind of convergence whatsoever, while the shape of the gamma responds vigorously to chaining yet without that complete convergence. An extremely compelling explanation for this dichotomy is found in the expressions for both averages, $\lambda*\Gamma(1+1/k)$ for the Weibull and $\lambda*k$ for the gamma. For the gamma the average is responding steadily to any increase in k, while the average of the Weibull quickly becomes progressively more and more indifferent to further increase in k beyond a certain point. The general rule for shape parameters can be then stated as follow: whenever derivative **d(average)/d(parameter)** continuously falls and then diminishes in the limit (hence parameter is deemed irrelevant to the expression of centrality in the limit) no chainability can be found whatsoever, otherwise, it responds vigorously to chaining but without that complete convergence. It is noted that the median here proves a better measure of centrality than the average in predicting and explaining chainability behavior, as might be expected perhaps for being a much more robust measure of centrality in general.

The rough conjecture regarding a general principle here can then be stated as follow: **A parameter that does <u>not</u> continuously involve itself in the expression of centrality is <u>not</u> chainable at all, and does not show even feeble convergence under chaining. 'Continuously involved' means that d(center)/d(parameter) never diminishes as parameter gets large.** The converse is not always being exactly observed, as there are cases where parameters that do involve themselves in an additive expression of centrality turned out not to be chainable at all (suggesting perhaps that an additive expression in expressions of centrality is weaker than the multiplicative manner!?) Three such cases come to mind with the generalized exponential $1/\lambda*\exp(-(X-\mu)/\lambda)$, Fisher-Tippett and the Normal. Their medians are respectively: $1.000*\mu+ 0.693*\lambda$, $1.000*\mu+ 0.367*\lambda$ and $1.000*\mu+0.000*\sigma$, emphasizing (with a factor of 1) the location parameter $\mu$ which responds vigorously (but not completely) to chaining, and de-emphasizing (with factors less than 1) scale which is in fact totally indifferent to chaining. A more appropriate view to be taken here perhaps is to observe the simple fact that basically both are involved in the expression of centrality and hence both are needed to be chained for a full and complete convergence to be observed, not just one of them. Cases where the two parameters are involved in a multiplicative expression of centrality (such as μ*λ) seem to lend each single parameter sufficient weight as to be chainable or almost so in and by itself.

**WHY should only those parameters that are involved in the expression of centrality be chainable? Conceptually the answer is quite straightforward: because the mere act of chaining such parameters in any way at a minimum yields a set of distributions with increasing (stretching) upper bounds – derived from the fact that that parameter affects location (centrality) a great deal - and hopefully with lower bounds staying fixed near 0 or 1 or any other low IPOT number, and as was seen with our simple averaging schemes earlier, this results in an overall density that is diminishing, having that one-sided tail to the right so typical of logarithmic distribution! Nothing of the sort could happen if parameter to be varied by chaining does not affect centrality in the least!**

Let us recall the short chain in figures 16 and 17, **exponential( *Uniform(0, 1)* ).** There, the parameter being chained is ρ scale, which is also the expression for the average, and varying ρ there via chaining means that we are directly stretching upper bounds, while lower bounds are fixed at 0 by default, that is, by virtue of it being the exponential. Visualization of this process supports the conceptual argument given above.

In cases where the scale parameter λ and X are of different dimensions, no chainability can be obtained. Hence X/log(λ), X/ln(λ), X/√λ, X/(Nth root of λ) in the PDF expression does not lent itself to chainability. On the other hand, powers of λ in the PDF such as $X/\lambda^2$, $X/\lambda^{3.5}$, etc. are chainable given that power is 1 or higher. This is so because any power transformation (of 1 or larger) of any logarithmic data or variable (or parameter!) is also logarithmic, hence chaining λ to a logarithmic distribution implies that $\lambda^2$ or $\lambda^{3.5}$ are also logarithmic and thus could be considered as a whole the actual parameter instead of λ and as such having the same dimension as X. On the other hand, log, square root, or any Nth root transformations do not result in any logarithmic inheritance hence un-chainable. Same argument applies to log(X)/λ ,√X/λ, Nth root(X)/λ , where X is of a different dimension than that of λ.

Here are two tables summarizing many simulation results. The designation "**BEN**" indicates an outright Benford behavior for the 2-sequence chain, namely chainability, "**NOT**" in capital letters indicates a total indifference to chaining, not even having any sort of improvement towards the logarithmic in LD distributions, while "**not**" in lower case letters indicates a lack of full logarithmic behavior yet having a definite and vigorous response to chaining resulting in a considerable movement towards the logarithmic distribution:

| DISTRIBUTION | PDF EXPRESSION | PARAM1 | PARAM2 | BOTH PARAM | MEAN |
|---|---|---|---|---|---|
| Weibull $(0, +\infty)$ | $(k/\lambda)(x/\lambda)^{(k-1)}e^{-(x/\lambda)^k}$ | λ BEN | k NOT | λ & k BEN | λ*Γ(1+1/k) |
| Rayleigh $(0, +\infty)$ | $\frac{x\exp\left(\frac{-x^2}{2\sigma^2}\right)}{\sigma^2}$ | BEN | | | $\sigma\sqrt{\frac{\pi}{2}}$ |
| Exp 1 $(0, +\infty)$ | ρ*exp(-ρx) | BEN | | | 1/ρ |
| Exp 2 $(0, +\infty)$ | (1/ρ)*exp(-x/ρ) | BEN | | | ρ |
| Exp 3 $(0, +\infty)$ | (1/√ρ)*exp(-x/√ρ) | not | | | √ρ |
| Exp 4 $(0, +\infty)$ | $(1/\rho^{7.5})*\exp(-x/\rho^{7.5})$ | BEN | | | $\rho^{7.5}$ |
| Exp 5 $(0, +\infty)$ | $(1/\rho^{8})*\exp(-x/\rho^{8})$ | BEN | | | $\rho^{8}$ |
| Exp 6 $(0, +\infty)$ | (1/log(ρ))*exp(-x/log(ρ)) | not | | | log(ρ) |
| Gamma $(0, +\infty)$ | $x^{k-1}\frac{\exp(-x/\theta)}{\Gamma(k)\,\theta^k}$ | θ BEN | k not | θ & k BEN | θk |
| Wald $(0, +\infty)$ | $\left[\frac{\lambda}{2\pi x^3}\right]^{1/2}\exp\frac{-\lambda(x-\mu)^2}{2\mu^2 x}$ | μ BEN | λ NOT | μ & λ BEN | μ |
| LogNormal $(0, +\infty)$ | $\frac{1}{x\sigma\sqrt{2\pi}}\exp\left(-\frac{[\ln(x)-\mu]^2}{2\sigma^2}\right)$ | σ NOT | μ NOT | σ & μ NOT | $\exp(\mu+\sigma^2/2)$ |
| Uniform a = 0 (0, b) or $(0, +\infty)$ | 1/b | BEN | | | b/2 |
| Gompertz $(0, +\infty)$ | $be^{-bx}e^{-\eta e^{-bx}}\left[1+\eta\left(1-e^{-bx}\right)\right]$ | B BEN | η noT | b & η BEN | long & complex expression |
| Nakagami $(0, +\infty)$ | $\frac{2\mu^\mu}{\Gamma(\mu)\omega^\mu}x^{2\mu-1}\exp\left(-\frac{\mu}{\omega}x^2\right)$ | ω not | μ NOT | ω & μ not | $\frac{\Gamma(\mu+\frac{1}{2})}{\Gamma(\mu)}\left(\frac{\omega}{\mu}\right)^{1/2}$ |
| Gupta Kundu $(0, +\infty)$ | $\alpha\lambda*\exp(-\lambda x)*(1-\exp(-\lambda x))^{(\alpha-1)}$ | λ BEN | α NOT | λ & α BEN | 1/λ*[ψ(α+1) - ψ(1)] ψ being the digamma |

FIG 31: Chainability of some distributions defined on (0, +∞)

| DISTRIBUTION | PDF EXPRESSION | PARAM 1 | PARAM 2 | BOTH PARAM | MEAN |
|---|---|---|---|---|---|
| Normal $(-\infty, +\infty)$ | $\frac{1}{\sigma\sqrt{2\pi}} \exp\left(-\frac{(x-\mu)^2}{2\sigma^2}\right)$ | $\sigma$ NOT | $\mu$ not | $\sigma$ & $\mu$ BEN | $\mu$ |
| Origin-centered Normal $(-\infty, +\infty)$ | $1/\sigma\sqrt{(2\pi)} * \exp(-x^2/2\sigma^2)$ | BEN | | | 0 |
| Fisher-Tippett $(-\infty, +\infty)$ | $(1/b)*e^{-(x-\mu)/\lambda} * \exp(-e^{-(x-\mu)/\lambda})$ | $\lambda$ NOT | $\mu$ not | $\lambda$ & $\mu$ BEN | $\mu + 0.57721*\lambda$ |
| Logistic $(-\infty, +\infty)$ | $\frac{e^{-(x-\mu)/s}}{s\left(1+e^{-(x-\mu)/s}\right)^2}$ | s NOT | $\mu$ not | s & $\mu$ BEN | $\mu$ |
| Cauchy-Lorentz $(-\infty, +\infty)$ | $\frac{1}{\pi\gamma\left[1+\left(\frac{x-x_0}{\gamma}\right)^2\right]}$ | $\gamma$ NOT | $X_D$ BEN- -almost | $\gamma$&$X_D$ BEN | infinite and hence not part of our conjecture |
| Pareto $[a, +\infty)$ | $(\theta/a)(x/a)^{-(\theta+1)}$ | a BEN | $\theta$ NOT | a & $\theta$ BEN | $a*(\theta/(\theta-1))$ if $\theta >= 1$<br>$\infty$ if $\theta < 1$ |
| Generalized exp1 $[\mu, +\infty)$ | $\rho*\exp(-\rho*(x-\mu))$ | $\rho$ NOT | $\mu$ not | $\rho$ & $\mu$ BEN | $\mu + 1/\rho$ |
| Generalized exp2 $[\mu, +\infty)$ | $(1/\rho)*\exp(-(x-\mu)/\rho)$ | $\rho$ NOT | $\mu$ not | $\rho$ & $\mu$ BEN | $\mu + \rho$ |

FIG 32: Chainability of some distributions defined on $(-\infty, +\infty)$ & $(\mu, +\infty)$

**Comments:**

For the single parameter distributions of the exponential and uniform on (0, b), the parameter is explicitly present in the expression for the center, hence full chainability.

The origin-center Normal ($\mu = 0$) should be viewed as the combination of two separate distributions, one on $(-\infty, 0)$ and the other on $(0, \infty)$, each with center expressed explicitly in terms of s.d. $\sigma$, hence chainability for each, and therefore chainability for the combination. This result is in spite of the fact that the overall mean here is a constant zero no matter what s.d happens to be, because the correct view to take here is that each side yields its own mean as a function of the s.d. yet they are both of opposites sign (but same absolute value) and only the aggregate mean is zero.

For the Wald, the expression for the center involves only location parameter $\mu$, omitting scale $\lambda$ altogether, hence $\mu$ is fully chainable while $\lambda$ is totally indifferent to chaining.

The case of the lognormal may at first appear a bit puzzling because there neither parameter is chainable in any way (nor when both are chained simultaneously) even though both are decisively involved in expressing the mean, and at least $\mu$ is also involved in expressing the median. But the absence of any chainability in the case of the lognormal springs from that severe dimensional mismatch there (between X and $\lambda$, and between X and $\mu$), as it uses **$(\ln(X)-\mu)/\lambda$** in its PDF expression.

The Nakagami suffers from a lack of dimensional compatibility between X and $\omega$, hence not chainable.

These last two distributions the lognormal and the Nakagami are the only ones found where even though ALL parameters are chained, it's futile and no convergence is seen.

For the Gupta Kundu distribution; shape parameter $\alpha$ is present in the expression for the mean $(1/\lambda)[\psi(\alpha+1) - \psi(1)]$, however it quickly converges to a number slightly less than 6 and so its d(mean)/d(parameter) quickly diminishes hence its total lack of chainability. Same with the expression of the median there, which is: $(-1/\lambda)*\ln[1- 1/(\alpha$ th root of 2)].

The mean and median for the Pareto are $a*\theta/(\theta-1)$ and $a*(\theta$ th root of 2) respectively whenever $\theta \geq 1$. Hence $\theta$ becomes progressively less relevant to the expression of centrality and as a consequence is not chainable at all. Parameter **a** on the other hand leads decisively to Benford under chaining as expected.

The center of the Normal is $\mu$, hence chaining $\mu$ alone is almost Benford, while chaining only $\sigma$ is not even close (unless $\sigma >> \mu$ , a situation where $\sigma$ now control overall location and spread on the x axis much more than $\mu$!). But the overall correct view is to require both parameters $\sigma$ and $\mu$ to be chained regardless of which is involved in the expression of centrality and which is not.

The first and second chain conjectures were given a rigorous proof by Steven Miller covering a wide range of many classes of classical distributions/parameters. Simulations clearly point to further applicability, much more than his relatively restricted cases.

Note that for the chain scheme, particular or arbitrary values used to conduct them (the 'last' parameter for the independent distribution at the 'bottom') become less relevant as we consider longer and longer chains. An infinite chain would completely and thoroughly remove any supposed dependency on the 'last' value whatsoever, while perfectly converging to the logarithmic at the same time.

## [27] LD-Invariance of Parametrical Transformation by Powers of Ten

For distributions with a single scale parameter of the forms $\lambda*f(X\lambda)$ and $(1/\lambda)*f(X/\lambda)$ in their p.d.f. such as the exponential for example, transformation of the single parameter via multiplication by an integral power of 10 (such as 10, 1000, 0.01, 1, and so forth) leaves LD distribution unaffected regardless whether distribution is logarithmic or not, and this fact has nothing to do with Benford's Law. Also note that all this can be extrapolated to any base. Here while operating under base 10, we use an integral powers of the base as a factor transforming the parameters, likewise for any other base.
For distributions with scale and location parameters of the form $(1/\lambda)*f((X-\mu)/\lambda)$ such as the Normal, Fisher-Tippett, Logistic, Cauchy–Lorentz and others defined over $(-\infty, +\infty)$, $(0,+\infty)$, or $(-\infty, 0)$, only a simultaneous transformation of both parameters using the same integral power of 10 yields such invariance. **Yet**, those of the form $(\lambda)*f(\lambda*(X-\mu))$ are **not** invariant! Theoretically, related distributions defined over $(\mu, +\infty)$, such as the generalized exponential $(1/\lambda)*\exp(-(X-\mu)/\lambda)$ say, could not be covered under the proof to be given below, as it depends on limits of integrations being 0 or $\infty$ and thus unaffected by multiplication (as in $0*R = 0$ and $\infty*R = \infty$). However, computer simulations confirm such invariance also for those defined over $(\mu, +\infty)$, and regardless of whether or not $\mu$ itself is an integral power of 10, hence this needs to be investigated further.

The quest for a general rule stating which distributions/parameters respond to chaining was roughly fulfilled, pointing to the expression of centrality as the predictor. Curiously, it is interesting to note that invariance and chainability seem to be connected in some ways. It is noted that for single-parameter distributions, invariance implies chainability. Also, for two-parameter distributions: **I)** if only one parameter is invariant, then chaining that invariant one guarantees the logarithmic, while the other non-invariant parameter is not chain-abe ; **II)** in the case where neither parameter is invariant alone, yet together (simultaneously) they are invariant, both are needed to be chained in order to obtain the logarithmic, while chaining only one is not sufficient, and this is quite typical, for example in cases such as the Normal, Fisher-Tippett, Logistic, generalized exponential, and Cauchy-Lorentz; **III)** in cases where invariance is not observed, not even when both are simultaneously transformed by the same factor, no chainability is found at all – not even if both parameters are chained simultaneously (such as in the Nakagami and the lognormal cases, suffering from severe dimensional mismatch).

Yet, it is not the case that chainability and invariance are truly correlated. Two counter examples preclude such a possibility. The first is the generalized exponential of the form (λ)*f(λ*(X-μ)), where no invariance is found whatsoever, not even when both are simultaneously transformed, yet chainability is found when both are being chained. The second counter example is the Wald distribution, where μ - the parameter that in some sense resembles location but which also pops out again in the PDF expression as something else - is not invariant, yet it is decisively chainable in its own right. Hence, at a minimum, the following weaker conjecture is asserted:

**Whenever LD-invariance occurs for the transformation of a given set of parameters either simultaneously or for each one individually, Benford is observed whenever all parameters of the said set are chained, that is, chainability.**

The converse is not true, **chainability does not imply invariance.**

Now, a formal proof for some of those invariance cases mentioned above follows (regarding 1st digit law only):

Given any proper density function pdf(x) = f(x) defined on (0, +∞), (-∞, 0) or (-∞, +∞) where all limits of integration for the entire range are either 0, +∞, or -∞. Let us examine how f(Qx) could serve also as a proper pdf where total area sums up to 1 by evaluating its entire area. We shall start with the (0, +∞) case, but the same results are easily obtained for the two other cases:

$$\int_0^{+\infty} f(Qx)dx$$

A change of variables $u = Qx$ and $du/dx = Q$ yields:

$$\int_{0*Q}^{+\infty*Q} f(u)\frac{du}{Q} = \frac{1}{Q}\int_0^{+\infty} f(u)du = \frac{1}{Q}*1 = \frac{1}{Q}$$

Hence Q*f(Qx) is the **only** proper probability density function of this form as its entire area sums up exactly to 1. Thus the general form of all such densities is

pdf(x)= k*f(kx).

Assume density function is of the form pdf(x)= $k_0$f($k_0$x) defined on the range (0, +∞ ): Prob( 1st digit = d | $k_0$ ) =

$$\sum_{int=-\infty}^{int=+\infty} \int_{(d)*10^{int}}^{(d+1)*10^{int}} k_0 f(k_0 x)dx$$

Where index **int** runs through **all** the integers, positive and negatives one, including zero.

Let us now examine 1st leading digits distribution when the parameter is transformed by an IPOT number, that is pdf(x)= $10^M k_0$*f($10^M k_0$x) where M is any integer:

$$\sum_{int=-\infty}^{int=+\infty} \int_{(d)*10^{int}}^{(d+1)*10^{int}} 10^M k_0 f(10^M k_0 x)dx$$

A Change of variables $u = x * 10^M$, du/dx = $10^M$ yields:

$$\sum_{int=-\infty}^{int=+\infty} \int_{(d)*10^{int}*10^M}^{(d+1)*10^{int}*10^M} 10^M k_0 f(k_0 u)\frac{du}{10^M}$$

or simply :

$$\sum_{int=-\infty}^{int=+\infty} \int_{(d)*10^{(int+M)}}^{(d+1)*10^{(int+M)}} k_0 f(k_0 u)du$$

Now, since index int goes over all possible integers to begin with, the introduction of some fixed M integer as addition into integration limits does not have any effect and thus can be omitted.

We are left with:

$$\sum_{int=-\infty}^{int=+\infty} \int_{(d)*10^{int}}^{(d+1)*10^{int}} k_0 f(k_0 u) du$$

Which yields the same (original) expression for the 1st leading digits probability of $k_0 f(k_0 x)$ and which completes the proof.

The same argument can be made for $f\left(\frac{x}{Q}\right)$ type of densities defined on (0, +∞), (-∞, 0) or (-∞, +∞). There, the general form is $\frac{1}{k} f\left(\frac{x}{k}\right)$.

Given any proper density function $pdf(x) = f\left(\frac{x-a}{b}\right)$ defined on (-∞, +∞). Let us examine how $f\left(\frac{x-Pa}{Qb}\right)$ could serve as a proper p.d.f. as well by evaluating its entire area:

$$\int_{-\infty}^{+\infty} f\left(\frac{x - Pa}{Qb}\right) dx$$

A change of variables $u = x/Q$ , or equivalently $x = Qu$, and $dx/du = Q$ yields:

$$\int_{\frac{-\infty}{Q}}^{\frac{+\infty}{Q}} f\left(\frac{Qu - \frac{PaQ}{Q}}{Qb}\right) * Q * du = Q * \int_{-\infty}^{+\infty} f\left(\frac{u - \frac{Pa}{Q}}{b}\right) du$$

Another change of variables $w = u - \frac{P}{Q}a + a$ , with $\frac{dw}{du} = 1$ or $dw = du$ yields:

$$Q * \int_{-\infty - \frac{P}{Q}a + a}^{+\infty - \frac{P}{Q}a + a} f\left(\frac{\mathrm{w} - \mathrm{a}}{\mathrm{b}}\right) \mathrm{dw} = Q * \int_{-\infty}^{+\infty} f\left(\frac{\mathrm{w} - \mathrm{a}}{\mathrm{b}}\right) \mathrm{dw}$$

Since the integral with variable w is equivalent to our original p.d.f its value is 1, we get:

$$= Q * 1 = Q$$

Hence $pdf = \frac{1}{\mathrm{Q}} * f\left(\frac{\mathrm{x} - \mathrm{Pa}}{\mathrm{Qb}}\right)$ is the only proper probability density function of this form as its entire area sums up exactly to 1. Thus the general form of all such densities is $\left(\frac{1}{b}\right) * f\left(\frac{\mathrm{x} - a}{b}\right)$, meaning that only **b** is needed to be introduced as some factor, not **a**. Note: a change of variables above focusing on parameter **a** as opposed to **b**, such as $\mathrm{u} = \mathrm{x}/\mathrm{P}$ could not lead to any such cancelation in the quotient, hence the dichotomy here between a and b.

Assume a density function defined on (-∞, +∞) of the form pdf(x) = $= \left(\frac{1}{b}\right) * f\left(\frac{\mathrm{x} - a}{b}\right)$ where each x appears as $\left(\frac{\mathrm{x} - a}{b}\right)$ everywhere in the pdf expression.

We would divide the entire interval of (-∞ , +∞ ) into (0, +∞ ) & (-∞ , 0) and in that order for the next definite integrals evaluations.

Prob( 1st digit = d | $a_0, b_0$) =

$$\sum_{int=-\infty}^{int=+\infty} \int_{(d)*10^{int}}^{(d+1)*10^{int}} \left(\frac{1}{b_0}\right) * f\left(\frac{x - a_0}{b_0}\right) dx \ +$$

$$\sum_{int=-\infty}^{int=+\infty} \int_{-(d+1)*10^{int}}^{-(d)*10^{int}} \left(\frac{1}{b_0}\right) * f\left(\frac{x - a_0}{b_0}\right) dx$$

Where **int** runs through all the integers.

Let us now examine 1st LD distribution for another density of the same type where both $a_0 \ and \ b_0$ are being simultaneously multiplied by the same factor $10^M$, M being any integer.

$$\sum_{int=-\infty}^{int=+\infty} \int_{(d)*10^{int}}^{(d+1)*10^{int}} \left(\frac{1}{10^M b_0}\right) * f\left(\frac{x - 10^M a_0}{10^M b_0}\right) dx \ +$$

$$\sum_{int=-\infty}^{int=+\infty} \int_{-(d+1)*10^{int}}^{-(d)*10^{int}} \left(\frac{1}{10^M b_0}\right) * f\left(\frac{x - 10^M a_0}{10^M b_0}\right) dx$$

A change of variables $u = x/10^M$ or equivalently $x = u10^M$, and dx/du = $10^M$ yields:

$$\sum_{int=-\infty}^{int=+\infty} \int_{(d)*10^{int}/10^M}^{(d+1)*10^{int}/10^M} \left(\frac{1}{10^M b_0}\right) f\left(\frac{u10^M - 10^M a_0}{10^M b_0}\right) 10^M du \ +$$

$$\sum_{int=-\infty}^{int=+\infty} \int_{-(d+1)*10^{int}/10^M}^{-(d)*10^{int}/10^M} \left(\frac{1}{10^M b_0}\right) f\left(\frac{u10^M - 10^M a_0}{10^M b_0}\right) 10^M du$$

$10^M$ nicely cancels out twice in the pfd expression. Also, since index **int** goes over all possible integers including negative ones, the introduction of $10^M$ into integration limits does not have any effect and thus can be omitted. We are left with:

$$\sum_{int=-\infty}^{int=+\infty} \int_{(d)*10^{int}}^{(d+1)*10^{int}} \left(\frac{1}{b_0}\right) * f\left(\frac{u-a_0}{b_0}\right) du \ +$$

$$\sum_{int=-\infty}^{int=+\infty} \int_{-(d+1)*10^{int}}^{-(d)*10^{int}} \left(\frac{1}{b_0}\right) * f\left(\frac{u-a_0}{b_0}\right) du$$

Which yields the same expression for the 1 st leading digits probability of the original pdf = $\left(\frac{1}{b_0}\right) * f\left(\frac{x-a_0}{b_0}\right)$ and which completes the proof.

The same argument can NOT be made for densities of the form b*f(b(x-a)), and no such invariance can be found.

It would be interesting to enlarge this proof to all higher order leading digits, that is, to show invariance under the generalized Benford statement.

## [28] Sample of Typical Data from a Variety of Internet Websites

As a crude test for the LogNormal-like conjecture presented in this article, a large sample of real everyday data was collected from The World Wide Web (Internet). Results mostly confirmed the conjecture that the totality of everyday data mimics the LogNormal distribution.

The data was collected in the spirit of Ted Hill's distribution of all distributions idea. Effort was put into obtaining only a relatively small number of values from each particular story or topic, while varying the topics as much as possible. In total 34,269 values were obtained from roughly 70 different sources or topics. Since some topics were related or similar in nature, this last value of 70 is realistically more like roughly 50 to 60. Year numbers such as 2001, 1976 etc. were omitted deliberately as their frequency seemed exaggerated. Numbers representing codes were also omitted of course. Negative numbers surprisingly were few and far between, representing less than 0.5% of overall data, and were also omitted.

Most of the sources were from the following 6 websites:

U.S. Department of Agriculture, Economic Research Service,
http://www.ers.usda.gov/ Topics included:

Animal Products
Countries & Regions
Crops
Diet, Health, & Safety
Farm Economy
Farm Practices & Management
Food & Nutrition Assistance
Food Sector
Natural Resources & Environment
Policy Topics
Research & Productivity
Rural Economy
Trade & International Markets

U.S Census Bureau. Mostly about population statistics as well as Employment, Health, and Business Dynamics statistics. http://www.census.gov/

Statistics Worldwide. Mostly regarding Business, Governmental, and some General International statistics. http://statisticsworldwide.com

The Scottish Government Data Center. Having a large variety of many different issues and topics. http://www.scotland.gov.uk/Disclaimers

National Collegiate Athletic Association, relating to Sport statistics.
http://www.ncaa.org/

The World Bank. Information regarding Health, Population, Nutrition, Life Expectancy, and Economics statistics. http://data.worldbank.org/topic/health

Leading digits came out very close the logarithmic. Here is the 1st order distribution: {0.288, 0.164, 0.124, 0.098, 0.083, 0.073, 0.061, 0.057, 0.053}.

Here is the histogram of related log of data (decimal/common). It strongly suggests that the totality of our everyday data could perhaps be simply the LogNormal, or that it is best approximated by a LogNormal-like distribution as was conjecture earlier. Notice the relatively huge range on the log(x)-axis here (about 12 units) as compared with roughly a modest 3-unit length encountered in the earlier simulations of related log densities. There, 2.5 to 3 units was the minimum range needed in order for the log density to bequeath approximately the logarithmic property to data.

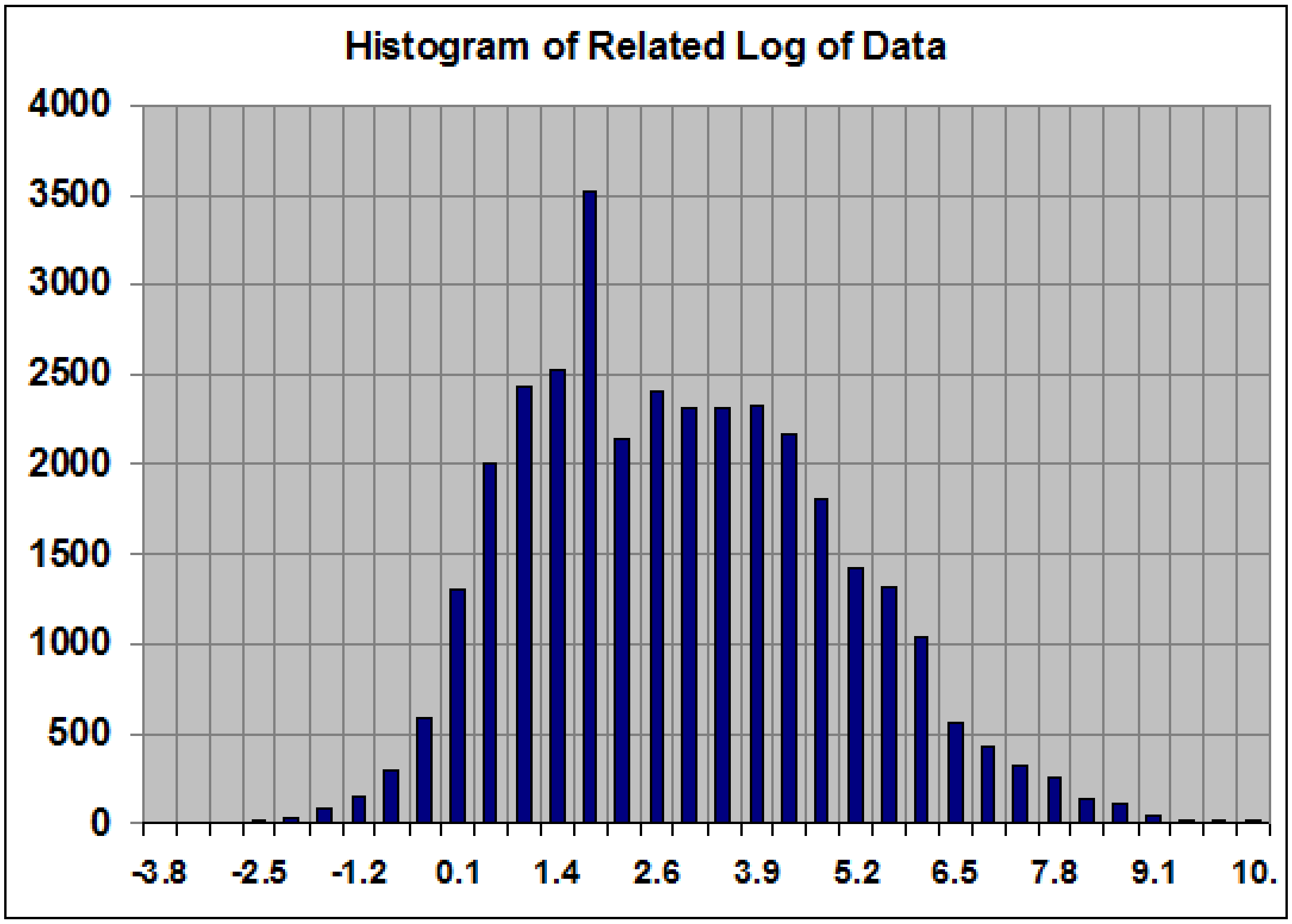

FIG 33: Density of Related Log of Sample Data Gathered from a Large Variety of Sources on the Internet

Also of note here is that for this dataset, both 0 and 1 serve as strong anchors for the data. In other words, that the area near the origin and 1 is where data congregates most, thining out on the x-axis onwards, having a one-sided tail to the right falling off strongly from that high concentration density curve near 0 and 1. Here there wasn't any temporary initial rise in the density similar to a lognormal with high shape parameter value. In this sense the data gathered here showed some similarity to the exponential distribution by consistently falling off from the very beginning at the origin.

Here are 9 charts of the histogrm of the data itself (from the variety of sources) when different regions come into focus. We note a sudden jump around the integers as well as for "integer-like" points such as 0.1, 0.2, 0.3, 10, 20, 30 and so forth.

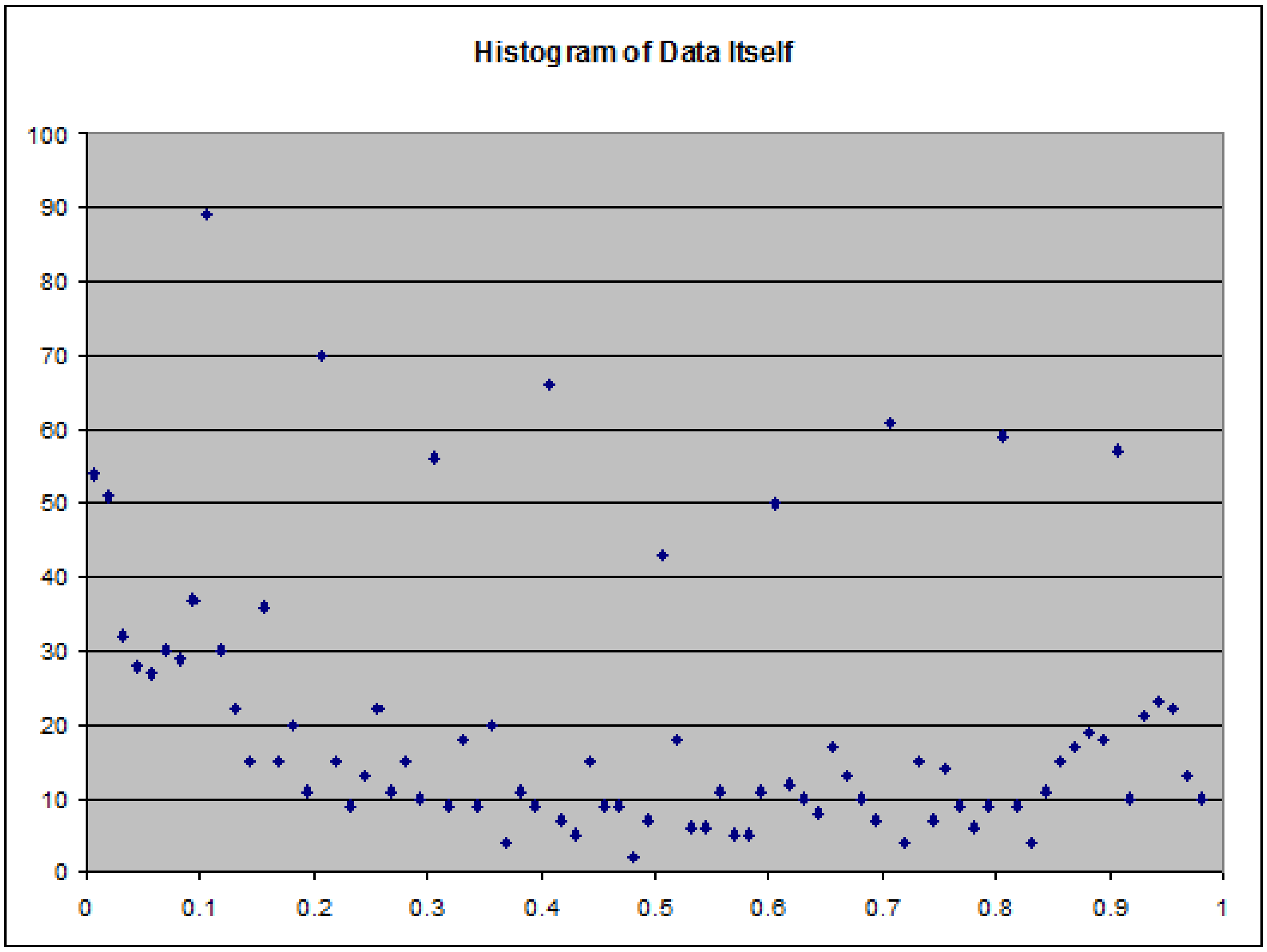

FIG 34s (begins)

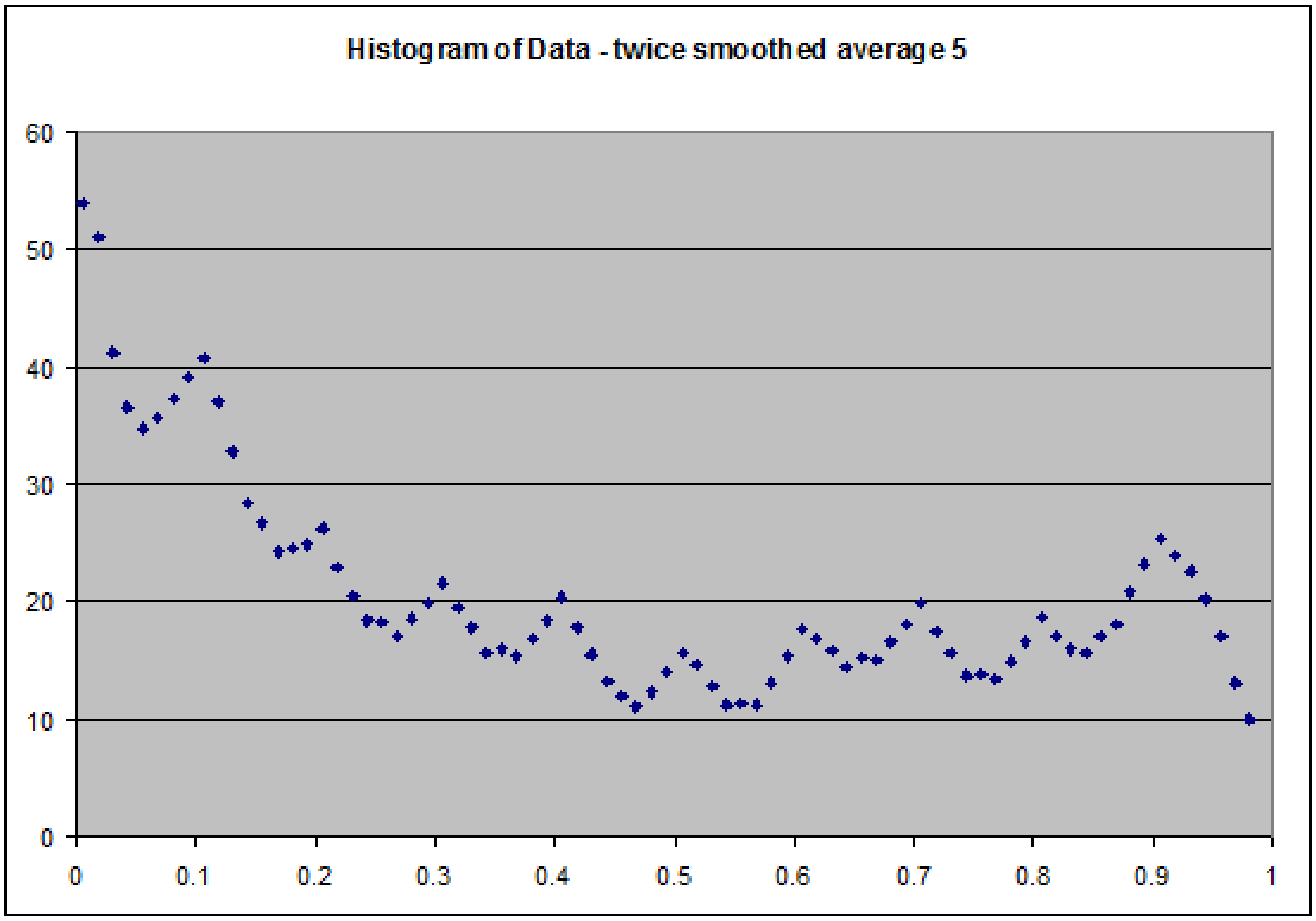
Histogram of Data - twice smoothed average 5
60
50
40
30
20
10
0
0
0.1
0.2
0.3
0.4
0.5
0.6
0.7
0.8
0.9
1

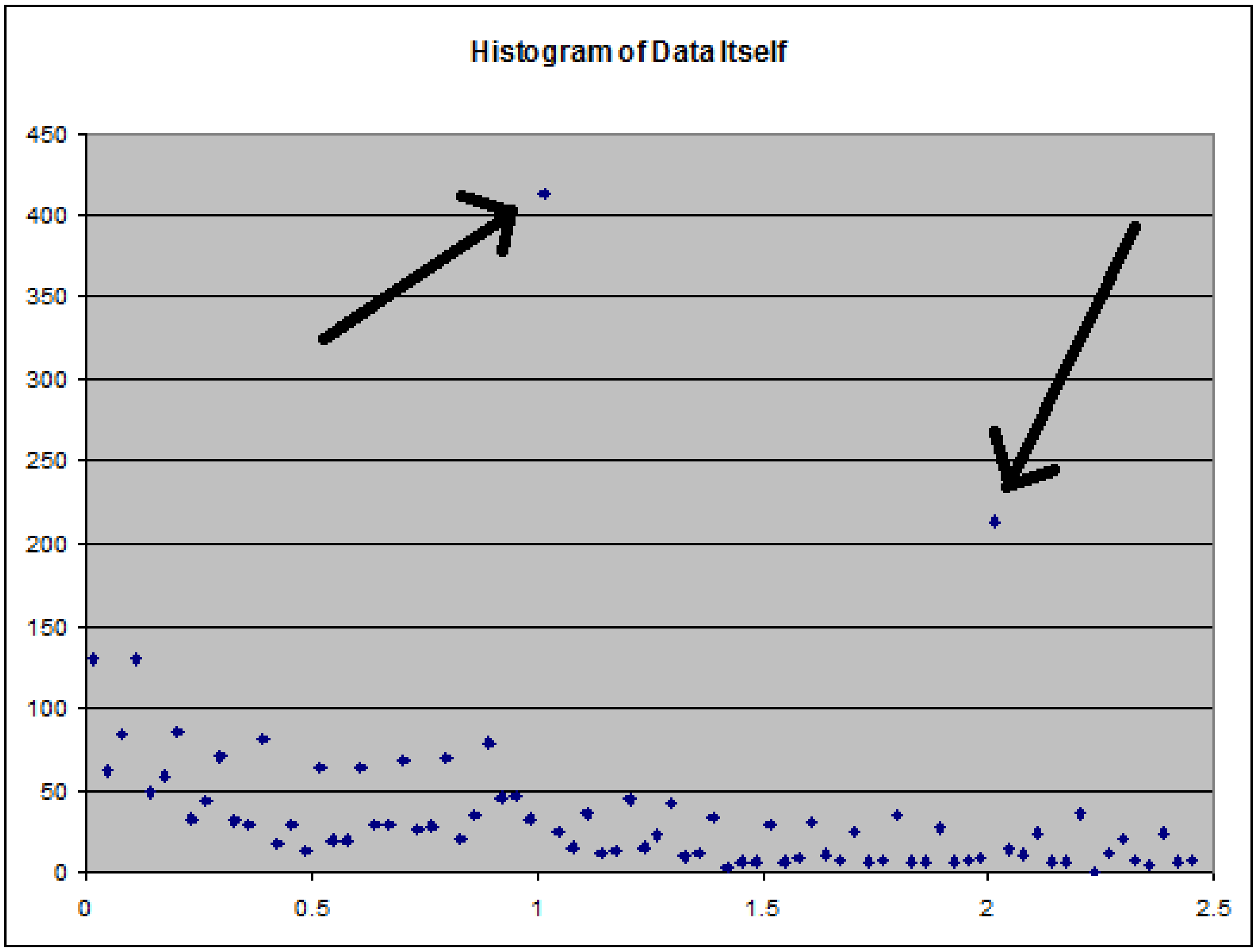

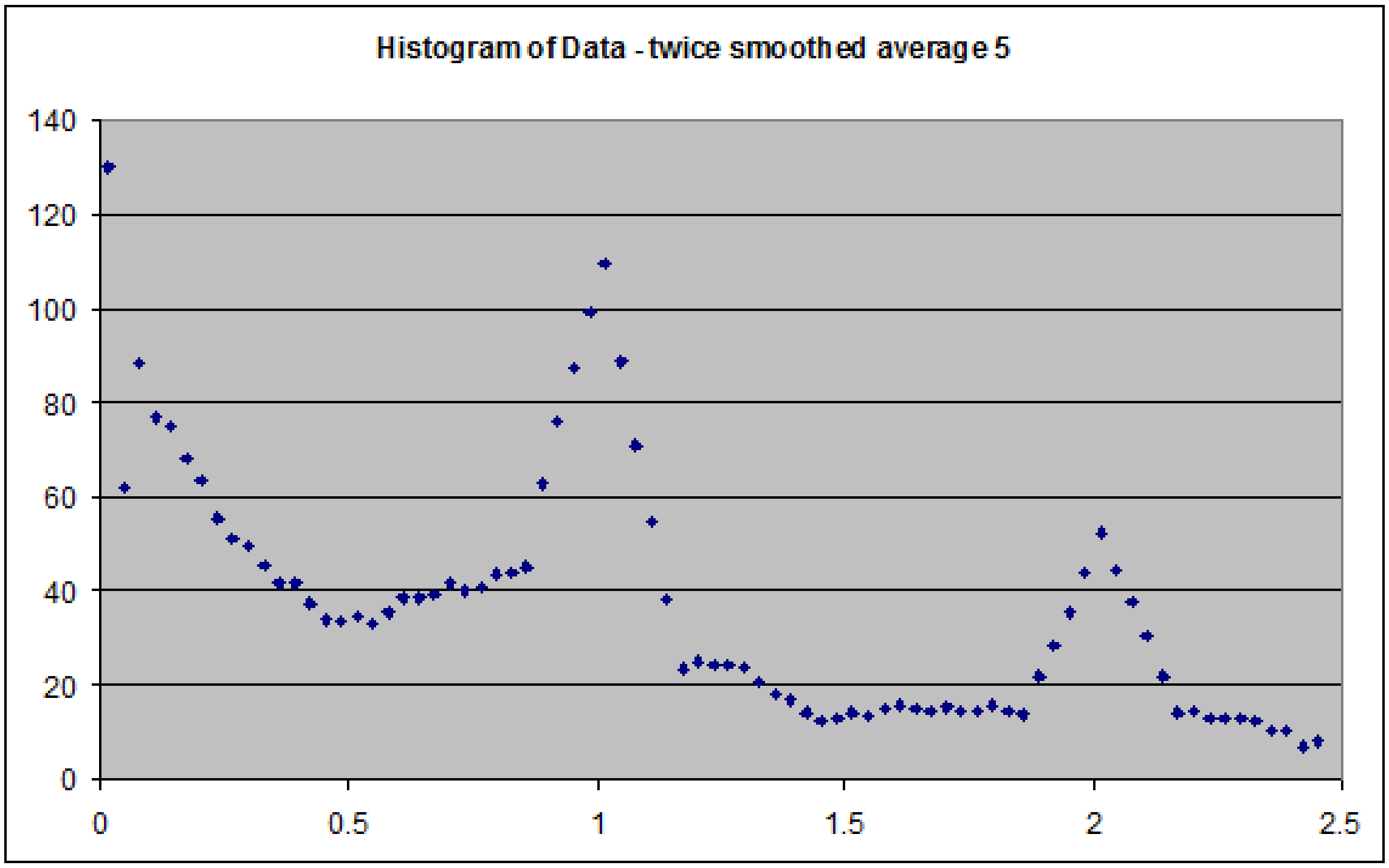
Histogram of Data - twice smoothed average 5
140
120
100
80
60
40
20
0
0
0.5
1
1.5
2
2.5

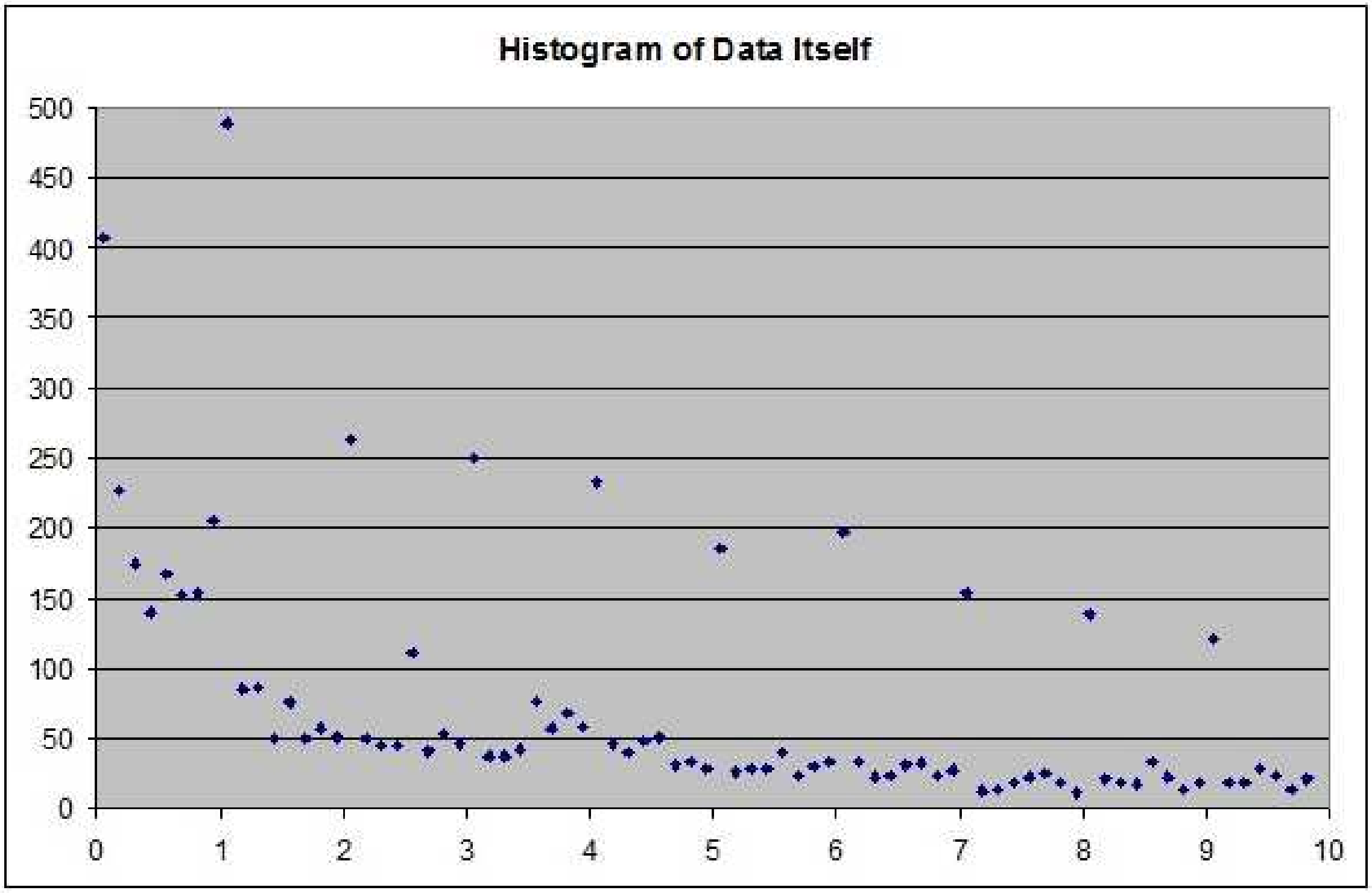

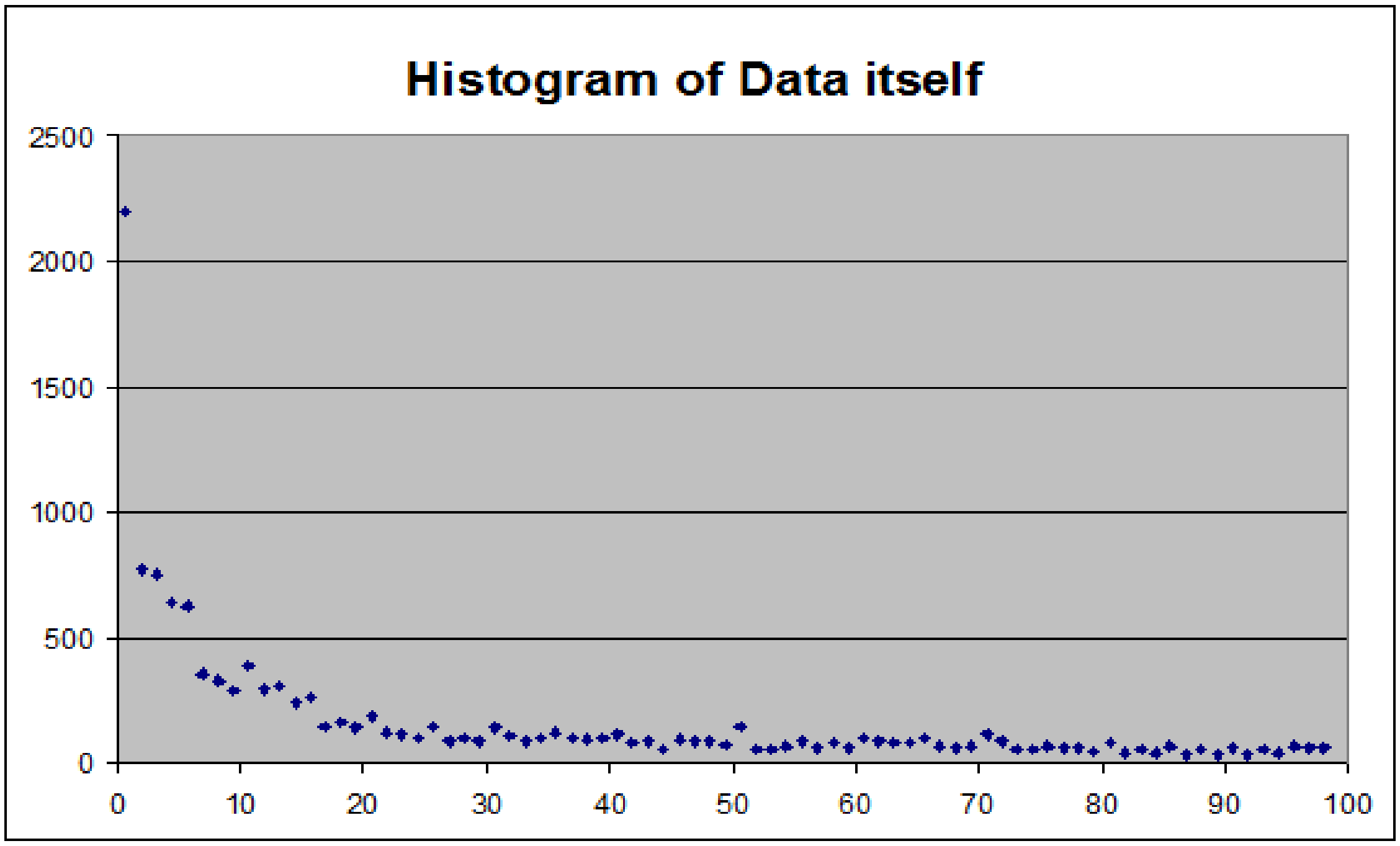

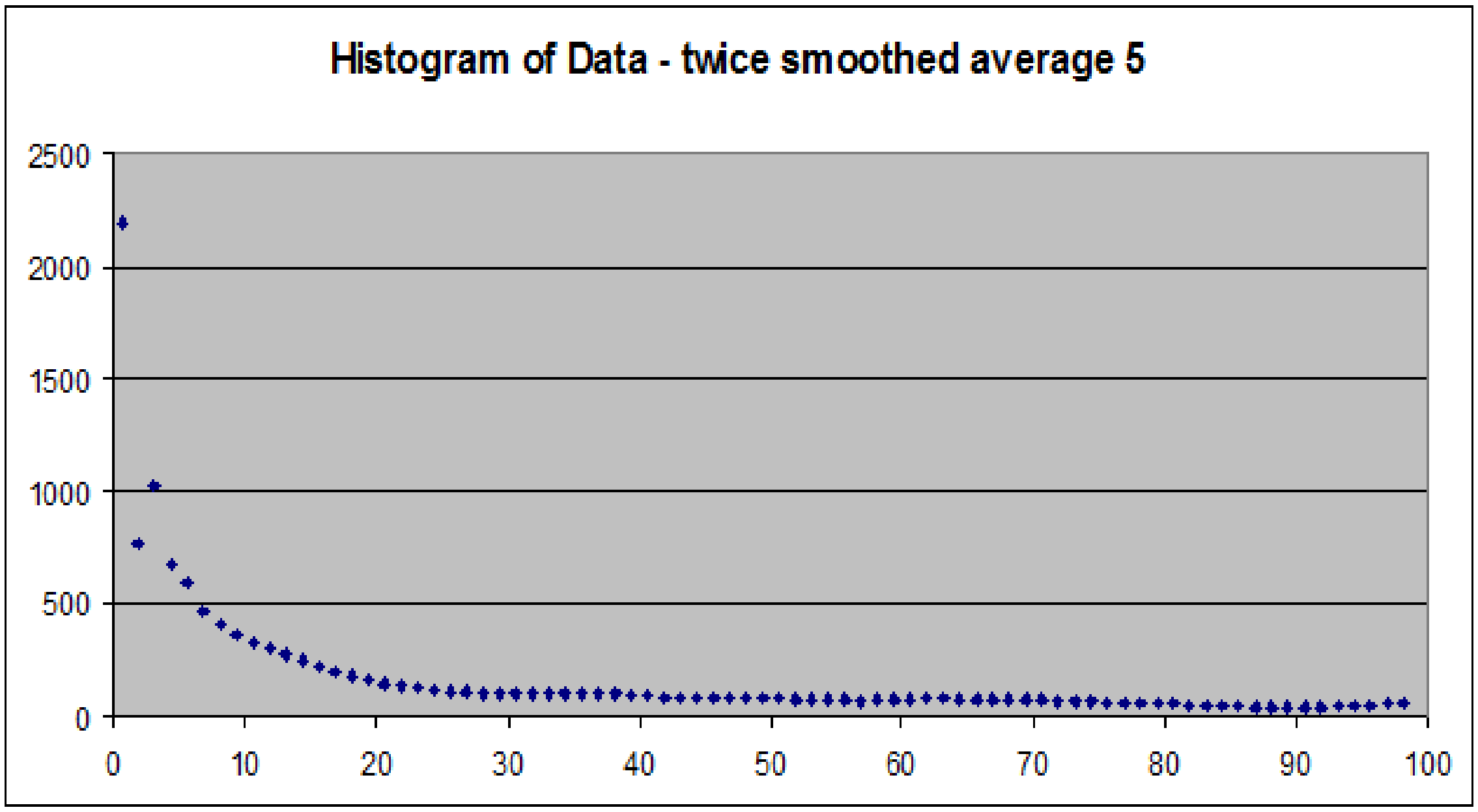
Histogram of Data - twice smoothed average 5
2500
2000
1500
1000
500
0
0
10
20
30
40
50
60
70
80
90
100

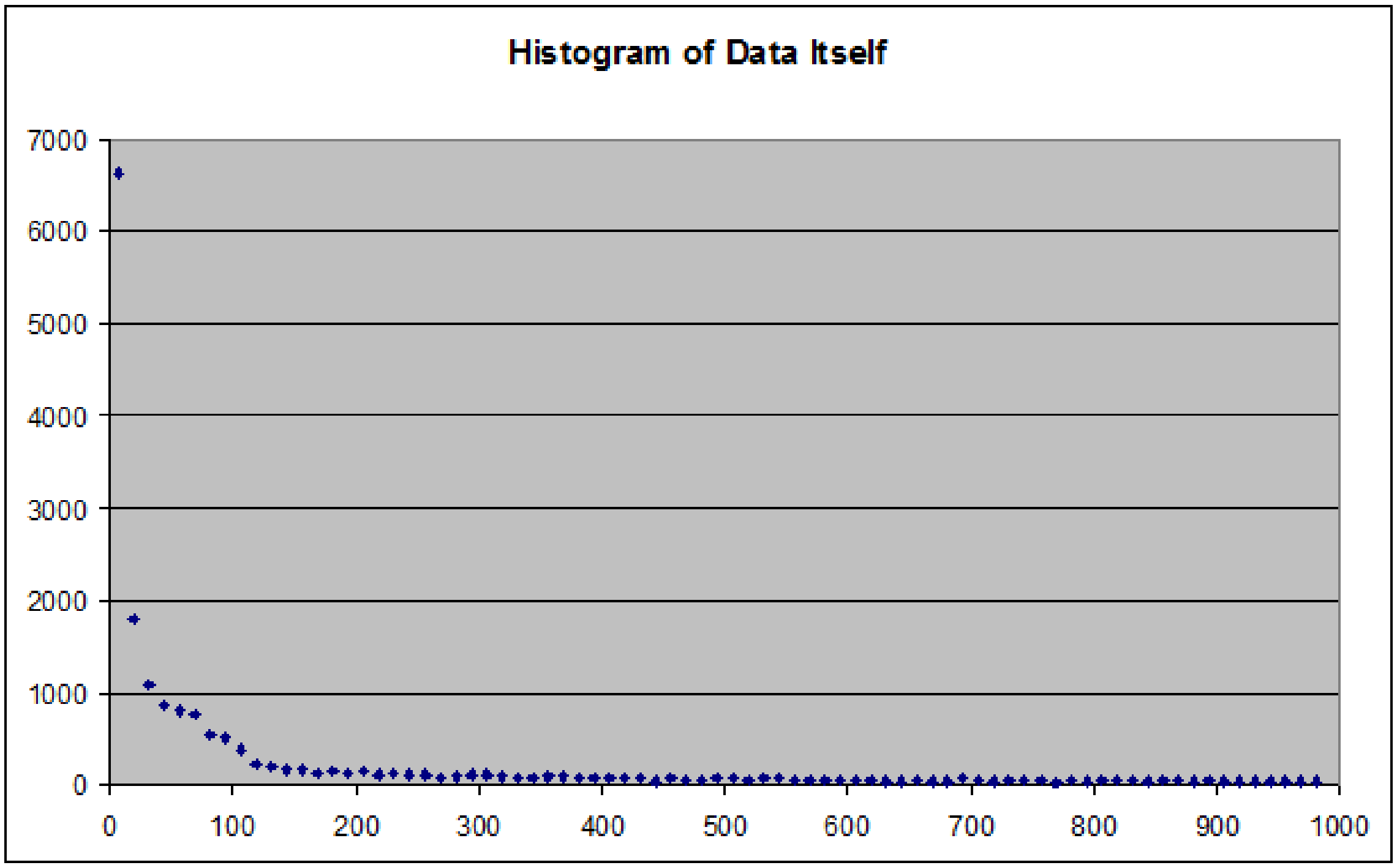

FIG 34s (end)

Here are some LD distributions, made more detailed interval by interval, of this data set from a large variety of sources: Between 0.1 and 1: { 0.178, 0.124, 0.102, 0.090, 0.079, 0.095, 0.094, 0.114, 0.126}. Between 1000 and 10,000 : {0.283, 0.186, 0.127, 0.092, 0.081, 0.067, 0.056, 0.056, 0.052}, while between $10^8$ and $10^9$ : {0.405, 0.160, 0.118, 0.095, 0.061, 0.084, 0.046, 0.027, 0.004}. Of note here is the strong tendency to increase skew-ness towards low digits as focus moves to the right.

The field of Benford's Law still does not have a rigorous mathematical proof/explanation for the phenomena where it occurs the most, namely for data in the physical world pertaining to a single issue, such as population data, length or rivers, flow of rivers, pulsars' rotation rate data, earthquake's depth, time between earthquake occurrences data, etc. These data types has nothing to do with mixtures of distribution, hence Hill's dist of dist proof can't explain away the Benford's phenomena manifested in those important cases - which thus remained (mathematically) mysterious for now. It is noted that very few real-life data types can be modeled on Hill's construct as mixtures of distributions (especially when the need to have a large variety of distribution forms - having numerous parameters - is noted.)

We shall also analyze one single-issue data type (non-Hill) relating to population for comparison. One almost perfectly such logarithmic data set is that of the US Census Bureau on population of all incorporated cities and towns in the USA. The US Census data on population of all 19,509 incorporated cities and towns in the USA starts at value 1, namely a single person living in an officially recognized town. Its top value is that of New York City with population of 8,391,881. The website link is at:
http://www.census.gov/popest/data/cities/totals/2009/files/SUB-EST2009-IP.csv
Vintage 2009 Incorporated Place and Minor Civil Division Population Dataset. This data set adheres to Benford's Law very closely. Its 2nd order digit distribution is also in close conformity with the law, as well as first-two digits and last-two digits combinations.
For example, 1st digits distribution is:
{29.4%, 18.1%, 12.0%, 9.5%, 8.0%, 7.0%, 6.0%, 5.3%, 4.6%}
While 2nd digits distribution is:
{11.9%, 11.4%, 11.3%, 10.5%, 10.2%, 9.4%, 9.6%, 8.7%, 8.8%, 8.1%}

The following chart below is the log histogram of population of USA incorporated cities and towns. While there is no good resemblance to either the Normal curve or the semi-circular one, the histogram starts and ends gradually on the log-axis itself without any abrupt spikes, and it clearly shows two distinct regions, the one on the left rising, and one on the right falling. Also the nice logarithmic behavior here can be clearly predicted and explained by the shape of the log histogram and the fact that the spread on the log-axis is wide enough, being around 4.5 to 5 units. LD-inflection point here is around log value of 2.9, which corresponds to the value of 794 for the data itself. This value of 794 is consistent with its digital development pattern as both 1st order as well as 2nd order mini LD distributions are less skewed and milder than the Benford condition to the left of this LD-inflection point, and more extremely skewed to the right of it, although it is s bit difficult to decipher all this due to the fact that it takes a whole sub-interval between IPOT values to be able to properly read local LD configuration.

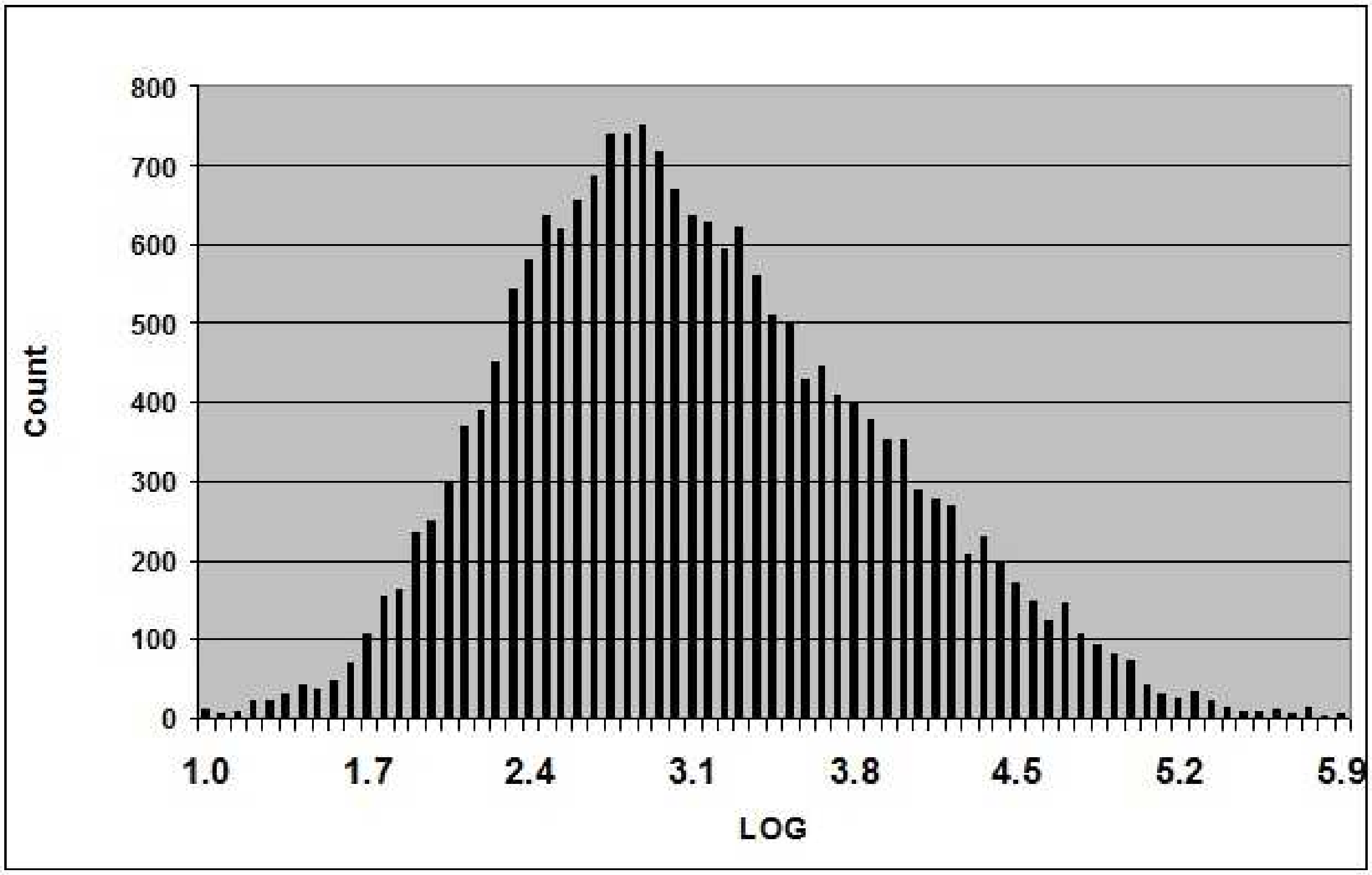

FIG 35: Density of Related Log of cities

Here are some detailed LD distributions, interval by interval, of the data on the population by cities and towns above: Between 10 and 100: {0.053, 0.081, 0.07, 0.092, 0.115, 0.139, 0.139, 0.17, 0.141}. Between 100 and 1000 {0.191, 0.174, 0.136, 0.115, 0.099, 0.088, 0.076, 0.061, 0.06}, while to the right between 100k and 1m:{0.629, 0.176, 0.06, 0.041, 0.03, 0.03, 0.015, 0.011, 0.007}. For both data sets above, the varied collection data set as well as the population data set, as we move the focus to higher values from the left regions towards the right regions, low digits are continuously and steadily obtaining more leadership from high digits. This steady evolution or development of digital proportions along the x-axis is what characterizes all random and statistical processes, and it is clearly demonstrated here, radically differentiating it from deterministic types of data (multiplication processes) which does not show any digital development whatsoever (as it's totally steady there). This feature of digital evolution or development is found in all random data, be it of Hill's type, namely a mixture of data/distributions, a chain of distribution, an averaging scheme, or a single-issue data type such as water flow in rivers, population, earthquake depth, time between successive earthquakes, pulsars rotation rates, and so forth. A markedly different situation appears in all simulations of exponential growth series, as no development whatsoever appears there; rather mini digit distributions of exponential growths are steady and constant everywhere throughout the entire range. We simulate one such (typical) growth and explore mini digital distributions, using 3% exponential growth series, starting at the initial value of 10 and considering the first 468 elements. Here are some mini LD distributions, interval by interval, of the data on this exponential growth series: Between 10 and 100: {0.308, 0.179, 0.115, 0.103, 0.077, 0.064, 0.064, 0.051, 0.038}. Between 1000 and 10,000: {0.308, 0.167, 0.128, 0.103, 0.077, 0.064, 0.064, 0.051, 0.038} and

between 1 million to 10 million:{0.295, 0.179, 0.128, 0.09, 0.09, 0.064, 0.051, 0.051, 0.051}, namely: totally stable inner digital distributions. Yet, in order to observe those development patterns (in the random case that is), it is crucial to thoughtfully and deliberately decide upon the exact nature of the partitioning of the entire range of data into smaller and adjacent sub-intervals, otherwise digital development patterns can not be observed. We seek an equality of digital opportunity between the various sub-intervals, not equality of length between them. It is only in order to obtain equal digital opportunity between the digits themselves within any single sub-interval that we require that each of the 9 digits should occupy the same distance/range in the usual sense within that given sub-interval. Thus, the requirement for the proper comparisons of digital conditions between sub-intervals necessitates constructing them in one very particular (unique) fashion by letting them stand between two adjacent integral powers of ten such as (0.1, 1], (1, 10], (10, 100], ($10^2$, $10^3$], and so forth. Sub-intervals constructed in any other way wouldn't do! The improper shorter sub-interval (10, 90] does not yield any numbers with first digit 9 (except for 90 itself). The same difficulty exists with the improper sub-interval (20, 100] where digit 1 doesn't have a chance to occur except just once. In contrast, the sub-interval (10, 200] is too long in our context because here digit 1 gets some undeserved advantage, dominating the portions (10, 20) as well as (100, 200], leaving by far less ground for digits 2 through 9. Different considerations also preclude choosing a longer sub-interval in this context. Choosing sub-intervals between two non-adjacent integral powers of ten such as (10, 1000] - which can be expressed as ($10^1$, $10^3$] – obscures pure comparisons. Such a choice could potentially confuse Leading Digits travail occurring on (10, 100] with a totally different one occurring on (100, 1000], thus these two sections should not be mixed. Aggregating Leading Digits here, instead of summarizing, would simply mask and obscure the more detailed forensic occurrences happening on each of the two separate sections.

Note that for example, ($10^{3.778}$, $10^{4.778}$], (2.5, 25], (5.97, 59.7], and so forth wouldn't do here either. For example on (20, 200], digit 1 dominates [100, 200), digit 2 dominates (20, 30), digit 3 dominates [30, 40), and so forth. Yet, for an interval such as (20, 200], digit 1 dominates ten times more range than any other digit, and this implies some bias. To show (to demonstrate) the severe distortion in observed digital development patterns under a misguided partition, we shall perform exactly one such erroneous partition on the US Census data on population of 19,509 incorporated cities and towns that was explored earlier, using the values 4, 4, 40, 400, and so forth, as border points: Here are some LD distributions, interval by interval, for such data and mistaken partition: Between 4 and 40: {0.236, 0.363, 0.316, 0.008, 0.008, 0.017, 0.008, 0.017, 0.025}. Between 400 and 4000: {0.299, 0.158, 0.093, 0.104, 0.089, 0.079, 0.069, 0.055, 0.054} while between 40k and 400k: {0.201, 0.056, 0.019, 0.205, 0.183, 0.121, 0.09, 0.072, 0.054}. No clear pattern emerges here, and there is even a slight if confused trend to lower skew-ness a bit as we more to the right. In any case, under the correct partition between AIPOT, this gradualism in leading digits development from the left to the right for all random data is much more prevalent and ubiquitous than Benford's Law itself. In other words, this tendency shows itself not only in Benford data but also in all other random/statistical data that has nothing to do with logarithmic behavior, such as payroll accounting data, US Census county area data (in miles squared, or kilometer squared), to name just a few examples. What accounts for this phenomenon? The answer may be that the very nature of random and statistical data in extreme generality is this log-gradualism and curviness, that it rarely starts nor ends abruptly around an initial/final value, and when counted and recorded in log sense shows a Normal-like or semi-circular-like behavior – be it Benford or not. That is, log is never flat/uniform, and it never abruptly starts nor ends out at a single point, but gradually builds up from the very log-axis itself to a certain plateau, and then dies out and falls all the way back to axis, which in turns implies digital development to some extend as argued earlier. For non-Benford data, deviation from the logarithmic is due typically to any one, some, or all of the following: (1) range is not wide enough to bequeath logarithmic behavior to data, (2) log curve is a bit too steep on one or two sides, (3) log curve is too asymmetric for the given range. What this assertion claims is that even for data that is not logarithmic at all, a LD-dichotomy between the 3 basic regions (left, center, right) is always still roughly observed. Thus also for numerous other types of data not obeying Benford's Law, the socialist left region of data shows roughly digital equality, the central region shows an approximate logarithmic behavior, and around the extreme far right we observe the most severe digital inequality. Two obvious applications here are immediate: the first is to forensically check non-logarithmic data for fraud, data types such as payroll and areas of counties say. That is, fraudsters who invent such data do not know anything about digital development and thus spread digits proportion equally everywhere, and their data does not show any digital development even though it's random by nature (or else their development is zigzagged and/or atypical). The statistician can easily detect such fraud. The second application is whenever the well-educated and sophisticated fraudster is well-aware of Benford's Law and all its higher orders proportions and invent/fake/concoct data according to the law, yet without any digital development, because even such smart fraudsters are unaware of the correct digital development in all real life random data, and most likely spread the Benford condition equally everywhere throughout the data (or end up with some zigzagged or meaningless digital development), all of which can be easily detected and checked by the statistician.

In order to measure digital development we first need to construct a singular numerical measure of digital skew-ness over and above the Benford condition on each of those mini sub-intervals. It is necessary then to decide upon what digits should be designated 'high' and what should be designated 'low', hopefully in a manner not considered arbitrary. An attractive approach is to designate 1, 2 as low digits, and 3, 4, 5, 6, 7, 8, 9 high digits, as the best equal division of digits in a probabilistic sense, as per Benford's Law itself. In other words, since about half of real life numbers are being led by digits 1, 2 (30.1% + 17.6% = 47.7%), and roughly the other half are being led by 3, 4, 5, 6, 7, 8, and 9 (12.5%+9.7%+7.9%+6.7%+5.8%+5.1%+4.6% = 52.3%), such a classification is natural. Hence one reasonable definition of skew-ness over and above the Benford condition is the sum of observed percent of numbers being led by digits 1 and 2 minus the Benford default sum of 47.7% for both digits, namely:
[ observed % of digit 1 + observed % of digit 2 ] – [ log(1+1/1) + log(1+1/2) ]
[ observed % of digit 1 + observed % of digit 2 ] – [ 47.7% ]. For the US cities population above, it yields the sequence of development from 1 to 10M as

{-25%, -34%, -11%, 9%, 19%, 33%, 30%}. An alternative measure could be the sum of (Oi-Bi)/Bi for digits 1 and 2, plus (Bi-Oi)/Bi for digits 3, 4, 5, 6, 7, 8, and 9, where Oi denotes the Actual proportion of numbers with first digit i within the specific sub-interval in question. Bi denotes the proportion of digit i according to Benford's Law. That is: $(O_1 - 0.301)/0.301 + (O_2 - 0.176)/0.176 + (0.125 - O_3)/0.125 + (0.097 - O_4)/0.097 + (0.079 - O_5)/0.079 + (0.067 - O_6)/0.067 + (0.058 - O_7)/0.058 + (0.051 - O_8)/0.051 + (0.046 - O_9)/0.046$. For the US incorporated cities population above, it yields the sequence of development from 1 to 10M as {-7, -8, -2, 2, 4, 6, 5}. By far the most elegant measure of skew-ness enabling us to express the whole story of digital development here would be Adrian Saville's regression line Y=b*X+a namely Observed = b * Benford + a where X represent the theoretical Benford proportions, and Y the actually observed proportions for the data set under consideration. Parameter b is the slope and parameter a is the intercept. The restriction that sum of all probabilities equals unity algebraically implies a one-to-one relationship between a and b, and thus we arrive at the calculated expression: intercept = ( 1 – slope )/9. It is therefore enough to focus merely on one parameter while the other can be thought of as totally redundant, and preferably the slope could serve as being more 'natural' as it expresses skew-ness. Saville's slope nicely measures skew-ness relative to the Benford condition, with 1 being exactly logarithmic, with >1 skewer than the logarithmic, and with <1 flatter and more equitable than the logarithmic. For the US cities population above, Saville's slopes yield the sequence of development from 1 to 10M as {-0.1, -0.4, 0.5, 1.3, 1.7, 2.3, 2.1}. In applying Saville's regression, there is no need to decide upon what digits exactly are defined as being high/low!

.

The importance of distinguishing between the random and the deterministic in the context of Benford's Law and leading digits in general should not be overlooked. Acknowledging the profound dichotomy between the two processes helps in shedding light on many aspects within the field. From its very inception, Leading Digits/Benford's Law has been suffering from this profound confusion and mixing of those two very different flavors within the discipline, causing many errors and leading research astray. An important example of such perilous omission is Pieter Allaart's summation test. The original insight and intuition here suggests that sums along digital lines may be all equal. For example, the sum of all numbers beginning with digit 1 was thought to equal the sum of all numbers beginning with digit 2, and so forth. It might be attractive to intuit that summing those very few numbers (of supposedly larger values) beginning with 9, should yield something similar if not equal to summing all those much more numerous numbers beginning with digit 1 (and having supposedly lower value.) Pieter Allaart has given a rigorous mathematically proof of the conjecture. The author of this article though observed that his proof is based on a restricted range for the data/distribution, spanning merely two adjacent integral powers of ten, such as 1 &10, having unity exponent difference. Is it extremely rare (if ever) that real-life data should fall into such 'small' restricted range, but if and when it does, such data (by default - automatically) is of the deterministic flavor representing exponential growth series, with the continuous density curve k/x necessarily being that unique curve representing it. Indirectly, by basing his proof on such a restricted interval, Allaart (passively) chose to consider only deterministic processes. As discussed earlier, actual statistical and random data typically falls over multiple integral log-distance, as well as showing a marked differentiation between the left, the center, and right regions, being that their log density is rising on the left and falling on the right. All this implies that low digits are relatively weaker on the left where values are small and matter less, and that low digits are relatively much stronger on the right (as compared with the logarithmic distribution) where values are large and matter more, hence that perfect trade-off as was the case for the mantissa itself does not apply here when summing actual amounts, and low digits should take, and do indeed take, higher overall sums than high digits. On the other hand, sums are equal along digital lines for the deterministic version of Benfordness where the log is uniform throughout, such as exponential series as well as k/x distribution, as shown by Allaart, but prevalence in real life data is extremely limited. The population data above on US cities and towns (with New York City excluded as an outlier) shows a decisive, highly significant and almost monotonic fall in sums along digital lines, Sums on digits 1 and 9 there differ by a factor of 5.9! Numerous empirical tests on real life random data show other common factors of sums along digit 1 and digit 9 being typically 5, or 6 or even 8! In general, if low digits are relatively weaker on the left where values are small and matter less, and if they are relatively stronger on the right where values are large and matter more, then low digits should take higher overall sum than high digits. Indeed all summation tests on actual statistical/random data relating to accounting, financial, population, governmental census, and other random types, show a strong and consistent bias towards higher sums for low digits. We shall now formally prove Allaart's statement in the 1st digit sense for the deterministic case where k/x on (a, b) serves as the exact density, and with a range such that a and b are two adjacent IPOT points of the particular (1, 10) case first, to be enlarged later to **any** other two IPOT values, adjacent or not. When a continuous random probability density distribution is considered, it is not possible to sum in the usual manner by adding discrete values, there are none, rather it is the **average** defined as ∫ x*f(x) dx, that can express a measure of addition. In the discrete case, we employ the idea of summation in the definition Average=[Sum]/[Number of Values] to express some measure of centrality within the entire data. To show **directly** why comparing averages for the case of the continuous probability density function is equivalent to comparing sums in the discrete case, one has to trace back the way one gets from densities to histograms, which is done by simply multiplying each tiny vertical rectangle inscribed within the density curve by the fixed value N representing the number of all the values within the data set fitting the density.

It will now be shown that each digit d has the same average on its sub-interval (d, d+1): Average for digit d = ∫ x*f(x) dx [from (d) to (d+1)] = ∫ x*[k/x] dx [from (d) to (d+1)]. The x term cancels out, and we are left with: ∫ k dx [from (d) to (d+1)] = k*[(d+1) – (d)] = k, namely the constant value of k for each of the 9 digits. Two essential features provide for this equality, one is the constant integrand term after the cancellation of the x term, and the other is the equality of length on the x-axis in the limits of integration for the various digits. Since those two features are present just the same for all other cases of adjacent IPOT points such as (100, 1000) and so forth, the proof can be enlarged to include them. For example, in the case of (100, 1000), the limits of integration for each digit d leading are from 100*d to 100*(d+1), which is also the same for all digits, and the value of the integral is also constant. In the case where k/x is defined over non-adjacent

IPOT values, such as (1, 100), the same equality holds, as more than one set of limits of integration are involved where the terms involving d cancel out. For (1, 100), the d-term cancellations are within (k+1) – (k) and 10*(k+1) – 10*(k), and together the differences are equal to the constant 11*k, and thus equality follows just the same. The proof is ***not*** valid though for cases where k/x is defined over an integral exponent difference of the non-IPOT type, such as (30, 3000), nor even for the cases where k/x is defined over a unity exponent difference such as (30, 300), since there the limits of integration are not of equal distance for each digit. For example, in the case of (30, 300) range, limits of integration for digits 1 & 2 span 100, while for digits 3 to 9 they span only 10, hence sums vary considerably. The outline of the above proof easily carries over to the same summation equality along first-two digits (FTD) line, where instead of integrating on the intervals (d, d+1), we integrate on much smaller sub-intervals. For example, for k/x defined on (1, 10), averages along FTD lines are calculated on (1.0, 1.1), (1.1, 1.2), (1.2, 1.3), … (9.8, 9.9), (9.9, 10), and with all 90 sub-intervals spanning the same length of limits of integration, resulting in the same uniformity of averages that was seen earlier.

Here are some of the histograms of the data itself (of population of US cities and towns incoprporated) when different regions come into focus. A temporary rise near the origin until around 90 is depicted, followed by an approximately flat region 90 to 150, while most of the region to the right of 150 shows a consistent fall in density.

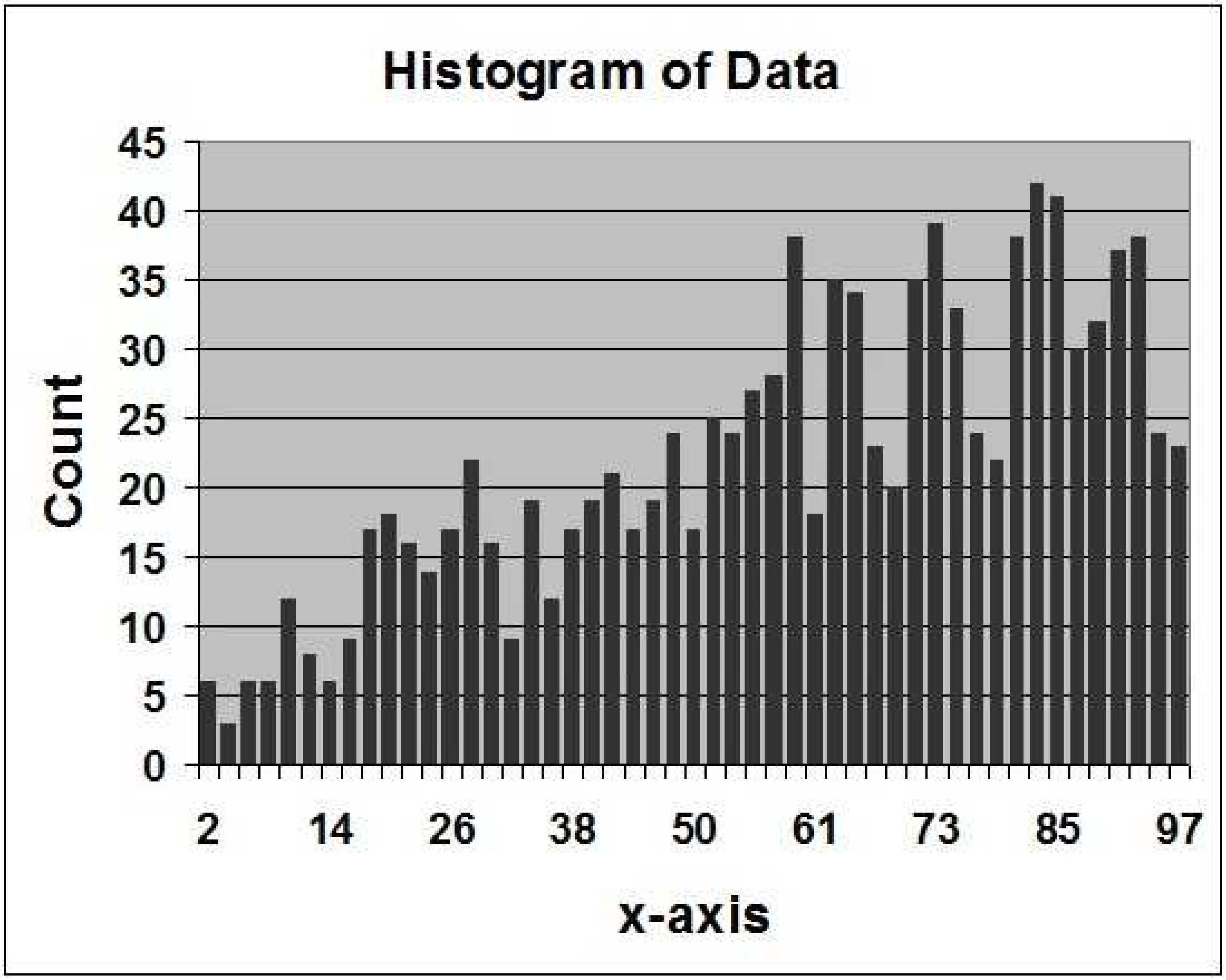

FIG 36s (begins)

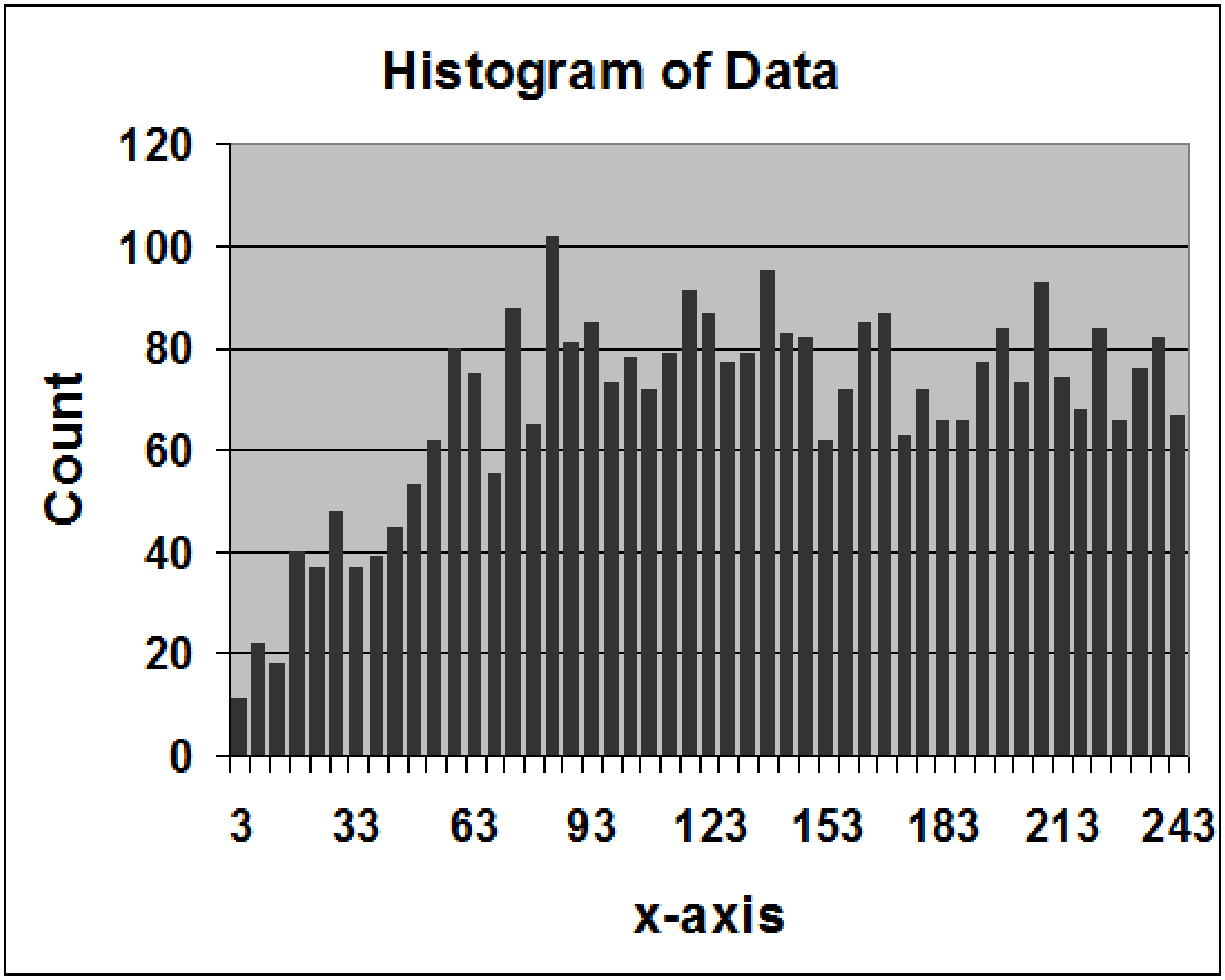

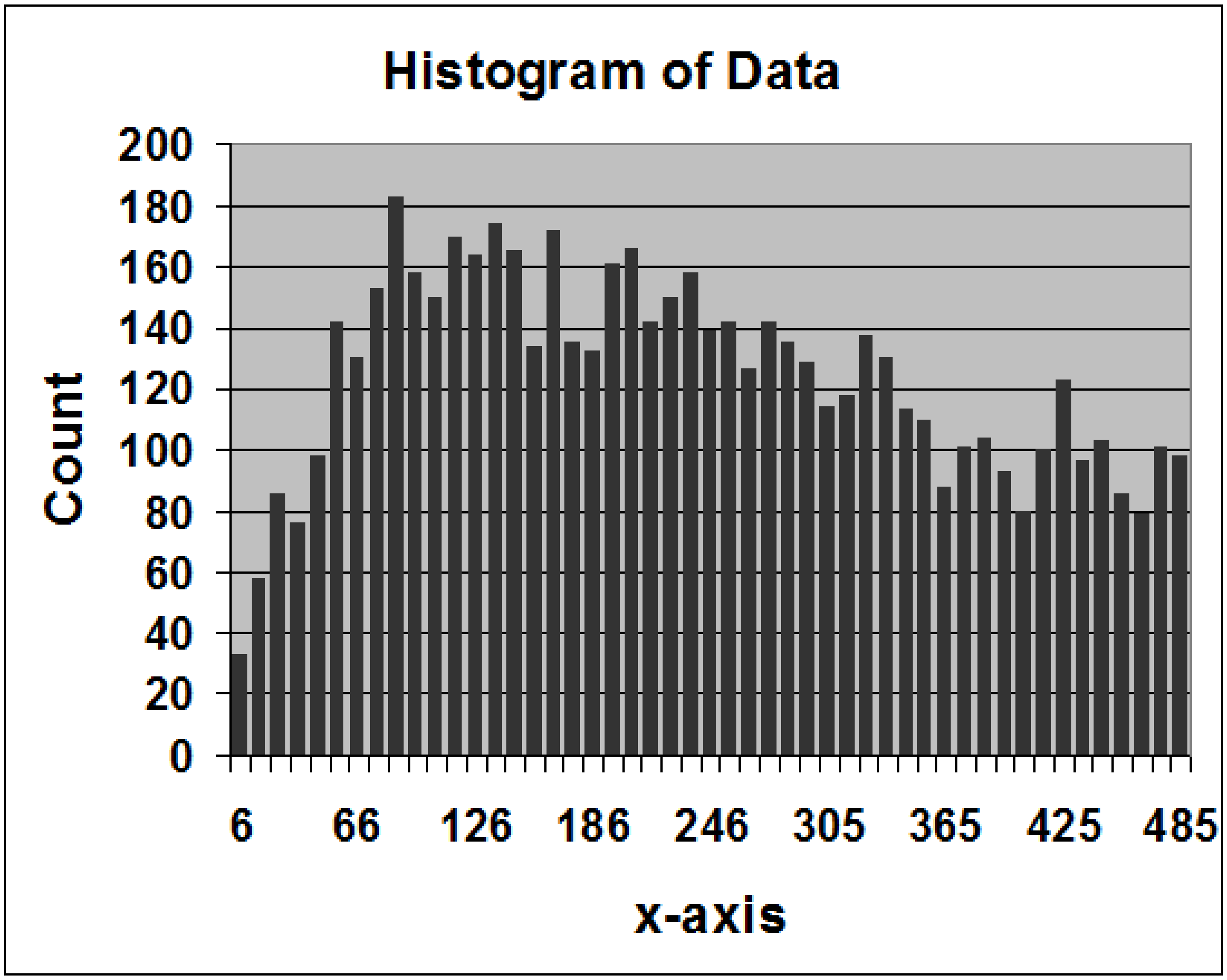

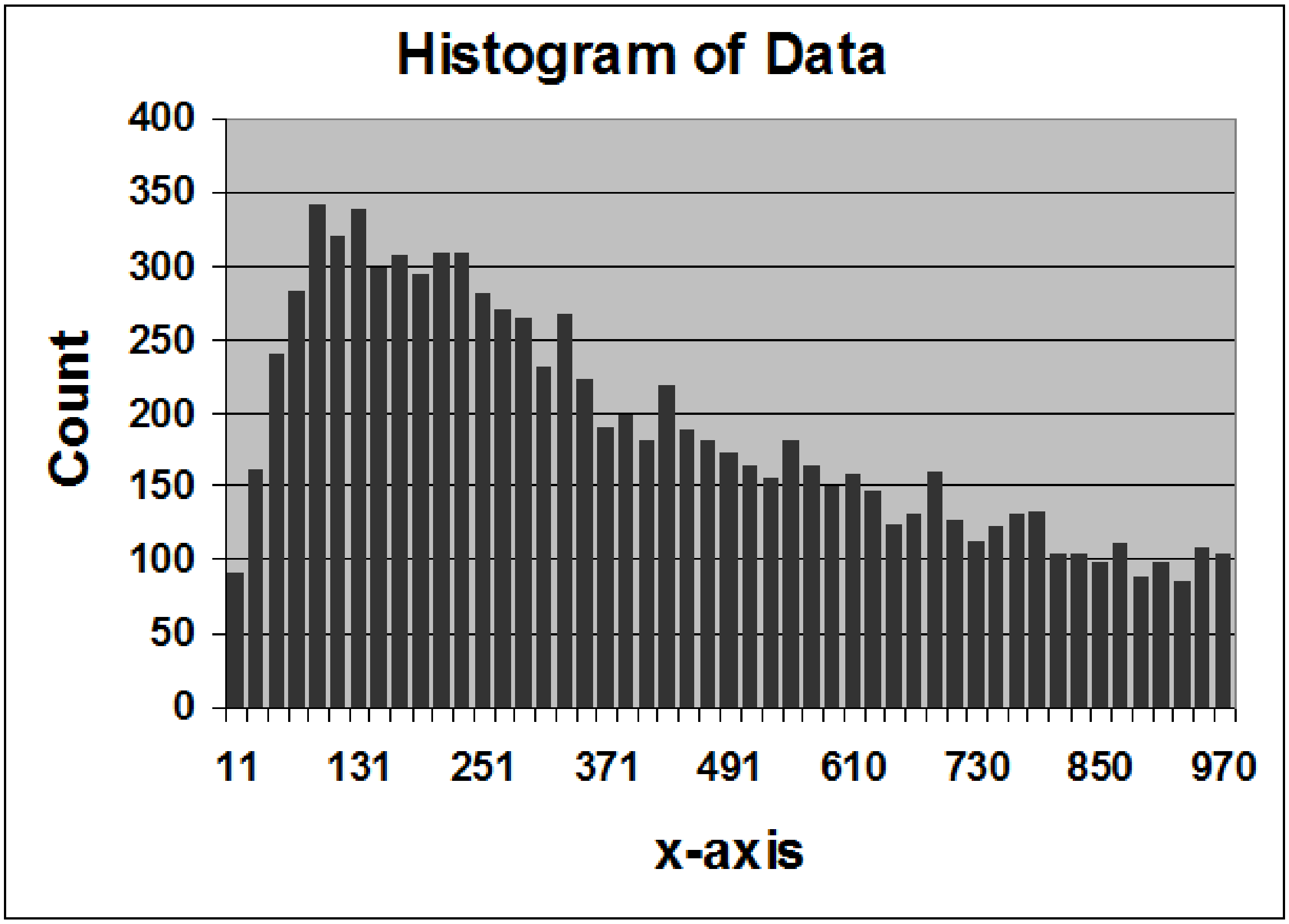

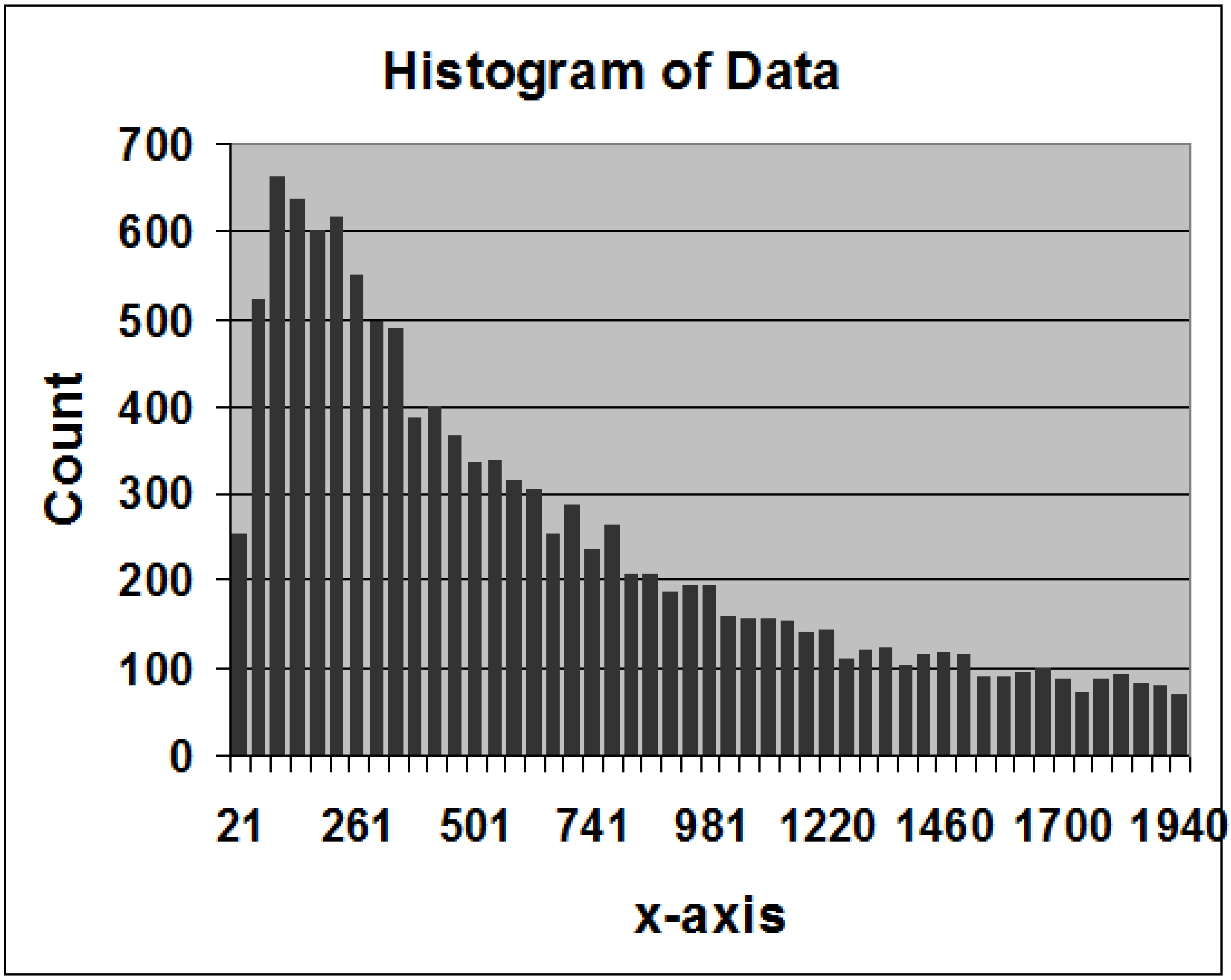

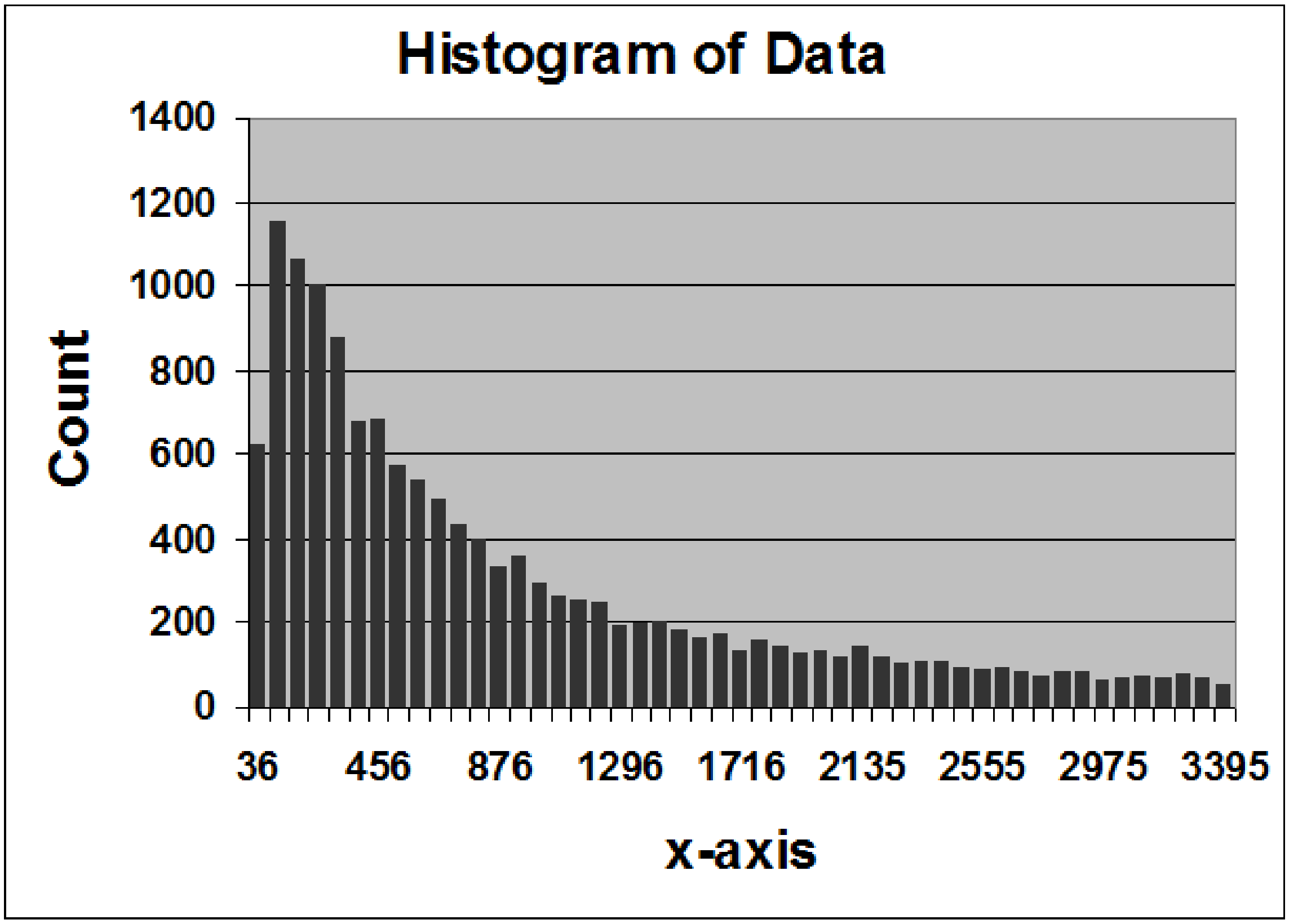

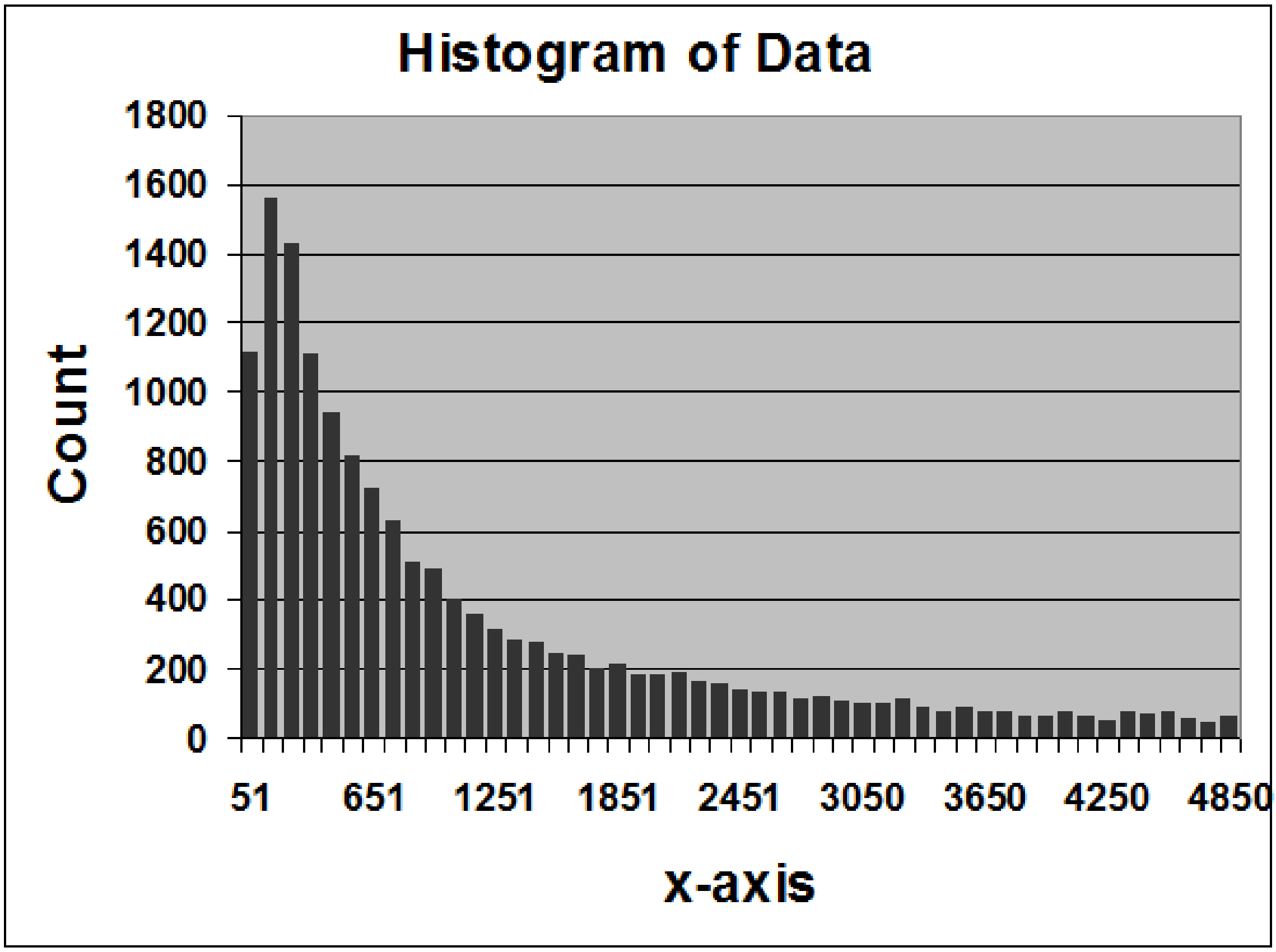

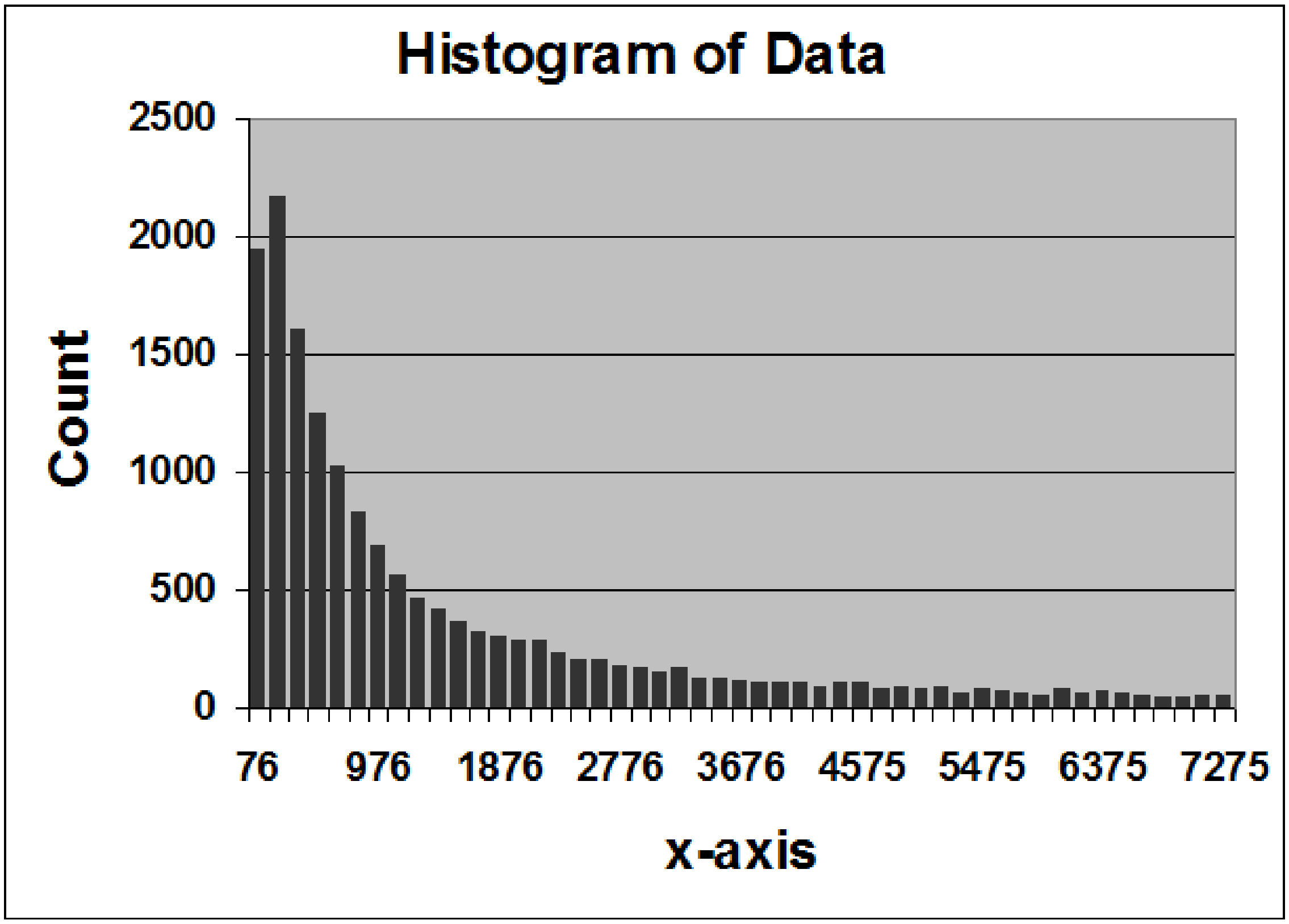

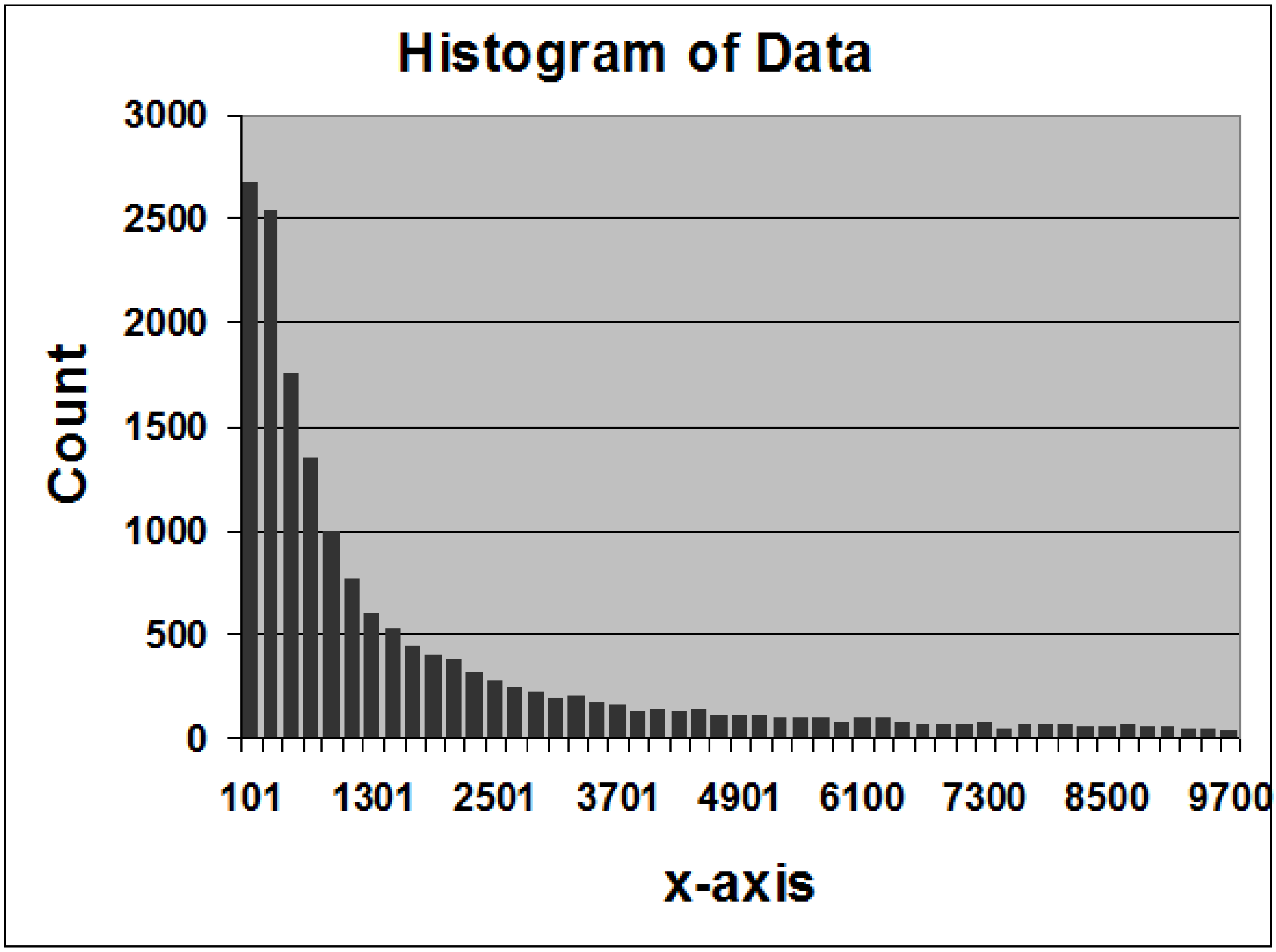

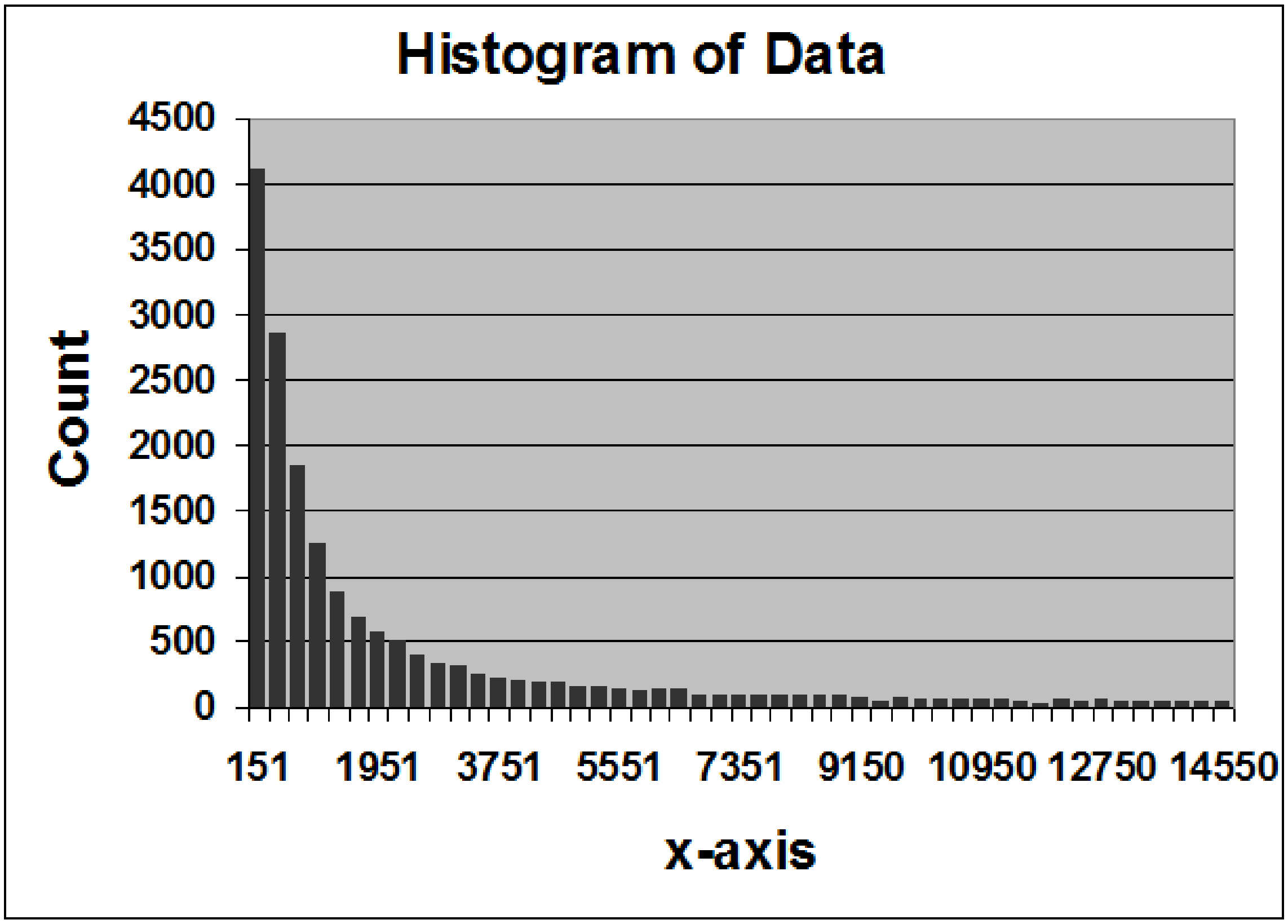

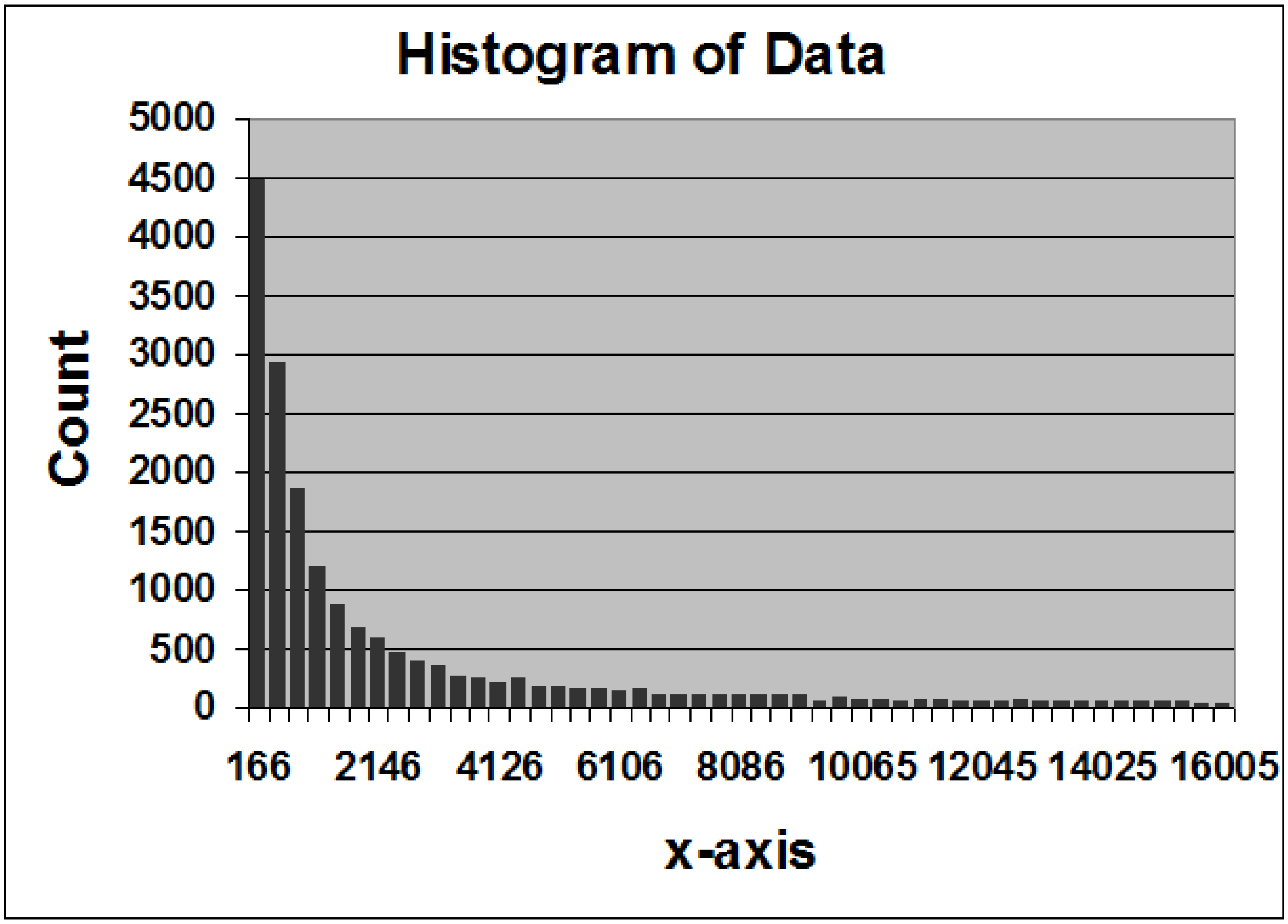

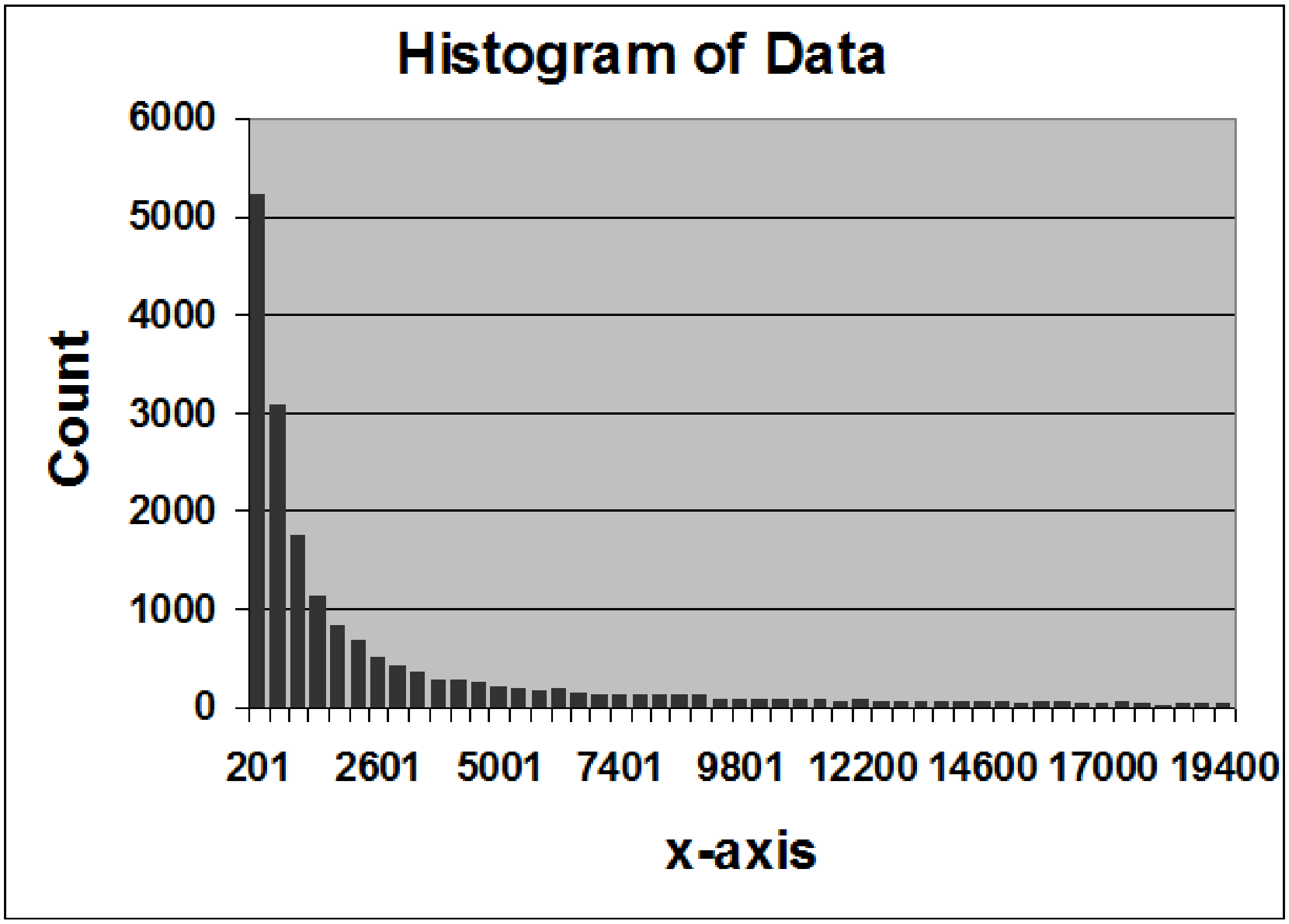

FIG 36s (ends)

## [29] Frank Benford's Non-Simple Averaging Scheme as Chain

Benford's own attempt in his 1938 article to establish the logarithmic distribution can now be viewed and explained in term of the 2nd chain conjecture. Benford employed a modified version of the simple averaging scheme to arrive exactly at the logarithmic distribution. His solution introduces a curious twist, putting more emphasis on shorter intervals than longer ones by having upper bounds vary exponentially. In other words, instead of having UB vary upwards by constantly adding 1 as have been done here earlier with all our averaging schemes including that of Flehinger, namely such as {99, 100, 101, … , 998, 999} and so forth, he varied them exponentially, such as say 2% growth as in {99, 101.0, 103.0, 105.1, 107.2 , … , 927.8, 946.4, 965.3, 984.6} and fractions ignored. His semi-logarithmic scale of P vs.UB and his calculations of the area under this curve to arrive at its average height are mathematically equivalent to the use of exponential growth series as upper bounds.

His scheme could then be viewed or interpreted simply as a 2-sequence chain **Uniform(0, *exponential growth*** ) where N is an integer. This very short chain is noted for its special feature of inserting a distribution that is logarithmic in its own right to serve parameter **b** of the uniform distribution there, although parameter **a** is being made fixed at 0, as opposed to being chained as well.

Benford's construction is an algorithm that blends two different structures and stands on two different legs. On one hand it is a collection of intervals, random uniform distributions in essence, and on the other hand it relies on detreministic multiplication processes because their upper bounds are derived from an exponential series.
Let us recap Benford's model and assign specific values for the model. Lower bounds for all intervals are fixed at the low value of 1. Upper bounds vary as an exponential having 2% growth, starting from 99 and ending just short of 999, hence upper bounds are {99, 100.9, 102.9, 105.0, 107.1 , … , 927.8, 946.3, 965.2, 984.5}. For the sake of facilitating calculations, only the integers are to be considered on each line, without effecting result at all almost. Although upper bound should vary from an IPOT number to another much higher such number as in a limit, let us focus on a narrower range of 99 to 999 without any loss of generality.

To illustrate the model more clearly, the following figure is added:

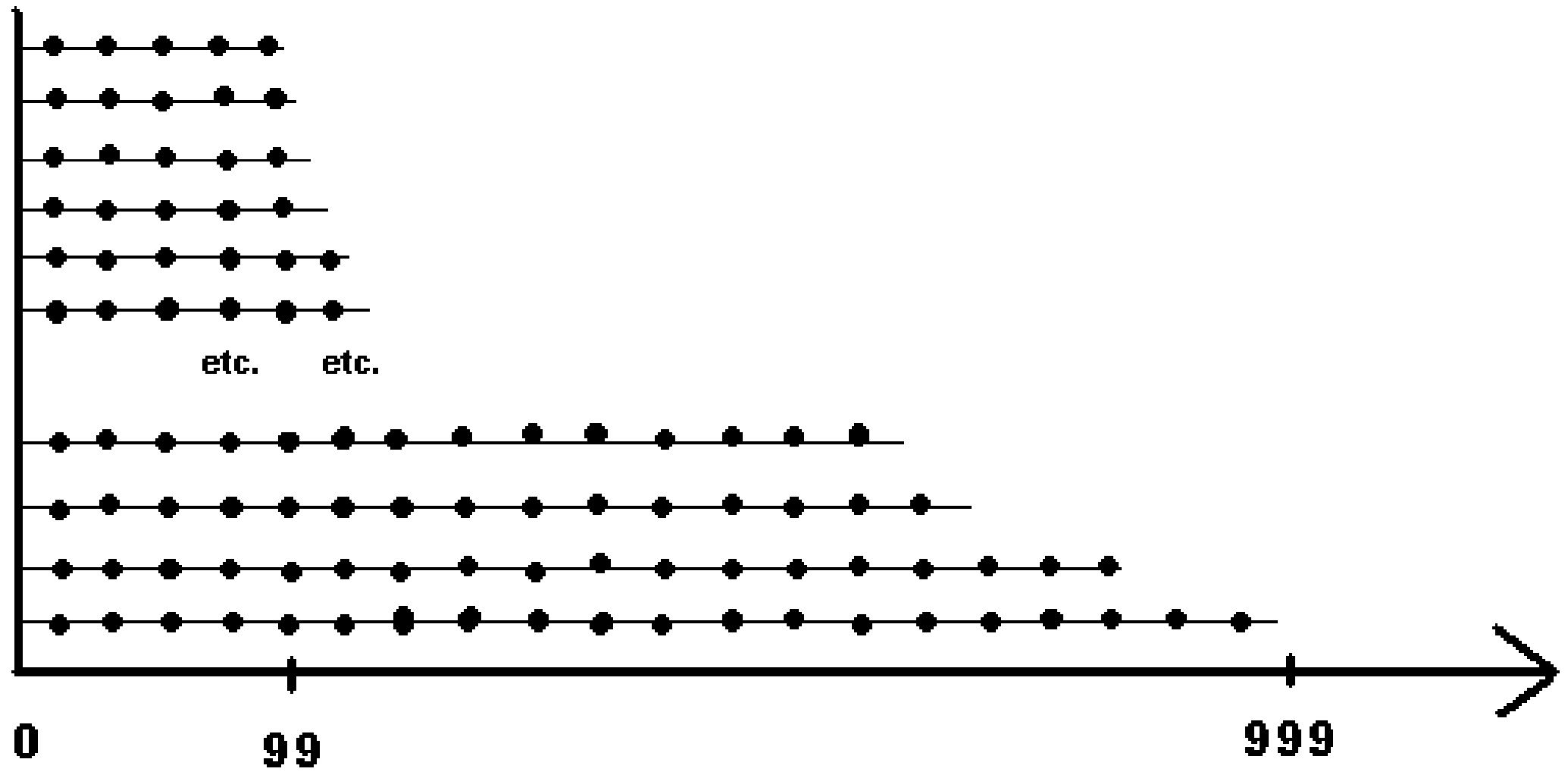

FIG 37 : Frank Benford's Attempt at Explaination, Successive Upper Bounds Grow Exponentially

We have depicted in the above figure as dots the integers placed on the continuous line, and showed only a fraction of them for lack of space. A line here corresponds to what we called an interval in the averaging schemes. As in our earlier effort to represent the simple averaging scheme as a single dataset, it is understood that dots on shorter lines are to be repeated, and that the shorter the line the more such duplications are needed, until all lines are of the common length of 999. This is so since all lines/intervals are accorded equal importance in the overall scheme - as seen in the act of taking the simple average of LD distributions of all the lines/intervals. With the minor issue of filling in at the extreme right side aside having a minute effect, the almost exact factor of duplication for each line is given by F=999/UB, where UB represent the value of x at end of the line/interval, or simply F=999/x. This factor F insures that all lines/intervals are now accorded equal importance (equal length). Hence the contribution of each line in generating dots (integers) for the system is the number of dots it carries up to its upper bound times the factor F, and this total per line is the same for all of them, namely 999 dots. An equivalent mirror image of this scheme is not to stretch out the lines, but instead to envision all the lines as having different densities of dots, with shorter ones denser, longer ones less dense, and all ending up with equal total numbers of dots due to this exact trade-off. The expression for the density (dots per unit length) for each line is then [UB*F]/UB, that is [UB*(999/UB)]/UB or simply 999/UB. For example, the shortest dense line with UB=99 has density of dots 999/99 or about 10. The longest diluted line with UB=999 has density of 999/999 or just 1.

As we move in the positive direction on the x-axis anywhere between 99 and 999, the number of lines per unit length are falling off as in the exponential series, namely at the rate of **-k/x**. This is so since discrete exponential growth series corresponds to the continuous k/x distribution via their common related uniform log as seemed earlier. Since each line carries a density of 999/UB dots per line or simply **999/x**, it follows that the decrease in the number of dots per unit length is expressed as **-k/x** times **999/x**, we get: $d(dots)/dx = (-k/x)*(999/x) = -Q/x^2$ , for Q being some positive constant.
Anti differentiating we get: dots per unit length = Q/x. We have thus established that the density is of the type k/x. Since this is considered over (99,999) or roughly $(10^2,10^3)$, exponent difference is integral (unity), namely that the range is bounded by two integral powers of ten, and by the proposition I in the earlier section it is logarithmic.

Benford's scheme is then the chain: **Uniform(0,** ***1/(x*ln10) dist over (10, 100)*** **)**.

Yet Benford's scheme can not constitute an explanation for the phenemena even though his abstract model is logarithmic, because it would be very hard to argue that typical everyday data follow this model. Hill's model on the other hand is perfectly suited for being a perfect representative and thus the explanation.

## [30] Direct Expression of 1st Digit Order & Significand

A useful expression yielding the first digit of a positive number X is given by:

$$1^{st}\ \text{Digit of } X = INT(\ X/10^{INT(LOG_{10}X)}\ )$$

The INT function refers to the integer on the (-∞, +∞) x-axis immediately to the left of the number R in question. If R is exactly an integer then INT(R) is that same integer. If R is positive then INT(R) is just the whole part excluding the fractional part. If $0 \le R < 1$ then INT(R) is 0. If R is negative then INT(R) then is the largest negative integer less than or equal to R.

Recall that $INT(LOG_{10}X)$ is the characteristic – the integral part of the log.
Also any X can be broken into $X=10^{Chracteristic}10^{Mantissa}$ , hence

$$X/10^{INT(LOG_{10}X)} = 10^{Chracteristic}10^{Mantissa}/10^{Chracteristic} = 10^{Mantissa} = \text{Significand}$$

And since significand in [1, 10) is that first part of scientific notation with the decimal location totally ignored, INT(Significand) yields the 1st leading digit.

# Miscellaneous Items:

## Programming code in C++:

### *SIMPLE AVERAGING SCHEME (city):*

```
#include <iostream>
#include <iomanip>
#include <cstdlib>
#include <ctime>
#include <cctype>
#include <cmath>

using std:: cout;
using std:: cin;
using std:: endl;
using std:: setprecision;
using std:: setw;

double leading(double);


int main()

{

        int LB = 1 ;     // Lower Bound for ALL intervals!    and    typically
assumed to be 1
        int minUB = 70 ;
   int maxUB = 800 ;
        double distribution[10]={0};
        double ALLdist[10]={0};
        int a;
```

```
//=============================================

for(int q=minUB ; q<(maxUB+1) ; q++)
{

        for(int j=LB ; j<(q+1) ; j++)  // j IS the integer!
        {
        a=0;
        a=a+leading(j);
   distribution[a]++;
        }

   for(int s=0 ; s<10 ; s++)
   distribution[s]=distribution[s]/q;

        for( s=0;s<10;s++)
        ALLdist[s]=ALLdist[s]+distribution[s];


        for( s=0 ; s<10 ; s++)
   distribution[s]=0;


}

//============================================


        cout<<endl;

        for(int j=0 ; j<10 ; j++)
        {   cout<<"% of "<<j<<" leading digit  ---> ";
            cout<<ALLdist[j]/(maxUB-minUB+1)<<endl;
        }

        return 0;
}
```

```
//-----------------------------------------------------


double leading(double x)
{
          double small,lead;

    for(int i=1  ; ; i++)
          {
                    if(x<10)
                    {small=x; break;}
                    else
                    x=(x/10);
          }

         lead=floor(small);
         return lead;
}//-----------------------------------------------------
```

## *AVERAGING OF AVERAGING SCHEME (Country):*

```
#include <iostream>
#include <iomanip>
#include <cstdlib>
#include <ctime>
#include <cctype>
#include <cmath>

using std:: cout;
using std:: cin;
using std:: endl;
using std:: setprecision;
using std:: setw;

double leading(double);




int main()

{



        int SUPERmax = 1000;
        int SUPERmin = 100;
  //int maxUB = 300 ;  ------------------> this is now to vary from SUPERmin  to
SUPERmax
        int minUB = 1 ;     // SAME for ALL cities! we have no flexibility here,
another much more sophisticated scheme would be needed to handle such
additional variability
        double LB = 1 ;        // Lower Bound for ALL intervals!    and typically
assumed to be 1
        int forview;
```

```
        double distribution[10]={0}; //street
        double ALLdist[10]={0};      //city
        double avgSUPER[10]={0};     //country
        int a;

for(double view=SUPERmin ; view<(SUPERmax+1) ; view++)  //        CITIES
view=longest street in the city
{


//===================================================

for(double q=minUB ; q<(view+1) ; q++)   //                    STREETS
q=street's length
{

        for(double j=LB ; j<(q+1) ; j++)  // j IS the integer!       HOUSES
j=house's number
        {
        a=0;
        a=a+leading(j);
   distribution[a]++;
        }

   for(int s=0 ; s<10 ; s++)
   distribution[s]=distribution[s]/(q-LB+1); // for a particular street q houses
long



        for( s=0;s<10;s++)
        ALLdist[s]=ALLdist[s]+distribution[s]; // for the particular city
(having view as the max UB) we add this street




        for( s=0 ; s<10 ; s++)
   distribution[s]=0;

}

//=============================================
```

```
   for(int s=0;s<10;s++)
         ALLdist[s]=ALLdist[s]/(view - minUB +1);   // averaging out al streets
for this particular city

         for( s=0;s<10;s++)
         avgSUPER[s]=avgSUPER[s]+ALLdist[s];   // and now we add this city
to the long country log/records

         for( s=0;s<10;s++)
         ALLdist[s]=0;

//       forview=0;  forview=forview + view;
//       if(forview%10==0)
                 cout<<view<<endl;   // clear out - empty all info - for the next
new generic city

}



         cout<<endl;

         for(int j=0 ; j<10 ; j++)
         {   cout<<"% of "<<j<<" leading digit  ---> ";
                 cout<<avgSUPER[j]/(SUPERmax-SUPERmin+1)<<endl;

         }

         return 0;
}
```

```
//-----------------------------------------------------
double leading(double x)
{
          double small,lead;

    for(int i=1  ; ; i++)
          {
                    if(x<10)
                    {small=x; break;}
                    else
                    x=(x/10);
          }

         lead=floor(small);
         return lead;
}//----------------------------------------------------


//double totalP;

/*
   // checking for 100% compliance
         cout<<endl<<q<<endl;
   totalP=0;
         for( s=0;s<10;s++)
         totalP=totalP+distribution[s];
         if(abs(totalP-10>0.001))
         {
         cout<<endl<<"not 100% "<<endl;
         return 0;
         }

  */
```

*AVERAGING OF AVERAGING OF AVERAGING SCHEME (Global):*

```
#include <iostream>
#include <iomanip>
#include <cstdlib>
#include <ctime>
#include <cctype>
#include <cmath>

using std:: cout;
using std:: cin;
using std:: endl;
using std:: setprecision;
using std:: setw;

double leading(double);



int main()


{

   double DUPERmax = 45;
        double DUPERmin = 40;
  //int SUPERmax = 1;-----------------> this is now to vary from DUPERmin to
DUPERmax
        double SUPERmin = 30;
  //int maxUB = 1 ;  -----------------> this is now to vary from SUPERmin  to
SUPERmax
        double minUB = 1 ;     // SAME for ALL cities! we have no flexibility
here, another much more sophisticated scheme would be needed to handle
such additional variability
        double LB = 1 ;        // Lower Bound for ALL intervals!    and typically
assumed to be 1
```

```
if((DUPERmax<DUPERmin)||(DUPERmin<SUPERmin)||(SUPERmin<minU
B)||(minUB<LB))
{
cout<<endl<<"wrong parameter input selection"<<endl;
return 0;
}


        double distribution[10]={0}; //street
        double ALLdist[10]={0};       //city
        double avgSUPER[10]={0};     //country
        double avgDUPER[10]={0};     //global

        int a;




for(int overall=DUPERmin ; overall<(DUPERmax+1) ; overall++)
{


for(double view=SUPERmin ; view<(overall+1) ; view++)  //        CITIES
view=logest street in the city
{


//=================================================

for(double q=minUB ; q<(view+1) ; q++)   //                    STREETS
q=street's length
{

        for(double j=LB ; j<(q+1) ; j++)  // j IS the integer!        HOUSES
j=house's number
        {
        a=0;
        a=a+leading(j);
   distribution[a]++;
        }
```

```
   for(int s=0 ; s<10 ; s++)
   distribution[s]=distribution[s]/(q-LB+1); // for a particular street q houses
long

         for( s=0;s<10;s++)
         ALLdist[s]=ALLdist[s]+distribution[s]; // for the particular city
(having view as the max UB) we add this street

         for( s=0 ; s<10 ; s++)
   distribution[s]=0;

}

//==============================================

   for(int s=0;s<10;s++)
         ALLdist[s]=ALLdist[s]/(view - minUB +1);   // averaging out al streets
for this particular city

         for( s=0;s<10;s++)
         avgSUPER[s]=avgSUPER[s]+ALLdist[s];   // and now we add this city
to the long country log/records

         for( s=0;s<10;s++)
         ALLdist[s]=0;   /*cout<<view<<endl;*/  // clear out - empty all info -
for the next new generic city

}
```

```
	for(int j=0 ; j<10 ; j++)
	avgSUPER[j]=avgSUPER[j]/(overall-SUPERmin+1);


	for(int s=0;s<10;s++)
	avgDUPER[s]=avgDUPER[s]+avgSUPER[s];


	for( s=0;s<10;s++)
	avgSUPER[s]=0;


 //	if(overall%50==0)
   cout<<overall<<" now, and will end at "<<DUPERmax<<endl;
//  if(overall>900)
//if(overall%10==0)
//   cout<<overall<<" now, and will end at "<<DUPERmax<<endl;



}




	cout<<endl;

	for(int j=0 ; j<10 ; j++)
	{   cout<<"% of "<<j<<" leading digit  ---> ";
		cout<<avgDUPER[j]/(DUPERmax-DUPERmin+1)<<endl;

	}

	return 0;
}
```

```
//-----------------------------------------------------
double leading(double x)
{
         double small,lead;

    for(int i=1  ; ; i++)
        {
                if(x<10)
                {small=x; break;}
                else
                x=(x/10);
        }

        lead=floor(small);
        return lead;
}//-----------------------------------------------------

//double totalP;

/*
   // checking for 100% compliance
        cout<<endl<<q<<endl;
   totalP=0;
        for( s=0;s<10;s++)
        totalP=totalP+distribution[s];
        if(abs(totalP-10>0.001))
        {
        cout<<endl<<"not 100% "<<endl;
        return 0;
        }
  */
```

# <u>Prime Numbers are NOT Benford:</u>

| | | | | | | | | | |
|---|---|---|---|---|---|---|---|---|---|
| 2 | 3 | 5 | 7 | 11 | 13 | 17 | 19 | 23 | 29 |
| 31 | 37 | 41 | 43 | 47 | 53 | 59 | 61 | 67 | 71 |
| 73 | 79 | 83 | 89 | 97 | 101 | 103 | 107 | 109 | 113 |
| 127 | 131 | 137 | 139 | 149 | 151 | 157 | 163 | 167 | 173 |
| 179 | 181 | 191 | 193 | 197 | 199 | 211 | 223 | 227 | 229 |
| 233 | 239 | 241 | 251 | 257 | 263 | 269 | 271 | 277 | 281 |
| 283 | 293 | 307 | 311 | 313 | 317 | 331 | 337 | 347 | 349 |
| 353 | 359 | 367 | 373 | 379 | 383 | 389 | 397 | 401 | 409 |
| 419 | 421 | 431 | 433 | 439 | 443 | 449 | 457 | 461 | 463 |
| 467 | 479 | 487 | 491 | 499 | 503 | 509 | 521 | 523 | 541 |
| 547 | 557 | 563 | 569 | 571 | 577 | 587 | 593 | 599 | 601 |
| 607 | 613 | 617 | 619 | 631 | 641 | 643 | 647 | 653 | 659 |
| 661 | 673 | 677 | 683 | 691 | 701 | 709 | 719 | 727 | 733 |
| 739 | 743 | 751 | 757 | 761 | 769 | 773 | 787 | 797 | 809 |
| 811 | 821 | 823 | 827 | 829 | 839 | 853 | 857 | 859 | 863 |
| 877 | 881 | 883 | 887 | 907 | 911 | 919 | 929 | 937 | 941 |
| 947 | 953 | 967 | 971 | 977 | 983 | 991 | 997 | 1009 | 1013 |
| 1019 | 1021 | 1031 | 1033 | 1039 | 1049 | 1051 | 1061 | 1063 | 1069 |
| 1087 | 1091 | 1093 | 1097 | 1103 | 1109 | 1117 | 1123 | 1129 | 1151 |
| 1153 | 1163 | 1171 | 1181 | 1187 | 1193 | 1201 | 1213 | 1217 | 1223 |
| 1229 | 1231 | 1237 | 1249 | 1259 | 1277 | 1279 | 1283 | 1289 | 1291 |
| 1297 | 1301 | 1303 | 1307 | 1319 | 1321 | 1327 | 1361 | 1367 | 1373 |
| 1381 | 1399 | 1409 | 1423 | 1427 | 1429 | 1433 | 1439 | 1447 | 1451 |
| 1453 | 1459 | 1471 | 1481 | 1483 | 1487 | 1489 | 1493 | 1499 | 1511 |
| 1523 | 1531 | 1543 | 1549 | 1553 | 1559 | 1567 | 1571 | 1579 | 1583 |
| 1597 | 1601 | 1607 | 1609 | 1613 | 1619 | 1621 | 1627 | 1637 | 1657 |
| 1663 | 1667 | 1669 | 1693 | 1697 | 1699 | 1709 | 1721 | 1723 | 1733 |
| 1741 | 1747 | 1753 | 1759 | 1777 | 1783 | 1787 | 1789 | 1801 | 1811 |
| 1823 | 1831 | 1847 | 1861 | 1867 | 1871 | 1873 | 1877 | 1879 | 1889 |
| 1901 | 1907 | 1913 | 1931 | 1933 | 1949 | 1951 | 1973 | 1979 | 1987 |
| 1993 | 1997 | 1999 | 2003 | 2011 | 2017 | 2027 | 2029 | 2039 | 2053 |
| 2063 | 2069 | 2081 | 2083 | 2087 | 2089 | 2099 | 2111 | 2113 | 2129 |
| 2131 | 2137 | 2141 | 2143 | 2153 | 2161 | 2179 | 2203 | 2207 | 2213 |
| 2221 | 2237 | 2239 | 2243 | 2251 | 2267 | 2269 | 2273 | 2281 | 2287 |
| 2293 | 2297 | 2309 | 2311 | 2333 | 2339 | 2341 | 2347 | 2351 | 2357 |
| 2371 | 2377 | 2381 | 2383 | 2389 | 2393 | 2399 | 2411 | 2417 | 2423 |
| 2437 | 2441 | 2447 | 2459 | 2467 | 2473 | 2477 | 2503 | 2521 | 2531 |
| 2539 | 2543 | 2549 | 2551 | 2557 | 2579 | 2591 | 2593 | 2609 | 2617 |
| 2621 | 2633 | 2647 | 2657 | 2659 | 2663 | 2671 | 2677 | 2683 | 2687 |
| 2689 | 2693 | 2699 | 2707 | 2711 | 2713 | 2719 | 2729 | 2731 | 2741 |
| 2749 | 2753 | 2767 | 2777 | 2789 | 2791 | 2797 | 2801 | 2803 | 2819 |
| 2833 | 2837 | 2843 | 2851 | 2857 | 2861 | 2879 | 2887 | 2897 | 2903 |
| 2909 | 2917 | 2927 | 2939 | 2953 | 2957 | 2963 | 2969 | 2971 | 2999 |
| 3001 | 3011 | 3019 | 3023 | 3037 | 3041 | 3049 | 3061 | 3067 | 3079 |
| 3083 | 3089 | 3109 | 3119 | 3121 | 3137 | 3163 | 3167 | 3169 | 3181 |
| 3187 | 3191 | 3203 | 3209 | 3217 | 3221 | 3229 | 3251 | 3253 | 3257 |
| 3259 | 3271 | 3299 | 3301 | 3307 | 3313 | 3319 | 3323 | 3329 | 3331 |
| 3343 | 3347 | 3359 | 3361 | 3371 | 3373 | 3389 | 3391 | 3407 | 3413 |
| 3433 | 3449 | 3457 | 3461 | 3463 | 3467 | 3469 | 3491 | 3499 | 3511 |
| 3517 | 3527 | 3529 | 3533 | 3539 | 3541 | 3547 | 3557 | 3559 | 3571 |
| 3581 | 3583 | 3593 | 3607 | 3613 | 3617 | 3623 | 3631 | 3637 | 3643 |
| 3659 | 3671 | 3673 | 3677 | 3691 | 3697 | 3701 | 3709 | 3719 | 3727 |
| 3733 | 3739 | 3761 | 3767 | 3769 | 3779 | 3793 | 3797 | 3803 | 3821 |
| 3823 | 3833 | 3847 | 3851 | 3853 | 3863 | 3877 | 3881 | 3889 | 3907 |
| 3911 | 3917 | 3919 | 3923 | 3929 | 3931 | 3943 | 3947 | 3967 | 3989 |
| 4001 | 4003 | 4007 | 4013 | 4019 | 4021 | 4027 | 4049 | 4051 | 4057 |
| 4073 | 4079 | 4091 | 4093 | 4099 | 4111 | 4127 | 4129 | 4133 | 4139 |
| 4153 | 4157 | 4159 | 4177 | 4201 | 4211 | 4217 | 4219 | 4229 | 4231 |
| 4241 | 4243 | 4253 | 4259 | 4261 | 4271 | 4273 | 4283 | 4289 | 4297 |
| 4327 | 4337 | 4339 | 4349 | 4357 | 4363 | 4373 | 4391 | 4397 | 4409 |
| 4421 | 4423 | 4441 | 4447 | 4451 | 4457 | 4463 | 4481 | 4483 | 4493 |
| 4507 | 4513 | 4517 | 4519 | 4523 | 4547 | 4549 | 4561 | 4567 | 4583 |
| 4591 | 4597 | 4603 | 4621 | 4637 | 4639 | 4643 | 4649 | 4651 | 4657 |
| 4663 | 4673 | 4679 | 4691 | 4703 | 4721 | 4723 | 4729 | 4733 | 4751 |
| 4759 | 4783 | 4787 | 4789 | 4793 | 4799 | 4801 | 4813 | 4817 | 4831 |
| 4861 | 4871 | 4877 | 4889 | 4903 | 4909 | 4919 | 4931 | 4933 | 4937 |
| 4943 | 4951 | 4957 | 4967 | 4969 | 4973 | 4987 | 4993 | 4999 | 5003 |
| 5009 | 5011 | 5021 | 5023 | 5039 | 5051 | 5059 | 5077 | 5081 | 5087 |
| 5099 | 5101 | 5107 | 5113 | 5119 | 5147 | 5153 | 5167 | 5171 | 5179 |
| 5189 | 5197 | 5209 | 5227 | 5231 | 5233 | 5237 | 5261 | 5273 | 5279 |
| 5281 | 5297 | 5303 | 5309 | 5323 | 5333 | 5347 | 5351 | 5381 | 5387 |
| 5393 | 5399 | 5407 | 5413 | 5417 | 5419 | 5431 | 5437 | 5441 | 5443 |
| 5449 | 5471 | 5477 | 5479 | 5483 | 5501 | 5503 | 5507 | 5519 | 5521 |
| 5527 | 5531 | 5557 | 5563 | 5569 | 5573 | 5581 | 5591 | 5623 | 5639 |
| 5641 | 5647 | 5651 | 5653 | 5657 | 5659 | 5669 | 5683 | 5689 | 5693 |
| 5701 | 5711 | 5717 | 5737 | 5741 | 5743 | 5749 | 5779 | 5783 | 5791 |
| 5801 | 5807 | 5813 | 5821 | 5827 | 5839 | 5843 | 5849 | 5851 | 5857 |
| 5861 | 5867 | 5869 | 5879 | 5881 | 5897 | 5903 | 5923 | 5927 | 5939 |
| 5953 | 5981 | 5987 | 6007 | 6011 | 6029 | 6037 | 6043 | 6047 | 6053 |
| 6067 | 6073 | 6079 | 6089 | 6091 | 6101 | 6113 | 6121 | 6131 | 6133 |
| 6143 | 6151 | 6163 | 6173 | 6197 | 6199 | 6203 | 6211 | 6217 | 6221 |
| 6229 | 6247 | 6257 | 6263 | 6269 | 6271 | 6277 | 6287 | 6299 | 6301 |
| 6311 | 6317 | 6323 | 6329 | 6337 | 6343 | 6353 | 6359 | 6361 | 6367 |
| 6373 | 6379 | 6389 | 6397 | 6421 | 6427 | 6449 | 6451 | 6469 | 6473 |
| 6481 | 6491 | 6521 | 6529 | 6547 | 6551 | 6553 | 6563 | 6569 | 6571 |
| 6577 | 6581 | 6599 | 6607 | 6619 | 6637 | 6653 | 6659 | 6661 | 6673 |
| 6679 | 6689 | 6691 | 6701 | 6703 | 6709 | 6719 | 6733 | 6737 | 6761 |
| 6763 | 6779 | 6781 | 6791 | 6793 | 6803 | 6823 | 6827 | 6829 | 6833 |
| 6841 | 6857 | 6863 | 6869 | 6871 | 6883 | 6899 | 6907 | 6911 | 6917 |
| 6947 | 6949 | 6959 | 6961 | 6967 | 6971 | 6977 | 6983 | 6991 | 6997 |
| 7001 | 7013 | 7019 | 7027 | 7039 | 7043 | 7057 | 7069 | 7079 | 7103 |
| 7109 | 7121 | 7127 | 7129 | 7151 | 7159 | 7177 | 7187 | 7193 | 7207 |
| 7211 | 7213 | 7219 | 7229 | 7237 | 7243 | 7247 | 7253 | 7283 | 7297 |
| 7307 | 7309 | 7321 | 7331 | 7333 | 7349 | 7351 | 7369 | 7393 | 7411 |
| 7417 | 7433 | 7451 | 7457 | 7459 | 7477 | 7481 | 7487 | 7489 | 7499 |
| 7507 | 7517 | 7523 | 7529 | 7537 | 7541 | 7547 | 7549 | 7559 | 7561 |
| 7573 | 7577 | 7583 | 7589 | 7591 | 7603 | 7607 | 7621 | 7639 | 7643 |
| 7649 | 7669 | 7673 | 7681 | 7687 | 7691 | 7699 | 7703 | 7717 | 7723 |
| 7727 | 7741 | 7753 | 7757 | 7759 | 7789 | 7793 | 7817 | 7823 | 7829 |
| 7841 | 7853 | 7867 | 7873 | 7877 | 7879 | 7883 | 7901 | 7907 | 7919 |
| 7927 | 7933 | 7937 | 7949 | 7951 | 7963 | 7993 | 8009 | 8011 | 8017 |
| 8039 | 8053 | 8059 | 8069 | 8081 | 8087 | 8089 | 8093 | 8101 | 8111 |
| 8117 | 8123 | 8147 | 8161 | 8167 | 8171 | 8179 | 8191 | 8209 | 8219 |
| 8221 | 8231 | 8233 | 8237 | 8243 | 8263 | 8269 | 8273 | 8287 | 8291 |
| 8293 | 8297 | 8311 | 8317 | 8329 | 8353 | 8363 | 8369 | 8377 | 8387 |
| 8389 | 8419 | 8423 | 8429 | 8431 | 8443 | 8447 | 8461 | 8467 | 8501 |
| 8513 | 8521 | 8527 | 8537 | 8539 | 8543 | 8563 | 8573 | 8581 | 8597 |
| 8599 | 8609 | 8623 | 8627 | 8629 | 8641 | 8647 | 8663 | 8669 | 8677 |
| 8681 | 8689 | 8693 | 8699 | 8707 | 8713 | 8719 | 8731 | 8737 | 8741 |
| 8747 | 8753 | 8761 | 8779 | 8783 | 8803 | 8807 | 8819 | 8821 | 8831 |
| 8837 | 8839 | 8849 | 8861 | 8863 | 8867 | 8887 | 8893 | 8923 | 8929 |
| 8933 | 8941 | 8951 | 8963 | 8969 | 8971 | 8999 | 9001 | 9007 | 9011 |
| 9013 | 9029 | 9041 | 9043 | 9049 | 9059 | 9067 | 9091 | 9103 | 9109 |
| 9127 | 9133 | 9137 | 9151 | 9157 | 9161 | 9173 | 9181 | 9187 | 9199 |
| 9203 | 9209 | 9221 | 9227 | 9239 | 9241 | 9257 | 9277 | 9281 | 9283 |
| 9293 | 9311 | 9319 | 9323 | 9337 | 9341 | 9343 | 9349 | 9371 | 9377 |
| 9391 | 9397 | 9403 | 9413 | 9419 | 9421 | 9431 | 9433 | 9437 | 9439 |
| 9461 | 9463 | 9467 | 9473 | 9479 | 9491 | 9497 | 9511 | 9521 | 9533 |
| 9539 | 9547 | 9551 | 9587 | 9601 | 9613 | 9619 | 9623 | 9629 | 9631 |
| 9643 | 9649 | 9661 | 9677 | 9679 | 9689 | 9697 | 9719 | 9721 | 9733 |
| 9739 | 9743 | 9749 | 9767 | 9769 | 9781 | 9787 | 9791 | 9803 | 9811 |
| 9817 | 9829 | 9833 | 9839 | 9851 | 9857 | 9859 | 9871 | 9883 | 9887 |
| 9901 | 9907 | 9923 | 9929 | 9931 | 9941 | 9949 | 9967 | 9973 | 10007 |
| 10009 | 10037 | 10039 | 10061 | 10067 | 10069 | 10079 | 10091 | 10093 | 10099 |
| 10103 | 10111 | 10133 | 10139 | 10141 | 10151 | 10159 | 10163 | 10169 | 10177 |
| 10181 | 10193 | 10211 | 10223 | 10243 | 10247 | 10253 | 10259 | 10267 | 10271 |
| 10273 | 10289 | 10301 | 10303 | 10313 | 10321 | 10331 | 10333 | 10337 | 10343 |
| 10357 | 10369 | 10391 | 10399 | 10427 | 10429 | 10433 | 10453 | 10457 | 10459 |
| 10463 | 10477 | 10487 | 10499 | 10501 | 10513 | 10529 | 10531 | 10559 | 10567 |
| 10589 | 10597 | 10601 | 10607 | 10613 | 10627 | 10631 | 10639 | 10651 | 10657 |
| 10663 | 10667 | 10687 | 10691 | 10709 | 10711 | 10723 | 10729 | 10733 | 10739 |
| 10753 | 10771 | 10781 | 10789 | 10799 | 10831 | 10837 | 10847 | 10853 | 10859 |
| 10861 | 10867 | 10883 | 10889 | 10891 | 10903 | 10909 | 10937 | 10939 | 10949 |
| 10957 | 10973 | 10979 | 10987 | 10993 | 11003 | 11027 | 11047 | 11057 | 11059 |
| 11069 | 11071 | 11083 | 11087 | 11093 | 11113 | 11117 | 11119 | 11131 | 11149 |
| 11159 | 11161 | 11171 | 11173 | 11177 | 11197 | 11213 | 11239 | 11243 | 11251 |
| 11257 | 11261 | 11273 | 11279 | 11287 | 11299 | 11311 | 11317 | 11321 | 11329 |
| 11351 | 11353 | 11369 | 11383 | 11393 | 11399 | 11411 | 11423 | 11437 | 11443 |
| 11447 | 11467 | 11471 | 11483 | 11489 | 11491 | 11497 | 11503 | 11519 | 11527 |
| 11549 | 11551 | 11579 | 11587 | 11593 | 11597 | 11617 | 11621 | 11633 | 11657 |
| 11677 | 11681 | 11689 | 11699 | 11701 | 11717 | 11719 | 11731 | 11743 | 11777 |
| 11779 | 11783 | 11789 | 11801 | 11807 | 11813 | 11821 | 11827 | 11831 | 11833 |
| 11839 | 11863 | 11867 | 11887 | 11897 | 11903 | 11909 | 11923 | 11927 | 11933 |
| 11939 | 11941 | 11953 | 11959 | 11969 | 11971 | 11981 | 11987 | 12007 | 12011 |
| 12037 | 12041 | 12043 | 12049 | 12071 | 12073 | 12097 | 12101 | 12107 | 12109 |
| 12113 | 12119 | 12143 | 12149 | 12157 | 12161 | 12163 | 12197 | 12203 | 12211 |
| 12227 | 12239 | 12241 | 12251 | 12253 | 12263 | 12269 | 12277 | 12281 | 12289 |
| 12301 | 12323 | 12329 | 12343 | 12347 | 12373 | 12377 | 12379 | 12391 | 12401 |
| 12409 | 12413 | 12421 | 12433 | 12437 | 12451 | 12457 | 12473 | 12479 | 12487 |
| 12491 | 12497 | 12503 | 12511 | 12517 | 12527 | 12539 | 12541 | 12547 | 12553 |
| 12569 | 12577 | 12583 | 12589 | 12601 | 12611 | 12613 | 12619 | 12637 | 12641 |
| 12647 | 12653 | 12659 | 12671 | 12689 | 12697 | 12703 | 12713 | 12721 | 12739 |
| 12743 | 12757 | 12763 | 12781 | 12791 | 12799 | 12809 | 12821 | 12823 | 12829 |
| 12841 | 12853 | 12889 | 12893 | 12899 | 12907 | 12911 | 12917 | 12919 | 12923 |
| 12941 | 12953 | 12959 | 12967 | 12973 | 12979 | 12983 | 13001 | 13003 | 13007 |
| 13009 | 13033 | 13037 | 13043 | 13049 | 13063 | 13093 | 13099 | 13103 | 13109 |
| 13121 | 13127 | 13147 | 13151 | 13159 | 13163 | 13171 | 13177 | 13183 | 13187 |
| 13217 | 13219 | 13229 | 13241 | 13249 | 13259 | 13267 | 13291 | 13297 | 13309 |
| 13313 | 13327 | 13331 | 13337 | 13339 | 13367 | 13381 | 13397 | 13399 | 13411 |
| 13417 | 13421 | 13441 | 13451 | 13457 | 13463 | 13469 | 13477 | 13487 | 13499 |
| 13513 | 13523 | 13537 | 13553 | 13567 | 13577 | 13591 | 13597 | 13613 | 13619 |
| 13627 | 13633 | 13649 | 13669 | 13679 | 13681 | 13687 | 13691 | 13693 | 13697 |
| 13709 | 13711 | 13721 | 13723 | 13729 | 13751 | 13757 | 13759 | 13763 | 13781 |
| 13789 | 13799 | 13807 | 13829 | 13831 | 13841 | 13859 | 13873 | 13877 | 13879 |
| 13883 | 13901 | 13903 | 13907 | 13913 | 13921 | 13931 | 13933 | 13963 | 13967 |
| 13997 | 13999 | 14009 | 14011 | 14029 | 14033 | 14051 | 14057 | 14071 | 14081 |
| 14083 | 14087 | 14107 | 14143 | 14149 | 14153 | 14159 | 14173 | 14177 | 14197 |
| 14207 | 14221 | 14243 | 14249 | 14251 | 14281 | 14293 | 14303 | 14321 | 14323 |
| 14327 | 14341 | 14347 | 14369 | 14387 | 14389 | 14401 | 14407 | 14411 | 14419 |
| 14423 | 14431 | 14437 | 14447 | 14449 | 14461 | 14479 | 14489 | 14503 | 14519 |
| 14533 | 14537 | 14543 | 14549 | 14551 | 14557 | 14561 | 14563 | 14591 | 14593 |
| 14621 | 14627 | 14629 | 14633 | 14639 | 14653 | 14657 | 14669 | 14683 | 14699 |
| 14713 | 14717 | 14723 | 14731 | 14737 | 14741 | 14747 | 14753 | 14759 | 14767 |
| 14771 | 14779 | 14783 | 14797 | 14813 | 14821 | 14827 | 14831 | 14843 | 14851 |
| 14867 | 14869 | 14879 | 14887 | 14891 | 14897 | 14923 | 14929 | 14939 | 14947 |
| 14951 | 14957 | 14969 | 14983 | 15013 | 15017 | 15031 | 15053 | 15061 | 15073 |
| 15077 | 15083 | 15091 | 15101 | 15107 | 15121 | 15131 | 15137 | 15139 | 15149 |
| 15161 | 15173 | 15187 | 15193 | 15199 | 15217 | 15227 | 15233 | 15241 | 15259 |
| 15263 | 15269 | 15271 | 15277 | 15287 | 15289 | 15299 | 15307 | 15313 | 15319 |
| 15329 | 15331 | 15349 | 15359 | 15361 | 15373 | 15377 | 15383 | 15391 | 15401 |
| 15413 | 15427 | 15439 | 15443 | 15451 | 15461 | 15467 | 15473 | 15493 | 15497 |
| 15511 | 15527 | 15541 | 15551 | 15559 | 15569 | 15581 | 15583 | 15601 | 15607 |
| 15619 | 15629 | 15641 | 15643 | 15647 | 15649 | 15661 | 15667 | 15671 | 15679 |
| 15683 | 15727 | 15731 | 15733 | 15737 | 15739 | 15749 | 15761 | 15767 | 15773 |
| 15787 | 15791 | 15797 | 15803 | 15809 | 15817 | 15823 | 15859 | 15877 | 15881 |
| 15887 | 15889 | 15901 | 15907 | 15913 | 15919 | 15923 | 15937 | 15959 | 15971 |
| 15973 | 15991 | 16001 | 16007 | 16033 | 16057 | 16061 | 16063 | 16067 | 16069 |
| 16073 | 16087 | 16091 | 16097 | 16103 | 16111 | 16127 | 16139 | 16141 | 16183 |
| 16187 | 16189 | 16193 | 16217 | 16223 | 16229 | 16231 | 16249 | 16253 | 16267 |
| 16273 | 16301 | 16319 | 16333 | 16339 | 16349 | 16361 | 16363 | 16369 | 16381 |
| 16411 | 16417 | 16421 | 16427 | 16433 | 16447 | 16451 | 16453 | 16477 | 16481 |
| 16487 | 16493 | 16519 | 16529 | 16547 | 16553 | 16561 | 16567 | 16573 | 16603 |
| 16607 | 16619 | 16631 | 16633 | 16649 | 16651 | 16657 | 16661 | 16673 | 16691 |
| 16693 | 16699 | 16703 | 16729 | 16741 | 16747 | 16759 | 16763 | 16787 | 16811 |
| 16823 | 16829 | 16831 | 16843 | 16871 | 16879 | 16883 | 16889 | 16901 | 16903 |
| 16921 | 16927 | 16931 | 16937 | 16943 | 16963 | 16979 | 16981 | 16987 | 16993 |
| 17011 | 17021 | 17027 | 17029 | 17033 | 17041 | 17047 | 17053 | 17077 | 17093 |
| 17099 | 17107 | 17117 | 17123 | 17137 | 17159 | 17167 | 17183 | 17189 | 17191 |
| 17203 | 17207 | 17209 | 17231 | 17239 | 17257 | 17291 | 17293 | 17299 | 17317 |
| 17321 | 17327 | 17333 | 17341 | 17351 | 17359 | 17377 | 17383 | 17387 | 17389 |
| 17393 | 17401 | 17417 | 17419 | 17431 | 17443 | 17449 | 17467 | 17471 | 17477 |
| 17483 | 17489 | 17491 | 17497 | 17509 | 17519 | 17539 | 17551 | 17569 | 17573 |
| 17579 | 17581 | 17597 | 17599 | 17609 | 17623 | 17627 | 17657 | 17659 | 17669 |
| 17681 | 17683 | 17707 | 17713 | 17729 | 17737 | 17747 | 17749 | 17761 | 17783 |
| 17789 | 17791 | 17807 | 17827 | 17837 | 17839 | 17851 | 17863 | 17881 | 17891 |
| 17903 | 17909 | 17911 | 17921 | 17923 | 17929 | 17939 | 17957 | 17959 | 17971 |
| 17977 | 17981 | 17987 | 17989 | 18013 | 18041 | 18043 | 18047 | 18049 | 18059 |
| 18061 | 18077 | 18089 | 18097 | 18119 | 18121 | 18127 | 18131 | 18133 | 18143 |
| 18149 | 18169 | 18181 | 18191 | 18199 | 18211 | 18217 | 18223 | 18229 | 18233 |
| 18251 | 18253 | 18257 | 18269 | 18287 | 18289 | 18301 | 18307 | 18311 | 18313 |
| 18329 | 18341 | 18353 | 18367 | 18371 | 18379 | 18397 | 18401 | 18413 | 18427 |
| 18433 | 18439 | 18443 | 18451 | 18457 | 18461 | 18481 | 18493 | 18503 | 18517 |
| 18521 | 18523 | 18539 | 18541 | 18553 | 18583 | 18587 | 18593 | 18617 | 18637 |
| 18661 | 18671 | 18679 | 18691 | 18701 | 18713 | 18719 | 18731 | 18743 | 18749 |
| 18757 | 18773 | 18787 | 18793 | 18797 | 18803 | 18839 | 18859 | 18869 | 18899 |
| 18911 | 18913 | 18917 | 18919 | 18947 | 18959 | 18973 | 18979 | 19001 | 19009 |
| 19013 | 19031 | 19037 | 19051 | 19069 | 19073 | 19079 | 19081 | 19087 | 19121 |
| 19139 | 19141 | 19157 | 19163 | 19181 | 19183 | 19207 | 19211 | 19213 | 19219 |
| 19231 | 19237 | 19249 | 19259 | 19267 | 19273 | 19289 | 19301 | 19309 | 19319 |
| 19333 | 19373 | 19379 | 19381 | 19387 | 19391 | 19403 | 19417 | 19421 | 19423 |
| 19427 | 19429 | 19433 | 19441 | 19447 | 19457 | 19463 | 19469 | 19471 | 19477 |
| 19483 | 19489 | 19501 | 19507 | 19531 | 19541 | 19543 | 19553 | 19559 | 19571 |
| 19577 | 19583 | 19597 | 19603 | 19609 | 19661 | 19681 | 19687 | 19697 | 19699 |
| 19709 | 19717 | 19727 | 19739 | 19751 | 19753 | 19759 | 19763 | 19777 | 19793 |
| 19801 | 19813 | 19819 | 19841 | 19843 | 19853 | 19861 | 19867 | 19889 | 19891 |
| 19913 | 19919 | 19927 | 19937 | 19949 | 19961 | 19963 | 19973 | 19979 | 19991 |
| 19993 | 19997 | 20011 | 20021 | 20023 | 20029 | 20047 | 20051 | 20063 | 20071 |
| 20089 | 20101 | 20107 | 20113 | 20117 | 20123 | 20129 | 20143 | 20147 | 20149 |
| 20161 | 20173 | 20177 | 20183 | 20201 | 20219 | 20231 | 20233 | 20249 | 20261 |
| 20269 | 20287 | 20297 | 20323 | 20327 | 20333 | 20341 | 20347 | 20353 | 20357 |
| 20359 | 20369 | 20389 | 20393 | 20399 | 20407 | 20411 | 20431 | 20441 | 20443 |
| 20477 | 20479 | 20483 | 20507 | 20509 | 20521 | 20533 | 20543 | 20549 | 20551 |
| 20563 | 20593 | 20599 | 20611 | 20627 | 20639 | 20641 | 20663 | 20681 | 20693 |
| 20707 | 20717 | 20719 | 20731 | 20743 | 20747 | 20749 | 20753 | 20759 | 20771 |
| 20773 | 20789 | 20807 | 20809 | 20849 | 20857 | 20873 | 20879 | 20887 | 20897 |
| 20899 | 20903 | 20921 | 20929 | 20939 | 20947 | 20959 | 20963 | 20981 | 20983 |
| 21001 | 21011 | 21013 | 21017 | 21019 | 21023 | 21031 | 21059 | 21061 | 21067 |
| 21089 | 21101 | 21107 | 21121 | 21139 | 21143 | 21149 | 21157 | 21163 | 21169 |
| 21179 | 21187 | 21191 | 21193 | 21211 | 21221 | 21227 | 21247 | 21269 | 21277 |
| 21283 | 21313 | 21317 | 21319 | 21323 | 21341 | 21347 | 21377 | 21379 | 21383 |
| 21391 | 21397 | 21401 | 21407 | 21419 | 21433 | 21467 | 21481 | 21487 | 21491 |
| 21493 | 21499 | 21503 | 21517 | 21521 | 21523 | 21529 | 21557 | 21559 | 21563 |
| 21569 | 21577 | 21587 | 21589 | 21599 | 21601 | 21611 | 21613 | 21617 | 21647 |
| 21649 | 21661 | 21673 | 21683 | 21701 | 21713 | 21727 | 21737 | 21739 | 21751 |
| 21757 | 21767 | 21773 | 21787 | 21799 | 21803 | 21817 | 21821 | 21839 | 21841 |
| 21851 | 21859 | 21863 | 21871 | 21881 | 21893 | 21911 | 21929 | 21937 | 21943 |
| 21961 | 21977 | 21991 | 21997 | 22003 | 22013 | 22027 | 22031 | 22037 | 22039 |
| 22051 | 22063 | 22067 | 22073 | 22079 | 22091 | 22093 | 22109 | 22111 | 22123 |
| 22129 | 22133 | 22147 | 22153 | 22157 | 22159 | 22171 | 22189 | 22193 | 22229 |
| 22247 | 22259 | 22271 | 22273 | 22277 | 22279 | 22283 | 22291 | 22303 | 22307 |
| 22343 | 22349 | 22367 | 22369 | 22381 | 22391 | 22397 | 22409 | 22433 | 22441 |
| 22447 | 22453 | 22469 | 22481 | 22483 | 22501 | 22511 | 22531 | 22541 | 22543 |
| 22549 | 22567 | 22571 | 22573 | 22613 | 22619 | 22621 | 22637 | 22639 | 22643 |
| 22651 | 22669 | 22679 | 22691 | 22697 | 22699 | 22709 | 22717 | 22721 | 22727 |
| 22739 | 22741 | 22751 | 22769 | 22777 | 22783 | 22787 | 22807 | 22811 | 22817 |
| 22853 | 22859 | 22861 | 22871 | 22877 | 22901 | 22907 | 22921 | 22937 | 22943 |
| 22961 | 22963 | 22973 | 22993 | 23003 | 23011 | 23017 | 23021 | 23027 | 23029 |
| 23039 | 23041 | 23053 | 23057 | 23059 | 23063 | 23071 | 23081 | 23087 | 23099 |
| 23117 | 23131 | 23143 | 23159 | 23167 | 23173 | 23189 | 23197 | 23201 | 23203 |
| 23209 | 23227 | 23251 | 23269 | 23279 | 23291 | 23293 | 23297 | 23311 | 23321 |
| 23327 | 23333 | 23339 | 23357 | 23369 | 23371 | 23399 | 23417 | 23431 | 23447 |
| 23459 | 23473 | 23497 | 23509 | 23531 | 23537 | 23539 | 23549 | 23557 | 23561 |
| 23563 | 23567 | 23581 | 23593 | 23599 | 23603 | 23609 | 23623 | 23627 | 23629 |
| 23633 | 23663 | 23669 | 23671 | 23677 | 23687 | 23689 | 23719 | 23741 | 23743 |
| 23747 | 23753 | 23761 | 23767 | 23773 | 23789 | 23801 | 23813 | 23819 | 23827 |
| 23831 | 23833 | 23857 | 23869 | 23873 | 23879 | 23887 | 23893 | 23899 | 23909 |
| 23911 | 23917 | 23929 | 23957 | 23971 | 23977 | 23981 | 23993 | 24001 | 24007 |
| 24019 | 24023 | 24029 | 24043 | 24049 | 24061 | 24071 | 24077 | 24083 | 24091 |
| 24097 | 24103 | 24107 | 24109 | 24113 | 24121 | 24133 | 24137 | 24151 | 24169 |
| 24179 | 24181 | 24197 | 24203 | 24223 | 24229 | 24239 | 24247 | 24251 | 24281 |
| 24317 | 24329 | 24337 | 24359 | 24371 | 24373 | 24379 | 24391 | 24407 | 24413 |
| 24419 | 24421 | 24439 | 24443 | 24469 | 24473 | 24481 | 24499 | 24509 | 24517 |
| 24527 | 24533 | 24547 | 24551 | 24571 | 24593 | 24611 | 24623 | 24631 | 24659 |
| 24671 | 24677 | 24683 | 24691 | 24697 | 24709 | 24733 | 24749 | 24763 | 24767 |
| 24781 | 24793 | 24799 | 24809 | 24821 | 24841 | 24847 | 24851 | 24859 | 24877 |
| 24889 | 24907 | 24917 | 24919 | 24923 | 24943 | 24953 | 24967 | 24971 | 24977 |
| 24979 | 24989 | 25013 | 25031 | 25033 | 25037 | 25057 | 25073 | 25087 | 25097 |
| 25111 | 25117 | 25121 | 25127 | 25147 | 25153 | 25163 | 25169 | 25171 | 25183 |
| 25189 | 25219 | 25229 | 25237 | 25243 | 25247 | 25253 | 25261 | 25301 | 25303 |
| 25307 | 25309 | 25321 | 25339 | 25343 | 25349 | 25357 | 25367 | 25373 | 25391 |
| 25409 | 25411 | 25423 | 25439 | 25447 | 25453 | 25457 | 25463 | 25469 | 25471 |
| 25523 | 25537 | 25541 | 25561 | 25577 | 25579 | 25583 | 25589 | 25601 | 25603 |
| 25609 | 25621 | 25633 | 25639 | 25643 | 25657 | 25667 | 25673 | 25679 | 25693 |
| 25703 | 25717 | 25733 | 25741 | 25747 | 25759 | 25763 | 25771 | 25793 | 25799 |
| 25801 | 25819 | 25841 | 25847 | 25849 | 25867 | 25873 | 25889 | 25903 | 25913 |
| 25919 | 25931 | 25933 | 25939 | 25943 | 25951 | 25969 | 25981 | 25997 | 25999 |
| 26003 | 26017 | 26021 | 26029 | 26041 | 26053 | 26083 | 26099 | 26107 | 26111 |
| 26113 | 26119 | 26141 | 26153 | 26161 | 26171 | 26177 | 26183 | 26189 | 26203 |
| 26209 | 26227 | 26237 | 26249 | 26251 | 26261 | 26263 | 26267 | 26293 | 26297 |
| 26309 | 26317 | 26321 | 26339 | 26347 | 26357 | 26371 | 26387 | 26393 | 26399 |
| 26407 | 26417 | 26423 | 26431 | 26437 | 26449 | 26459 | 26479 | 26489 | 26497 |
| 26501 | 26513 | 26539 | 26557 | 26561 | 26573 | 26591 | 26597 | 26627 | 26633 |
| 26641 | 26647 | 26669 | 26681 | 26683 | 26687 | 26693 | 26699 | 26701 | 26711 |
| 26713 | 26717 | 26723 | 26729 | 26731 | 26737 | 26759 | 26777 | 26783 | 26801 |
| 26813 | 26821 | 26833 | 26839 | 26849 | 26861 | 26863 | 26879 | 26881 | 26891 |
| 26893 | 26903 | 26921 | 26927 | 26947 | 26951 | 26953 | 26959 | 26981 | 26987 |
| 26993 | 27011 | 27017 | 27031 | 27043 | 27059 | 27061 | 27067 | 27073 | 27077 |
| 27091 | 27103 | 27107 | 27109 | 27127 | 27143 | 27179 | 27191 | 27197 | 27211 |
| 27239 | 27241 | 27253 | 27259 | 27271 | 27277 | 27281 | 27283 | 27299 | 27329 |
| 27337 | 27361 | 27367 | 27397 | 27407 | 27409 | 27427 | 27431 | 27437 | 27449 |
| 27457 | 27479 | 27481 | 27487 | 27509 | 27527 | 27529 | 27539 | 27541 | 27551 |
| 27581 | 27583 | 27611 | 27617 | 27631 | 27647 | 27653 | 27673 | 27689 | 27691 |
| 27697 | 27701 | 27733 | 27737 | 27739 | 27743 | 27749 | 27751 | 27763 | 27767 |
| 27773 | 27779 | 27791 | 27793 | 27799 | 27803 | 27809 | 27817 | 27823 | 27827 |
| 27847 | 27851 | 27883 | 27893 | 27901 | 27917 | 27919 | 27941 | 27943 | 27947 |
| 27953 | 27961 | 27967 | 27983 | 27997 | 28001 | 28019 | 28027 | 28031 | 28051 |
| 28057 | 28069 | 28081 | 28087 | 28097 | 28099 | 28109 | 28111 | 28123 | 28151 |
| 28163 | 28181 | 28183 | 28201 | 28211 | 28219 | 28229 | 28277 | 28279 | 28283 |
| 28289 | 28297 | 28307 | 28309 | 28319 | 28349 | 28351 | 28387 | 28393 | 28403 |
| 28409 | 28411 | 28429 | 28433 | 28439 | 28447 | 28463 | 28477 | 28493 | 28499 |
| 28513 | 28517 | 28537 | 28541 | 28547 | 28549 | 28559 | 28571 | 28573 | 28579 |
| 28591 | 28597 | 28603 | 28607 | 28619 | 28621 | 28627 | 28631 | 28643 | 28649 |
| 28657 | 28661 | 28663 | 28669 | 28687 | 28697 | 28703 | 28711 | 28723 | 28729 |
| 28751 | 28753 | 28759 | 28771 | 28789 | 28793 | 28807 | 28813 | 28817 | 28837 |
| 28843 | 28859 | 28867 | 28871 | 28879 | 28901 | 28909 | 28921 | 28927 | 28933 |
| 28949 | 28961 | 28979 | 29009 | 29017 | 29021 | 29023 | 29027 | 29033 | 29059 |
| 29063 | 29077 | 29101 | 29123 | 29129 | 29131 | 29137 | 29147 | 29153 | 29167 |
| 29173 | 29179 | 29191 | 29201 | 29207 | 29209 | 29221 | 29231 | 29243 | 29251 |
| 29269 | 29287 | 29297 | 29303 | 29311 | 29327 | 29333 | 29339 | 29347 | 29363 |
| 29383 | 29387 | 29389 | 29399 | 29401 | 29411 | 29423 | 29429 | 29437 | 29443 |
| 29453 | 29473 | 29483 | 29501 | 29527 | 29531 | 29537 | 29567 | 29569 | 29573 |
| 29581 | 29587 | 29599 | 29611 | 29629 | 29633 | 29641 | 29663 | 29669 | 29671 |
| 29683 | 29717 | 29723 | 29741 | 29753 | 29759 | 29761 | 29789 | 29803 | 29819 |
| 29833 | 29837 | 29851 | 29863 | 29867 | 29873 | 29879 | 29881 | 29917 | 29921 |
| 29927 | 29947 | 29959 | 29983 | 29989 | 30011 | 30013 | 30029 | 30047 | 30059 |
| 30071 | 30089 | 30091 | 30097 | 30103 | 30109 | 30113 | 30119 | 30133 | 30137 |
| 30139 | 30161 | 30169 | 30181 | 30187 | 30197 | 30203 | 30211 | 30223 | 30241 |
| 30253 | 30259 | 30269 | 30271 | 30293 | 30307 | 30313 | 30319 | 30323 | 30341 |
| 30347 | 30367 | 30389 | 30391 | 30403 | 30427 | 30431 | 30449 | 30467 | 30469 |
| 30491 | 30493 | 30497 | 30509 | 30517 | 30529 | 30539 | 30553 | 30557 | 30559 |
| 30577 | 30593 | 30631 | 30637 | 30643 | 30649 | 30661 | 30671 | 30677 | 30689 |
| 30697 | 30703 | 30707 | 30713 | 30727 | 30757 | 30763 | 30773 | 30781 | 30803 |
| 30809 | 30817 | 30829 | 30839 | 30841 | 30851 | 30853 | 30859 | 30869 | 30871 |
| 30881 | 30893 | 30911 | 30931 | 30937 | 30941 | 30949 | 30971 | 30977 | 30983 |
| 31013 | 31019 | 31033 | 31039 | 31051 | 31063 | 31069 | 31079 | 31081 | 31091 |
| 31121 | 31123 | 31139 | 31147 | 31151 | 31153 | 31159 | 31177 | 31181 | 31183 |
| 31189 | 31193 | 31219 | 31223 | 31231 | 31237 | 31247 | 31249 | 31253 | 31259 |
| 31267 | 31271 | 31277 | 31307 | 31319 | 31321 | 31327 | 31333 | 31337 | 31357 |
| 31379 | 31387 | 31391 | 31393 | 31397 | 31469 | 31477 | 31481 | 31489 | 31511 |
| 31513 | 31517 | 31531 | 31541 | 31543 | 31547 | 31567 | 31573 | 31583 | 31601 |
| 31607 | 31627 | 31643 | 31649 | 31657 | 31663 | 31667 | 31687 | 31699 | 31721 |
| 31723 | 31727 | 31729 | 31741 | 31751 | 31769 | 31771 | 31793 | 31799 | 31817 |
| 31847 | 31849 | 31859 | 31873 | 31883 | 31891 | 31907 | 31957 | 31963 | 31973 |
| 31981 | 31991 | 32003 | 32009 | 32027 | 32029 | 32051 | 32057 | 32059 | 32063 |
| 32069 | 32077 | 32083 | 32089 | 32099 | 32117 | 32119 | 32141 | 32143 | 32159 |
| 32173 | 32183 | 32189 | 32191 | 32203 | 32213 | 32233 | 32237 | 32251 | 32257 |
| 32261 | 32297 | 32299 | 32303 | 32309 | 32321 | 32323 | 32327 | 32341 | 32353 |
| 32359 | 32363 | 32369 | 32371 | 32377 | 32381 | 32401 | 32411 | 32413 | 32423 |
| 32429 | 32441 | 32443 | 32467 | 32479 | 32491 | 32497 | 32503 | 32507 | 32531 |
| 32533 | 32537 | 32561 | 32563 | 32569 | 32573 | 32579 | 32587 | 32603 | 32609 |
| 32611 | 32621 | 32633 | 32647 | 32653 | 32687 | 32693 | 32707 | 32713 | 32717 |
| 32719 | 32749 | 32771 | 32779 | 32783 | 32789 | 32797 | 32801 | 32803 | 32831 |
| 32833 | 32839 | 32843 | 32869 | 32887 | 32909 | 32911 | 32917 | 32933 | 32939 |
| 32941 | 32957 | 32969 | 32971 | 32983 | 32987 | 32993 | 32999 | 33013 | 33023 |
| 33029 | 33037 | 33049 | 33053 | 33071 | 33073 | 33083 | 33091 | 33107 | 33113 |
| 33119 | 33149 | 33151 | 33161 | 33179 | 33181 | 33191 | 33199 | 33203 | 33211 |
| 33223 | 33247 | 33287 | 33289 | 33301 | 33311 | 33317 | 33329 | 33331 | 33343 |
| 33347 | 33349 | 33353 | 33359 | 33377 | 33391 | 33403 | 33409 | 33413 | 33427 |
| 33457 | 33461 | 33469 | 33479 | 33487 | 33493 | 33503 | 33521 | 33529 | 33533 |
| 33547 | 33563 | 33569 | 33577 | 33581 | 33587 | 33589 | 33599 | 33601 | 33613 |
| 33617 | 33619 | 33623 | 33629 | 33637 | 33641 | 33647 | 33679 | 33703 | 33713 |
| 33721 | 33739 | 33749 | 33751 | 33757 | 33767 | 33769 | 33773 | 33791 | 33797 |
| 33809 | 33811 | 33827 | 33829 | 33851 | 33857 | 33863 | 33871 | 33889 | 33893 |
| 33911 | 33923 | 33931 | 33937 | 33941 | 33961 | 33967 | 33997 | 34019 | 34031 |
| 34033 | 34039 | 34057 | 34061 | 34123 | 34127 | 34129 | 34141 | 34147 | 34157 |
| 34159 | 34171 | 34183 | 34211 | 34213 | 34217 | 34231 | 34253 | 34259 | 34261 |
| 34267 | 34273 | 34283 | 34297 | 34301 | 34303 | 34313 | 34319 | 34327 | 34337 |
| 34351 | 34361 | 34367 | 34369 | 34381 | 34403 | 34421 | 34429 | 34439 | 34457 |
| 34469 | 34471 | 34483 | 34487 | 34499 | 34501 | 34511 | 34513 | 34519 | 34537 |
| 34543 | 34549 | 34583 | 34589 | 34591 | 34603 | 34607 | 34613 | 34631 | 34649 |
| 34651 | 34667 | 34673 | 34679 | 34687 | 34693 | 34703 | 34721 | 34729 | 34739 |
| 34747 | 34757 | 34759 | 34763 | 34781 | 34807 | 34819 | 34841 | 34843 | 34847 |
| 34849 | 34871 | 34877 | 34883 | 34897 | 34913 | 34919 | 34939 | 34949 | 34961 |
| 34963 | 34981 | 35023 | 35027 | 35051 | 35053 | 35059 | 35069 | 35081 | 35083 |
| 35089 | 35099 | 35107 | 35111 | 35117 | 35129 | 35141 | 35149 | 35153 | 35159 |
| 35171 | 35201 | 35221 | 35227 | 35251 | 35257 | 35267 | 35279 | 35281 | 35291 |
| 35311 | 35317 | 35323 | 35327 | 35339 | 35353 | 35363 | 35381 | 35393 | 35401 |
| 35407 | 35419 | 35423 | 35437 | 35447 | 35449 | 35461 | 35491 | 35507 | 35509 |
| 35521 | 35527 | 35531 | 35533 | 35537 | 35543 | 35569 | 35573 | 35591 | 35593 |
| 35597 | 35603 | 35617 | 35671 | 35677 | 35729 | 35731 | 35747 | 35753 | 35759 |
| 35771 | 35797 | 35801 | 35803 | 35809 | 35831 | 35837 | 35839 | 35851 | 35863 |
| 35869 | 35879 | 35897 | 35899 | 35911 | 35923 | 35933 | 35951 | 35963 | 35969 |
| 35977 | 35983 | 35993 | 35999 | 36007 | 36011 | 36013 | 36017 | 36037 | 36061 |
| 36067 | 36073 | 36083 | 36097 | 36107 | 36109 | 36131 | 36137 | 36151 | 36161 |
| 36187 | 36191 | 36209 | 36217 | 36229 | 36241 | 36251 | 36263 | 36269 | 36277 |
| 36293 | 36299 | 36307 | 36313 | 36319 | 36341 | 36343 | 36353 | 36373 | 36383 |
| 36389 | 36433 | 36451 | 36457 | 36467 | 36469 | 36473 | 36479 | 36493 | 36497 |
| 36523 | 36527 | 36529 | 36541 | 36551 | 36559 | 36563 | 36571 | 36583 | 36587 |
| 36599 | 36607 | 36629 | 36637 | 36643 | 36653 | 36671 | 36677 | 36683 | 36691 |
| 36697 | 36709 | 36713 | 36721 | 36739 | 36749 | 36761 | 36767 | 36779 | 36781 |
| 36787 | 36791 | 36793 | 36809 | 36821 | 36833 | 36847 | 36857 | 36871 | 36877 |
| 36887 | 36899 | 36901 | 36913 | 36919 | 36923 | 36929 | 36931 | 36943 | 36947 |
| 36973 | 36979 | 36997 | 37003 | 37013 | 37019 | 37021 | 37039 | 37049 | 37057 |
| 37061 | 37087 | 37097 | 37117 | 37123 | 37139 | 37159 | 37171 | 37181 | 37189 |
| 37199 | 37201 | 37217 | 37223 | 37243 | 37253 | 37273 | 37277 | 37307 | 37309 |
| 37313 | 37321 | 37337 | 37339 | 37357 | 37361 | 37363 | 37369 | 37379 | 37397 |
| 37409 | 37423 | 37441 | 37447 | 37463 | 37483 | 37489 | 37493 | 37501 | 37507 |
| 37511 | 37517 | 37529 | 37537 | 37547 | 37549 | 37561 | 37567 | 37571 | 37573 |
| 37579 | 37589 | 37591 | 37607 | 37619 | 37633 | 37643 | 37649 | 37657 | 37663 |
| 37691 | 37693 | 37699 | 37717 | 37747 | 37781 | 37783 | 37799 | 37811 | 37813 |
| 37831 | 37847 | 37853 | 37861 | 37871 | 37879 | 37889 | 37897 | 37907 | 37951 |
| 37957 | 37963 | 37967 | 37987 | 37991 | 37993 | 37997 | 38011 | 38039 | 38047 |
| 38053 | 38069 | 38083 | 38113 | 38119 | 38149 | 38153 | 38167 | 38177 | 38183 |
| 38189 | 38197 | 38201 | 38219 | 38231 | 38237 | 38239 | 38261 | 38273 | 38281 |
| 38287 | 38299 | 38303 | 38317 | 38321 | 38327 | 38329 | 38333 | 38351 | 38371 |
| 38377 | 38393 | 38431 | 38447 | 38449 | 38453 | 38459 | 38461 | 38501 | 38543 |
| 38557 | 38561 | 38567 | 38569 | 38593 | 38603 | 38609 | 38611 | 38629 | 38639 |
| 38651 | 38653 | 38669 | 38671 | 38677 | 38693 | 38699 | 38707 | 38711 | 38713 |
| 38723 | 38729 | 38737 | 38747 | 38749 | 38767 | 38783 | 38791 | 38803 | 38821 |
| 38833 | 38839 | 38851 | 38861 | 38867 | 38873 | 38891 | 38903 | 38917 | 38921 |
| 38923 | 38933 | 38953 | 38959 | 38971 | 38977 | 38993 | 39019 | 39023 | 39041 |
| 39043 | 39047 | 39079 | 39089 | 39097 | 39103 | 39107 | 39113 | 39119 | 39133 |
| 39139 | 39157 | 39161 | 39163 | 39181 | 39191 | 39199 | 39209 | 39217 | 39227 |
| 39229 | 39233 | 39239 | 39241 | 39251 | 39293 | 39301 | 39313 | 39317 | 39323 |
| 39341 | 39343 | 39359 | 39367 | 39371 | 39373 | 39383 | 39397 | 39409 | 39419 |
| 39439 | 39443 | 39451 | 39461 | 39499 | 39503 | 39509 | 39511 | 39521 | 39541 |
| 39551 | 39563 | 39569 | 39581 | 39607 | 39619 | 39623 | 39631 | 39659 | 39667 |
| 39671 | 39679 | 39703 | 39709 | 39719 | 39727 | 39733 | 39749 | 39761 | 39769 |
| 39779 | 39791 | 39799 | 39821 | 39827 | 39829 | 39839 | 39841 | 39847 | 39857 |
| 39863 | 39869 | 39877 | 39883 | 39887 | 39901 | 39929 | 39937 | 39953 | 39971 |
| 39979 | 39983 | 39989 | 40009 | 40013 | 40031 | 40037 | 40039 | 40063 | 40087 |
| 40093 | 40099 | 40111 | 40123 | 40127 | 40129 | 40151 | 40153 | 40163 | 40169 |
| 40177 | 40189 | 40193 | 40213 | 40231 | 40237 | 40241 | 40253 | 40277 | 40283 |
| 40289 | 40343 | 40351 | 40357 | 40361 | 40387 | 40423 | 40427 | 40429 | 40433 |
| 40459 | 40471 | 40483 | 40487 | 40493 | 40499 | 40507 | 40519 | 40529 | 40531 |
| 40543 | 40559 | 40577 | 40583 | 40591 | 40597 | 40609 | 40627 | 40637 | 40639 |
| 40693 | 40697 | 40699 | 40709 | 40739 | 40751 | 40759 | 40763 | 40771 | 40787 |
| 40801 | 40813 | 40819 | 40823 | 40829 | 40841 | 40847 | 40849 | 40853 | 40867 |
| 40879 | 40883 | 40897 | 40903 | 40927 | 40933 | 40939 | 40949 | 40961 | 40973 |
| 40993 | 41011 | 41017 | 41023 | 41039 | 41047 | 41051 | 41057 | 41077 | 41081 |
| 41113 | 41117 | 41131 | 41141 | 41143 | 41149 | 41161 | 41177 | 41179 | 41183 |
| 41189 | 41201 | 41203 | 41213 | 41221 | 41227 | 41231 | 41233 | 41243 | 41257 |
| 41263 | 41269 | 41281 | 41299 | 41333 | 41341 | 41351 | 41357 | 41381 | 41387 |
| 41389 | 41399 | 41411 | 41413 | 41443 | 41453 | 41467 | 41479 | 41491 | 41507 |
| 41513 | 41519 | 41521 | 41539 | 41543 | 41549 | 41579 | 41593 | 41597 | 41603 |
| 41609 | 41611 | 41617 | 41621 | 41627 | 41641 | 41647 | 41651 | 41659 | 41669 |
| 41681 | 41687 | 41719 | 41729 | 41737 | 41759 | 41761 | 41771 | 41777 | 41801 |
| 41809 | 41813 | 41843 | 41849 | 41851 | 41863 | 41879 | 41887 | 41893 | 41897 |
| 41903 | 41911 | 41927 | 41941 | 41947 | 41953 | 41957 | 41959 | 41969 | 41981 |
| 41983 | 41999 | 42013 | 42017 | 42019 | 42023 | 42043 | 42061 | 42071 | 42073 |
| 42083 | 42089 | 42101 | 42131 | 42139 | 42157 | 42169 | 42179 | 42181 | 42187 |
| 42193 | 42197 | 42209 | 42221 | 42223 | 42227 | 42239 | 42257 | 42281 | 42283 |
| 42293 | 42299 | 42307 | 42323 | 42331 | 42337 | 42349 | 42359 | 42373 | 42379 |
| 42391 | 42397 | 42403 | 42407 | 42409 | 42433 | 42437 | 42443 | 42451 | 42457 |
| 42461 | 42463 | 42467 | 42473 | 42487 | 42491 | 42499 | 42509 | 42533 | 42557 |
| 42569 | 42571 | 42577 | 42589 | 42611 | 42641 | 42643 | 42649 | 42667 | 42677 |
| 42683 | 42689 | 42697 | 42701 | 42703 | 42709 | 42719 | 42727 | 42737 | 42743 |
| 42751 | 42767 | 42773 | 42787 | 42793 | 42797 | 42821 | 42829 | 42839 | 42841 |
| 42853 | 42859 | 42863 | 42899 | 42901 | 42923 | 42929 | 42937 | 42943 | 42953 |
| 42961 | 42967 | 42979 | 42989 | 43003 | 43013 | 43019 | 43037 | 43049 | 43051 |
| 43063 | 43067 | 43093 | 43103 | 43117 | 43133 | 43151 | 43159 | 43177 | 43189 |
| 43201 | 43207 | 43223 | 43237 | 43261 | 43271 | 43283 | 43291 | 43313 | 43319 |
| 43321 | 43331 | 43391 | 43397 | 43399 | 43403 | 43411 | 43427 | 43441 | 43451 |
| 43457 | 43481 | 43487 | 43499 | 43517 | 43541 | 43543 | 43573 | 43577 | 43579 |
| 43591 | 43597 | 43607 | 43609 | 43613 | 43627 | 43633 | 43649 | 43651 | 43661 |
| 43669 | 43691 | 43711 | 43717 | 43721 | 43753 | 43759 | 43777 | 43781 | 43783 |
| 43787 | 43789 | 43793 | 43801 | 43853 | 43867 | 43889 | 43891 | 43913 | 43933 |
| 43943 | 43951 | 43961 | 43963 | 43969 | 43973 | 43987 | 43991 | 43997 | 44017 |
| 44021 | 44027 | 44029 | 44041 | 44053 | 44059 | 44071 | 44087 | 44089 | 44101 |
| 44111 | 44119 | 44123 | 44129 | 44131 | 44159 | 44171 | 44179 | 44189 | 44201 |
| 44203 | 44207 | 44221 | 44249 | 44257 | 44263 | 44267 | 44269 | 44273 | 44279 |
| 44281 | 44293 | 44351 | 44357 | 44371 | 44381 | 44383 | 44389 | 44417 | 44449 |
| 44453 | 44483 | 44491 | 44497 | 44501 | 44507 | 44519 | 44531 | 44533 | 44537 |
| 44543 | 44549 | 44563 | 44579 | 44587 | 44617 | 44621 | 44623 | 44633 | 44641 |
| 44647 | 44651 | 44657 | 44683 | 44687 | 44699 | 44701 | 44711 | 44729 | 44741 |
| 44753 | 44771 | 44773 | 44777 | 44789 | 44797 | 44809 | 44819 | 44839 | 44843 |
| 44851 | 44867 | 44879 | 44887 | 44893 | 44909 | 44917 | 44927 | 44939 | 44953 |
| 44959 | 44963 | 44971 | 44983 | 44987 | 45007 | 45013 | 45053 | 45061 | 45077 |
| 45083 | 45119 | 45121 | 45127 | 45131 | 45137 | 45139 | 45161 | 45179 | 45181 |
| 45191 | 45197 | 45233 | 45247 | 45259 | 45263 | 45281 | 45289 | 45293 | 45307 |
| 45317 | 45319 | 45329 | 45337 | 45341 | 45343 | 45361 | 45377 | 45389 | 45403 |
| 45413 | 45427 | 45433 | 45439 | 45481 | 45491 | 45497 | 45503 | 45523 | 45533 |
| 45541 | 45553 | 45557 | 45569 | 45587 | 45589 | 45599 | 45613 | 45631 | 45641 |
| 45659 | 45667 | 45673 | 45677 | 45691 | 45697 | 45707 | 45737 | 45751 | 45757 |
| 45763 | 45767 | 45779 | 45817 | 45821 | 45823 | 45827 | 45833 | 45841 | 45853 |
| 45863 | 45869 | 45887 | 45893 | 45943 | 45949 | 45953 | 45959 | 45971 | 45979 |
| 45989 | 46021 | 46027 | 46049 | 46051 | 46061 | 46073 | 46091 | 46093 | 46099 |
| 46103 | 46133 | 46141 | 46147 | 46153 | 46171 | 46181 | 46183 | 46187 | 46199 |
| 46219 | 46229 | 46237 | 46261 | 46271 | 46273 | 46279 | 46301 | 46307 | 46309 |
| 46327 | 46337 | 46349 | 46351 | 46381 | 46399 | 46411 | 46439 | 46441 | 46447 |
| 46451 | 46457 | 46471 | 46477 | 46489 | 46499 | 46507 | 46511 | 46523 | 46549 |
| 46559 | 46567 | 46573 | 46589 | 46591 | 46601 | 46619 | 46633 | 46639 | 46643 |
| 46649 | 46663 | 46679 | 46681 | 46687 | 46691 | 46703 | 46723 | 46727 | 46747 |
| 46751 | 46757 | 46769 | 46771 | 46807 | 46811 | 46817 | 46819 | 46829 | 46831 |
| 46853 | 46861 | 46867 | 46877 | 46889 | 46901 | 46919 | 46933 | 46957 | 46993 |
| 46997 | 47017 | 47041 | 47051 | 47057 | 47059 | 47087 | 47093 | 47111 | 47119 |
| 47123 | 47129 | 47137 | 47143 | 47147 | 47149 | 47161 | 47189 | 47207 | 47221 |
| 47237 | 47251 | 47269 | 47279 | 47287 | 47293 | 47297 | 47303 | 47309 | 47317 |
| 47339 | 47351 | 47353 | 47363 | 47381 | 47387 | 47389 | 47407 | 47417 | 47419 |
| 47431 | 47441 | 47459 | 47491 | 47497 | 47501 | 47507 | 47513 | 47521 | 47527 |
| 47533 | 47543 | 47563 | 47569 | 47581 | 47591 | 47599 | 47609 | 47623 | 47629 |
| 47639 | 47653 | 47657 | 47659 | 47681 | 47699 | 47701 | 47711 | 47713 | 47717 |
| 47737 | 47741 | 47743 | 47777 | 47779 | 47791 | 47797 | 47807 | 47809 | 47819 |
| 47837 | 47843 | 47857 | 47869 | 47881 | 47903 | 47911 | 47917 | 47933 | 47939 |
| 47947 | 47951 | 47963 | 47969 | 47977 | 47981 | 48017 | 48023 | 48029 | 48049 |
| 48073 | 48079 | 48091 | 48109 | 48119 | 48121 | 48131 | 48157 | 48163 | 48179 |
| 48187 | 48193 | 48197 | 48221 | 48239 | 48247 | 48259 | 48271 | 48281 | 48299 |
| 48311 | 48313 | 48337 | 48341 | 48353 | 48371 | 48383 | 48397 | 48407 | 48409 |
| 48413 | 48437 | 48449 | 48463 | 48473 | 48479 | 48481 | 48487 | 48491 | 48497 |
| 48523 | 48527 | 48533 | 48539 | 48541 | 48563 | 48571 | 48589 | 48593 | 48611 |
| 48619 | 48623 | 48647 | 48649 | 48661 | 48673 | 48677 | 48679 | 48731 | 48733 |
| 48751 | 48757 | 48761 | 48767 | 48779 | 48781 | 48787 | 48799 | 48809 | 48817 |
| 48821 | 48823 | 48847 | 48857 | 48859 | 48869 | 48871 | 48883 | 48889 | 48907 |
| 48947 | 48953 | 48973 | 48989 | 48991 | 49003 | 49009 | 49019 | 49031 | 49033 |
| 49037 | 49043 | 49057 | 49069 | 49081 | 49103 | 49109 | 49117 | 49121 | 49123 |
| 49139 | 49157 | 49169 | 49171 | 49177 | 49193 | 49199 | 49201 | 49207 | 49211 |
| 49223 | 49253 | 49261 | 49277 | 49279 | 49297 | 49307 | 49331 | 49333 | 49339 |
| 49363 | 49367 | 49369 | 49391 | 49393 | 49409 | 49411 | 49417 | 49429 | 49433 |
| 49451 | 49459 | 49463 | 49477 | 49481 | 49499 | 49523 | 49529 | 49531 | 49537 |
| 49547 | 49549 | 49559 | 49597 | 49603 | 49613 | 49627 | 49633 | 49639 | 49663 |
| 49667 | 49669 | 49681 | 49697 | 49711 | 49727 | 49739 | 49741 | 49747 | 49757 |
| 49783 | 49787 | 49789 | 49801 | 49807 | 49811 | 49823 | 49831 | 49843 | 49853 |
| 49871 | 49877 | 49891 | 49919 | 49921 | 49927 | 49937 | 49939 | 49943 | 49957 |
| 49991 | 49993 | 49999 | 50021 | 50023 | 50033 | 50047 | 50051 | 50053 | 50069 |
| 50077 | 50087 | 50093 | 50101 | 50111 | 50119 | 50123 | 50129 | 50131 | 50147 |
| 50153 | 50159 | 50177 | 50207 | 50221 | 50227 | 50231 | 50261 | 50263 | 50273 |
| 50287 | 50291 | 50311 | 50321 | 50329 | 50333 | 50341 | 50359 | 50363 | 50377 |
| 50383 | 50387 | 50411 | 50417 | 50423 | 50441 | 50459 | 50461 | 50497 | 50503 |
| 50513 | 50527 | 50539 | 50543 | 50549 | 50551 | 50581 | 50587 | 50591 | 50593 |
| 50599 | 50627 | 50647 | 50651 | 50671 | 50683 | 50707 | 50723 | 50741 | 50753 |
| 50767 | 50773 | 50777 | 50789 | 50821 | 50833 | 50839 | 50849 | 50857 | 50867 |
| 50873 | 50891 | 50893 | 50909 | 50923 | 50929 | 50951 | 50957 | 50969 | 50971 |
| 50989 | 50993 | 51001 | 51031 | 51043 | 51047 | 51059 | 51061 | 51071 | 51109 |
| 51131 | 51133 | 51137 | 51151 | 51157 | 51169 | 51193 | 51197 | 51199 | 51203 |
| 51217 | 51229 | 51239 | 51241 | 51257 | 51263 | 51283 | 51287 | 51307 | 51329 |
| 51341 | 51343 | 51347 | 51349 | 51361 | 51383 | 51407 | 51413 | 51419 | 51421 |
| 51427 | 51431 | 51437 | 51439 | 51449 | 51461 | 51473 | 51479 | 51481 | 51487 |
| 51503 | 51511 | 51517 | 51521 | 51539 | 51551 | 51563 | 51577 | 51581 | 51593 |
| 51599 | 51607 | 51613 | 51631 | 51637 | 51647 | 51659 | 51673 | 51679 | 51683 |
| 51691 | 51713 | 51719 | 51721 | 51749 | 51767 | 51769 | 51787 | 51797 | 51803 |
| 51817 | 51827 | 51829 | 51839 | 51853 | 51859 | 51869 | 51871 | 51893 | 51899 |
| 51907 | 51913 | 51929 | 51941 | 51949 | 51971 | 51973 | 51977 | 51991 | 52009 |
| 52021 | 52027 | 52051 | 52057 | 52067 | 52069 | 52081 | 52103 | 52121 | 52127 |
| 52147 | 52153 | 52163 | 52177 | 52181 | 52183 | 52189 | 52201 | 52223 | 52237 |
| 52249 | 52253 | 52259 | 52267 | 52289 | 52291 | 52301 | 52313 | 52321 | 52361 |
| 52363 | 52369 | 52379 | 52387 | 52391 | 52433 | 52453 | 52457 | 52489 | 52501 |
| 52511 | 52517 | 52529 | 52541 | 52543 | 52553 | 52561 | 52567 | 52571 | 52579 |
| 52583 | 52609 | 52627 | 52631 | 52639 | 52667 | 52673 | 52691 | 52697 | 52709 |
| 52711 | 52721 | 52727 | 52733 | 52747 | 52757 | 52769 | 52783 | 52807 | 52813 |
| 52817 | 52837 | 52859 | 52861 | 52879 | 52883 | 52889 | 52901 | 52903 | 52919 |
| 52937 | 52951 | 52957 | 52963 | 52967 | 52973 | 52981 | 52999 | 53003 | 53017 |
| 53047 | 53051 | 53069 | 53077 | 53087 | 53089 | 53093 | 53101 | 53113 | 53117 |
| 53129 | 53147 | 53149 | 53161 | 53171 | 53173 | 53189 | 53197 | 53201 | 53231 |
| 53233 | 53239 | 53267 | 53269 | 53279 | 53281 | 53299 | 53309 | 53323 | 53327 |
| 53353 | 53359 | 53377 | 53381 | 53401 | 53407 | 53411 | 53419 | 53437 | 53441 |
| 53453 | 53479 | 53503 | 53507 | 53527 | 53549 | 53551 | 53569 | 53591 | 53593 |
| 53597 | 53609 | 53611 | 53617 | 53623 | 53629 | 53633 | 53639 | 53653 | 53657 |
| 53681 | 53693 | 53699 | 53717 | 53719 | 53731 | 53759 | 53773 | 53777 | 53783 |
| 53791 | 53813 | 53819 | 53831 | 53849 | 53857 | 53861 | 53881 | 53887 | 53891 |
| 53897 | 53899 | 53917 | 53923 | 53927 | 53939 | 53951 | 53959 | 53987 | 53993 |
| 54001 | 54011 | 54013 | 54037 | 54049 | 54059 | 54083 | 54091 | 54101 | 54121 |
| 54133 | 54139 | 54151 | 54163 | 54167 | 54181 | 54193 | 54217 | 54251 | 54269 |
| 54277 | 54287 | 54293 | 54311 | 54319 | 54323 | 54331 | 54347 | 54361 | 54367 |
| 54371 | 54377 | 54401 | 54403 | 54409 | 54413 | 54419 | 54421 | 54437 | 54443 |
| 54449 | 54469 | 54493 | 54497 | 54499 | 54503 | 54517 | 54521 | 54539 | 54541 |
| 54547 | 54559 | 54563 | 54577 | 54581 | 54583 | 54601 | 54617 | 54623 | 54629 |
| 54631 | 54647 | 54667 | 54673 | 54679 | 54709 | 54713 | 54721 | 54727 | 54751 |
| 54767 | 54773 | 54779 | 54787 | 54799 | 54829 | 54833 | 54851 | 54869 | 54877 |
| 54881 | 54907 | 54917 | 54919 | 54941 | 54949 | 54959 | 54973 | 54979 | 54983 |
| 55001 | 55009 | 55021 | 55049 | 55051 | 55057 | 55061 | 55073 | 55079 | 55103 |
| 55109 | 55117 | 55127 | 55147 | 55163 | 55171 | 55201 | 55207 | 55213 | 55217 |
| 55219 | 55229 | 55243 | 55249 | 55259 | 55291 | 55313 | 55331 | 55333 | 55337 |
| 55339 | 55343 | 55351 | 55373 | 55381 | 55399 | 55411 | 55439 | 55441 | 55457 |
| 55469 | 55487 | 55501 | 55511 | 55529 | 55541 | 55547 | 55579 | 55589 | 55603 |
| 55609 | 55619 | 55621 | 55631 | 55633 | 55639 | 55661 | 55663 | 55667 | 55673 |
| 55681 | 55691 | 55697 | 55711 | 55717 | 55721 | 55733 | 55763 | 55787 | 55793 |
| 55799 | 55807 | 55813 | 55817 | 55819 | 55823 | 55829 | 55837 | 55843 | 55849 |
| 55871 | 55889 | 55897 | 55901 | 55903 | 55921 | 55927 | 55931 | 55933 | 55949 |
| 55967 | 55987 | 55997 | 56003 | 56009 | 56039 | 56041 | 56053 | 56081 | 56087 |
| 56093 | 56099 | 56101 | 56113 | 56123 | 56131 | 56149 | 56167 | 56171 | 56179 |
| 56197 | 56207 | 56209 | 56237 | 56239 | 56249 | 56263 | 56267 | 56269 | 56299 |
| 56311 | 56333 | 56359 | 56369 | 56377 | 56383 | 56393 | 56401 | 56417 | 56431 |
| 56437 | 56443 | 56453 | 56467 | 56473 | 56477 | 56479 | 56489 | 56501 | 56503 |
| 56509 | 56519 | 56527 | 56531 | 56533 | 56543 | 56569 | 56591 | 56597 | 56599 |
| 56611 | 56629 | 56633 | 56659 | 56663 | 56671 | 56681 | 56687 | 56701 | 56711 |
| 56713 | 56731 | 56737 | 56747 | 56767 | 56773 | 56779 | 56783 | 56807 | 56809 |
| 56813 | 56821 | 56827 | 56843 | 56857 | 56873 | 56891 | 56893 | 56897 | 56909 |
| 56911 | 56921 | 56923 | 56929 | 56941 | 56951 | 56957 | 56963 | 56983 | 56989 |
| 56993 | 56999 | 57037 | 57041 | 57047 | 57059 | 57073 | 57077 | 57089 | 57097 |
| 57107 | 57119 | 57131 | 57139 | 57143 | 57149 | 57163 | 57173 | 57179 | 57191 |
| 57193 | 57203 | 57221 | 57223 | 57241 | 57251 | 57259 | 57269 | 57271 | 57283 |
| 57287 | 57301 | 57329 | 57331 | 57347 | 57349 | 57367 | 57373 | 57383 | 57389 |
| 57397 | 57413 | 57427 | 57457 | 57467 | 57487 | 57493 | 57503 | 57527 | 57529 |
| 57557 | 57559 | 57571 | 57587 | 57593 | 57601 | 57637 | 57641 | 57649 | 57653 |
| 57667 | 57679 | 57689 | 57697 | 57709 | 57713 | 57719 | 57727 | 57731 | 57737 |
| 57751 | 57773 | 57781 | 57787 | 57791 | 57793 | 57803 | 57809 | 57829 | 57839 |
| 57847 | 57853 | 57859 | 57881 | 57899 | 57901 | 57917 | 57923 | 57943 | 57947 |
| 57973 | 57977 | 57991 | 58013 | 58027 | 58031 | 58043 | 58049 | 58057 | 58061 |
| 58067 | 58073 | 58099 | 58109 | 58111 | 58129 | 58147 | 58151 | 58153 | 58169 |
| 58171 | 58189 | 58193 | 58199 | 58207 | 58211 | 58217 | 58229 | 58231 | 58237 |
| 58243 | 58271 | 58309 | 58313 | 58321 | 58337 | 58363 | 58367 | 58369 | 58379 |
| 58391 | 58393 | 58403 | 58411 | 58417 | 58427 | 58439 | 58441 | 58451 | 58453 |
| 58477 | 58481 | 58511 | 58537 | 58543 | 58549 | 58567 | 58573 | 58579 | 58601 |
| 58603 | 58613 | 58631 | 58657 | 58661 | 58679 | 58687 | 58693 | 58699 | 58711 |
| 58727 | 58733 | 58741 | 58757 | 58763 | 58771 | 58787 | 58789 | 58831 | 58889 |
| 58897 | 58901 | 58907 | 58909 | 58913 | 58921 | 58937 | 58943 | 58963 | 58967 |
| 58979 | 58991 | 58997 | 59009 | 59011 | 59021 | 59023 | 59029 | 59051 | 59053 |
| 59063 | 59069 | 59077 | 59083 | 59093 | 59107 | 59113 | 59119 | 59123 | 59141 |
| 59149 | 59159 | 59167 | 59183 | 59197 | 59207 | 59209 | 59219 | 59221 | 59233 |
| 59239 | 59243 | 59263 | 59273 | 59281 | 59333 | 59341 | 59351 | 59357 | 59359 |

**Considering the first 1000 primes (that is all the primes up to 7919) gives this unsatisfactory result:**

| Digit | Count | Primes | Benford | chi-sqr |
|---|---|---|---|---|
| 1 | 160 | 16.0% | 30.1% | 66.07 |
| 2 | 160 | 16.0% | 17.6% | 1.47 |
| 3 | 130 | 13.0% | 12.5% | 0.21 |
| 4 | 130 | 13.0% | 9.7% | 11.30 |
| 5 | 130 | 13.0% | 7.9% | 32.62 |
| 6 | 140 | 14.0% | 6.7% | 79.72 |
| 7 | 120 | 12.0% | 5.8% | 66.30 |
| 8 | 20 | 2.0% | 5.1% | 18.97 |
| 9 | 10 | 1.0% | 4.6% | 27.94 |
| | 1000 | 100.0% | 100.0% | 304.60 |

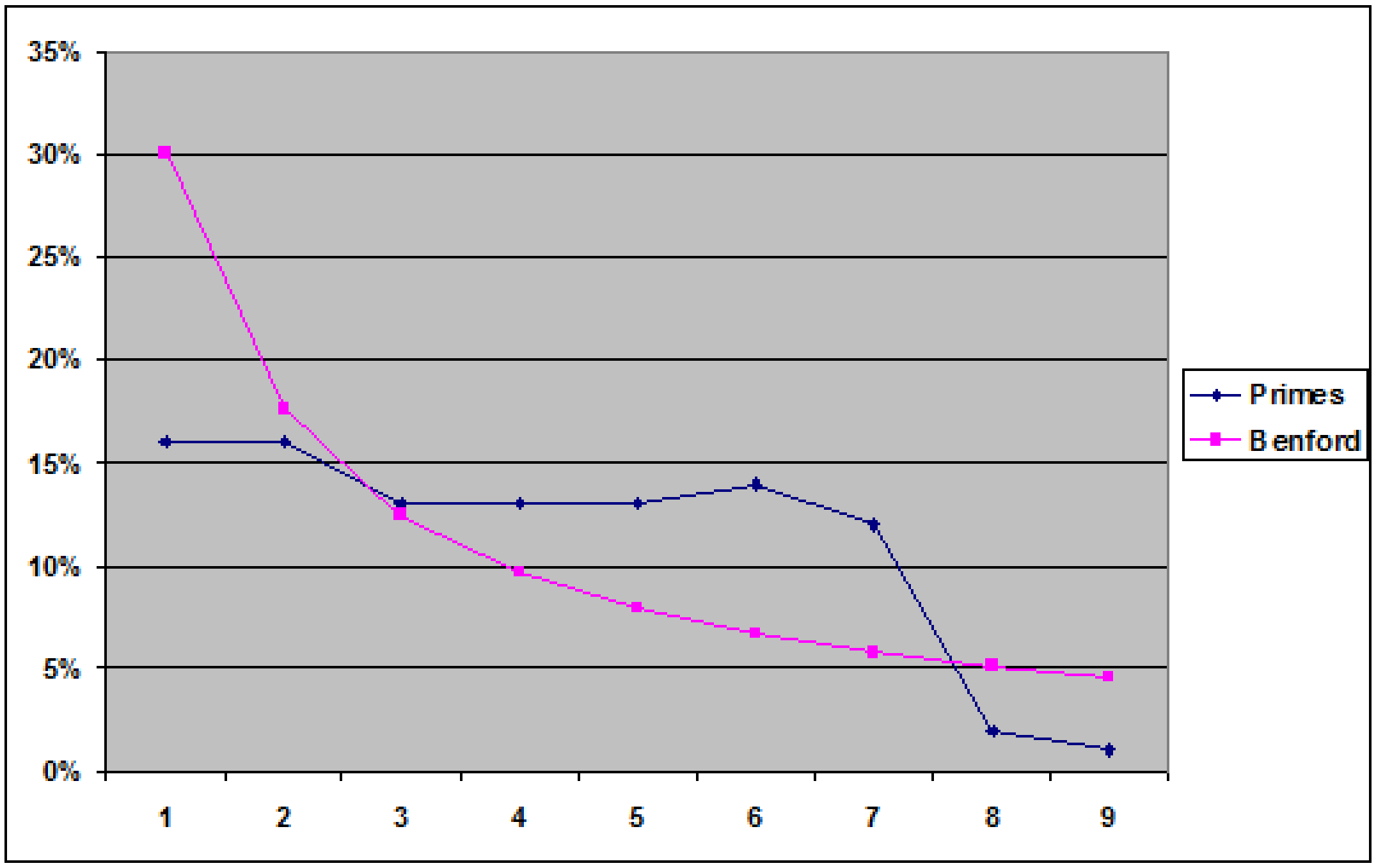

**Considering the first 3500 primes (that is all the primes up to 32609) gives slightly better result:**

| Digit | count | LD here | Benford |
|---|---|---|---|
| 1 | 1193 | 34.1% | 30.1% |
| 2 | 1129 | 32.3% | 17.6% |
| 3 | 394 | 11.3% | 12.5% |
| 4 | 139 | 4.0% | 9.7% |
| 5 | 131 | 3.7% | 7.9% |
| 6 | 135 | 3.9% | 6.7% |
| 7 | 125 | 3.6% | 5.8% |
| 8 | 127 | 3.6% | 5.1% |
| 9 | 127 | 3.6% | 4.6% |

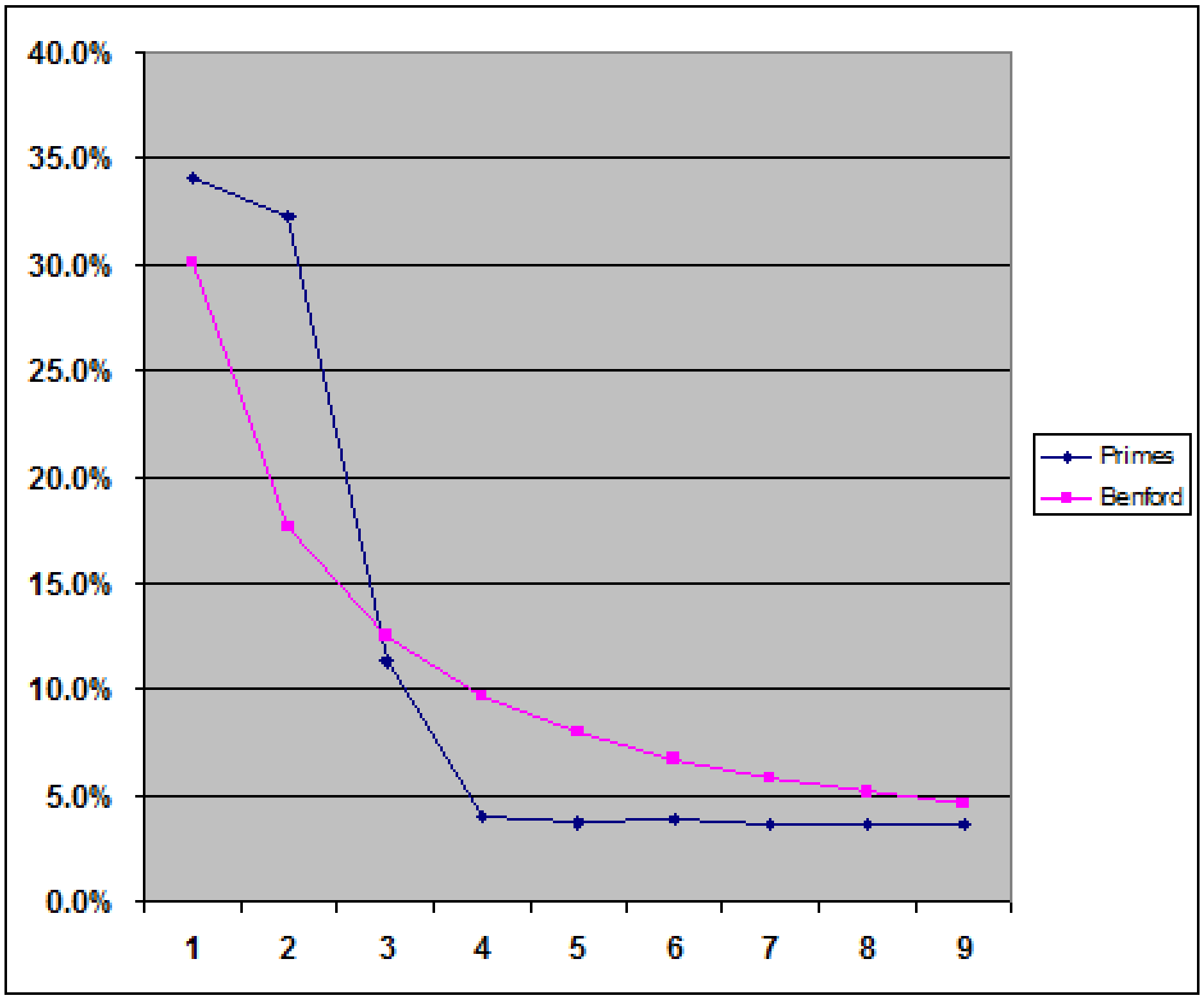

**Considering the first 6000 primes (that is all the primes up to 59359) gives a different result:**

| Digit | count | Primes | Benford |
|---|---|---|---|
| 1 | 1193 | 19.9% | 30.1% |
| 2 | 1129 | 18.8% | 17.6% |
| 3 | 1097 | 18.3% | 12.5% |
| 4 | 1069 | 17.8% | 9.7% |
| 5 | 998 | 16.6% | 7.9% |
| 6 | 135 | 2.3% | 6.7% |
| 7 | 125 | 2.1% | 5.8% |
| 8 | 127 | 2.1% | 5.1% |
| 9 | 127 | 2.1% | 4.6% |

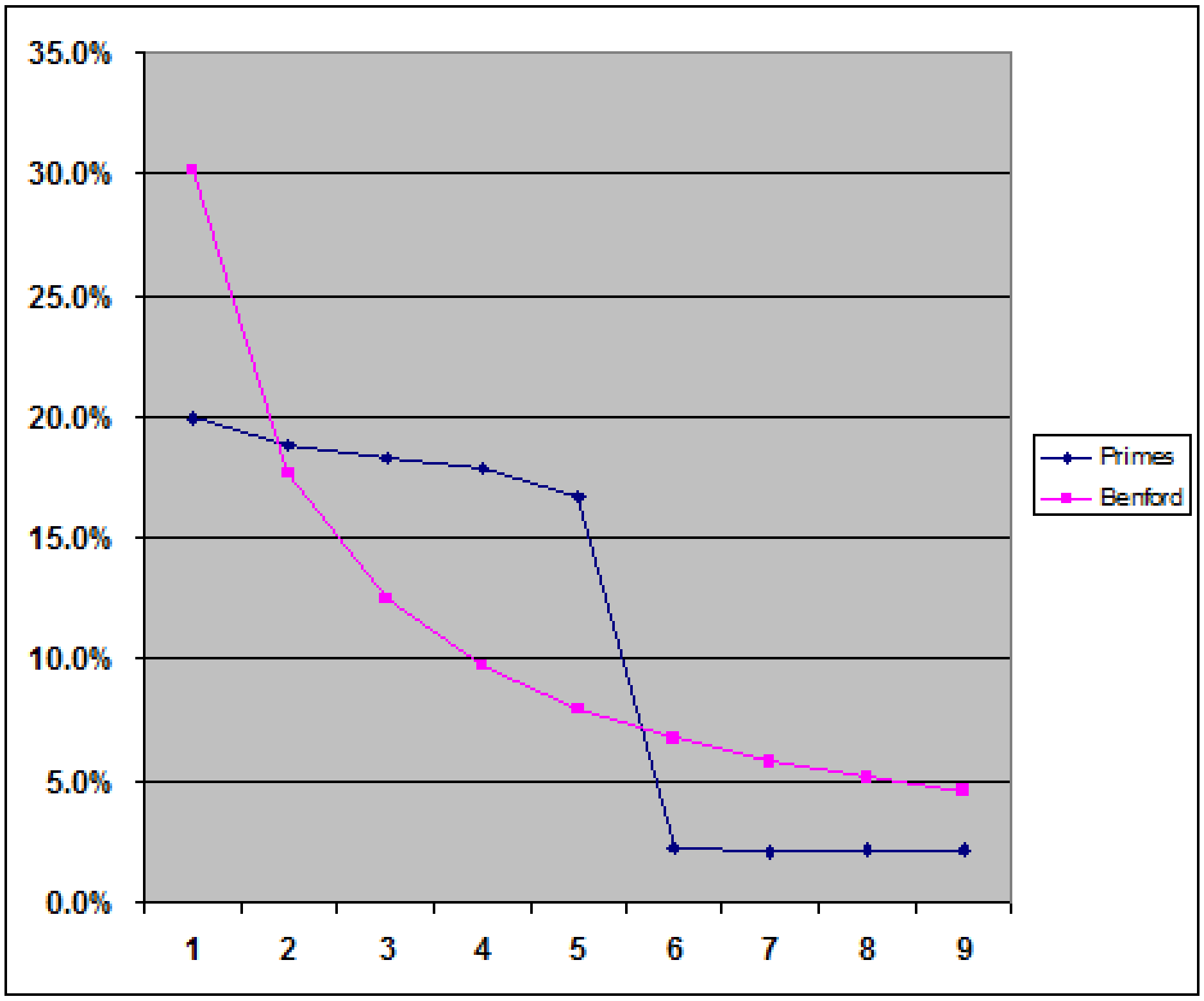

**Yet, digital distributions of these 3 collections of prime numbers are close to being monotonically decreasing.**

**Alex Ely Kossovsky**
**akossovs@yahoo.com**